\documentclass{memo-l}
\usepackage{CHQgreek}
\usepackage{CHQdefs}
\usepackage{latexsym}
\usepackage{amssymb}
\usepackage{diagrams}

\begin{document}
\frontmatter
\title{Homogeneous Spaces, Tits Buildings, and Isoparametric Hypersurfaces}
\author{Linus Kramer}
\address{Mathematisches Institut,
Universit\"at W\"urzburg,
Am Hubland,
D~97074~W\"urzburg,
Germany}
\email{kramer@mathematik.uni-wuerzburg.de}
\thanks{Supported by a Heisenberg Fellowship by the Deutsche
Forschungsgemeinschaft.}
\date{February 22, 2001}

\begin{abstract}
We classify 1-connected compact homogeneous spaces which have the same
rational cohomology as a product of spheres $\SS^{n_1}\times\SS^{n_2}$,
with $3\leq n_1\leq n_2$ and $n_2$ odd. As an application, we classify
compact generalized quadrangles (buildings of type $C_2)$ which
admit a point transitive automorphism group, and isoparametric
hypersurfaces which admit a transitive isometry group on one focal
manifold.
\end{abstract}

\maketitle
\setcounter{page}{6}

\tableofcontents

%                                                                      %
%                                                                      %
%                   Compact homogeneous quadrangles                    %
%                                                                      %
%                          Linus Kramer                                %
%                                                                      %
%                           Memoirs AMS                                %
%                                                                      %
%                         Wuerzburg 2000                               %
%                                                                      %
%                                                                      %
%                          CHQ-pre.tex                                 %
%                                                                      %
%                                                                      %
%                                                                      %
%%%%%%%%%%%%%%%%%%%%%%%%%%%%%%%%%%%%%%%%%%%%%%%%%%%%%%%%%%%%%%%%%%%%%%%%
\chapter*{Introduction}

The classification 
of compact Lie groups acting transitively on 
spheres due to Montgomery, Samelson and Borel 
was one of the main achievements 
in the early theory of compact transitive Lie transformation groups.
Later, Hsiang and Su classified compact transitive groups on 
(sufficiently highly connected) Stiefel manifolds. Their results were
extended by Scheerer, Schneider and other authors.
Given a compact 1-connected manifold $X$, it is in general difficult 
to classify all compact Lie groups which act transitively on $X$.
The problem becomes more complicated if only certain homotopy invariants
of $X$ such as the cohomology ring are known.
If the Euler characteristic of $X$ is positive, then 
results of Borel, De Siebenthal and Wang can be used. However,
if the Euler characteristic is 0, there is no general classification method.
\begin{center}
${*}\qquad{*}\qquad{*}$
\end{center}
In this book we classify all 1-connected
homogeneous spaces $G/H$ of compact Lie groups
which
have the same rational cohomology as a product
of spheres
\[
\SS^{n_1}\times\SS^{n_2},
\]
with $3\leq n_1\leq n_2$ and $n_2$ odd. Note that this implies that
the Euler characteristic of $G/H$ is 0.
Examples of such spaces are --- besides products of spheres ---
Stiefel manifolds of orthonormal 2-frames in real, complex, or
quaternionic vector spaces; another class of examples are certain
homogeneous sphere bundles. The following theorem is a direct consequence
of this rational classification.

\medskip
\noindent\textbf{Theorem}
{\em Let $X=G/H$ be a 1-connected compact homogeneous space of a compact
connected Lie group $G$. Assume that $G$ acts effectively and 
contains no normal transitive
subgroup, and that $X$ has the same integral cohomology as
a product of spheres
\[
\SS^{n_1}\times\SS^{n_2}
\]
with $3\leq n_1\leq n_2$ and $n_2$ odd.
There are the following possibilities for $G/H$ and the numbers
$(n_1,n_2)$.

\noindent
(1) Stiefel manifolds
\[
\begin{array}{cl}
\SO(2n)/\SO(2n-2)=V_2(\RR^{2n}) & (2n-2,2n-1) \\
\SU(n)/\SU(n-2)=V_2(\CC^n)      & (2n-3,2n-1) \\ 
\Sp(n)/\Sp(n-2)=V_2(\HH^n)      & (4n-5,4n-1).
\end{array}
\]

\noindent
(2) Certain homogeneous sphere bundles
\[
\begin{array}{cl}
\Sp(n)\times\Sp(2)/\Sp(n-1)\cdot\Sp(1) & (7,4n-1) \\
\Sp(n)\times\SU(3)/\Sp(n-1)\cdot\Sp(1) & (5,4n-1) \\ 
\Sp(n)\times\Sp(2)/\Sp(n-1)\cdot\Sp(1)\cdot\Sp(1) & (4,4n-1) . 
\end{array}
\]

\noindent
(3) Products of homogeneous spheres
\[
K_1/H_1\times K_2/H_2=\SS^{n_1}\times\SS^{n_2}
\]
where $K_1/H_1$ is one of the spaces
\begin{gather*}
\SO(n)/\SO(n-1)=\SS^{n-1}\\
\SU(n)/\SU(n-1)=\SS^{2n-1}\\
\Sp(n)/\Sp(n-1)=\SS^{4n-1}\\
\G_2/\SU(3)=\SS^6\\
\Spin(7)/\G_2=\SS^7\\
\Spin(9)/\Spin(7)=\SS^{15}
\end{gather*}
and $K_2/H_2$ is one of the spaces
\begin{gather*}
\SO(2n)/\SO(2n-1)=\SS^{2n-1}\\
\SU(n)/\SU(n-1)=\SS^{2n-1}\\
\Sp(n)/\Sp(n-1)=\SS^{4n-1}\\
\Spin(7)/\G_2=\SS^7\\
\Spin(9)/\Spin(7)=\SS^{15}
\end{gather*}

\noindent
(4) Some sporadic spaces
\[
\begin{array}{cl}
\E_6/\Ffour & (9,17) \\
\Spin(10)/\Spin(7) & (9,15) \\
\Spin(9)/\G_2=V_2(\OO^2) & (7,15) \\
\Spin(8)/\G_2=\SS^7\times\SS^7 & (7,7) \\
\SU(6)/\Sp(3)=\SU(5)/\Sp(2) & (5,9) \\
\Spin(10)/\SU(5)=\Spin(9)/\SU(4) & (6,15) \\
\Spin(7)/\SU(3)=V_2(\RR^8) & (6,7) \\
\Sp(3)/\Sp(1)\times\Sp(1) & (4,11) \\
\Sp(3)/\Sp(1)\times{}^\HH\rho_{3\lambda_1}(\Sp(1)) & (4,11) \\
\SU(5)/\SU(3)\times\SU(2) & (4,9).
\end{array}
\]}

\noindent
The proof proceeds as follows.
We show first that $G/H$ has the same rational homotopy groups
as the product $\SS^{n_1}\times\SS^{n_2}$. This homotopy theoretic
result follows from a generalization of a theorem by Cartan and Serre.
The rational homotopy groups of a compact Lie group $G$ can be determined
explicitly; they depend only on the Dynkin diagram
of $G$ and the rank of the central torus. 
In particular, a compact connected Lie group has the same rational
homotopy groups as a product of odd-dimensional spheres.
It follows in our situation that $\rk(G)-\rk(H)\in\{1,2\}$, depending
on whether $n_1$ is even or odd.
Using this fact, we determine the rational Leray-Serre
spectral sequence of the principal bundle
\[
H\too G\too G/H.
\]
If $n_1$ is odd, then the spectral sequence collapses, and
$G$ has the same rational cohomology as $H\times\SS^{n_1}\times\SS^{n_2}$. 
If $n_1$ is even, then there are non-zero differentials in the
spectral sequence and the situation is more complicated.
In both cases we obtain a relation 
between the numbers $(n_1,n_2)$ and the degrees of the primitive elements
in the rational cohomology of $G$ and $H$.

It follows also that $G$ has a semisimple normal transitive subgroup 
$K\unlhd G$ 
which has at most two almost simple factors. By a general result
about transitive permutation groups it suffices to determine
the pair $(K,K\cap H)\subseteq (G,H)$ in order to determine all
possibilities for the larger group $G$. Replacing $G$ by $K$ we may
therefore assume that $G$ is almost simple
or that $G$ is semisimple with two almost simple factors. 
The condition $n_1\geq 3$ guarantees then that $H$ is also semisimple. 
We determine all such pairs
$(G,H)$ with the right rational cohomology and all possible 
embeddings $H\rInto G$ using representation theory.
Thus, we obtain an explicit classification of these homogeneous spaces
together with the transitive group actions.

In the special case that $G/H$ has the same integral cohomology as
$\SS^{n_1}\times\SS^{n_2}$ we obtain the list of homogeneous 
spaces given in the theorem above.
In the course of the proof we determine also all compact connected Lie groups
which act transitively on 1-connected rational homology spheres;
in particular, we reprove the
classification of transitive actions on spheres
and on spaces which have the same homology as the Stiefel manifolds
$V_2(\RR^{2n+1})$.
\begin{center}
${*}\qquad{*}\qquad{*}$
\end{center}
We apply our result to a problem in submanifold geometry.
A closed hypersurface in a sphere is called isoparametric if
its principal curvatures are constant.
Hsiang and Lawson classified all isoparametric hypersurfaces which
admit a transitive group of isometries; these homogeneous
isoparametric hypersurfaces arise as principal orbits of
isotropy representations
of non-compact symmetric spaces of rank 2. 

By a result of M\"unzner, the number of distinct principal curvatures
of an isoparametric hypersurface is $g=1,2,3,4,6$; the hypersurfaces
with $g=1,2,3$ have been classified in the 30s by Segre and Cartan.
Some hypersurfaces with $g=6$ were classified by Dorfmeister and Neher;
the full classification for $g=6$ still seems to be an open problem.
The case of hypersurfaces with $g=4$ is much more difficult.
Takagi proved uniqueness for isoparametric hypersurfaces with $g=4$
distinct principal curvatures and multiplicities $(1,k)$. On the other
hand, Ferus, Karcher and M\"unzner showed that there are many
non-homogeneous isoparametric hypersurfaces with $g=4$ distinct principal
curvatures, and Stolz recently obtained sharp number theoretic restrictions
on the possible dimensions of such hypersurfaces.
Some of the known inhomogeneous examples have homogeneous focal manifolds,
so the question arises if the classification by Hsiang and Lawson can be
generalized to transitive actions on focal manifolds. 
In view of the classification by Segre, Cartan, Dorfmeister and Neher,
the cases $g=1,2,3$ are less interesting. If $g=4$ and if 
the multiplicities $(m_1,m_2)$ of the isoparametric hypersurface
are large enough, then the focal manifolds have the same integral
cohomology as a product of spheres.
We apply our result about homogeneous spaces to this problem and
obtain a complete classification.

\medskip
\noindent
\textbf{Theorem}
{\em Let $M$ be an isoparametric hypersurface with 4 distinct principal
curvatures. Suppose that the isometry group of $M$ acts transitively
on one of the focal manifolds, and that this focal manifold is
2-connected. Then either the hypersurface
itself is homogeneous (and explicitly known),
or it is of Clifford type with
multiplicities $(8,7)$ or $(3,4k-4)$.}

\medskip\noindent
Recently, Wolfrom showed in his Ph.D. Thesis that the theorem above
holds also if one drops the assumption that the focal manifold is
$2$-connected, and proved a similar result for $g=6$. The
final result is as follows.

\medskip\noindent
\textbf{Theorem}
{\em Let $M$ be an isoparametric hypersurface, and
suppose that the isometry group of $M$ acts transitively
on one of the focal manifolds. Then either the hypersurface
itself is homogeneous (and explicitly known),
or it is of Clifford type with
multiplicities $(8,7)$ or $(3,4k-4)$.}

\medskip\noindent
This theorem gives in particular a new, independent proof for
the classification by Hsiang and Lawson.
\begin{center}
${*}\qquad{*}\qquad{*}$
\end{center}
Another application is in topological geometry. Every isotropic 
simple algebraic group, in particular every non-compact real simple
Lie group, gives rise to a spherical
Tits building. These buildings are 
characterized by the so-called Moufang condition.
In the case of a non-compact simple Lie group, the building 
inherits a compact
topology from the group action. These compact buildings are closely
related to symmetric spaces, F\"urstenberg boundaries and isoparametric 
submanifolds. They play also a r\^ole in the theory of Hadamard spaces
and rigidity results.

Tits classified all irreducible
spherical buildings of rank at least 3 by showing that
they automatically satisfy the Moufang property. 
In contrast to this generalized polygons, i.e.~spherical
buildings of rank 2, need not
be Moufang, and there is no way to classify them without further
assumptions. 

In view of the examples above, it is natural to 
consider compact generalized polygons
and to try to classify them in terms of their automorphism
groups. A result of Knarr and the author states that such a building is
of type $A_2$, $C_2$, or $G_2$. Topologically, a compact
generalized polygon looks very similar to an isoparametric foliation
with $g=3,4,6$ distinct principal curvatures, respectively; in particular,
the cohomology of these spaces can be determined.
This is the analogue of M\"unzner's theorem mentioned above. 

An $A_2$-building is the same as a projective
plane; all compact homogeneous projective planes have been classified
by L\"owen and Salzmann. The compact homogeneous $G_2$-buildings
have been classified by the author. The remaining cases are
the generalized quadrangles, i.e.~the buildings of type $C_2$.
As in the case of isoparametric hypersurfaces, this is much
more involved. The easiest case here are the quadrangles with Euler
characteristic 4 which have been classified by the author. The remaining
case, namely compact quadrangles of Euler characteristic 0, is 
very interesting, since the
inhomogeneous isoparametric hypersurfaces discovered by 
Ferus, Karcher and M\"unzner are examples of such
(non-Moufang) quadrangles.

By general arguments, a transitive automorphism group on 
a 1-connected compact quadrangle contains a compact
transitive Lie subgroup (the automorphism group of a compact
building is in general not compact).
Thus, we can apply our result to classify transitive
actions of compact Lie groups on compact quadrangles. We obtain a list of 
all possible transitive actions. 
For three infinite series we show that
the group action determines the quadrangle up to isomorphism.
In fact, we can (to some extent) use the
same geometric methods and arguments both for isoparametric hypersurfaces
with $g=4$ distinct principal curvatures and for compact quadrangles.

Nevertheless, the classification of the compact quadrangles is much more
difficult than the classification of isoparametric hypersurfaces.
This is due to the fact that an isoparametric hypersurface sits inside some
Euclidean space on which the group acts linearly.
Even though this surrounding space also exists for generalized quadrangles,
it has no natural Euclidean structure, and it is not
a priori clear that the group action has to be linear. Therefore
representation theory has to be replaced by arguments about
compact transformation groups acting on locally compact spaces.
The result is as follows.

\medskip
\noindent
\textbf{Theorem}
{\em Let $\frak G$ be a compact connected quadrangle. Assume that
the point space is 9-connected, and that the automorphism group
is point transitive. Then $\frak G$ is a Moufang quadrangle
(in fact the dual of 
a classical quadrangle associated to a hermitian form
over $\RR$, $\CC$ or $\HH$).}

\medskip\noindent
At this point, we should mention the following related results by
Grundh\"ofer, Knarr and the author.

\medskip
\noindent
\textbf{Theorem}
{\em Let $\Delta$ be a compact connected irreducible
spherical building of rank at least 2.
Assume that the automorphism group
is chamber transitive. Then $\Delta$ is the Moufang building
associated to a simple non-compact Lie group $G$. If
$H\subseteq\Aut(\Delta)$
is a connected chamber transitive subgroup, then either $H=G$,
or $H$ is compact.}

\medskip\noindent
All compact connected chamber transitive groups in the theorem
above were determined by Eschenburg and Heintze.
Combining these results, we have the following theorem.

\medskip
\noindent
\textbf{Theorem}
{\em Let $\Delta$ be a compact connected irreducible 
spherical building of rank 
$k\geq 2$. Assume that the automorphism group
is transitive on one type of vertices of $\Delta$.
If the building is of type $C_2$
assume in addition that either 
(1) the vertex space in question is 9-connected, or 
(2) that the two vertex sets have the same dimension, or
(3) that the action is chamber transitive.
Then $\Delta$ is a Moufang building associated to a simple non-compact
real Lie group of real rank $k$.}

\medskip\noindent
In the $C_2$-case, the assumption on the connectivity cannot be 
dropped completely,
since there are counterexamples which are not highly connected.
In view of the results in the present book, the following conjecture
for the $C_2$-case is very natural.

\medskip\noindent
\textbf{Conjecture}
{\em Let $\frak G$ be a compact connected generalized quadrangle.
Assume that the automorphism group is point or line transitive. The either
$\frak G$ is a Moufang quadrangle, or $\frak G$ is of Clifford type
with topological parameters $(3,4k)$ or $(7,8)$.}

\medskip\noindent
Finally, I should mention here the following new and beautiful result by
Immervoll: every isoparametric hypersurface with $g=4$ distinct principal
curvatures is a $C_2$-building.
\begin{center}
${*}\qquad{*}\qquad{*}$
\end{center}
The material is organized as follows. In the first chapter we collect
some well-known facts about the algebraic topology of fibrations,
spectral sequences, and Eilenberg-MacLane spaces. This chapter
should be accessible for any reader with basic knowledge about
algebraic topology. In Chapter 2 we prove an extension of the
Cartan-Serre theorem about rational homotopy groups by a standard
homotopy theoretic method. This is applied in Chapter 3 in order to
determine the rational Leray-Serre spectral
sequence associated to the transitive group action. 

All facts about representations of compact Lie groups on real, 
complex or quaternionic vector spaces which are used in the
classification are presented in Chapter 4, which
also contains tables about compact almost simple
Lie groups and their low-dimensional representations. These facts
may be of some independent interest.

Our classification of homogeneous spaces is carried out in
Chapters 5 and 6. The classification is stated at the end of
Chapter 3. Readers who are willing to accept this result without
entering into the details can skip Chapters 5 and 6. However,
the representation theory of Chapter 4 is used again in the last two
chapters.

In Chapter 7 we classify transitive actions of compact Lie groups
on compact quadrangles. The homogeneous focal manifolds of isoparametric
hypersurfaces with $g=4$ distinct principal curvatures are
classified in Chapter 8, together with all possible transitive group 
actions. Here, the reader is assumed to be familiar either with
topological geometry or with isoparametric hypersurfaces. 

The logical dependencies of the chapters are as follows:
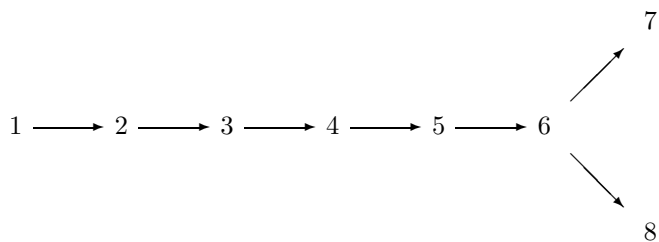
\begin{diagram}[size=2em]
&&&&&&&&&&&& 7 \\
&&&&&&&&&&& \ruTo \\
 1 & \rTo & 2 & \rTo & 3 & \rTo & 4 & \rTo & 5 & \rTo & 6 \\
&&&&&&&&&&& \rdTo \\
&&&&&&&&&&&& 8
\end{diagram} 
In particular, Chapter 7 (topological geometry) and Chapter 8
(submanifold geometry) are in principle independent of each other.
Nevertheless, the subjects of these two chapters,
isoparametric  hypersurfaces and compact polygons, share many
geometric properties which I tried to emphasize. Often, the proofs
in Chapters 7 and 8 are quite similar and geometric.
Besides the global properties of isoparametric hypersurfaces,
very little differential geometry is needed in the classification.
I hope that the present book will be useful both for
differential geometers and for topological geometers, and that it
helps to broaden the bridge between the two fields.
\begin{center}
${*}\qquad{*}\qquad{*}$
\end{center}
The book is a revised, corrected and expanded version of my
{\sl Habilitationsschrift}. I would like to thank 
my teachers Theo Grundh\"ofer and Reiner Salzmann for their
constant interest and support.
Harald Biller, Oliver Bletz, Norbert Knarr, Gerhard R\"ohrle,
Stephan Stolz, Markus Stroppel, Hendrik Van Maldeghem, 
and Martin Wolfrom made helpful suggestions or spotted errors.
I used a computer program by Richard B\"odi  and Michael Joswig
to check the tables for the representations of simple Lie groups.
Robert Bryant and Friedrich Knop helped me with a 
question about a certain representation.
The commutative diagrams in the original manuscript
were drawn with Paul Taylors diagrams \TeX-package.
Finally, I would like to thank my wife, Katrin Tent, not only
for reading the manuscript. Without her support, this book would not
have been possible.
 
\bigskip
\hfill Gerbrunn, February 2001\par
\hfill Linus Kramer

\vfill
\begin{raggedleft}
\textit{Man vergilt seinem Lehrer schlecht,\\
wenn man immer nur der Sch\"uler bleibt.\\
F.~N.\\}
\end{raggedleft}

\mainmatter
%%%%%%%%%%%%%%%%%%%%%%%%%%%%%%%%%%%%%%%%%%%%%%%%%%%%%%%%%%%%%%%%%%%%%%%%
%                                                                      %
%                                                                      %
%                   Compact homogeneous quadrangles                    %
%                                                                      %
%                          Linus Kramer                                %
%                                                                      %
%                           Memoirs AMS                                %
%                                                                      %
%                         Wuerzburg 2000                               %
%                                                                      %
%                                                                      %
%                            CHQ1.tex                                  %
%                                                                      %
%                                                                      %
%                                                                      %
%%%%%%%%%%%%%%%%%%%%%%%%%%%%%%%%%%%%%%%%%%%%%%%%%%%%%%%%%%%%%%%%%%%%%%%%
\chapter{The Leray-Serre spectral sequence}

Our main tool from algebraic topology is the Leray-Serre spectral
sequence. It  relates the cohomology rings of the fibre and the
base of a fibration with the cohomology of the total space.
Although spectral sequences are standard devices in topology,
they tend to be somewhat intimidating to non-specialists.
The aim of this
chapter is to give a basic introduction to the relevant
notions and techniques.

Let us consider
a specific example, the Leray-Serre spectral sequence with
field coefficients.
Let $K$ be a field, and let $F,B$ be topological
spaces. The K\"unneth Theorem asserts that the $K$-cohomology of
the product $E=F\times B$ is given by 
\[
\bfH^k(E;K)\isom\bigoplus_{i+j=k}\bfH^i(F;K)\otimes\bfH^j(B;K).
\]
It is convenient to visualize the $K$-modules
$\bfE^{i,j}=\bfH^i(F;K)\otimes\bfH^j(B;K)$ as distributed
on the lattice $\ZZ^2\SUB\RR^2$: attached to the point
$(i,j)$ is the vector space $\bfE^{ij}$. Then
$\bfH^k(E;K)$ is obtained by adding up all vector spaces
along the line $i+j=k$.

Now suppose that $E$ is not a product, but the total space of
a ($K$-simple)
fibre bundle $F\too E\too B$. Then the cohomology of $E$
is 'smaller' than the cohomology of the trivial bundle $F\times B$.
The recipe to obtain the $K$-modules $\bfH^k(E;K)$ is as follows.
Start with the collection of $K$-modules
$\bfE^{i,j}_2=\bfH^i(F;K)\otimes\bfH^j(B;K)$ as before.
There exists a collection of maps, the differentials, denoted
$d_2:E_2^{i,j}\too E_2^{i+2,j-1}$.
These differentials should be visualized as arrows going from
$(i,j)$ to $(i+2,j-1)$.
\begin{diagram}[size=2em]
\bfE_2^{i,j} \\
& \rdTo(4,2)^{d_2} \\ 
&&&& \qquad\bfE_2^{i+2,j-1}
\end{diagram}
They satisfy the relation $d_2\circ d_2=0$,
so we can take their cohomology (kernel mod image) at each
point $(i,j)$. Call the resulting $K$-module $\bfE_3^{i,j}$.
Note that this is a quotient of a submodule of $\bfE_2^{i,j}$,
so its dimension is smaller.
Again, these $K$-modules should be viewed as distributed in the
plane. There is another differential $d_3$, this time from
$(i,j)$ to $(i+3,j-2)$. Now this process is iterated
{\sl ad infinitum}. The arrows $d_2,d_3,d_4,\ldots$
become longer, and their slope approaches $-1$.
In the limit, one obtains a collection of $K$-modules
denoted $\bfE_\infty^{i,j}$.
Similarly as in the K\"unneth Theorem one has
\[
\bfH^k(E;K)\isom\bigoplus_{i+j=k}\bfE^{i,j}_\infty.
\]

In fact, one has not to go to infinity in this situation.
Since $\bfE_2^{i,j}=0$ for $i<0$ or $j<0$, the arrows $d_r$
starting or ending at $(i,j)$ are trivial maps for $r$
large enough (e.g. $r>\max\{i,j\}$). 
Thus, the modules $\bfE^{i,j}_r$ become stationary after
some time.

However, there is one big
problem: in general, no information is given about the arrows
$d_2,d_3,\ldots$! Thus, it seems to be impossible to determine
$\bfE_3^{i,j},\bfE_4^{i,j},\ldots$.
Here, the multiplicative structure of the cohomology becomes
important. The arrows $d_r$ act as derivations, and this
makes it (often) feasible to determine all terms in the
spectral sequence.

The material of this chapter can be found in
McCleary \cite{McCl85}, Borel \cite{BorelThesis,BorelLNM}, 
Spanier \cite{Spa66}, Whitehead \cite{Whi78}, and
Fomenko-Fuchs-Gutenmacher \cite{FFG86}.

\emph{Throughout this chapter, $R$ is a principal ideal domain.}

\section{Additive structure}

\begin{Num}\textsc{Graded and bigraded modules}\psn
A \emph{graded $R$-module} is a direct sum 
\[
\bfM=\bfM^\bullet=\bigoplus_{i\in\ZZ} \bfM^i
\]
of $R$-modules $\bfM^i$, indexed by the integers. 
Similarly, a \emph{bigraded $R$-module} is a direct sum
\[
\bfM=\bfM^{\bullet\bullet}=
\bigoplus_{i,j\in\ZZ} \bfM^{i,j}
\]
of $R$-modules $\bfM^{i,j}$, indexed by pairs of integers. 
The elements of $\bfM^i$ or $\bfM^{i,j}$ are called
\emph{homogeneous} of \emph{degree} $i$ or \emph{bidegree} $(i,j)$,
respectively. A graded or bigraded module is of \emph{finite type}
if the $\bfM^i$ or $\bfM^{i,j}$ are finitely generated.
A submodule $\bfN\SUB\bfM^\bullet$ is graded if
$\bfN=\bigoplus_{i\in\ZZ}\bfN^i$, where
$\bfN^i=\bfN\cap\bfM^i$.
\end{Num}

\begin{Num}\textsc{Total gradings and tensor products}\psn
Associated to a bigraded module $\bfM^{\bullet\bullet}$
is the graded module $\bfM^\bullet=
\Tot(\bfM)$ which is graded by the total
degree,
\[
\Tot(\bfM)^i=
\bigoplus_{j+k=i}\bfM^{j,k}.
\]
A typical example for a bigraded module is a tensor product
of graded modules. Put
\[
(\bfM\otimes\bfN)^{i,j}=\bfM^i\otimes\bfN^j.
\]
Then the corresponding graded module is
\[
(\bfM\otimes\bfN)^i=\bigoplus_{j+k=i}(
\bfM^j\otimes\bfN^k).
\]
\end{Num}
Another example is obtained from filtrations.

\begin{Num}\textsc{Filtrations and associated gradings}\psn
A \emph{filtration} of a module $\bfM$ is
a collection of submodules 
\[
\cdots
\SUB\bfF^i\bfM\SUB\bfF^{i-1}\bfM\SUB
\bfF^{i-2}\bfM\SUB\bfF^{i-3}\bfM\SUB
\cdots
\]
The filtration is \emph{bounded} if $\bfF^0\bfM=\bfM$,
and \emph{convergent} if $\bigcap_{i\geq 0}\bfF^i\bfM=0$.
If $\bfM=\bfM^\bullet$ is graded, then one requires
that the submodules in the filtration are graded.
An example of a filtration is the following. Let $X$ be a
CW complex, and let $X^{(k)}$ denote its $k$-skeleton.
Then
\[
\bfF^i\bfH^j(X;R)=
\mathrm{ker}[\bfH^j(X^{(i-1)};R)\oot\bfH^j(X;R)]
\]
defines a filtration of the cohomology module of $X$.

Associated to a (bounded and convergent) filtration is
the bigraded module $\bfG(\bfM)$ which is defined
by
\[
\bfG(\bfM)^{i,j}=
\frac{\bfF^{i}\bfM^{i+j}}
{\bfF^{i+1}\bfM^{i+j}}
\]
In general, it can be difficult to recover the graded module $\bfM$ 
from $\bfG(\bfM)$ (this is a problem about module extensions). 
However, if $\bfM^\bullet$ is
of finite type, and if $R$ is a field, then
\[
\Tot(\bfG(\bfM))^i=
\bigoplus_{j\in\ZZ}\frac{\bfF^{j}\bfM^i}
{\bfF^{j+1}\bfM^i}\isom\bfM^i.
\]
\end{Num}

\begin{Num}\textsc{Differential graded modules}\psn
A map $f:\bfM^\bullet\too\bfN^\bullet$ of degree $r$ 
between graded modules is an $R$-linear
map which increases degrees by $r$,
i.e.~$f(\bfM^i)\SUB\bfN^{i+r}$. A differential is a map
$d:\bfM^\bullet\too\bfM^\bullet$ of degree 1, with
$d^2=0$. The \emph{cohomology} of $(\bfM,d)$ is the graded module
\[
\bfH^i(\bfM)=\frac{\mathrm{ker}[\bfM^i\TOO d\bfM^{i+1}]}
{\mathrm{im}[\bfM^{i-1}\TOO d\bfM^i]}.
\]
An $R$-linear map $f:\bfM^{\bullet\bullet}\too
\bfN^{\bullet\bullet}$  has bidegree $(r,s)$ if
$f(\bfM^{i,j})\SUB \bfN^{i+r,j+s}$. A 
\emph{differential} $d$ of bidegree $(r,1-r)$ is a map
$d:\bfM\too\bfM$ of bidegree $(r,1-r)$ with $d^2=0$.
The \emph{cohomology} of $(\bfM,d)$ is the bigraded module
\[
\bfH^{i,j}(\bfM)=\frac{\mathrm{ker}[\bfM^{i,j}\TOO d\bfM^{i+r,j+1-r}]}
{\mathrm{im}[\bfM^{i-r,j-1+r}\TOO d\bfM^{i,j}]}.
\]
Note that $d$ is a differential on $\Tot(\bfM)$, and that
$\Tot(\bfH^{\bullet\bullet}(\bfM))=
\bfH^\bullet(\Tot(\bfM))$. The pair $(\bfM,d)$ is called a
\emph{differential (bi)graded module}.
\end{Num}

\begin{Num}\textsc{Spectral sequences}\psn
An $\bfE_2$-\emph{spectral sequence} is a collection of
differential bigraded modules $\bfE_r$, indexed by $r=2,3,4,\cdots$,
endowed with differentials $d_r$ of bidegree 
$(r,1-r)$, such that 
\[
\bfE_{r+1}\isom\bfH(\bfE_r).
\]
The spectral sequence \emph{converges} if for every 
pair $(i,j)$ there exists a number $n=n_{i,j}$ such the two maps 
\begin{diagram}
\ & \\
& \rdTo^{d_r}(2,1) &
 \bfE^{i,j}_r \\
&&& \rdTo^{d_r}(2,1) 
\end{diagram}
are trivial for all $r\geq n$. The resulting module 
$\bfE_n^{i,j}\isom\bfE_{n+1}^{i,j}\isom\bfE_{n+2}^{i,j}\isom
\bfE_{n+3}^{i,j}\isom\cdots$ is denoted by
$\bfE_\infty^{i,j}$, and one says that the spectral sequence
\emph{converges} to the bigraded module $\bfE_\infty^{\bullet\bullet}$.
The spectral sequence \emph{collapses} if all differentials
$d_r$ vanish, and it collapses at $\bfE_n$ if all differentials vanish 
for $r\geq n$, in which case $\bfE_n\isom\bfE_\infty$.
Note also that if $\bfE_n^{i,j}=0$, then $\bfE_r^{i,j}=0$ for all
$r\geq n$.

Each term $\bfE_r$ can be visualized as the grid $\ZZ^2\SUB\RR^2$. The 
point with coordinates $(i,j)$ represents $\bfE_r^{i,j}$, and the 
differentials are arrows between elements of the grid
pointing $r$ steps to the right- and $r-1$ steps downwards.
As $r$ increases, the arrows become longer and their slope approaches
$-1$.
\begin{diagram}
\bfE_r^{i,j} && \rLine^{r} && {} \\
& \rdTo(4,2)^{d_r} &&& \dLine_{r-1}\\ 
&&&& \qquad\bfE_r^{i+r,j+1-r}
\end{diagram}
Typically, the $\bfE_2$-term contains lots of zeros. The first differentials
are zero until the arrows become long enough to reach from one non-zero
entry to another one. If the region in $\RR^2$
containing the non-zero terms is
bounded, the arrows become too long after some time,
and the spectral sequence collapses.
\end{Num}

\begin{Lem}
If $R$ is a field, and if $\dim(\bfE_n)<\infty$, then
\[
\dim(\bfE_\infty)\leq\dim(\bfE_n).
\]
Equality holds if and only if the spectral sequence collapses
at $n$.

\begin{proof}
We have 
\[
\dim(\bfE_r)=\dim(\mathrm{im}(d_r))
+\dim(\mathrm{ker}(d_r)),
\]
whence 
\[
\dim(\bfE_{r+1})=\dim(\bfE_r)-2\dim(\mathrm{im}(d_r)).
\]
\end{proof}
\end{Lem}

\begin{Thm}[The Leray-Serre spectral sequence]
\label{LeraySerreSS}
\ \psn Let
\begin{diagram}
F & \rTo & E \\
&& \dTo \\
&& B
\end{diagram}
be a fibration over a path-connected space $B$.
The fundamental group $\pi_1(B)$ acts on the fibre
$F$ and hence on the cohomology $\bfH^\bullet(F;R)$. If this action
is trivial (e.g. if $B$ is 1-connected), then the fibration is
called \emph{$R$-simple}. Suppose that this is the case.
Then there is an $\bfE_2$-spectral sequence which converges to
the bigraded module associated to some 
(bounded and convergent) filtration of $\bfH^\bullet(E;R)$,
with $\bfE_2^{i,j}\isom\bfH^i(B;\bfH^j(F;R))$.

\begin{proof} See Spanier \cite{Spa66} 9.4.9.
\end{proof}
\end{Thm}
Note that
$\bfE_r^{i,j}=0$ if $i<0$ or $j<0$. Such a \emph{first-quadrant spectral
sequence} is always convergent because the arrows eventually stick
out of the first quadrant.

Put $\bfF^i\bfH^j=\bfF^i\bfH^j(E;R)$. There are short exact sequences
\begin{diagram}
&&0 & \rTo & \bfF^i\bfH^i & \rTo^{\isom} &\bfE^{i,0}_\infty & \rTo  & 0 \\
0 & \rTo &
\bfF^i\bfH^i & \rTo & \bfF^{i-1}\bfH^i & \rTo &\bfE^{i-1,1}_\infty 
&\rTo & 0 \\
0& \rTo &
\bfF^{i-1}\bfH^i & \rTo & \bfF^{i-2}\bfH^i & \rTo & \bfE^{i-2,2}_\infty 
&\rTo & 0 \\
&&&&\dDots\\
0& \rTo &
\bfF^{1}\bfH^i & \rTo & \bfH^i & \rTo & \bfE^{0,i}_\infty &\rTo & 0 
\end{diagram}
We consider the 'edges' of the first quadrant, i.e.~the $X$- and $Y$-axis.
Note that 
\[
\bfE^{0,i}_\infty\isom\frac{\bfF^0\bfH^i}{\bfF^1\bfH^i}
\]
is a quotient of $\bfF^0\bfH^i=\bfH^i(E;R)$, and that 
\[
\bfE^{0,i}_\infty\SUB\bfE_2^{0,i}
\]
because all arrows which come from the left are zero.
Since the base $B$ is assumed to be path connected,
$\bfH^0(B;\bfH^\bullet(F;R))\isom\bfH^\bullet(F;R)$ and
the following diagram commutes, cp.~Spanier~\cite{Spa66}~9.5.
\begin{diagram}
\bfH^i(F;R) && \lTo && \bfH^i(E;R) \\
\dTo^\isom &&& & \dTo \\
\bfE_2^{0,i} && \lInto & & \bfE_\infty^{0,i}
\end{diagram}
The projection $E\too B$ can also be interpreted in terms of the 
Leray-Serre spectral sequence. Note that there is a surjection 
\[
\bfE_2^{i,0}\too\bfE_\infty^{i,0}
\]
because all arrows starting on the $X$-axis are zero.
Moreover,
\[
\bfE^{i,0}_\infty\isom\frac{\bfF^i\bfH^i}{\bfF^{i+1}\bfH^i}
=\bfF^i\bfH^i\SUB\bfH^i.
\]
If the fibre $F$ is path connected, then
$\bfH^\bullet(B;\bfH^0(F;R))\isom\bfH^\bullet(B;R)$ and the diagram
\begin{diagram}
\bfH^i(E;R) && \lTo && \bfH^i(B;R) \\
\uInto && && \uTo_\isom \\
\bfE_\infty^{i,0} && \lTo && \bfE_2^{i,0}
\end{diagram}
commutes, cp.~Spanier~\cite{Spa66}~9.5.

\begin{Num}\textsc{Example}
\label{Acyc}
Here is a very simple example. Suppose that the fibre $F$ is
$R$-acyclic, $\widetilde\bfH^\bullet(F;R)=0$.
Then all non-zero terms of the spectral sequence
are contained in the $X$-axis, the spectral sequence collapses, and
$\bfH^\bullet(E;R)\oot\bfH^\bullet(B;R)$ is an isomorphism. A similar
result holds if $B$ is $R$-acyclic; in this case
$\bfH^\bullet(F;R)\OOT\isom\bfH^\bullet(E;R)$.
\end{Num}

Note that \ref{LeraySerreSS}
says nothing about the differentials of the 
spectral sequence. Thus, it is almost impossible to determine
$\bfE_\infty$ from $\bfE_2$, except for very special situations.
An important additional ingredient is the multiplicative structure 
of the cohomology rings.

%\newpage

\section{Multiplicative structure}

\begin{Num}\textsc{Graded algebras}\psn
A \emph{graded $R$-algebra}
is a graded $R$-module $\bfM^\bullet$ with
an $R$-linear map
\[
\Tot(\bfM\otimes\bfM)\too\bfM
\]
of degree 0. Thus, $\bfM$ is an $R$-algebra, with the extra condition
that $\bfM^i\cdot\bfM^j\SUB\bfM^{i+j}$. If there is a unit element
$1$, then we require that $1\in\bfM^0$. The multiplication is
\emph{graded commutative} if 
\[
m_1\cdot m_2=(-1)^{\deg(m_1)\deg(m_2)}m_2\cdot m_1
\]
holds for all homogeneous elements.
If $\bfM$, $\bfN$ are graded algebras, then
$\Tot(\bfM\otimes\bfN)$ is again a graded algebra if we define
\[
(m_1\otimes n_1)\cdot(m_2\otimes n_2)=(-1)^{\deg(n_1)\deg(m_2)}
(m_1\cdot m_2)\otimes(n_1\cdot n_2).
\]
If $\bfM,\bfN$ are graded commutative, then so is their tensor product.
In the case of a filtration of a graded algebra $\bfM$ we require that
\[
\bfF^i\bfM^j\cdot\bfF^k\bfM^l\SUB\bfF^{i+k}\bfM^{j+l}.
\]
The corresponding notions in the bigraded case are very similar.
Here, we require that
\[
\bfM^{i,j}\cdot\bfM^{k,l}\SUB\bfM^{i+k,j+l};
\]
the multiplication
is \emph{bigraded commutative} if it is graded commutative with
respect to the total degree.
\end{Num}

\begin{Num}\textsc{Differential graded algebras}\psn
If $d$ is a differential on a graded algebra $\bfM$, and if $d$ 
satisfies the Leibniz rule
\[
d(m_1\cdot m_2)=d(m_1)\cdot m_2+(-1)^{\deg(m_1)}m_1\cdot d(m_2),
\]
then $(\bfM,d)$ is called a \emph{differential graded algebra}.
Similarly, a \emph{differential bigraded algebra} is a bigraded algebra
$\bfM$ with a differential of bidegree $(r,1-r)$, such that
$\Tot(\bfM)$ is a differential graded algebra, i.e.
\[
d(m_1\cdot m_2)=d(m_1)\cdot m_2+(-1)^{i+j}m_1\cdot d(m_2),
\]
where $\deg(m_1)=(i,j)$.
\end{Num}
Now we can sharpen the statement about the Leray-Serre spectral
sequence.
The $\bfE_2$-term $\bfE_2^\bubu\isom\bfH^\bullet(B;\bfH^\bullet(F;R))$
is a differential bigraded algebra. If $R$ is a
field, then there is an isomorphism of (bi)graded algebras
\[
\bfH^\bullet(B;\bfH^\bullet(F;R))\isom
\bfH^\bullet(B;R)\otimes\bfH^\bullet(F;R).
\]
Note that the multiplication on the right-hand side is
\[
(b_1\otimes f_1)\cdot(b_2\otimes f_2)=
(-1)^{\deg(f_1)\deg(b_2)}(b_1\smile b_2)\otimes(f_1\smile f_2).
\]
The cohomology of a differential graded algebra is again a
graded algebra. The differentials in the Leray-Serre spectral sequence
are compatible with this algebra structure,
i.e.~the $\bfE_r$  are differential bigraded algebras, and
$\bfE_\infty$ is a graded algebra. Moreover, the corresponding
filtration on $\bfH^\bullet(E;R)$ is compatible with the
cup-product.

The following examples illustrate the multiplicative properties;
they will also be useful later. 

\begin{Lem}
\label{LemKpin1}
Suppose that
$\bfH^\bullet(F;\QQ)\isom\EA_\QQ(u)$ is an exterior algebra
on one homogeneous generator $u$ of odd degree $n$. Suppose that
\begin{diagram}
F & \rTo & E \\
&& \dTo \\
&& B
\end{diagram}
is a $\QQ$-simple
fibration, and that $E$ is $\QQ$-acyclic. Then
$\bfH^\bullet(B;\QQ)\isom\QQ[a]$ is a polynomial algebra on one 
generator $a$ of degree $n+1$.

\begin{proof}
The only non-zero terms in $\bfE_2$ (and hence in $\bfE_r$, $2\leq r\leq
\infty$) are contained in the two
horizontal strips $\bfE_2^{\bullet,0}$ and $\bfE_2^{\bullet,n}$.
Thus, the only possibly non-zero differential is $d_{n+1}$, and
$\bfE_{n+2}\isom\bfE_\infty$. Therefore, the sequences
\[
0\too\bfE^{i,n}_{n+1}\too\bfE^{i+n+1,0}_{n+1}\too 0
\]
are exact for all $i\in\ZZ$. But 
\[
\bfE^{i,n}_{n+1}\isom\bfE^{i,n}_2\isom\bfH^i(B;\QQ)
\]
and
\[
\bfE^{i+n+1,0}_{n+1}\isom\bfE^{i+n+1,0}_2\isom\bfH^{i+n+1}(B;\QQ).
\]
This shows already that
\[
\bfH^i(B;\QQ)\isom\begin{cases}
\QQ & \text{ for }i\equiv 0\pmod{n+1} \\
0 & \text { else.}
\end{cases}
\]
Let $d_{n+1}(1\otimes u)=a\otimes 1$, for some 
$a\in\bfH^{n+1}(B;\QQ)$. We claim that $a^k$ spans
$\bfH^{k(n+1)}(B;\QQ)$, for all $k\geq 0$. This is true for
$k=1$, so we proceed by induction. Suppose that $a^k$ spans
$\bfH^{k(n+1)}(B;\QQ)$. Then $a^k\otimes u$ spans 
$\bfE_2^{k(n+1),n}\isom\bfE_{n+1}^{k(n+1),n}$. Now 
\[
d_{n+1}(a^k\otimes u)=d_{n+1}((a^k\otimes 1)\cdot(1\otimes u))=
(a^k\otimes 1)\cdot(a\otimes1)=a^{k+1}\otimes 1.
\]
This element spans $\bfE_{n+1}^{(k+1)(n+1),0}$, and the claim
follows.
\end{proof}
\end{Lem}
Suppose now that $\bfH^\bullet(F;\QQ)$ is a polynomial algebra
in an element $a$ of even degree. Here, the multiplicative 
structure is more important.

\begin{Lem}
\label{LemKpin2}
Suppose that
$\bfH^\bullet(F;\QQ)\isom\QQ[a]$ is a polynomial algebra
in a homogeneous generator $a$ of even degree $n$. Suppose that
\begin{diagram}
F & \rTo & E \\
&& \dTo \\
&& B
\end{diagram}
is a $\QQ$-simple fibration, and that $E$ is $\QQ$-acyclic. Then
$\bfH^\bullet(B;\QQ)\isom\EA_\QQ(u)$ is an exterior algebra on one 
generator $u$ of degree $n+1$.

\begin{proof}
The non-zero terms of $\bfE_2\isom\bfH^\bullet(B;\QQ)\otimes\QQ[a]$
are contained in the horizontal strips $\bfE_2^{\bullet,kn}$, for
$k\geq 0$. Thus, the differentials $d_2,\ldots d_n$ vanish
and $\bfE_{n+1}\isom\bfE_2$. Now $\bfE_2^{i,0}\isom\bfH^i(B;\QQ)=0$
for $1<i<n+1$, and the sequence
\[
0\too\bfE_{n+1}^{0,n}\TOO{d_{n+1}}\bfE_{n+1}^{n+1,0}\too0
\]
is exact. Put $d_{n+1}(1\otimes a)=u\otimes 1$. By the Leibniz rule,
$d_{n+1}(1\otimes a^k)=ku\otimes a^{k-1}$, whence
\[
d_{n+1}(u\otimes a^k)=0
\]
for $k\geq 1$.
The bigraded subalgebra $\EA_\QQ(u)\otimes\QQ[a]\SUB\bfE_{n+1}$ 
generated by
$a$ and $u$ is thus closed under the differential $d_{n+1}$, and
its cohomology is trivial, as is easily seen (using the divisibility
of $\QQ$). 
Suppose $\bfH^\bullet(B;\QQ)\neq\EA_\QQ(u)$. Then there is a
minimal number $i\geq 1$ such that
$\bfH^i(B;\QQ)\setminus\EA_\QQ(u)\neq 0$.
But now $\bfE_{n+2}^{j,k}=0$ for $0<j<i$ and $k\geq 0$.
Thus, all arrows that start or end at $\bfE_{n+2}^{i,0}
\isom\bfH^i(B;\QQ)$ are zero. This contradicts the fact that
$\bfE_\infty^{i,0}=0$.
\end{proof}
\end{Lem}
The two lemmata have an immediate application. Let $\pi$ be an
abelian group.
Recall that an Eilenberg-MacLane space of type $(\pi,n)$ is a space
$K(\pi,n)$ with the property that
\[
\pi_k(K(\pi,n))\isom\begin{cases}
\pi & \text{ for }k=n \\
0 & \text{ else.}
\end{cases}
\]
Let $PK(\pi,n)$ denote the path space. The map that sends a
path to its endpoint is a fibration whose fibre (over the
base point) is the loop space $\Omega K(\pi,n)$,
\begin{diagram}
\Omega K(\pi,n) & \rTo & PK(\pi,n) \\
&& \dTo \\
&& K(\pi,n)\rlap{.}
\end{diagram}
The path space is contractible and hence acyclic. The long exact 
homotopy sequence of this fibration shows then that $\Omega K(\pi,n)$ 
is a space of type $K(\pi,n-1)$. Now $\SS^1$ is clearly a space
of type $(\ZZ,1)$. Combining the two lemmata we see the following.

\begin{Prop}
\label{CoKpin}
Let $K(\ZZ,n)$ be a space of type $(\ZZ,n)$. If $n$ is
odd, then the rational cohomology 
\[
\bfH^\bullet(K(\ZZ,n);\QQ)\isom\EA(u)
\]
is an exterior algebra on one homogeneous generator of degree
$n$. If $n$ is even, then
\[
\bfH^\bullet(K(\ZZ,n);\QQ)\isom\QQ[a]
\]
is a polynomial algebra on one homogeneous generator $a$ of degree
$n$.

\begin{proof}
The 1-sphere $\SS^1$ is clearly a space of type $(\ZZ,1)$. 
The proof is a straight-forward induction based on the two lemmata
and the path-space fibration.
\end{proof}
\end{Prop}
If $\pi$ is a finite cyclic group, then a space of type
$(\pi,n)$ is $\QQ$-acyclic, cp.~Spanier \cite{Spa66} 9.5.6. 
This, combined with the result above, yields the following.

\begin{Prop}
\label{H(Kpin)}
Let $\pi$ be a finitely generated abelian group of rank $k$.
The rational cohomology of a space of type $(\pi,n)$ is a free 
graded commutative algebra on $k$ homogeneous generators of
degree $n$. Thus, if $n$ is odd, then
\[
\bfH^\bullet(K(\pi,n);\QQ)\isom\EA(u_1,\ldots,u_k),
\]
where $\deg(u_1)=\ldots=\deg(u_k)=n$,
and if $n$ is even then
\[
\bfH^\bullet(K(\pi,n);\QQ)\isom\QQ[a_1,\ldots,a_k],
\]
where $\deg(a_1)=\ldots=\deg(a_k)=n$.

\begin{proof}
Let $\pi=C_1\oplus\cdots\oplus C_r$ be a direct sum
of cyclic groups. Any two Eilenberg-MacLane spaces of type
$(\pi,n)$ are (weakly) homotopy equivalent; thus 
\[
K(\pi,n)\homot K(C_1,n)\times\cdots\times K(C_r,n),
\]
and the claim follows from the K\"unneth Theorem and \ref{CoKpin}.
\end{proof}
\end{Prop}

%\newpage

\section{Notes on collapsing}

We collect some criteria which ensure that the Leray-Serre
spectral sequence collapses.
\begin{Num}
\label{TrivialOnY-Axis}
Let
\begin{diagram}
F & \rTo^\iota & E \\
&&\dTo_p \\
&& B
\end{diagram}
be an $R$-simple fibration. Assume that either 
$\bfH^\bullet(F;R)$ or $\bfH^\bullet(B;R)$ is a free $R$-module.
Then the $\bfE_2$-term of the Leray-Serre spectral sequence is 
\[
\bfE_2\isom\bfH^\bullet(B;R)\otimes\bfH^\bullet(F;R).
\]
Suppose that the first differentials $d_2,d_3,\ldots,d_{k-1}$ vanish,
$\bfE_2\isom\bfE_3\isom\ldots\isom\bfE_k$. Because of the Leibniz
rule, the differential $d_k$ is completely determined by
its restriction to the $Y$-axis $\bfE^{0,\bullet}_k$:
\begin{align*}
d_k(b\otimes f)&=d_k((b\otimes1)\cdot(1\otimes f)) \\
&={\underbrace{d_k(b\otimes1)}_0}
\cdot(1\otimes f)+(-1)^{\deg(b)}(b\otimes1)\cdot d_k(1\otimes f)\\
&=(-1)^{\deg(b)}(b\otimes1)\cdot d_k(1\otimes f)
\end{align*}
In particular, if $d_k$ is trivial on $\bfE^{0,\bullet}_k$, then
$d_k=0$.
\end{Num}

\begin{Def}
If the map 
\[
\iota^\bullet:\bfH^\bullet(E;R)\too\bfH^\bullet(F;R)
\]
is surjective, then $F$ is called \emph{totally non-homologous to 0} in
$E$.
\end{Def}

\begin{Thm}[Leray-Hirsch]
\label{LHThm}
Suppose that 
\begin{diagram}
F & \rTo^\iota & E \\
&& \dTo_p \\
&& B
\end{diagram}
is an $R$-simple fibration, and that $F$ and $B$ are of finite type.
Suppose moreover that $\bfH^\bullet(F)$ is a finitely generated
free $R$-module. If $F$ is totally 
non-homologous to 0, then the Leray-Serre spectral sequence collapses,
$p^\bullet$ is an injection, and there is an additive isomorphism
\[
\bfH^\bullet(E)\isom\bfH^\bullet(B)\otimes\bfH^\bullet(F).
\]
The kernel of the map $\bfH^\bullet(F;R)\oot\bfH^\bullet(E;R)$ is
the ideal generated by 
\[
\mathrm{im}(p^\bullet)\cap\bigoplus_{i\geq 1}\bfH^i(E;R).
\]

\begin{proof}
See Spanier \cite{Spa66} 5.7.9 and Mimura-Toda \cite{MiTo91} III 
Thm.~4.2.
\end{proof}
\end{Thm}

%%%%%%%%%%%%%%%%%%%%%%%%%%%%%%%%%%%%%%%%%%%%%%%%%%%%%%%%%%%%%%%%%%%%%%%%
%                                                                      %
%                                                                      %
%                   Compact homogeneous quadrangles                    %
%                                                                      %
%                          Linus Kramer                                %
%                                                                      %
%                           Memoirs AMS                                %
%                                                                      %
%                         Wuerzburg 2000                               %
%                                                                      %
%                                                                      %
%                            CHQ2.tex                                  %
%                                                                      %
%                                                                      %
%                                                                      %
%%%%%%%%%%%%%%%%%%%%%%%%%%%%%%%%%%%%%%%%%%%%%%%%%%%%%%%%%%%%%%%%%%%%%%%%
\chapter{Ranks of homotopy groups}

Suppose that $X$ is a 1-connected space which has the same integral 
(co-)ho\-mo\-lo\-gy as a product of spheres. The aim of this chapter 
is to determine the ranks of the homotopy groups of $X$ in terms
of the rational cohomology of $X$. The idea is to
generalize a result of Cartan-Serre \cite{CarSer52} who determined the
rational homotopy groups of spaces whose rational cohomology ring
is freely generated. 
The note of Cartan and Serre appeared in the
Comptes Rendus, and they gave only a short hint how to prove
their theorem; a proof can be found in Fomenko-Fuchs-Gutenmacher
\cite{FFG86}. We prove here a more general version of the 
Cartan-Serre Theorem.

Everything in this chapter could also be done in the
framework of rational homotopy theory and minimal models,
see e.g.~F\'elix-Halperin-Thomas \cite{FHT}.
I chose a more traditional homotopy theoretic approach.

\section{The Whitehead tower}

Let $X$ be a 1-connected space, and let $\pi_k=\pi_k(X)$ be the
$k$th homotopy group of $X$. The idea is the following.
There is a sequence of spaces and fibrations
\[
\cdots
\too X\bra{k+1}\too X\bra{k}\too X\bra{k-1} \too\cdots\too X\bra{1}=X
\]
with the following properties.
\begin{description}
\item[\textbf{WT$_1$}] Each $X\bra{k}$ is $(k-1)$-connected.
\item[\textbf{WT$_2$}] The map $X\bra{k}\too X$ induces an
isomorphism on all homotopy groups
of degree at least $k$.
\item[\textbf{WT$_3$}] The homotopy fibre of the map
$X\bra{k+1}\too X\bra{k}$ is an
Eilenberg-MacLane space of type $(\pi_k,k-1)$.
\end{description}
This is sometimes called the \emph{Whitehead tower} or 
\emph{upside-down Postnikov tower} of $X$. It is rather
easy to compute the rank of $\pi_k$ in terms of the cohomology group
$\bfH^k(X\bra{k};\QQ)$:
\[
\rk(\pi_k)=\rk(\pi_k(X\bra{k}))=\dim_\QQ\bfH^k(X\bra{k};\QQ).
\]
Thus, if there is a way to calculate the rational cohomology of the
spaces $X\bra{k}$ from the knowledge of 
$\bfH^\bullet(X;\QQ)$, then one can determine
the ranks of the homotopy groups. 

The Whitehead tower can be constructed as follows.
Suppose we have already constructed $X\bra{k}$. Then
$X\bra{k}$ is $(k-1)$-connected. There exists an Eilenberg-MacLane
space $K(\pi_k,k)$ and a map $\theta:X\bra{k}\too K(\pi_k,k)$ 
that induces an isomorphism between $\pi_k(X\bra{k})$ and
$\pi_k(K(\pi_k,k))$ (one way to obtain this space is to kill all
higher-dimensional homotopy groups of $X\bra{k}$, cp.~e.g.
Whitehead \cite{Whi78} V.2.4).
Consider the path space fibration
\begin{diagram}
\Omega K(\pi_k,k) &  \rTo & & PK(\pi_k,k) \\
                &  &     & \dTo \\
                &  &     & K(\pi_k,k) .
\end{diagram}
Since the path space of any $0$-connected space is contractible,
$\Omega K(\pi,k)$ is an 
Eilenberg-MacLane space of type $(\pi,k-1)$.
We pull the path-space fibration $PK(\pi_k,k)\too K(\pi_k,k)$ back
via $\theta$ and denote the resulting total space by $X\bra{k+1}$,
\begin{diagram}
\Omega K(\pi_k,k) & & \rTo & & X\bra{k+1} & & \rTo & & PK(\pi_k,k) \\
             & & & &  \dTo & & & & \dTo \\
             & & &   & X\bra{k}   & & \rTo^\theta & & K(\pi_k,k).
\end{diagram}
The interesting part of the exact homotopy sequences of these
fibrations is
\begin{diagram}
& \rTo &\pi_k(K(\pi_k,k)) & & \rTo^\partial_\isom & & 
\pi_{k-1}(K(\pi_k,k-1)) & \rTo & \\
& & \dTo_\isom & & & & \dEq \\
& \rTo & \pi_k(X\bra{k}) & & \rTo^\partial & & \pi_{k-1}(K(\pi_k,k-1))
& \rTo & .
\end{diagram}
It shows that $\pi_k(X\bra{k+1})\isom 0\isom \pi_{k-1}(X\bra{k+1})$;
thus, $X\bra{k+1}$ has the claimed properties \textbf{WT$_1$, WT$_2$, and
WT$_3$}. 

Suppose that $\pi_k$ is \emph{finite}. 
Then $K(\pi_k,k-1)$ is $\QQ$-acyclic, and 
$\bfH^\bullet(X\bra{k+1};\QQ)\OOT\isom\bfH^\bullet(X\bra{k};\QQ)$
is an isomorphism, cp.~\ref{Acyc}. So, finite homotopy groups
in low degrees do not matter.

\begin{Lem}
\label{starteq}
Let $\pi_k$ be the first infinite homotopy group of $X$.
Then there is a chain of isomorphisms
\[
\bfH^\bullet(X\bra{k};\QQ)\OOT\isom
\bfH^\bullet(X\bra{k-1};\QQ)\OOT\isom
\bfH^\bullet(X\bra{k-2};\QQ)\OOT\isom \cdots \OOT\isom
\bfH^\bullet(X;\QQ)
.
\]
\qed
\end{Lem}
By Hurewicz' Theorem we have $\pi_k\isom\bfH_k(X\bra{k})$. Since
\[
\bfH^k(X\bra{k};\QQ)\isom\Hom(\bfH_k(X\bra{k}),\QQ)\isom\Hom(\pi_k,\QQ)
\isom \QQ^{\rk(\pi_k)},
\]
we have the following result.

\begin{Lem}
\label{startfin}
Let $\bfH^k(X;\QQ)$ be the first non-trivial rational cohomology group 
of $X$, for $k\geq 1$. Then 
\[
\rk(\pi_k)=\dim_\QQ\bfH^k(X;\QQ),
\]
and all groups $\pi_j$ are finite, for $0\leq j<k$
(note that we assumed that $\pi_0=\pi_1=0$).\qed
\end{Lem}
So one can tell the rank of the first infinite homotopy group
right away from the rational cohomology of $X$.
Suppose now that $\rk(\pi_k)=r$. The question is if
we can compute the cohomology of $X\bra{k+1}$. In the diagram
\begin{diagram}
\Omega K(\pi_k,k) & & \rTo & & X\bra{k+1} & & \rTo & & PK(\pi_k,k) \\
             & & & &  \dTo & & & & \dTo \\
             & & &   & X\bra{k}   & & \rTo^\theta & & K(\pi_k,k),
\end{diagram}
we consider first the fibration on the right
\begin{diagram}
\Omega K(\pi_k,k) &  \rTo & & PK(\pi_k,k) \\
                &  &     & \dTo \\
                &  &     & K(\pi_k,k)\rlap{.}
\end{diagram}
This is very similar to what we did in 
\ref{LemKpin1}, \ref{LemKpin2}.
Recall that $\Omega K(\pi_k,k)$ is a space of type $(\pi_k,k-1)$,
\[
\Omega K(\pi_k,k)\homot K(\pi_k,k-1).
\]
Consider the Leray-Serre spectral sequence of this
fibration. The total space $PK(\pi_k,k)$ is acyclic, hence
the $\bfE_\infty$-term is trivial. Since we are working over
the field $\QQ$, the $\bfE_2$-term is given by
\[
\bfE_2^{s,t}\isom\bfH^s(K(\pi,k);\QQ)\otimes\bfH^t(K(\pi,k-1);\QQ).
\]
Let $k\geq 2$ be \emph{even}. Then $\bfH^\bullet(K(\pi,k-1);\QQ)$
is an exterior algebra. Thus, the only non-trivial terms in $\bfE_2$ 
(and hence in $\bfE_n$, $2\leq n\leq \infty$) 
are contained in the
horizontal strips 
$\bfE_2^{\bullet,0}$, $\bfE_2^{\bullet,k-1},\ldots,
\bfE_2^{\bullet,r(k-1)}$.
Therefore, the first possibly non-trivial differential 
is $d_k$, and
\[
\bfE_2\isom \bfE_3\isom\cdots\isom \bfE_k.
\]
The $\bfE_k$-term looks as follows.
\begin{center}
\begin{picture}(200,160)(10,-20)
\put(0,40){\vector(0,1){40}}
\put(80,0){\vector(1,0){40}}
\put(-10,60){\makebox(0,0)[b]{$t$}}
\put(100,-10){\makebox(0,0)[b]{$s$}}
\put(10,10){\framebox(160,5)[b]{}}
\put(10,60){\framebox(160,5)[b]{}}
\put(10,110){\framebox(160,5)[b]{}}
\put(190,10){\makebox(0,0)[bl]{$t=0$}}
\put(190,60){\makebox(0,0)[bl]{$t=k-1$}}
\put(190,110){\makebox(0,0)[bl]{$t=2(k-1)$}}
\put(60,62){\vector(3,-2){75}}
\put(115,35){\makebox(0,0)[b]{$d_k$}}
\put(80,112){\vector(3,-2){75}}
\put(135,85){\makebox(0,0)[b]{$d_k$}}
\end{picture}
\end{center}
Since $\bfE_\infty^{0,k-1}=0$, the sequence
\[
0\too\bfE_n^{0,k-1}\TOO{d_k}\bfE_k^{k,0}\too0
\]
is exact. Let $w_1,\ldots,w_r$ be a basis for
$\bfH^{k-1}(K(\pi,k-1);\QQ)$. We define a $\QQ$-linear map $\tau$ by
\[
d_k(1\otimes w_i)=\tau(w_i)\otimes 1.
\]
The elements $\tau(w_1),\ldots,\tau(w_r)$ span $\bfH^k(K(\pi,k);\QQ)$,
therefore,
$\bfH^\bullet(K(\pi,k);\QQ)$ is a (free graded) polynomial algebra in these
elements,
\[
\bfH^\bullet(K(\pi,k);\QQ)=\QQ[b_1,\ldots,b_r]=
\QQ[\tau(w_1),\ldots,\tau(w_r)].
\]
We claim that $\bfE_{k+1}$ is trivial. Put
\[
\bfC_i^\bubu=\QQ[\tau(w_i)]\otimes\EA_\QQ(w_i).
\]
This is a bigraded subalgebra of $\bfE_k$ which is closed under
the differential $d_k$, and it is easy to check
that the cohomology of each $\bfC_i$ is trivial,
$\bfH(\bfC_i)=0$. But
\[
\mathrm{Tot}(\bfE_k)=\mathrm{Tot}(\bfC_1\otimes\cdots\otimes\bfC_r).
\]
By the K\"unneth Theorem, 
\[
\bfH(\bfE_k)\isom\bfH(\bfC_1)\otimes\cdots\otimes
\bfH(\bfC_r)
\]
is trivial, and thus
$\bfE_{k+1}=\bfE_{k+2}=\cdots=\bfE_\infty$.

If $k$ is \emph{odd}, then the non-trivial terms in $\bfE_2$ are
in the vertical strips 
$\bfE_2^{0,\bullet}$, 
$\bfE_2^{k,\bullet},\ldots,\bfE_2^{kr,\bullet}$. 
Again, 
\[
\bfE_2\isom \bfE_3\isom\cdots\isom \bfE_k 
\]
and the first possibly non-trivial differential is $d_k$.
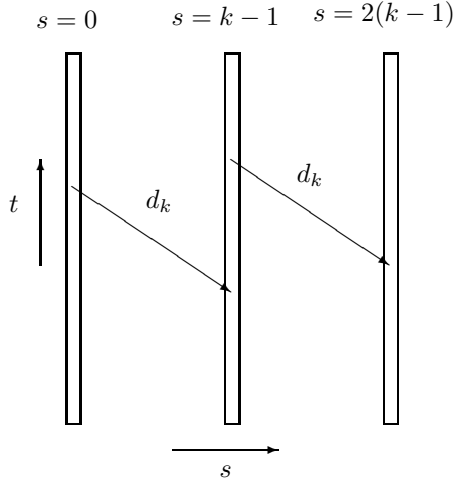
\begin{center}
\begin{picture}(100,200)(10,-20)
\put(0,70){\vector(0,1){40}}
\put(50,0){\vector(1,0){40}}
\put(-10,90){\makebox(0,0)[b]{$t$}}
\put(70,-10){\makebox(0,0)[b]{$s$}}
\put(10,10){\framebox(5,140)[b]{}}
\put(70,10){\framebox(5,140)[b]{}}
\put(130,10){\framebox(5,140)[b]{}}
\put(10,160){\makebox(0,0)[b]{$s=0$}}
\put(70,160){\makebox(0,0)[b]{$s=k-1$}}
\put(130,160){\makebox(0,0)[b]{$s=2(k-1)$}}
\put(12,100){\vector(3,-2){60}}
\put(45,90){\makebox(0,0)[b]{$d_k$}}
\put(72,110){\vector(3,-2){60}}
\put(102,100){\makebox(0,0)[b]{$d_k$}}
\end{picture}
\end{center}
As before the sequence
\[
0\too1\otimes\bfH^{k-1}(K(\pi,k-1);\QQ)\TOO{d_k}\bfH^k(K(\pi,k);\QQ)\otimes 1
\too0
\]
is exact, and we may put
\[
d_k(1\otimes b_i)=\tau(b_i)\otimes 1.
\]
It follows that
\[
\bfH^\bullet(K(\pi,k);\QQ)=\EA_\QQ(w_1,\ldots,w_r)=\EA_\QQ(\tau(b_1),\ldots,
\tau(b_r))
\]
is an exterior algebra generated by the $\tau(b_i)$.
By the same reasoning as above, $\bfE_{k+1}\isom\bfE_{k+2}\isom
\cdots\isom\bfE_\infty$ is trivial.

%\newpage

\section{The Cartan-Serre Theorem}

Now we get back to our real problem, the fibration
\begin{diagram}
\Omega K(\pi_k,k) & & \rTo & & X\bra{k+1} \\
             & & & &  \dTo \\
             & & &   & X\bra{k} & & \rTo^\theta && K(\pi_k,k)\rlap{.}
\end{diagram}
The map in cohomology $\theta^\bullet$ induces a map between the
two spectral sequences. The image of $\theta^\bullet$ is precisely
the subalgebra of $\bfH^\bullet(X\bra{k};\QQ)$ generated by 
$\bfH^k(X\bra{k};\QQ)$.

Suppose that $k$ is even. Again, the non-zero terms in $\bfE_2$
are contained in the horizontal strips $\bfE_2^{\bullet,j(k-1)}$,
for $0\leq j\leq r$, and the same is true for $\bfE_i$, $i\geq 2$.
Therefore, the first possibly non-trivial differential again
is $d_k$, and 
\[
\bfE_2\isom \bfE_3\isom\cdots\isom \bfE_k.
\]
Moreover, we know $d_k$ on the image of $\theta^\bullet$.
Let $\bfA^\bullet=\theta^\bullet(\bfH^\bullet(K(\pi_k,k);\QQ))$ 
denote this image.
Suppose that the following holds: the cohomology algebra 
$\bfH^\bullet(X\bra{k};\QQ)$
decomposes as
\[
\bfH^\bullet(X\bra{k};\QQ)\isom\mathrm{Tot}
(\bfA^\bullet\otimes\bfB^\bullet),
\]
for some subalgebra $\bfB^\bullet\SUB\bfH^\bullet(X\bra{k};\QQ)$.
This holds for example if $\bfH^\bullet(X\bra{k};\QQ)$ is a
free graded commutative algebra. Endow the complex
$\bfB^\bullet$ with the trivial (zero) differential.
Then there is an isomorphism of differential graded algebras
\[
\bfE_k^\bubu\isom\mathrm{Tot}(\bfA^\bullet\otimes\bfH^\bullet(K(\pi,k-1);\QQ))
\otimes\bfB^\bullet,
\]
and the K\"unneth Theorem yields
\[
\mathrm{Tot}(\bfE_{k+1})=\mathrm{Tot}(\bfH(\bfE_k^\bubu))=
\mathrm{Tot}(\bfH(\bfA^\bullet\otimes\bfH^\bullet(K(\pi,k-1);\QQ))
\otimes\bfB^\bullet).
\]
If the differential graded algebra 
$\bfA^\bullet\otimes\bfH^\bullet(K(\pi,k-1);\QQ)$ is
acyclic, then we obtain the classical Cartan-Serre Theorem.

\begin{Thm}[Cartan-Serre]
Suppose that $\bfH^\bullet(X;\QQ)$ is a free graded anti-commutative
algebra on homogeneous generators $a_1,\ldots,a_r$ of even degrees
and homogeneous generators $u_1,\ldots,u_s$ of odd degrees
\[
\bfH^\bullet(X;\QQ)\isom \QQ[a_1,\ldots,a_r]\otimes 
\EA_\QQ(u_1,\ldots,u_s).
\]
Let
\[
r_k=|\{i|\ \deg(a_i)=k\}|+|\{i|\ \deg(u_i)=k\}|
\]
denote the number of homogeneous generators of degree $k$.
Then 
\[
\rk(\pi_k(X))=r_k
\]
for all $k$. In particular, only finitely many homotopy groups of
$X$ are infinite, cp.~Cartan-Serre \cite{CarSer52} Prop. 3.

\begin{proof}
Let $k$ be the minimum of the degrees of the homogeneous generators.
By Lemma \ref{starteq} we have an isomorphism
$\bfH^\bullet(X\bra{k};\QQ)\oot\bfH^\bullet(X;\QQ)$, 
and for $1\leq i\leq k$ the groups $\pi_i$ are finite. Suppose that
$k$ is even and that $a_1,\ldots,a_m$ span $\bfH^k(X\bra{k};\QQ)$.
By Lemma \ref{startfin} $\rk(\pi_k)=m$.
We may put 
\[
\bfB^\bullet=\QQ[a_{m+1},\ldots,a_r]\otimes\EA_\QQ(u_1,\ldots,u_s).
\]
Then $\bfH^\bullet(K(\pi_k,k);\QQ)$ maps
isomorphically onto 
\[
\bfA^\bullet=\QQ[a_1,\ldots,a_m]\SUB\bfH^\bullet(X\bra{k};\QQ).
\]
In the spectral sequence we have 
\begin{align*}
\bfE_{k+1}^\bullet & \isom
\bfH([\bfA^\bullet\otimes\bfH^\bullet(K(\pi,k-1);\QQ)]\otimes\bfB^\bullet) \\
&\isom 
\bfH(\bfA^\bullet\otimes\bfH^\bullet(K(\pi,k-1);\QQ))
\otimes\bfH(\bfB^\bullet)) \\
&\isom
\bfH(\bfB^\bullet) \\
& \isom \bfB^\bullet,
\end{align*}
where all non-zero terms are contained in the strip
$\bfE^{\bullet,0}_{k+1}$. Thus, 
\[
\bfH^\bullet(X\bra{k+1};\QQ)\isom\QQ[a_{m+1},\ldots,a_r]\otimes
\EA_\QQ(u_1,\ldots,u_s)
\]
is isomorphic to $\bfH^\bullet(X\bra{k};\QQ)$ factored by the
ideal generated by $\bfH^k(X\bra{k};\QQ)$. The result follows now
by induction on the number of generators, starting with the case
where $X$ is $\QQ$-acyclic. The case when $k$ is odd is similar.
\end{proof}
\end{Thm}
This implies of course Serre's famous result: the only infinite homotopy 
group of $\SS^{2n+1}$ is $\pi_{2n+1}$. However, we need also to
consider spaces which have the same cohomology as a product of an
odd- and and an even-dimensional sphere. This is not covered by the
Cartan-Serre theorem.
Thus, suppose that the rational cohomology of $X$ is of the form
\[
\bfH^\bullet(X;\QQ)\isom\QQ[a]/(a^m)\otimes\EA_\QQ(u_1,\ldots,u_r).
\]
By the reduction process of the Whitehead tower,
there is no loss of generality in assuming that $k=\deg(a)<\deg(u_i)$
for $i=1,\ldots,r$. Then we may put 
\[
\bfA^\bullet=\QQ[a]/(a^m)
\]
and
\[
\bfB^\bullet=\EA_\QQ(u_1,\ldots,u_r).
\]
Let $d_k(1\otimes v)=a\otimes 1$. 
The cohomology of the differential bigraded algebra
\[
\QQ[a]/(a^m)\otimes\EA_\QQ(v)\SUB\bfE_k
\]
is $\EA_\QQ(w)$, for an element $w$ of bidegree $(k(m-1),k-1)$
(and total degree $km-1$).
Therefore $\bfE_{k+1}\isom\EA_\QQ(w,u_1,\ldots,u_r)$.
At this stage the spectral sequence collapses, and thus
\[
\bfE_\infty\isom\bfH^\bullet(X\bra{k+1};\QQ)\isom
\EA_\QQ(w,u_1,\ldots,u_r).
\]

\begin{Thm}
\label{PiEvenOdd}
Let $X$ be a 1-connected space whose rational cohomology is of the form
\[
\bfH^\bullet(X;\QQ)\isom\QQ[a]/(a^m)\otimes\EA_\QQ(u_1,\ldots,u_r),
\]
where $\deg(a)$ is even and the degrees of the $u_i$ are odd.
Let $r_i=|\{j|\ \deg(u_j)=i\}|$.
Then
\[
\rk(\pi_k)=
\begin{cases}
1 & \text{ for }k=\deg(a) \\
r_k & \text{ for }k\neq\deg(a)\text{ and }k\neq m\cdot\deg(a)-1 \\
r_k+1 & \text{ for }k=m\cdot\deg(a)-1.
\end{cases}
\]
\qed
\end{Thm}
In particular, if $X$ has the same cohomology as $\SS^{2n}$,
i.e.~if $\bfH^\bullet(X;\QQ)\isom\QQ[a]/(a^2)$,
then the only infinite homotopy groups are $\pi_{2n}$ and
$\pi_{4n-1}$.
This is the even-dimensional version of Serre's finiteness result
for homotopy groups of spheres.

%%%%%%%%%%%%%%%%%%%%%%%%%%%%%%%%%%%%%%%%%%%%%%%%%%%%%%%%%%%%%%%%%%%%%%%%
%                                                                      %
%                                                                      %
%                   Compact homogeneous quadrangles                    %
%                                                                      %
%                          Linus Kramer                                %
%                                                                      %
%                           Memoirs AMS                                %
%                                                                      %
%                         Wuerzburg 2000                               %
%                                                                      %
%                                                                      %
%                            CHQ3.tex                                  %
%                                                                      %
%                                                                      %
%                                                                      %
%%%%%%%%%%%%%%%%%%%%%%%%%%%%%%%%%%%%%%%%%%%%%%%%%%%%%%%%%%%%%%%%%%%%%%%%
\chapter{Some homogeneous spaces}

In this chapter we analyze the algebraic topology of
certain homogeneous spaces, using the homotopy theoretic
results of the previous chapter. The spaces we are
interested in have the same cohomology as a product
of two spheres. Although some of the theorems are
tailored for this special situation, it should be emphasized
that the methods developed here can be adapted to a much
larger class of homogeneous spaces. 

The first section contains various results about
the structure and the topology of compact Lie groups.
In the second section we determine the rational Leray-Serre
spectral sequence of the principal bundle
\begin{diagram}
H & \rTo & G & \rTo & G/H
\end{diagram}
where $G/H$ has the same cohomology as $\SS^{n_1}\times\SS^{n_2}$,
where $2\leq n_1\leq n_2$ and $n_2$ is odd.

If $n_1$ is odd, then the spectral sequence collapses. This
result is also contained in Onishchik \cite{Oni94}. However,
Onishchik uses real cohomology of Lie algebras instead of
fibre-bundle techniques and spectral sequences, hence his proof 
is quite different. The result can also be found in a
paper of Hsiang-Su \cite{HsSu68}
about transitive action on Stiefel manifolds. 
However, Hsiang-Su were obviously not aware of the
Cartan-Serre Theorem, because they included a superfluous condition
on the rational homotopy groups.

We obtain a relation between degrees of the primitive elements 
in the cohomology of $G$ and $H$ and the numbers $n_1$ and $n_2$.
The second important result is that $G$ has a semisimple normal transitive 
subgroup with at most two almost simple factors. 

In the last section we state the classification of all
compact 1-connected homogeneous spaces which have the same integral
cohomology as a product $\SS^{n_1}\times\SS^{n_2}$, for
$3\leq n_1\leq n_2$,
$n_2$ odd. The actual classification is carried out in Chapters 5 and 6,
using the results of this chapter and the representation theory developed
in the next chapter. This classification is the first
main result of this book. However, as I mentioned before, the reader
who has a different classification problem in mind should have no problems
to adapt the techniques and results of Chapter 3 and 4 for his
own purposes.

%\newpage

\section{Structure of compact Lie groups}

We call a compact connected Lie group \emph{almost simple} if it
has no normal closed subgroup of positive dimension. Such a
group is connected with finite center, and its
Lie algebra is simple.
The compact simple Lie algebras are the series
\begin{align*}
\fa_n&=\su_{n+1},\ n\geq 1, \\
\fb_n&=\so_{2n+1},\ n\geq 2, \\
\fc_n&=\sp_n,\ n\geq 3, \\
\fd_n&=\so_{2n},\ n\geq 4, 
\end{align*}
and the five exceptional compact Lie algebras 
\[
\fg_2,\ \ff_4,\ \fe_6,\ \fe_7, \text{ and }
\fe_8.
\]
There are isomorphisms 
\begin{gather*}
\su_2\isom\so_3\isom\sp_1 \\
\so_4\isom\su_2\oplus\su_2 \\
\so_5\isom\sp_2 \\
\so_6\isom\su_4.
\end{gather*}
The corresponding classical groups are
the complex unitary groups 
\[
\SU(n+1)=\SU_{n+1}\CC,
\]
the orthogonal groups in odd dimensions 
\[
\SO(2n+1)=\SO_{2n+1}\RR,
\] 
the quaternion unitary groups
\[
\Sp(n)=\U_n\HH,
\]
and the orthogonal groups in even dimensions 
\[
\SO(2n)=\SO_{2n}\RR.
\]
The groups $\SU(n)$ and $\Sp(n)$ are simply connected;
the universal coverings of the orthogonal groups $\SO(n)$,
$n\geq 3$, are the groups $\Spin(n)$.
A compact Lie group of type $\fg_2$, $\ff_4$ or $\fe_8$ is 
automatically simple and simply connected, and we denote these
unique groups by $\G_2$, $\Ffour$, and $\E_8$, respectively. 
The simply connected groups of type
$\fe_6,\fe_7$ have centers $\ZZ/3$ and $\ZZ/2$, respectively.
We denote these groups by
$\E_6$ and $\E_7$. 

A prefix P denotes the corresponding
simple group, e.g. $\PSO(6)$ or $\mathrm{PE}_6$.
These groups and the compact 1-torus 
\[
\TT=\SO(2)
\]
are the building
blocks of all compact connected Lie groups.

\begin{Thm}[Structure of compact connected Lie groups]\ \psn
Let $G$ be a compact connected Lie group. Then there exist
simply connected almost simple compact Lie groups
$G_1,\ldots,G_r$ and a torus $\TT^s$ and a surjection
\[
f:G_1\times\cdots\times G_r\times\TT^s\too G
\]
with finite kernel.

\begin{proof}
See e.g.~Hofmann-Morris \cite{HoMo} Ch.~6, Thm.~6.19, or
Salzmann \emph{et al.} \cite{CPP95} 93.11.
\end{proof}
\end{Thm}
A Lie group $G$ is an H-space. Therefore the rational cohomology
of $G$ is a finite dimensional associative Hopf algebra,
\[
\bfH^\bullet(G;\QQ)\isom\EA_\QQ(u_1,\ldots,u_r),
\]
cp.~Spanier \cite{Spa66} 5.8.13, generated by primitive elements
$\{u_1,\ldots,u_r\}$, i.e.~with comultiplication
\[
u_i\mapstoo u_i\otimes1+1\otimes u_i,\quad\text{ for }i=1,\ldots,r,
\]
cp.~e.g. Whitehead \cite{Whi78} III.8.
Since $G$ is an H-space, the fundamental group $\pi_1(G)$ acts trivially
on the homotopy groups $\pi_r(G)$, for $r\geq 1$, 
cp.~Whitehead \cite{Whi78} III.4.18. To determine the rational cohomology
of $G$, it suffices to consider the
case where $s=0$, since $G$ is topologically a product of a $\TT^s$
and a semisimple group. So assume that $s=0$. The
map $f$ described above induces then an isomorphism
\[
\bfH^\bullet(G_1\times\cdots\times G_r;\QQ)
\oot\bfH^\bullet(G;\QQ).
\]
(To see the last implication, consider the rational Leray-Serre
spectral sequence of the $\QQ$-simple fibration 
\[
G_1\times \cdots\times G_r\too G\too K(\pi_1,1).
\]
The base 
$K(\pi_1,1)$ is $\QQ$-acyclic, since $\pi_1=\pi_1(G)$ is finite
--- we assumed that $s=0$.)
Thus, we see that the Cartan-Serre theorem applies in particular
to compact connected Lie groups; in particular, we have
an isomorphism
\[
\bfH^\bullet(G_1\times\cdots\times G_r\times\TT^s;\QQ)
\lTo^{\cong}\bfH^\bullet(G;\QQ)
\]
(now $s$ is again arbitrary).

\begin{Num}
Put $m_i=\deg(u_i)$. These numbers are known for all compact almost
simple Lie groups; they are as follows.
\begin{center}
\label{TypeTable}
\begin{tabular}{ll}
Type & $(m_1,\ldots,m_r)$ \\ \hline
$\fa_n$ & $(3,5,7,\ldots,2n+1)$ \\
$\fb_n$ & $(3,7,11,\ldots,4n-1)$ \\
$\fc_n$ & $(3,7,11,\ldots,4n-1)$ \\
$\fd_n$ & $(3,7,11,\ldots,4n-5,2n-1)$ \\
$\fe_6$ & $(3,9,11,15,17,23)$ \\
$\fe_7$ & $(3,11,15,19,23,27,35)$ \\
$\fe_8$ & $(3,15,23,27,35,39,47,59)$ \\
$\ff_4$ & $(3,11,15,23)$ \\
$\fg_2$ & $(3,11)$ 
\end{tabular}
\end{center}
cp.~Mimura-Toda \cite{MiTo91} III.6.5, VI.5.10 or Mimura
\cite{Mi95} Thm.~2.2.
In particular, we see the following: If $G$ is a compact connected
Lie group, then the torus factor has dimension 
$\dim_\QQ\bfH^1(G;\QQ)$, and if $G$ has no torus factors, then
the number of almost simple
factors is $\dim_\QQ\bfH^3(G;\QQ)$. In fact, it is known
that $\pi_3(G)\isom\ZZ$ for all almost simple compact Lie groups,
cp.~Mimura \cite{Mi95} Thm.~3.9.
Let $P_G\SUB\bfH^\bullet(G;\QQ)$ denote the (graded) vector
space spanned by the primitive
homogeneous generators of the cohomology.
We recall the following facts.

\begin{itemize}
\item The dimension of $P_G$ is the rank of $G$.
\item The dimension of $P_G^1$ is the rank of the central torus of
$G$.
\item The dimension of $P_G^3$ is the number of almost simple
factors of $G$.
\item If $G$ is almost simple and of rank at least 2, then 
$G$ is exceptional if and only if $P_G^5=P_G^7=0$ (if the rank of
$G$ is at least 3 then it suffices that $P_G^7=0$).
\end{itemize}

These facts follow essentially from the structure of the rational
cohomology of compact almost simple Lie groups given by the
table above.
\end{Num}

\begin{Thm}[Borel]
\label{BorelHom1}
Let $G$ be a compact connected Lie group, and let $i:H\SUB G$ be a
closed connected subgroup.
Then $i^\bullet(P_G)$ is contained in $P_H$; we denote the
restriction of $i^\bullet$ to $P_G$ by 
\[
P_i:P_G\too P_H.
\]
The image $i^\bullet(\bfH^\bullet(G;\QQ))$ is an exterior algebra
(and a Hopf algebra) generated by $i^\bullet(P_G)$.
The kernel of $i^\bullet$ is the ideal generated by the homogeneous
elements contained in it.
\begin{proof}
See Borel \cite{BorelThesis} \S21.
\end{proof}
\end{Thm}
Late, we will need the following result of Borel about compact
transformation groups.

\begin{Thm}[Borel]
\label{BorelHom2}
Suppose that $H\SUB G$ is a closed connected subgroup
of a compact Lie group $G$. Consider the principal bundle
\begin{diagram}
H & \rTo & G \\
&& \dTo_p \\
&& G/H
\end{diagram}
The image $p^\bullet(\bfH^\bullet(G/H;\QQ))$ is an exterior
algebra (and a Hopf algebra) generated by
$P_G\cap p^\bullet(\bfH^\bullet(G/H;\QQ))$.
\begin{proof}
See Borel \cite{BorelThesis} \S21.
\end{proof}
\end{Thm}

%\newpage

\section{Certain homogeneous spaces}

\begin{Num}\textsc{Irreducible actions}\psn
\label{DefIrreducibleAction}
We begin with some general observations about homogeneous spaces.
Let $X$ be a homogeneous space of a compact connected
Lie group $K$. The action of 
$K$ on $X$ is called \emph{irreducible} if $K$ has no proper
normal transitive subgroup. There exists always a 
normal connected transitive  subgroup $G\SUB K$ such that the action of 
$G$ on $X$ is irreducible. 
\end{Num}
Let $L$ be a normal complement of $G$, i.e.~$K=G\cdot L$.
Then 
\[
L\SUB\Cen_{\mathrm{Sym}(X)}(G)=C
\]
is contained in the group $C$ 
of all permutations of $X$ which centralize $G$. This group $C$ can
be recovered from $(G,X)$ as follows. Let $H=G_x$ be the stabilizer
of an element $x\in X$ and put $X=G/H$. We define an action of the
normalizer $N=\mathrm{Nor}_{G}(H)$ on $X=G/H$ by 
\[
n\cdot (gH)=gn^{-1}H;
\]
this action centralizes the action of $G$ on $X$, and 
the kernel of this action is $H$. Thus we have an injection
$N/H\rInto C$. We claim that this injection is an isomorphism
\[
C\isom N/H.
\]
Let $c\in C$ and note that $H=G_x$ fixes $c(x)$, because $C$ centralizes
the action of $H$. Choose $n\in G$ such that $n^{-1}(x)=c(x)$.
For $h\in H$ we have $hn^{-1}(x)=n^{-1}(x)$, so $n\in N$.
Now let $y\in X$ be an arbitrary element. Choose $g\in G$ such that
$g(x)=y$. Then $c(y)=cg(x)=gc(x)=gn^{-1}(x)$, and
this is precisely the $N$-action on $X=G/H$
described above. Thus $C\isom N/H$. Note that $N/H$ is locally
isomorphic to $\Cen_G(H)/\Cen(H)$.
Therefore we have the following result.

\begin{Prop}
Let $X$ be a compact space. To classify all compact connected groups
which act transitively and effectively on $X$, it suffices to classify 
all pairs $(G,H)$, where $G$ is compact and connected, and where
$G$ acts irreducibly on $X=G/H$, and then to determine the 
$G$-normalizer (or centralizer) of the isotropy group $H$.
\qed
\end{Prop}
Now we consider a compact 1-connected homogeneous space $G/H$ 
which has the same integral
cohomology as a product of spheres $\SS^{n_1}\times\SS^{n_2}$,
where $n_1\leq n_2$, and $n_2$ is odd. Thus, the cohomology of
$G/H$ is
\[
\bfH^\bullet(G/H)\isom\EA_\ZZ(u,v) 
\quad\quad\deg(u)=n_1,\ \deg(v)=n_2
\]
if $n_1$ is odd, and
\[
\bfH^\bullet(G/H)\isom\ZZ[a]/(a^2)\otimes\EA_\ZZ(u) 
\quad\quad\deg(a)=n_1,\ \deg(u)=n_2
\]
if $n_1$ is even.
We may assume that $G$ is connected, because $G/H$ is connected, and
the evaluation map is open.

\begin{Thm}[Onishchik]
\label{Case(I)General}
Let $G/H$ be a 1-connected homogeneous space of a
compact connected Lie group $G$. Suppose that the rational 
cohomology of $G/H$ is an exterior algebra
\[
\bfH^\bullet(G/H;\QQ)\isom\EA_\QQ(w_1,\ldots,w_r)
\]
in homogeneous generators of odd degrees.
Let $P_{G/H}\SUB\bfH^\bullet(G/H;\QQ)$
denote the graded vector space spanned by the
homogeneous generators $w_1,\ldots,w_r$.
Then the spectral sequence of the fibration $G\too G/H$ collapses,
and there is an exact sequence
\[
0\oot P_H\oot P_G\oot P_{G/H}\oot0.
\]
In particular,
\[
\rk(G)-\rk(H)=r,
\]
and we have the following diagram:
\begin{diagram}
\bfH^\bullet(H;\QQ) & \lTo & \bfH^\bullet(G;\QQ) & \lTo &
\bfH^\bullet(G/H;\QQ) \\
\dTo_\isom && \dTo_\isom && \dTo_\isom \\
\EA_\QQ(u_1,\ldots,u_s) & \lTo &
\EA_\QQ(u_1,\ldots,u_s,w_1,\ldots,w_r) & \lTo &
\EA_\QQ(w_1,\ldots,w_r)\rlap{.}
\end{diagram}

\begin{proof}
Consider the long exact homotopy sequence of the fibration
$H\too G\too G/H$, tensored with $\QQ$. 
All even-dimensional homotopy groups of $G$, $H$ and $G/H$ are finite.
The field $\QQ$ is a flat $\ZZ$-module, hence
there are short exact sequences
\[
0\too\pi_k(H)\otimes \QQ\too
\pi_k(G)\otimes \QQ\too
\pi_k(G/H)\otimes \QQ \too 0
\]
for all $k$ (the groups are trivial if $k$ is even).
This proves that $\rk(G)=\rk(H)+r$.

The rational spectral sequence of the fibration $H\too G\too G/H$ has
\linebreak
$\bfH^\bullet(G/H;\QQ)\otimes\bfH^\bullet(H;\QQ)$ as its $\bfE_2$-term.
Therefore 
\[
\dim\bfE_2=\dim\bfH^\bullet(G/H;\QQ)\cdot
\dim\bfH^\bullet(H;\QQ)=2^r\cdot 2^{\rk(H)}.
\]
But 
\[
\dim\bfE_\infty=\dim\bfH^\bullet(G;\QQ)=
2^{\rk(G)}=2^{r+\rk(H)}=\dim\bfE_2.
\]
Thus, the spectral sequence collapses. This implies that
the map $H\too G$ induces a surjection in cohomology, and that
$G\too G/H$ induces an injection. The claim follows from 
\ref{BorelHom1} and \ref{BorelHom2} and the Leray-Hirsch
Theorem \ref{LHThm}.
\end{proof}
\end{Thm}
The result above follows also from Theorem 1 on p.~216 in Onishchik 
\cite{Oni94} (for cohomology with real coefficients).

\begin{Num}\textsc{How to use this}\psn
\label{HowTo(I)}
Let $G/H$ be a homogeneous space as in the theorem above.
Let 
\[
(m_1^G,\ldots,m_{r+s}^G)
\]
denote the degrees of the homogeneous generators
of the rational cohomology of $G$, and let
\[
(m_1^H,\ldots,m_s^H)
\]
denote the corresponding numbers for $H$. 
Then, after a suitable permutation
of the indices, $m_i^G=m_i^H$ for $i=1,\ldots,s$, so
$(m_1^H,\ldots,m_s^H)$ has to be a subsequence of 
$(m_1^G,\ldots,m_{r+s}^G)$, and the sequence which remains after deleting
all entries of the subsequence gives the degrees of the homogeneous
generators of the cohomology of $G/H$. This is essentially the method
that we are going to use in the case that $G$ is an almost simple
group. For example, if $G=\SU(5)$ and $H=\SU(3)$, then we obtain
\[
\begin{array}{ll}
m^G: & (\not\!3,\not\!5,7,9) \\
m^H: & (3,5) \\ 
\hline
 & (7,9)
\end{array}
\]
If $G$ is not almost simple, we need an additional reduction method,
i.e.~a bound on the number of almost simple factors (see 
\ref{LengthAtMost2} below).
\end{Num}
Now we assume that the cohomology of $X=G/H$ is
\[
\bfH^\bullet(X;\QQ)\isom \QQ[a]/(a^2)\otimes\bigwedge(u)
\]
where $\deg(a)=n_1$ is even, $\deg(u)=n_2$ is odd, and $n_1<n_2$.
There are short exact sequences
\[
0\too\pi_k(H)\otimes\QQ\too\pi_k(G)\otimes\QQ\too\pi_k(G/H)\otimes\QQ\too 0
\]
for $k\neq n_1-1$, and an exact sequence
\[
0\too\pi_{n_1}(G/H)\otimes\QQ\too\pi_{n_1-1}(H)\otimes\QQ
\too\pi_{n_1-1}(G)\otimes \QQ\too0,
\]
since $\pi_{n_1-1}(G/H)$ is finite.
This shows that $\rk(G)-\rk(H)=1$. The rational
Leray-Serre spectral sequence
of the fibration $H\too G\too G/H$ is slightly more complicated
because there are non-trivial differentials. We have
\[
\dim\bfE_2=\dim\bfH^\bullet(G/H;\QQ)\cdot\dim\bfH^\bullet(H;\QQ)
=4\cdot 2^{\rk(H)},
\]
and 
\[
\dim\bfE_\infty=\dim\bfH(G;\QQ)=2^{\rk(G)}=2\cdot2^{\rk(H)}.
\]
This shows already that the spectral sequence does not collapse.
The terms $\bfE_2^{s,t}$ in the strip $0<s<n_1$ are zero, because
$\bfH^s(G/H;\QQ)=0$ for these values of $s$. Therefore
$d_k$ is trivial on $\bfE_k^{0,\bullet}$ for $k<n_1$ and
this implies by \ref{TrivialOnY-Axis}
that $d_k$ is trivial for $k=2,\ldots,n_1-1$
on all terms of $\bfE_k$. 
Recall that the following diagram commutes.
\begin{diagram}
\bfH^\bullet(G/H;\QQ) & \rTo  &&&&& \bfH^\bullet(G;\QQ) \\
\dTo^\isom & &&&&& \dEq \\
\bfH^\bullet(G/H;\QQ)\otimes1 & \rEq & \bfE_2^{\bullet,0} & \rTo &
\bfE_\infty^{\bullet,0} & \rInto & \bfH^\bullet(G;\QQ)
\end{diagram}
Now we use the fact that the image of the cohomology of
$G/H$ is generated by the primitive elements contained in it.
A primitive element has odd degree. The possibly non-zero terms of odd 
degree in $\bfE_\infty^{\bullet,0}$ are 
$\bfH^{n_2}(G/H;\QQ)$ and $\bfH^{n_1+n_2}(G/H;\QQ)$. 
They cannot generate an element
of degree $n_1$, hence $\bfE_\infty^{n_1,0}=0$. 
In the spectral sequence,
the only possibly non-zero arrow which ends at $\bfE_k^{n_1,0}$ is
\[
d_{n_1}:\bfE_{n_1}^{0,n_1-1}\too\bfE_{n_1}^{n_1,0}.
\]
Thus this map is surjective.
Pick an element $w_1\in\bfH^{n_1-1}(H;\QQ)$ with
$d_{n_1}(1\otimes w_1)=a\otimes 1$. All primitive elements in 
$\bfH^\bullet(H;\QQ)$ of degree less than $n_1-1$ are mapped to
zero. Therefore, $w_1$ is not a sum of products of primitive
elements of lower degree, i.e.~$w_1$ itself is a primitive
element. Put
\[
\bfH^\bullet(H;\QQ)=\EA_\QQ(w_1,w_2,\ldots w_r),
\]
where $\rk(H)=r$. We claim that the elements $w_2,\ldots,w_r$ can be chosen
in such a way that the following holds.
\begin{itemize}
\item[(1)] $d_{n_1}(1\otimes w_i)=0$ for $i=2,\ldots,r$
\item[(2)] The spectral sequence collapses at 
$n_1+1$, i.e.~$\bfE_{n_1+1}=\bfE_\infty$.
\end{itemize}
For the second claim, it suffices to show that $\dim\bfE_{n_1+1}
\leq 2^{r+1}$, since we know already that
$\dim\bfE_{n_1+1}\geq \dim\bfE_\infty=2^{r+1}$, and the equality
$\dim\bfE_{n_1+1}=\dim\bfE_\infty$
forces that $d_k=0$ for all $k\geq n_1+1$. Since $\bfE_{n_1+1}$ 
is the cohomology of $\bfE_{n_1}$, we have
\[
\dim\bfE_{n_1+1}=\dim\bfE_{n_1}-2\dim(\mathrm{im}(d_{n_1})).
\]
So we need to show that $\dim(\mathrm{im}(d_{n_1}))\geq 2^r$ because then
$\dim\bfE_{n_1+1}\leq 4\cdot 2^r-2\cdot 2^r=2^{r+1}$.
Define a linear map $\phi$ by $d_{n_1}(1\otimes z)=a\otimes\phi(z)$.
Then 
\begin{align*}
d_{n_1}(1\otimes w_1z) & =d_{n_1}(1\otimes w_1)(1\otimes z)-
(1\otimes w_1)d_{n_1}(1\otimes z) \\
& =(a\otimes 1)(1\otimes z)-
(1\otimes w_1)(a\otimes\phi(z)) \\
& =a\otimes(z- w_1\phi(z)).
\end{align*}
Now $w_1(z- w_1\phi(z))=w_1z\neq 0$, provided that
$z\in\EA_\QQ(w_2,\ldots,w_r)$. Thus, $d_{n_1}$ is injective
on the $2^{r-1}$-dimensional space 
$1\otimes \left(w_1\EA_\QQ(w_2,\ldots,w_r)\right)$.
For $u\otimes z\in\bfE_2^{n_2,\bullet}$ we have similarly
\begin{align*}
d_{n_1}(u\otimes z)&
=d_{n_1}(u\otimes 1)(1\otimes z)-(u\otimes 1)d_{n_1}(1\otimes z) \\
&=-(u\otimes 1)(a\otimes \phi(z)) \\
&=-(ua)\otimes\phi(z).
\end{align*}
Therefore, $d_{n_1}$ is also injective on the $2^{r-1}$-dimensional
subspace 
\[
\textstyle
u\otimes\left(w_1\EA_\QQ(w_2,\ldots,w_r)\right).
\]
The images of $1\otimes \left(w_1\EA_\QQ(w_2,\ldots,w_r)\right)$
and $u\otimes \left(w_1\EA_\QQ(w_2,\ldots,w_r)\right)$ have
trivial intersection (because of the grading), therefore
the image of $d_{n_1}$ has dimension 
\[
\dim(\mathrm{im}(d_{n_1}))=2^r.
\]
We have proved claim (2). Moreover, we know the cohomology of the
map $G\TOO p G/H$, and we use this to prove claim (1).

\begin{Lem}
The map $\bfH^\bullet(G;\QQ)\oot\bfH^\bullet(G/H;\QQ)$ is given by the
diagram
\begin{diagram}
\bfH^\bullet(G;\QQ) && \lTo^{p^\bullet} && \bfH^\bullet(G/H;\QQ) \\
\dEq & && & \dEq \\
\EA_\QQ(u,v_1,\ldots,v_r) && \lTo && \QQ[a]/(a^2)\otimes\EA_\QQ(u) ,
\end{diagram}
where $p^\bullet(u)=u$ and $p^\bullet(a)=0$.
\qed
\end{Lem}
We still have not yet determined the differential $d_{n_1}$. 
We apply \ref{BorelHom1} to the inclusion $h:H\SUB G$.
Recall that the following diagram commutes.
\begin{diagram}
\bfH^\bullet(G;\QQ) &&\rTo^{h^\bullet} && \bfH^\bullet(H;\QQ) \\
\dTo &&&& \dEq \\
\bfE_\infty^{0,\bullet} & \rInto & \bfE_2^{0,\bullet} & \rTo^{\isom\ } &
\bfH^\bullet(H;\QQ) 
\end{diagram}
Thus, we can identify the image of $h^\bullet$ with
$\bfE_\infty^{0,\bullet}\SUB\bfE_2^{0,\bullet}$. We know that
$\dim\bfE_\infty^{0,\bullet}=2^{r-1}$. Therefore
$\dim(\mathrm{im}(h^\bullet)\cap P_H)=r-1$, and we can
find primitive elements $w_2,\ldots w_r\in P_H$ which generate 
$\mathrm{im}(h^\bullet)\isom\bfE_\infty^{0,\bullet}$. This implies 
that $d_{n_1}(1\otimes w_i)=0$ for $i=2,\ldots,r$.
\begin{Lem}
There exist primitive elements  $w_1,\ldots w_r$ which generate
\linebreak
$\bfH^\bullet(H;\QQ)$ such that 
\[
d_{n_1}(1\otimes w_1)=a\otimes 1,\quad\text{ and }\quad
d_{n_1}(1\otimes w_i)=0\text{ for }i=2,\ldots,r.
\]
\qed
\end{Lem}
Consider the map $P_G\too P_H$. Its image has dimension $(r-1)$, thus
its kernel has dimension 2.
We may choose an $(r-1)$-dimensional subspace of $P_G$ which
maps isomorphically onto $h^\bullet(P_G)$.
The image of $u$ is contained in $\bfE_\infty^{n_1,0}$ and thus
in the kernel of $h^\bullet$. We choose one more primitive
element $v$ such that $u,v$ span the kernel of $P_G\too P_H$.
Note that $\deg(v)=2n_1-1$ by \ref{PiEvenOdd}.

\begin{Thm}
\label{Case(II)General}
Let $X=G/H$ be a 1-connected homogeneous space of a
compact connected Lie group $G$. Suppose that the rational 
cohomology of $G/H$ is of the form
\[
\bfH^\bullet(G/H;\QQ)\isom\QQ[a]/(a^2)\otimes\EA_\QQ(u),
\]
where $\deg(a)$ is even and $\deg(u)$ is odd, and $\deg(a)<\deg(u)$.
Then $\rk(G)-\rk(H)=1$, and
the following diagram gives the induced maps in cohomology,
\begin{diagram}
\bfH^\bullet(H;\QQ) & \lTo^{h^\bullet} &
\bfH^\bullet(G;\QQ) & \lTo^{p^\bullet} &
\bfH^\bullet(G/H;\QQ) \\
\dEq && \dEq && \dEq \\
\EA_\QQ(w_1,\ldots,w_r) & \lTo &
\EA_\QQ(w_2,\ldots,w_r,u,v) & \lTo &
\QQ[a]/(a^2)\otimes\EA_\QQ(u)
\end{diagram}
where $p^\bullet(a)=0$,  and $h^\bullet(u)=0=h^\bullet(v)$.
Moreover, $\deg(v)=2\deg(a)-1$.
\end{Thm}

\begin{Num}\textsc{How to use this}\psn
\label{HowTo(II)}
Let $G/H$ be a homogeneous space as in the theorem above.
Let 
\[
(m_1^G,\ldots,m_{s+1}^G)
\]
denote the degrees of the homogeneous 
generators of the rational cohomology of $G$, and let
\[
(m_1^H,\ldots,m_s^H)
\]
denote the corresponding numbers for $H$. 
Then, after a suitable permutation
of the indices, $m_i^G=m_i^H$ for $i=1,\ldots,s-1$, and
$m_s^G=2m_s^H+1$ so
$(m_1^H,\ldots,m_{s-1}^H,2m_s^H+1)$ has to be a subsequence of 
$(m_1^G,\ldots,m_{s+1}^G)$. The homogeneous generators of
the cohomology of $G/H$ then have degrees $n_1=m_s^H+1$ and 
$n_2=m_{s+1}^G$.
Here is an example: $G=\SU(5)$ and $H=\SU(3)\times\SU(2)$.
\[
\begin{array}{llll}
m^G : &&& (\not\!3,\not\!5,\not\!7,9) \\
m^H : & (3,5,3)  & \leadsto & (3,5,\fbox{7}) \\ \hline
&&& (\fbox{4},9)
\end{array}
\]
\end{Num}

We want to show that there exists a transitive normal semisimple 
subgroup of $G$ with at most 2 almost simple factors. To prove
this, we use the following result.

\begin{Lem}[Hsiang-Su, Onishchik]
\label{Length2}
Let $G/H$ be a homogeneous space of a compact connected Lie group
$G$, and let $N\SUB G$ be a closed connected
normal subgroup. Let $\eta:H\SUB G$ and
$\nu:N\SUB G$ be the inclusion maps. Let $P_\eta$ and $P_\nu$
denote the induced maps on the vector spaces of primitive elements.
If $\mathrm{ker}(P_\eta)$ injects
into $P_N$ under $P_\nu$, then $N$ acts transitively on $X$.

\begin{proof}
Passing to a suitable compact covering, we may assume that 
$G\isom N\times G/N$. 
The maps $N\too G\too G/N$ yield a short exact
sequence 
\[
0\oot P_N\oot P_G\oot P_{G/N}\oot 0,
\]
and we obtain a diagram

%\newpage
\begin{diagram}
&&&& 0 \\
&&&& \uTo \\
&&&& P_N \\
&&&& \uTo & \luTo^\phi \\
0 & \lTo & \mathrm{im}(P_\eta) & \lTo & P_G & \lTo & \mathrm{ker}
(P_\eta) &
\lTo & 0 . \\
&&&\luTo_\psi& \uTo^{P_\nu} \\
&&&& P_{G/N} \\
&&&& \uTo \\
&&&& 0 \\
\end{diagram}
The row and the column are exact, hence $\phi$ is injective 
if and only if
$\psi$ is injective. By our assumption $\phi$ is injective, hence
the composite $H\rInto G\too G/N$ induces an injection 
\[
P_H\oot P_{G/N}.
\]
Therefore, the map 
$\bfH^\bullet(H;\QQ)\oot\bfH^\bullet(G/N;\QQ)$ is also an injection.
We factor it through the image $HN/N\cong H/H\cap N$ of $H$ in $G/N$ 
to obtain
\[
\bfH^\bullet(H;\QQ)\oot
\bfH^\bullet(HN/N;\QQ)\oot
\bfH^\bullet(G/N;\QQ).
\]
In particular, the map on the right must be an injection.
But $HN/N\SUB G/N$; therefore, the map on the right is
an isomorphism, and $HN/N=G/N$. Thus $HN=G$, and $N$ acts transitively on
$G/H$.
\end{proof}
\end{Lem}
A similar result is proved in Hsiang-Su 
\cite{HsSu68} Prop.1.4 and in Onishchik \cite{Oni94} \S17.1.

\begin{Prop}
\label{LengthAtMost2}
Let $G/H$ be a homogeneous space whose cohomology is either of the
form
\[
\EA_\ZZ(u,v),
\]
where $\deg(u),\deg(v)\geq 3$ are odd, or of the form
\[
\ZZ[a]/(a^2)\otimes\EA_\ZZ(u).
\]
where $\deg(a)$ is even and $\deg(u)$ is odd, and $\deg(u)>\deg(a)\geq 4$.
There exists a normal transitive semisimple subgroup $N\SUB G$
with at most two almost simple factors.

\begin{proof}
Passing to a suitable covering, we may assume that $G$ is 
a product of torus groups and simply connected almost simple groups,
\[
G=G_1\times\cdots\times G_s.
\]
Then $P_G=P_{G_1}\oplus\cdots\oplus P_{G_s}$. The kernel
of $P_\eta:P_G\too P_H$ is 2-dimensional. The inclusion
$G_i\rInto G$ induces the projection $P_G\too P_{G_i}$.
There exist numbers $i,j$ such that $\mathrm{ker}(P_\eta)$ injects into
$P_{G_i}\oplus P_{G_j}$ under the map induced from
$G_i\times G_j\rInto G$.
Therefore, $N=G_i\times G_j$ acts transitively on $G/H$ by 
\ref{Length2}.
\end{proof}
\end{Prop}

%\newpage

\section{The integral classification}

We use the following terminology. The point space of
the $n$-dimensional projective geometry over the skew field
$\FF$ is denoted by 
\[
\FF\mathrm{P}^n,
\]
for $\FF=\RR,\CC,\HH$.
The point space
of the projective Cayley plane is denoted by $\OO\mathrm{P}^2$.
The \emph{Stiefel manifold} of orthonormal $k$-frames in $\FF^n$ is
\[
V_k(\FF^n)=\{(v_,\ldots,v_k)\in\FF^{kn}|\ (v_i|v_j)=\delta_{i,j}\},
\]
where $(x|y)=\sum_{i=1}^n \bar x_iy_i$, for
$\FF=\RR,\CC,\HH,\OO$ (for $\FF=\OO$ we restrict $k$ to $k=1,2$ since
$V_k(\OO^n)$ is a singular algebraic variety and not a manifold for
$k>2$).
We define the \emph{oriented
Grassmann manifolds} as
\[
\widetilde{G}_k(\RR^n)=\SO(n)/\SO(k)\times\SO(n-k)
\text{ and }
\widetilde{G}_k(\CC^n)=\SU(n)/\SU(k)\times\SU(n-k).
\]

In Chapter 5 and 6 we prove the following main result.

\begin{Thm}
\label{MainTheorem}
Let $X=G/H$ be a 1-connected homogeneous space of a compact connected Lie 
group $G$. Suppose that the action of $G$ on $X$ is irreducible,
cp.~\ref{DefIrreducibleAction}. 

\textbf{(A)}
Assume that $X$ has the same rational cohomology as
a product of spheres
\[
\SS^{n_1}\times\SS^{n_2},
\]
where $n_2\geq n_1\geq 3$ and $n_2$ is odd.
Then $(G,H)$ is one of the pairs 
which we discuss in chapters 5-6. 

\textbf{(B)}
Suppose that $X$ has the same integral cohomology as
a product of spheres
\[
\SS^{n_1}\times\SS^{n_2},
\]
where $n_2\geq n_1\geq 3$ and $n_2$ is odd. There are the following
possibilities (and no others).

\textbf{(B1)}
If $n_1$ is odd, then $(G,H)$ is one of the pairs in
\ref{Case(I)Simple}, \ref{Case(I)SemiSimple}, or \ref{SplitCase}. 
More precisely we have the following spaces.
\begin{description}
\item[\fbox{Stiefel manifolds}]
\begin{gather*}
\SU(n)/\SU(n-2)=V_2(\CC^n),\quad n\geq 3 \\
\Sp(n)/\Sp(n-2)=V_2(\HH^n),\quad n\geq 2
\end{gather*}
\item[\fbox{Homogeneous sphere bundles}]
\begin{gather*}
\Sp(n)\times\SU(3)/\Sp(n-1)\cdot\Sp(1) \\
\Sp(n)\times\Sp(2)/\Sp(n-1)\cdot\Sp(1)
\end{gather*}
\item[\fbox{Products of homogeneous spheres}]
\[
K_1/H_1\times K_2/H_2
\]
where $K_1/H_1$ and $K_2/H_2$ is one of the spaces
\begin{gather*}
\SO(2n)/\SO(2n-1)=\SS^{2n-1} \\
\SU(n)/\SU(n-1)=\SS^{2n-1} \\
\Sp(n)/\Sp(n-1)=\SS^{4n-1} \\
\Spin(9)/\Spin(7)=\SS^{15} \\
\Spin(7)/\G_2=\SS^7
\end{gather*}
\item[\fbox{Sporadic spaces}]
\begin{gather*}
\E_6/\Ffour \\
\Spin(10)/\Spin(7) \\
\Spin(9)/\G_2=V_2(\OO^2) \\
\Spin(8)/\G_2=\SS^7\times\SS^7 \\
\SU(6)/\Sp(3)=\SU(5)/\Sp(2) 
\end{gather*}
\end{description}
\textbf{(B2)}
If $n_1\geq 4$ is even, then $(G,H)$ is one of the pairs in 
\ref{Case(II)Simple}, \ref{Case(II)SemiSimple}, or \ref{SplitCase}.
Thus, $G/H$ is one of the following spaces.
\begin{description}
\item[\fbox{Stiefel manifolds}]
\[
\SO(2n)/\SO(2n-2)=V_2(\RR^{2n}),\quad n\geq 3 
\]
\item[\fbox{Homogeneous sphere bundles}]
\[
\Sp(n)\times\Sp(2)/\Sp(n-1)\cdot\Sp(1)\cdot\Sp(1) 
\]
\item[\fbox{Products of homogeneous spheres}]
\[
K_1/H_1\times K_2/H_2
\]
where $K_1/H_1$ is one of the spaces
\begin{gather*}
\SO(2n+1)/\SO(2n)=\SS^{2n} \\
\G_2/\SU(3)=\SS^6
\end{gather*}
and $K_2/H_2$ is one of the spaces
\begin{gather*}
\SO(2n)/\SO(2n-1)=\SS^{2n-1} \\
\SU(n)/\SU(n-1)=\SS^{2n-1} \\
\Sp(n)/\Sp(n-1)=\SS^{4n-1} \\
\Spin(9)/\Spin(7)=\SS^{15} \\
\Spin(7)/\G_2=\SS^7
\end{gather*}
\item[\fbox{Sporadic spaces}]
\begin{gather*}
\Spin(10)/\SU(5)=\Spin(9)/\SU(4) \\
\Spin(7)/\SU(3)=V_2(\RR^8) \\
\Sp(3)/\Sp(1)\times\Sp(1) \\
\Sp(3)/\Sp(1)\times{}^\HH\rho_{3\lambda_1}(\Sp(1)) \\
\SU(5)/\SU(3)\times\SU(2)
\end{gather*}
(see Chapter 4 for the definition of
$^\HH\rho_{3\lambda_1}:\Sp(1)\too\Sp(3)$).
\end{description}
\textbf{(C)}
Suppose that $G/H$ is an $(m-1)$-connected rational cohomology 
$(2m+1)$-sphere, with $\pi_{m-1}(G/H)\isom\ZZ/2$, and that $m\geq 3$.
Then $m=2n-1$ is odd, $G/H=V_2(\RR^{2n+1})$ is a real Stiefel manifold,
and $G/H$ is one of the pairs
\begin{gather*}
\SO(2n+1)/\SO(2n-1)=V_2(\RR^{2n+1}),\quad n\geq 2 \\
\G_2/\SU(2)=V_2(\RR^7),
\end{gather*}
cp.~\ref{HomogeneousRealStiefel}.
\qed
\end{Thm}

\newpage

It is interesting to picture the distribution of the possible values
for $(n_1,n_2)$. Of course, one has to disregard the products of
homogeneous spheres, since for them every pair $(n_1,n_2)$ is possible.
The diagram below displays the numbers $(n_1,n_2-n_1)$.
The three horizontal infinite series are the 
Stiefel manifolds and the three vertical series are the homogeneous
sphere bundles; they are marked as white circles.
The sporadic cases are marked as black dots. A double circle
$\circ\!\!\circ$ or double dot indicates that there are two
irreducible group actions which yield the same numbers $(n_1,n_2-n_1)$.
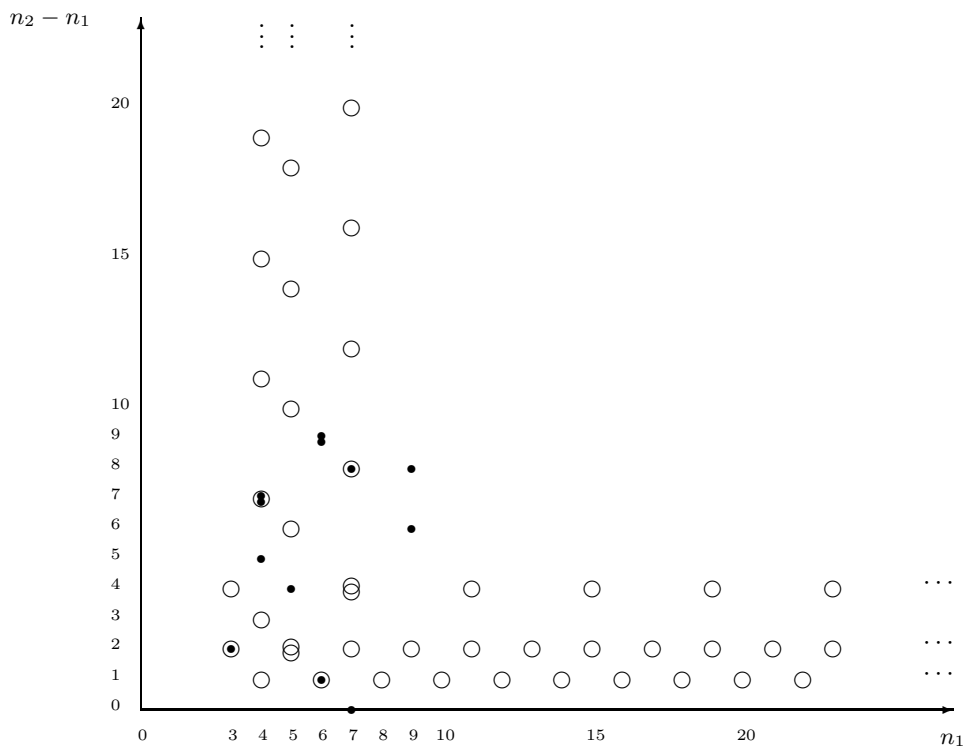
\begin{center}
\unitlength=4mm
\begin{picture}(30,30)(-2,-4)
\put(0,0){\vector(1,0){27}}
\put(27,-1){\makebox(0,0){\small$n_1$}}
\put(0,0){\vector(0,1){23}}
\put(-3,23){\makebox(0,0){\small$n_2-n_1$}}
\multiput(4,3)(0,4){5}{\circle{0.5}}
\put(5,1.9){\circle{0.5}}
\multiput(5,6)(0,4){4}{\circle{0.5}}
\put(7,3.9){\circle{0.5}}
\multiput(7,8)(0,4){4}{\circle{0.5}}
\put(3,4){\circle{0.5}}
\put(7,4.1){\circle{0.5}}
\multiput(11,4)(4,0){4}{\circle{0.5}}
\put(3,2){\circle{0.5}}
\put(5,2.1){\circle{0.5}}
\multiput(7,2)(2,0){9}{\circle{0.5}}
\multiput(4,1)(2,0){10}{\circle{0.5}}
\put(26,4){$\cdots$}
\put(26,2){$\cdots$}
\put(26,1){$\cdots$}
\put(3.9,22){$\vdots$}
\put(4.9,22){$\vdots$}
\put(6.9,22){$\vdots$}
\put(9,6){\circle*{0.3}}
\put(9,8){\circle*{0.3}}
\put(7,8){\circle*{0.3}}
\put(7,0){\circle*{0.3}}
\put(5,4){\circle*{0.3}}
\put(3,2){\circle*{0.3}}
\put(6,8.9){\circle*{0.3}}
\put(6,9.1){\circle*{0.3}}
\put(6,1){\circle*{0.3}}
\put(4,6.9){\circle*{0.3}}
\put(4,7.1){\circle*{0.3}}
\put(4,5){\circle*{0.3}}
\put(-1,0){\tiny$0$}
\put(-1,1){\tiny$1$}
\put(-1,2){\tiny$2$}
\put(-1,3){\tiny$3$}
\put(-1,4){\tiny$4$}
\put(-1,5){\tiny$5$}
\put(-1,6){\tiny$6$}
\put(-1,7){\tiny$7$}
\put(-1,8){\tiny$8$}
\put(-1,9){\tiny$9$}
\put(-1,10){\tiny$10$}
\put(-1,15){\tiny$15$}
\put(-1,20){\tiny$20$}
\put(-0.1,-1){\tiny$0$}
\put(2.9,-1){\tiny$3$}
\put(3.9,-1){\tiny$4$}
\put(4.9,-1){\tiny$5$}
\put(5.9,-1){\tiny$6$}
\put(6.9,-1){\tiny$7$}
\put(7.9,-1){\tiny$8$}
\put(8.9,-1){\tiny$9$}
\put(9.8,-1){\tiny$10$}
\put(14.8,-1){\tiny$15$}
\put(19.8,-1){\tiny$20$}
\end{picture}
\end{center}
An immediate consequence is the following result.

\begin{Cor}
\label{HighlyCon}
Let $X=G/H$ be a homogeneous space as in Theorem \ref{MainTheorem} (B).
Assume that $X$ is $9$-connected, or equivalently, that
$n_1>9$. Then $X$ is either a product of homogeneous spheres, or
$X$ is a Stiefel manifold.
\qed
\end{Cor}

Onishchik \cite{Oni94}
defines the \emph{rank} of a homogeneous space $G/H$ as
\[
\rk(G/H)=\dim(\mathrm{ker}[P_H\oot P_G]).
\]
The homogeneous spaces $G/H$ classified in Chapter 5 and 6
are therefore homogeneous of rank 1 or 2.
Onishchik's book \cite{Oni94}
contains tables (on p.~265 and p.~270)
of homogeneous spaces
of rank 1 (on p.~265) and rank 2 (on p.~270), provided that $G$ is 
almost simple (and for the case of rank 2, if $H$ is also simple).
Our classification shows that the table on p.~270 \emph{loc.~cit.}
is incomplete; the
spaces $\Spin(11)/\Spin(7)$ and $\Spin(10)/\Spin(7)$ are missing,
cp.~\ref{Spin11} and \ref{Spin10}.

%%%%%%%%%%%%%%%%%%%%%%%%%%%%%%%%%%%%%%%%%%%%%%%%%%%%%%%%%%%%%%%%%%%%%%%%
%                                                                      %
%                                                                      %
%                   Compact homogeneous quadrangles                    %
%                                                                      %
%                          Linus Kramer                                %
%                                                                      %
%                           Memoirs AMS                                %
%                                                                      %
%                         Wuerzburg 2000                               %
%                                                                      %
%                                                                      %
%                            CHQ4.tex                                  %
%                                                                      %
%                                                                      %
%                                                                      %
%%%%%%%%%%%%%%%%%%%%%%%%%%%%%%%%%%%%%%%%%%%%%%%%%%%%%%%%%%%%%%%%%%%%%%%%
\chapter{Representations of compact Lie groups}

In this chapter we collect the necessary facts about 
(finite dimensional) representations
of compact semisimple Lie groups. First we describe the irreducible
representations on real, complex and quaternionic vector spaces, and
we explain
how the problem of classifying certain subgroups of compact classical
groups can be reduced to representation theory. 

The next section gives tables of the fundamental representations of
all compact almost simple Lie groups, and also of all low-dimensional
irreducible representations. The last section introduces an algebraic
invariant, the Dynkin index, of a homomorphism between almost
simple compact Lie groups. 

I believe that these results and tables are of some independent
interest for geometers.

%Everything in this chapter is certainly known to experts.
The representation theory is based on Tits \cite{Tits71},
\cite{TitsTabellen}, 
and on the appendices in Onishchik-Vinberg \cite{OniVin90}.
The section about the Dynkin index is taken from Onishchik 
\cite{Oni94} p. 58--61.

Throughout this chapter $G$ is a compact semisimple Lie group.
Such a group may be viewed as an algebraic semisimple $\RR$-group, 
and all results about rational representations of algebraic
groups apply.

\section{The classification of irreducible representations}

The material in this section is taken from Tits \cite{TitsTabellen}
and Tits \cite{Tits71}, see in particular \cite{Tits71} 7.2 and 8.1.1.

Let $D,E$ be central finite dimensional 
(skew) fields over $\RR$, i.e.~$D,E\in\{\RR,\CC,\HH\}$, with
$D\subseteq E$. Let $X$ be a finite dimensional $E$-module. 
The ring of all $E$-linear
endomorphisms of $X$ is denoted by 
\[
\End_E(X).
\]
The group of invertible
endomorphisms is denoted by $\GL_E(X)=\End_E(X)^\times$.

We may also view $X$ as a $D$-module by \emph{restriction of scalars}; 
this module is denoted by 
\[
\rest^E_D(X),
\]
and we put $\End_D(X)=\End_D(\rest^E_D(X))$ for short. Clearly,
$\End_E(X)\SUB\End_D(X)$.
Conversely, if $Y$ is a 
$D$-module, then we can consider the $E$-module 
\[
\tens^D_E(Y)=Y\otimes_DE
\]
obtained by \emph{extension of scalars}. 

A continuous group homomorphism $\rho:G\too \GL_E(X)$ is called
an $E$-re\-pre\-sen\-ta\-t\-ion of $G$ on $X$; we call $X$ a $G$-$E$-module.
A homomorphism between two $G$-$E$-modules is a $E$-linear map which
is $G$-equivariant. In particular, we have the ring
\[
\End^G_E(X)=\Cen_{\End_E(X)}(\rho(G))
\]
of all $G$-equivariant endomorphisms of $X$.
We can apply the above functors $\rest$ and $\tens$ to $G$-$E$-modules;
then $\rest^E_D(X)$ is a $G$-$D$-module, and if $Y$ is a
$G$-$D$-module, then $\tens^D_E(Y)$ is a $G$-$E$-module.

A non-zero $G$-$E$-module is called ($E$-)\emph{simple} if it has no proper
non-zero
$G$-$E$-submodule, and \emph{semisimple} if it decomposes into a direct sum
of simple $G$-$E$-modules. Since $G$ is compact,
every $G$-$E$-module is semisimple.

If two $G$-$E$-modules are isomorphic, then the corresponding
representations are called \emph{equivalent}.
\emph{Note that equivalence depends on the ground field $E$.}
Let $G,G'$ be two compact semisimple
Lie groups, with $E$-modules $X,X'$. If there exists an isomorphism
$f:G\too G'$ and an $E$-linear isomorphism  $F:X\too X'$ such that 
$f(g)F(x)=F(gx)$, then the modules are called \emph{quasi-isomorphic}, 
and the representations of $G$ and $G'$ are called 
\emph{quasi-equivalent}.

We recall the following facts from complex representation theory.
Let $\Lambda$ denote the \emph{weight lattice} of $G$,
and let $\Lambda_{+}$
denote the set of \emph{dominant weights}. The weight lattice is generated
by \emph{fundamental weights} $\lambda_1,\ldots,\lambda_k$, where $k$ is
the rank of $G$,
\[
\Lambda=\bigoplus_{i=1}^k\lambda_i\ZZ.
\]
The dominant weights are the weights with non-negative coefficients,
\[
\Lambda_{+}=\bigoplus_{i=1}^k\lambda_i\NN_0.
\]

\begin{Num}\textsc{Classification of simple $G$-$\CC$-modules}\psn
To each dominant weight $\lambda\in\Lambda_{+}$
one can associate a $\CC$-module $X_\lambda$ and an
irreducible representation
$\rho_\lambda:G\too\End_\CC(X_\lambda)$. Conversely, every simple
$G$-$\CC$-module is $\CC$-isomorphic to a unique module of this type.
Thus, $\Lambda_{+}$ corresponds bijectively to the isomorphism classes
of irreducible $G$-$\CC$-modules.
The representations $\rho_{\lambda_i}$ corresponding to the fundamental
weights are called the \emph{fundamental representations}.
\end{Num}
To classify the real and quaternionic $G$-modules one needs
some additional structure.

\begin{Num}\textsc{The action of the Galois group}\psn
The Galois group $\Gamma=\mathrm{Gal}(\CC/\RR)$ acts on the collection of
all $G$-$\CC$-modules as follows. If $\rho:G\too\End_\CC(X)$ is
a representation, then $\bar\rho:G\too\End_\CC(X)$ is given by
applying complex conjugation to each entry of the matrix $\rho(g)$
(with respect to some basis of $X$). There is a corresponding action
of $\Gamma$ on the semigroup $\Lambda_{+}$, i.e.~$\overline{\rho_\lambda}=
\rho_{\bar\lambda}$. This action is additive, hence it suffices
to know $\bar\lambda_i$ for $i=1,\ldots,k$.
\end{Num}
Let $\lambda\in\Lambda_{+}$, and let $X=X_\lambda$ be the corresponding
simple $G$-$\CC$-module. 
A \emph{real structure} on $X$ is a conjugate-linear
involution $\theta\in\End_\RR^G(X)$, i.e.
\[
\theta^2=1\neq\theta\quad\text{ and }\quad\theta(xa)=\theta(x)\bar a
\]
for all $x\in X$ and $a\in\CC$.
If such a real structure exists, then 
$^\RR X=\mathrm{Fix}(\theta,X_\lambda)$ is a real $G$-invariant subspace. 
Moreover, $X=\tens^\RR_\CC(^\RR X)$.
We denote the representation of $G$ on $^\RR X$ by
\[
{}^\RR\rho_\lambda:G\too\End_\RR(^\RR X).
\]
If no such real structure exists on $X$, then we put
$^\RR X=\rest^\CC_\RR(X)$ and
\[
{}^\RR\rho_\lambda=\rest^\CC_\RR(\rho_\lambda):G\too\End_\RR(^\RR X).
\]

\begin{Num}\textsc{Classification of simple $G$-$\RR$-modules}\psn
Every simple $G$-$\RR$-module $Y$ is isomorphic to a module $^\RR X$
obtained by the method described above.
Moreover, if $^\RR\rho$ and $^\RR\sigma$ are $\RR$-equivalent 
irreducible representations, then $\rho=\sigma$ or $\rho=\bar\sigma$.
Thus, the $\RR$-irreducible $G$-$\RR$-modules correspond bijectively
to the orbit space $\Lambda_{+}/\Gamma$.
\end{Num}
A \emph{quaternionic structure} on a $G$-$\CC$-module $X$ is a semilinear
map $\theta\in\End_\RR^G(X)$ with
\[
\theta^2=-1\quad\text{ and }\quad\theta(xa)=\theta(x)\bar a
\]
for all $x\in X$ and $a\in \CC$.
If such a quaternionic structure exists, then the ring generated 
over $\RR$ by 
$\theta$ and $x\mapstoo xi$ is clearly isomorphic to $\HH$.
We denote the corresponding $G$-$\HH$-module by $^\HH X$.
Thus, $X=\rest^\HH_\CC(^\HH X)$, and 
\[
^\HH\rho:G\too\End_\HH(^\HH X)
\]
is defined by $\rho=\rest^\HH_\CC(^\HH\rho)$.
If no such quaternionic structure exists, then we put
$^\HH X=\tens^\CC_\HH(X)$ and $^\HH\rho=\tens^\CC_\HH(\rho)$.

\begin{Num}\textsc{Classification of simple $G$-$\HH$-modules}\psn
Every simple $G$-$\HH$-module $W$ is isomorphic to a module $^\HH Y$
obtained by the method described above.
Moreover, if $^\HH\rho$ and $^\HH\sigma$ are $\HH$-equivalent 
irreducible representations, then $\rho=\sigma$ or $\rho=\bar\sigma$.
Thus, the $\HH$-irreducible $G$-$\HH$-modules correspond bijectively
to the orbit space $\Lambda_{+}/\Gamma$.
\end{Num}

\begin{Num}\textsc{Existence of real and quaternionic structures}\psn
Let $\lambda\in\Lambda_{+}$. A real or quaternionic structure exists
on $X_\lambda$ if and only if 
$\lambda\in\mathrm{Fix}(\Gamma,\Lambda_{+})$. 
There is a map 
\[
\beta:\mathrm{Fix}(\Gamma,\Lambda_{+})\too\mathrm{Br}(\RR)
\]
from the $\Gamma$-fixed dominant weights into the Brauer group 
$\mathrm{Br}(\RR)=\{\RR,\HH\}\isom\ZZ/2$
such that $X_\lambda$ admits a real (resp. a quaternionic) structure
if and only if $\beta(\lambda)=\RR$ (resp. $\beta(\lambda)=\HH$).
This map is additive in the sense that
$\beta(\lambda+\mu)=\beta(\lambda)\beta(\mu)$.

For an arbitrary $\Gamma$-invariant dominant weight $\lambda=\bar\lambda$,
the value $\beta(\lambda)$ can be determined as follows.
Let $\lambda=\sum n_i\lambda_i\in\mathrm{Fix}(\Gamma,\Lambda_{+})$ and
put $Q(\lambda)=\{i|\ \lambda_i=\bar\lambda_i,\ \beta(\lambda_i)=\HH\}$.
Then $\beta(\lambda)=\HH$ if and only if
$\sum_{i\in Q(\lambda)}n_i\equiv 1\pmod2$.
Thus, it suffices to know $\beta(\lambda_i)$ for the
$\Gamma$-invariant fundamental weights.
\end{Num}

%\newpage

\section{Subgroups of classical groups}

We are not really interested in homomorphism of $G$ into the linear group
$\GL_D(X)$, but rather into certain compact subgroups of $\GL_D(X)$.
The Cartan decomposition of the groups $\GL_n\RR$, $\GL_n\CC$ and
$\GL_n\HH$ leads to the following result.
\begin{Lem}
\label{PolarTrick}
Let $(G,K)$ denote one of the pairs $(\GL_n\RR,\O(n))$,
$(\GL_n\CC,\U(n))$, or $(\GL_n\HH,\Sp(n))$.
Let  
\[
\rho,\sigma:H\too K\SUB G
\]
be homomorphisms which are conjugate by an element $g\in G$,
\[
g\rho(h)=\sigma(h)g
\]
for all $h\in H$. Then there exists an element $k\in K$ such that
\[
k\rho(h)=\sigma(h)k.
\]

\begin{proof}
We use the Cartan decomposition $G=KP$. Note that $P$ is a $K$-invariant
subset of $G$. We decompose $g=kp$ with $k\in K$, $p\in P$.
Then 
\[
\underbrace{\sigma(h)k}_{\in K}p=kp\rho(h)=
\underbrace{k\rho(h)}_{\in K}
\underbrace{\rho(h)^{-1}p\rho(h)}_{\in P}.
\]
It follows from the uniqueness of the Cartan decomposition that 
\[
\sigma(h)k=k\rho(h)\text{ and }
p=\rho(h)^{-1}p\rho(h).
\]
\end{proof}
\end{Lem}
This can be improved.
\begin{Lem}
Let $\rho,\sigma:H\too\SU(n)$ be homomorphisms. If $\rho$ and $\sigma$
are conjugate by an element $u\in \U(n)$, then they are conjugate by an 
element in $\SU(n)$.

Similarly, suppose that $\rho,\sigma:H\too\SO(n)$ are homomorphisms which
 are conjugate by an element $u\in\O(n)$. If $n$ is odd, then they are
conjugate by an element of $\SO(n)$.

\begin{proof}
We may multiply $u$ by a number $z$ such that
$z^n=\det(u)^{-1}$. Then $\det(zu)=1$.
\end{proof}
\end{Lem}
This has the following consequence. In order to determine all
subgroups of a given type of a compact classical group $G$,
\emph{up to automorphisms of $G$}, it suffices to classify all 
representations of groups of this given type on the natural
$G$-module. If $G=\SO(n)$ is an orthogonal group, then the same method
gives all subgroups of a given type in the universal covering
$\Spin(n)$ (by considering the connected component of the lifted subgroup).
The only thing that one has to check in this case is whether two
different subgroups of $\SO(n)$ become equivalent in $\Spin(n)$.
This happens indeed: for example,
all subgroups of type $\fb_3$ in $\Spin(8)$
are equivalent under the automorphism group of $\Spin(8)$, although
there are different (not quasi-isomorphic) representations 
of $\Spin(7)$ on $\RR^8$.

%\newpage

\section{Useful formulas}

The results in this chapter follow from the tables in Tits 
\cite{TitsTabellen}, cp.~also Onish\-chik-Vinberg \cite{OniVin90}
in particular p. 299--305, Salzmann \emph{et al.}
\cite{CPP95} p. 616--630. B\"odi-Joswig \cite{BoediJoswig93}
contains an algorithm which determines all irreducible 
representations of a 
given simple group (up to a certain dimension). I frequently used the
computer implementation of this algorithm by the authors.

We use the following fact.
Let $\lambda,\mu\in\Lambda_{+}$ be dominant weights of an almost simple
compact Lie group. Suppose that 
\[
\mu=\lambda+n_1\lambda_1+\ldots+ n_k\lambda_k,
\]
with $n_1,\ldots,n_k\geq 0$, and that not all $n_i$ are zero. Then 
\[
\dim(\rho_\lambda)<\dim(\rho_\mu).
\]
Note also that 
\[
\dim_\RR(^\RR\rho)\in\{\dim_\CC(\rho),2\dim_\CC(\rho))\}
\]
and that
\[
\dim_\HH(^\HH\rho)\in\left\{\frac12\dim_\CC(\rho),\dim_\CC(\rho)\right\}.
\]
Suppose that $G$ is simply connected. There exist homomorphisms
$e_i:\Cen(G)\too\CC^{*}$ such that
\[
\rho_{m_1\lambda_1+\cdots m_k\lambda_k}(z)=e_1(z)^{m_1}\cdots e_k(z)^{m_k}
\]
for all $z\in \Cen(G)$. Thus, it suffices to know the $e_i$ for $i=1,\ldots,k$
to determine the kernel of a representation.
Note also that $\ker(\rho)=\ker({}^\RR\rho)=\ker({}^\HH\rho)$.
We denote the $k$th exterior power of a vector space $X$ by
\[
\EA^k X,
\]
and the $k$th symmetric power by 
\[
\mathrm{S}^kX.
\]

\newpage

\begin{Num}\textsc{Fundamental weights for $\SU(2)$}\psn
\label{A1Modules}%
The weight lattice has one generator, $\lambda_1$, and
$\bar\lambda_1=\lambda_1$. The center is 
\[
\Cen(\SU(2))\isom\bra{z|\ z^2=1}\isom\ZZ/2
\]
and $e_1(z)=-1$. Moreover, $\beta(\lambda_1)=\HH$, hence
\[
\beta(k\lambda_1)=
\begin{cases} 
\RR & \text{ if $k\equiv 0\pmod 2$} \\
\HH & \text{ if $k\equiv 1\pmod 2$.}
\end{cases}
\]
The dimension of $\rho_{k\lambda_1}$ is
\[
\dim(\rho_{k\lambda_1})=k+1
\]
and thus
\[
\dim(^\RR\rho_{k\lambda_1})=
\begin{cases}
k+1 & \text{ if $k\equiv 0\pmod 2$} \\
2(k+1) &  \text{ if $k\equiv 1\pmod 2$.} \\
\end{cases}
\]
In fact 
\[
X_{k\lambda_1}=\mathrm{S}^k\CC^2.
\]
\end{Num}

\newpage

\begin{Num}\textsc{Fundamental weights for $\SU(n+1)$}\psn
We label the Dynkin diagram as follows, for $n\geq 2$.
\begin{diagram}[abut]
\overset{1}\bullet & \rLine & 
\overset{2}\bullet & \rLine & 
\overset{3}\bullet & \rDots & 
\overset{\llap{$\scriptstyle n-1$}}\bullet & \rLine & 
\overset{n}\bullet
\end{diagram}
The center of the universal covering is 
\[
\Cen(\SU(n+1))\isom\bra{z|\ z^{n+1}=1}\isom\ZZ/(n+1),
\]
and $e_i(z)=z^i$.
The dimensions of the fundamental representations are
\[
\dim(\rho_{\lambda_i})=\binom{n+1}{i}
\]
for $i=1,\ldots,n$.
The natural module is $R=X_{\lambda_1}=\CC^{n+1}$, and
\[
X_{\lambda_i}=\EA^iR
\qquad
X_{k\lambda_1}=\mathrm{S}^kR
\qquad
X_{\lambda_1+\lambda_n}=\mathrm{Ad}.
\]
The Galois group $\Gamma$ 
acts as $\bar\lambda_i=\lambda_{n+1-i}$ for $i=1,\ldots,n$. 
Moreover
\[
\beta(\lambda_{\frac{n+1}2})=\RR \text{ for $n+1\equiv 0\pmod 2$.}
\]
\end{Num}

%\newpage

\begin{Num}\textsc{Low-dimensional simple $\SU(n+1)$-modules}\psn
\label{AnModules}%
Let $X$ be a simple $\SU(n+1)$-$\CC$-module, with $\dim(X)\leq 4(n+1)$, 
and let $G\SUB\GL(X)$ denote the represented group. If
$n\geq 9$, then $(G,X)$ is quasi-isomorphic to $(\SU(n+1),\CC^{n+1})$. 
For $2\leq n\leq 8$, there are other low-dimensional simple modules.
\begin{center}
\begin{tabular}{llllll}
$G$ & weight & $X$ & $\dim(X)$ & $^\RR X$ & $\dim(^\RR X)$ \\
\hline
$\SU(n+1)$ & $\lambda_1$ or $\lambda_n$ 
 & $\CC^{n+1}$ & $n+1$  & $\CC^{n+1}$ & $2(n+1)$ \\
\hfill$n\geq 2$ &             & \\
\hline
$\SU(9)$  & $\lambda_2$ or  $\lambda_7$ 
 & $\EA^2\CC^9$        & $36$ & $\EA^2\CC^9$ & $72$ \\
\hline
$\SU(8)/\bra{z^4}$  & $\lambda_2$ or  $\lambda_6$ 
 & $\EA^2\CC^8$        & $28$ & $\EA^2\CC^8$ & $56$ \\
%$\PSU(8)$  & $\lambda_1+\lambda_7$ 
% & $\frak{sl}_8\CC$ & $63$ & $\frak{su}_8\CC$ & $63$ \\
\hline
$\SU(7)$  & $\lambda_2$ or  $\lambda_5$ 
 & $\EA^2\CC^7$        & $21$ & $\EA^2\CC^7$ & $42$ \\
$\PSU(7)$  & $\lambda_1+\lambda_6$ 
% & $\frak{sl}_7\CC$ & $48$ & $\frak{su}_7\CC$ & $48$ \\
%$\SU(7)$  & $2\lambda_1$ or $2\lambda_6$ 
 &  & $28$ & & $56$ \\
%$\SU(7)$  & $\lambda_3$ or $\lambda_4$ 
% & $\EA^3\CC^7$        & $35$ & $\EA^3\CC^7$ & $70$  \\
\hline
$\SU(6)/\bra{z^3}$  & $\lambda_2$ or  $\lambda_4$ 
 & $\EA^2\CC^6$        & $15$ & $\EA^2\CC^6$ & $30$ \\
%$\PSU(6)$  & $\lambda_1+\lambda_5$ 
% & $\frak{sl}_6\CC$ & $35$ & $\frak{su}_6\CC$ & $35$ \\
$\SU(6)/\bra{z^2}$  & $\lambda_3$ 
 & $\EA^3\CC^6$ & $20$ & $\EA^3\CC^6$ 
& $40$ \\
$\SU(6)/\bra{z^3}$  & $2\lambda_1$ or $2\lambda_5$ 
 & S$^2\CC^6$        & $21$ & S$^2\CC^6$ & $42$  \\
\hline
$\SU(5)$  & $\lambda_2$ or  $\lambda_3$ 
 & $\EA^2\CC^5$        & $10$ & $\EA^2\CC^5$ & $20$ \\
%$\PSU(5)$  & $\lambda_1+\lambda_4$ 
% & $\frak{sl}_5\CC$ & $24$ & $\frak{su}_5\CC$ & $24$ \\
$\SU(5)$  & $2\lambda_1$ or $2\lambda_3$
 &  & $15$ & & $30$ \\
\hline
$\SO(6)$  & $\lambda_2$ 
 & $\CC^6$           & $6$ & $\RR^6$ & $6$ \\
$\PSU(4)$ & $\lambda_1+\lambda_3$ 
 & $\frak{sl}_4\CC $ & $15$ & $\frak{su}_4\CC$ & $15$ \\
$\SO(6)$  & $2\lambda_1$ or $2\lambda_3$ 
 & S$^2\CC^4$        & $10$ & S$^2\CC^4$ & $20$  \\
%$\PSU(4)$  & $2\lambda_2$  
% & & $20$ & & $20$   \\
\hline
$\SU(3)$  & $2\lambda_1$ or $2\lambda_2$ 
 & S$^2\CC^3$        & $6$ & S$^2\CC^3$  & $12$ \\
$\PSU(3)$ & $\lambda_1+\lambda_2$ 
 & $\frak{sl}_3\CC$  & $8$ & $\frak{su}_3\CC$ & $8$ \\
$\SU(3)$  & $3\lambda_1$ or $3\lambda_2$ 
 & S$^3\CC^3$        & $10$ & S$^3\CC^3$ & $20$  
\end{tabular}
\end{center}

\begin{proof}
We have $\dim(\rho_{\lambda_2})= n(n+1)/2>4n+4$ for $n\geq 9$.
Now $\dim(\rho_{\lambda_i})\geq\dim(\rho_{\lambda_2})$ for $i=2,\ldots,n-2$, 
and $\dim(\rho_{2\lambda_1})=\binom{n+2}2>4(n+1)$ for $n\geq7$.
Therefore $\dim(\rho_\lambda)>\dim(\rho_{\lambda_1})$ for
$n\geq 9$ and $\lambda\neq\lambda_1,\lambda_n$.
\end{proof}
\end{Num}

\newpage

\begin{Num}\textsc{Fundamental weights for $\Spin(2n+1)$}\psn
We label the Dynkin diagram as follows, for $n\geq 3$.
\begin{diagram}[abut]
\overset{1}\bullet & \rLine & 
\overset{2}\bullet & \rLine & 
\overset{3}\bullet & \rDots & 
\overset{\llap{$\scriptstyle n-1$}}\bullet & \rImplies & 
\overset{n}\bullet
\end{diagram}
The center of the universal covering is 
\[
\Cen(\Spin(2n+1))=\bra{z|\ z^2=1}\isom\ZZ/2,
\]
$e_i(z)=1$ for $i=1,\ldots,n-1$, and $e_n(z)=-1$.
The dimensions of the fundamental representations are
\[
\dim(\rho_{\lambda_i})=\binom{2n+1}{i}
\]
for $i=1,\ldots,n-1$, and 
\[
\dim(\rho_{\lambda_n})=2^n.
\]
The natural module is $R=X_{\lambda_1}=\CC^{2n+1}$, and
\[
X_{\lambda_i}=\EA^iR
\quad(1\leq i\leq n-1)\qquad
X_{2\lambda_n}=\EA^n R
\qquad
X_{\lambda_2}=\EA^2 R=\mathrm{Ad}\text{ if }n\geq 3,
\]
$\mathrm{Ad}=X_{2\lambda_2}$ for $n=2$.
The Galois group $\Gamma$ acts as
$\lambda_i=\bar\lambda_i$ for $i=1,\ldots,n$. Moreover
\[
\beta(\lambda_i)=\RR \text{ for $i=1,\ldots,n-1$}
\]
and
\[
\beta(\lambda_n)=\begin{cases}
\RR & \text{for $n\equiv 0,3 \pmod 4$} \\
\HH & \text{for $n\equiv 1,2 \pmod 4$.}
\end{cases}
\]
\end{Num}

%\newpage

\begin{Num}\textsc{Low-dimensional simple $\Spin(2n+1)$-modules}\psn
\label{BnModules}%
Let $X$ be a simple $\Spin(2n+1)$-$\CC$-module, with $\dim(X)\leq 4(2n+1)$, 
and let $G\SUB\GL(X)$ denote the represented group. If
$n\geq 6$, then $(G,X)$ is quasi-isomorphic to $(\SO(2n+1),\CC^{2n+1})$. 
For $2\leq n\leq 5$, we have additional simple modules.
\begin{center}
\begin{tabular}{llllll}
$G$ & weight & $X$ & $\dim(X)$ & $^\RR X$ & $\dim(^\RR X)$ \\
\hline
$\SO(2n+1)$ & $\lambda_1$ & $\CC^{2n+1}$ & $2n+1$ & $\RR^{2n+1}$ & $2n+1$ \\
\hfill$n\geq 2$ &             & \\
\hline
$\Spin(11)$  & $\lambda_5$ & & $32$ &  & $64$ \\
\hline
$\Spin(9)$  & $\lambda_4$ & $\OO^2\otimes_\RR\CC$ & $16$ & $\OO^2$ & $16$ \\
$\SO(9)$    & $\lambda_2$ & $\frak{\so}_9\CC$ & $36$ & $\frak{so}_9$ & $36$ \\
\hline
$\Spin(7)$  & $\lambda_3$ & $\CC^8$ & $8$ & $\RR^8$  & $8$ \\
$\SO(7)$  & $\lambda_2$ & $\frak{\so}_7\CC$ & $21$ & $\frak{so}_7$ & $21$ \\
$\SO(7)$  & $2\lambda_1$ & & $27$ && $27$ \\
\hline
$\Sp(2)$    & $\lambda_2$ & $\HH^2$    & $4$ & $\HH^2$ & $8$  \\
$\SO(5)$    & $2\lambda_2$& $\frak{so}_5\CC$ & $10$     & $\frak{so}_5$ & $10$ \\
$\SO(5)$    & $2\lambda_1$&  & $14$     &  & $14$ \\
$\SO(5)$    & $\lambda_1+\lambda_2$ &  & $16$  &  & $32$ \\
$\SO(5)$    & $3\lambda_1$ &  & $20$  &  & $40$ 
\end{tabular}
\end{center}

\begin{proof}
We have $\dim(\rho_{\lambda_2})=n(2n+1)>4(2n+1)$ for $n\geq 5$, and
$\dim(\rho_{\lambda_i})>\dim(\rho_{\lambda_2})$ for
$i=1,\ldots,n-1$. Moreover, $2^n>4(2n+1)$ for $n\geq 6$.
\end{proof}
\end{Num}

\newpage

\begin{Num}\textsc{Fundamental weights for $\Sp(n)$}\psn
We label the Dynkin diagram as follows, for $n\geq 3$.
\begin{diagram}[abut]
\overset{1}\bullet & \rLine & 
\overset{2}\bullet & \rLine & 
\overset{3}\bullet & \rDots & 
\overset{\llap{$\scriptstyle n-1$}}\bullet & \lImplies & 
\overset{n}\bullet
\end{diagram}
The center of the universal covering is 
\[
\Cen(\Sp(n))=\bra{z|\ z^2=1}\isom\ZZ/2,
\]
and $e_i(z)=(-1)^{i}$.
The dimensions of the fundamental representations are
\[
\dim(\rho_1)=2n
\]
and
\[
\dim(\rho_{\lambda_i})=\binom{2n}{i}-\binom{2n}{i-2}
=\binom{2n}{i-2}\frac{(2n+1)(2n-2i+2)}{(2n-i+1)(2n-i+2)}
\]
for $i=2,\ldots,n$.
The natural module is $R=X_{\lambda_1}=\CC^{2n}=\HH^n$, and
\[
X_{k\lambda_1}=\mathrm{S}^kR
\qquad
X_{2\lambda_1}=\mathrm{Ad}.
\]
The Galois group $\Gamma$ acts as
$\lambda_i=\bar\lambda_i$ for $i=1,\ldots,n$. Moreover
\[
\beta(\lambda_i)=
\begin{cases}
\RR &\text{ for $i\equiv 0\pmod2$}\\
\HH &\text{ for $i\equiv 1\pmod2$}
\end{cases}
\]
\end{Num}

%\newpage

\begin{Num}\textsc{Low-dimensional simple $\Sp(n)$-modules}\psn
\label{CnModules}%
Let $X$ be a simple $\Sp(n)$-$\CC$-module, with $\dim(X)\leq 4\cdot 2n$, 
and let $G\SUB\GL(X)$ denote the represented group. If
$n\geq 5$, then $(G,X)$ is quasi-isomorphic to
\linebreak $(\Sp(n),\HH^n)$. 
For $2\leq n\leq 4$, we have additional simple modules.
\begin{center}
\begin{tabular}{llllllll}
$G$ & weight & $X$ & $\dim(X)$ & $^\RR X$ & $\dim({}^\RR X)$ & $^\HH X$ &
$\dim({}^\HH X)$ \\
\hline
$\Sp(n)$ & $\lambda_1$ & $\HH^n$ & $2n$ & $\HH^n$ & $4n$ & $\HH^n$ & $n$ \\
\hfill$n\geq 3$ &             & \\
\hline
$\PSp(4)$   & $\lambda_2$ && 27 && 27 && 27 \\
\hline
$\PSp(3)$   & $\lambda_2$ && 14 && 14 && 14 \\
$\Sp(3)$    & $\lambda_3$ && 14 && 28 && 7 \\
$\PSp(3)$   & $2\lambda_1$ & $\frak{sp}_6\CC$ & 21 & $\frak{sp}_3$ & 21 &
 & 21 \\
\end{tabular}
\end{center}

\begin{proof}
We have $\dim(\rho_{\lambda_2})=\binom{2n}{2}-1>8n$ for $n\geq 5$.
An unpleasant calculation shows that $\dim(\rho_{\lambda_k})\geq
\dim(\rho_{\lambda_2})$ for $k=3,\ldots,n$.
\end{proof}
\end{Num}

\newpage

\begin{Num}\textsc{Fundamental weights for $\Spin(2n)$}\psn
We label the Dynkin diagram as follows, for $n\geq 4$.
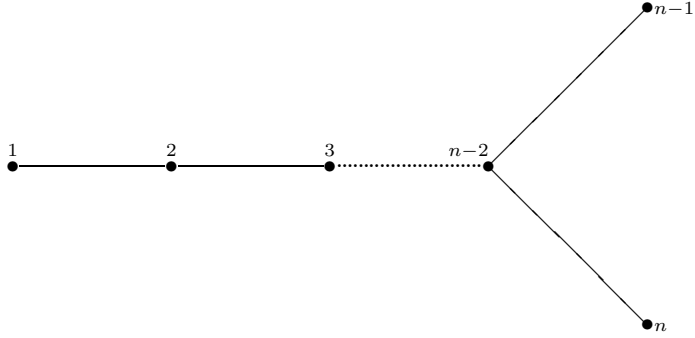
\begin{diagram}[abut]
&&&&&&&&\bullet\rlap{$\scriptstyle n-1$} \\
&&&&&&&\ruLine \\
\overset{1}\bullet & \rLine & 
\overset{2}\bullet & \rLine & 
\overset{3}\bullet & \rDots & 
\overset{\llap{$\scriptstyle n-2$}{}}\bullet \\
&&&&&&&\rdLine \\
&&&&&&&&\bullet\rlap{$\scriptstyle n$}
\end{diagram}
The center of the universal covering is as follows.

If $n$ is \emph{odd}, then
\[
\Cen(\Spin(2n))=\bra{z|\ z^4=1}\isom \ZZ/4,
\]
$e_i(z)=(-1)^i$ for $i=1,\ldots,n-2$, and $e_{n-1}(z)=e_n(z)=\sqrt{-1}$.

If $n$ is \emph{even}, then
\[
\Cen(\Spin(2n))=\bra{z,z'|\ z^2={z'}^2=1,\ zz'=z'z}\isom\ZZ/2\oplus\ZZ/2,
\]
$e_i(z)=e_i(z')=(-1)^i$ for $i=1,\ldots,n-2$, and
$e_{n-1}(z)=e_n(z')=1$, $e_{n-1}(z')=e_n(z)=-1$.

The dimensions of the fundamental representations are
\[
\dim(\rho_{\lambda_i})=\binom{2n}{i}
\]
for $i=1,\ldots,n-2$, and
\[
\dim(\rho_{\lambda_{n-1}})=\dim(\rho_{\lambda_n})=2^{n-1}.
\]
The natural module is $R=X_{\lambda_1}=\CC^{2n}$, and
\[
X_{\lambda_i}=\EA^iR
\quad
(1\leq i\leq n-2)
\qquad
X_{\lambda_2}=\EA^2 R=\mathrm{Ad},
\]
$X_{\lambda_{n-1}+\lambda_n}=\EA^{n-1}R$.
The Galois group $\Gamma$ acts as
$\lambda_i=\bar\lambda_i$ for $i=1,\ldots,n-2$, and
\[
\bar\lambda_{n-1}=
\begin{cases}
\lambda_n     & \text{ for $n\equiv1\pmod 2$} \\
\lambda_{n-1} & \text{ for $n\equiv0\pmod2$.}
\end{cases}
\]
Moreover
\[
\beta(\lambda_i)=\RR
\]
for $i-1,\ldots,n-2$, and
\[
\beta(\lambda_{n-1})=\beta(\lambda_n)=
\begin{cases}
\RR &\text{ for $n\equiv 0\pmod4$}\\
\HH &\text{ for $n\equiv 2\pmod4$.}
\end{cases}
\]
\end{Num}

\begin{Num}\textsc{Low-dimensional simple $\Spin(2n)$-modules}\psn
\label{DnModules}%
Let $X$ be a simple $\Spin(2n)$-$\CC$-module, with $\dim(X)\leq 4\cdot 2n$, 
and let $G\SUB\GL(X)$ denote the represented group. If
$n\geq 7$, then $(G,X)$ is quasi-isomorphic to $(\SO(2n),\CC^n)$. 
For $n=4,5,6$ there are additional simple modules.
\begin{center}
\begin{tabular}{llllll}
$G$ & weight & $X$ & $\dim(X)$ & $^\RR X$ & $\dim({}^\RR X)$ \\
\hline
$\SO(2n)$ & $\lambda_1$ & $\CC^{2n}$ & $2n$ & $\RR^{2n}$ & $2n$ \\
\hfill$n\geq 5$ &             & \\
\hline
$\Spin(12)$ & $\lambda_5$ or $\lambda_6$ 
 & $\HH^{16}$ & 16 & $\HH^{16}$ & 64 \\
\hline
$\Spin(10)$ & $\lambda_4$ or $\lambda_5$ 
 & $\CC^{16}$ & 16 & $\CC^{16}$ & 32 \\
\hline
$\SO(8)$ & $\lambda_1$, $\lambda_3$ or $\lambda_4$ 
 & $\CC^8$ & 8 & $\RR^8$ & 8 \\
$\PSO(8)$ & $\lambda_2$ & $\frak{so}_8\CC$ & 28 & $\frak{so}_8$ & 28 
\end{tabular}
\end{center}

\begin{proof}
We have $\dim(\rho_{\lambda_2})=\binom{2n}2>8n$ for $n\geq 5$, and
$\dim(\rho_{\lambda_i})>\dim(\rho_{\lambda_2})$ for $i=3,\ldots,n-2$.
Moreover, $2^{n-1}>8n$ for $n\geq 7$.
\end{proof}
\end{Num}

\newpage

\begin{Num}\textsc{Fundamental weights for $\E_6$}\psn
We label the Dynkin diagram as follows.
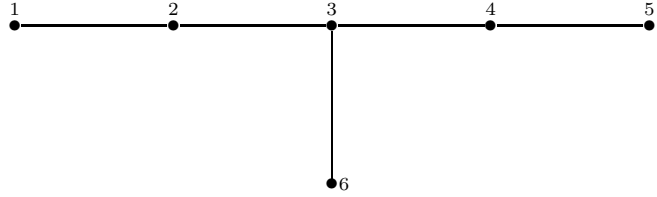
\begin{diagram}[abut]
\overset{1}\bullet & \rLine &
\overset{2}\bullet & \rLine & 
\overset{3}\bullet & \rLine &
\overset{4}\bullet & \rLine &
\overset{5}\bullet  \\
&&&& \dLine \\
&&&& \bullet\rlap{$\scriptstyle 6$} \\
\end{diagram}
The center of the universal covering is 
\[
\Cen(\E_6)=\bra{z|\ z^3=1}\isom\ZZ/3,
\]
and $e_i(z)=\exp(2\pi\sqrt{-1}/3)$.
The dimensions of the fundamental representations are
\[
\dim(\rho_{\lambda_i})=
\begin{cases}
27 & \text{ for $i=1,5$} \\
351 & \text{ for $i=2,4$} \\
2\, 925 & \text{ for $i=3$} \\
78 & \text{ for $i=6$.}
\end{cases}
\]
The natural module is $X_{\lambda_1}$, and
$X_{\lambda_6}=\mathrm{Ad}$.

The Galois group $\Gamma$ acts as 
$\lambda_i=\bar\lambda_{6-i}$ for $i=1,2$, $\bar\lambda_3=\lambda_3$,
and $\bar\lambda_6=\lambda_6$.
Moreover, $\beta(\lambda_3)=\beta(\lambda_6)=\RR$.
\end{Num}

\begin{Num}\textsc{Low-dimensional simple $\E_6$-modules}\psn
\label{E6Modules}%
The $\E_6$-modules of dimension at most 108 are the following.
\begin{center}
\begin{tabular}{llllll}
$G$ & weight & $X$ & $\dim(X)$ &$^\RR X$ & $\dim({}^\RR X)$ \\
\hline
$\E_6$ & $\lambda_1$ or $\lambda_5$ && 27 && 54  \\
$\mathrm{PE}_6$ & $\lambda_6$ & $\frak{e}_6\CC$ & 78 & $\frak{e}_6$ & 78
\end{tabular}
\end{center}
\qed
\end{Num}

\newpage

\begin{Num}\textsc{Fundamental weights for $\E_7$}\psn
We label the Dynkin diagram as follows.
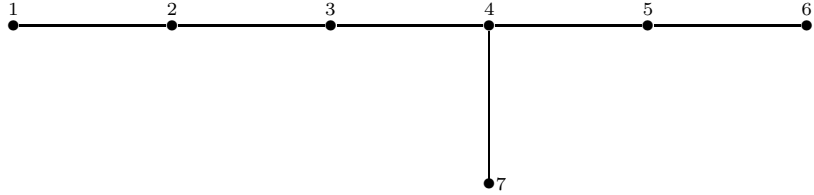
\begin{diagram}[abut]
\overset{1}\bullet & \rLine &
\overset{2}\bullet & \rLine & 
\overset{3}\bullet & \rLine &
\overset{4}\bullet & \rLine &
\overset{5}\bullet & \rLine &
\overset{6}\bullet  \\
&&&&&& \dLine \\
&&&&&& \bullet\rlap{$\scriptstyle 7$} \\
\end{diagram}
The center of the universal covering is 
\[
\Cen(\E_7)=\bra{z|\ z^2=1}\isom\ZZ/2,
\]
and $e_1(z)=e_3(z)=e_7(z)=-1$, $e_2(z)=e_4(z)=e_5(z)=e_6(z)=1$.
The dimensions of the fundamental representations are
\[
\dim(\rho_{\lambda_i})=
\begin{cases}
56 & \text{ for $i=1$} \\
1\, 539 & \text{ for $i=2$} \\
27\, 664 & \text{ for $i=3$} \\
365\, 750 & \text{ for $i=4$} \\
8\, 645  & \text{ for $i=5$} \\
133 & \text{ for $i=6$} \\
912 & \text{ for $i=7$.} \\
\end{cases}
\]
The natural module is $X_{\lambda_1}$, and
$X_{\lambda_6}=\mathrm{Ad}$.

The Galois group $\Gamma$ acts as
$\lambda_i=\bar\lambda_i$ for $i=1,\ldots,7$, and
\[
\beta(\lambda_i)=
\begin{cases}
\RR &\text{ for $i=2,4,5,6$} \\
\HH &\text{ for $i=1,3,7$.}
\end{cases}
\]
\end{Num}

\begin{Num}\textsc{Low-dimensional simple $\E_7$-modules}\psn
\label{E7Modules}%
The $\E_6$-modules of dimension at most 224 are the following.
\begin{center}
\begin{tabular}{llllll}
$G$ & weight & $X$ & $\dim(x)$ & $^\RR X$ & $\dim({}^\RR X)$ \\
\hline
$\E_7$ & $\lambda_1$ && 56  && 112  \\
$\mathrm{PE}_7$ & $\lambda_6$ & $\frak{e}_7\CC$ & 133 & $\frak{e}_7$ & 133
\end{tabular}
\end{center}
\qed
\end{Num}

\newpage

\begin{Num}\textsc{Fundamental weights for $\E_8$}\psn
We label the Dynkin diagram as follows.
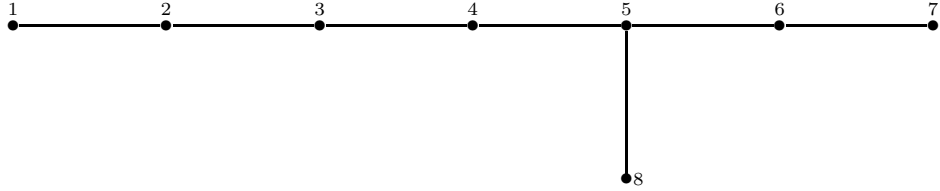
\begin{diagram}[size=2.9em,abut]
\overset{1}\bullet & \rLine &
\overset{2}\bullet & \rLine & 
\overset{3}\bullet & \rLine &
\overset{4}\bullet & \rLine &
\overset{5}\bullet & \rLine &
\overset{6}\bullet & \rLine &
\overset{7}\bullet  \\
&&&&&&&& \dLine \\
&&&&&&&& \bullet\rlap{$\scriptstyle 8$} \\
\end{diagram}
The group is simply connected and simple.
The dimensions of the fundamental representations are
\[
\dim(\rho_{\lambda_i})=
\begin{cases}
248 & \text{ for $i=1$} \\
30\, 380 & \text{ for $i=2$} \\
2\, 450\, 240 & \text{ for $i=3$} \\
146\, 325\, 270 & \text{ for $i=4$} \\
6\, 899\, 079\, 264  & \text{ for $i=5$} \\
6\, 696\, 000  & \text{ for $i=6$} \\
3\, 875 & \text{ for $i=7$} \\
147\, 250 & \text{ for $i=8$.} \\
\end{cases}
\]
The natural module is $X_{\lambda_1}=\mathrm{Ad}$.

The Galois group acts as
$\lambda_i=\bar\lambda_i$ for $i=1,\ldots,8$, and
$\beta(\lambda_i)=\RR$ for $i=1,\ldots,8$.
\end{Num}

\begin{Num}\textsc{Low-dimensional simple $\E_8$-modules}\psn
\label{E8Modules}%
The $\E_6$-modules of dimension at most 992 are the following.
\begin{center}
\begin{tabular}{llllll}
$G$ & weight & $X$ & $\dim(X)$ & $^\RR X$ & $\dim({}^\RR X)$ \\
\hline
$\mathrm{E}_8$ & $\lambda_1$ & $\frak{e}_8\CC$ & 248 & $\frak{e}_8$ & 248
\end{tabular}
\end{center}
\qed
\end{Num}

\newpage

\begin{Num}\textsc{Fundamental weights for $\Ffour$}\psn
We label the Dynkin diagram as follows.
\begin{diagram}[abut]
\overset{1}\bullet & \rLine &
\overset{2}\bullet & \lImplies & 
\overset{3}\bullet & \rLine &
\overset{4}\bullet \\
\end{diagram}
The group is simple and simply connected.
The dimensions of the fundamental representations are
\[
\dim(\rho_{\lambda_i})=
\begin{cases}
26 & \text{ for $i=1$} \\
273 & \text{ for $i=2$} \\
1\, 274 & \text{ for $i=3$} \\
52 & \text{ for $i=4$.}
\end{cases}
\]
The natural module is $X_{\lambda_1}$, and
$X_{\lambda_4}=\mathrm{Ad}$.
The Galois group $\Gamma$ acts as $\lambda_i=\bar\lambda_i$,
and $\beta(\lambda_i)=\RR$ for $i=1,\ldots,4$.
\end{Num}

\begin{Num}\textsc{Low-dimensional simple $\Ffour$-modules}\psn
\label{F4Modules}%
The $\Ffour$-modules of dimension at most 104 are the following.
\begin{center}
\begin{tabular}{llllll}
$G$ & weight & $X$ & $\dim(X)$ & $^\RR X$ & $\dim({}^\RR X)$ \\
\hline
$\Ffour$ & $\lambda_1$ &  & 26 &  & 26 \\
$\Ffour$ & $\lambda_4$ & $\frak{f}_4\CC$ & 52 & $\frak{f}_4$ & 52
\end{tabular}
\end{center}
\qed
\end{Num}

\newpage

\begin{Num}\textsc{Fundamental weights for $\G_2$}\psn
We label the Dynkin diagram as follows.
\newarrow{Drei}3333{>}
\begin{diagram}[abut]
\overset{1}\bullet & \lImplies & 
\overset{2}\bullet
\end{diagram}
\vspace{-28.3pt}
\begin{diagram}[abut]
\bullet  & \lLine & 
\bullet 
\end{diagram}
The group is simple and simply connected.
The dimensions of the fundamental representations are
\[
\dim(\rho_{\lambda_1})=7
\]
and
\[
\dim(\rho_{\lambda_2})=14.
\]
The natural module is $X_{\lambda_1}$, and
$X_{\lambda_2}=\mathrm{Ad}$.
The Galois group acts as $\lambda_i=\bar\lambda_i$ for $i=1,2$, and
\[
\beta(\lambda_i)=\RR
\]
for $i=1,2$.
\end{Num}

\begin{Num}\textsc{Low-dimensional simple $\G_2$-modules}\psn
\label{G2Modules}%
The $G_2$-modules of dimension at most $28$ are the following.
\begin{center}
\begin{tabular}{llllll}
$G$ & weight & $X$ & $\dim(X)$ & $^\RR X$ & $\dim({}^\RR X)$ \\
\hline
$\G_2$ & $\lambda_1$ & $\CC^7$ & 7  & $\RR^7$ & 7   \\
$\G_2$ & $\lambda_2$ & $\frak{g}_2\CC$ & 14 & $\frak{g}_2$ & 14 \\
$\G_2$ & $2\lambda_1$ &  & 27 &      & 27   \\
\end{tabular}
\end{center}

\qed
\end{Num}

\newpage

\begin{Num}\textsc{Certain subgroups of exceptional Lie groups}\psn
\label{SubgroupsInExceptionalGroups}%
Borel-De Siebenthal \cite{BorelSieb49}
determined the maximal connected subgroups of
maximal rank for all compact almost simple Lie groups.
Dynkin \cite{DynkinAMST57}
and Seitz \cite{Seitz91} determined the maximal connected subgroups
of the exceptional complex almost simple Lie groups
which have strictly smaller rank. 
The complexification is a functorial equivalence between
compact connected Lie groups and reductive algebraic $\CC$-groups.
Therefore, the results of Dynkin and Seitz
apply to compact almost simple Lie
groups. We obtain in particular the following 
result: The two standard inclusions
\[
\G_2\SUB\Ffour\text{ and }\Ffour\SUB\E_6
\]
are unique up to conjugation, cp.~Seitz \cite{Seitz91} Table 1,
p.~193.

The reasoning here is as follows.
Let $\G_2^\CC$ and $\Ffour^\CC$ denote the
corresponding complex Lie groups. An inclusion $\G_2\rInto\Ffour$
yields (via complexification) an inclusion
$\G_2^\CC\rInto^\phi\Ffour^\CC$.
By the result of Dynkin and Seitz, this inclusion is conjugate
to the standard inclusion $\G_2^\CC\rInto^s\Ffour^\CC$ by an
element $g\in\Ffour^\CC$. Since $\G_2\SUB\G_2^\CC$ is a maximal compact
subgroup, we can assume that $g$ transforms
$\phi(\G_2)$ into $s(\G_2)$ (all maximal compact subgroups in $\G_2^\CC$
are conjugate). Now write $g=kp$ as in \ref{PolarTrick}, with
$k\in\Ffour$. It follows by a similar argument as in \ref{PolarTrick}
that $\phi(\G_2)$ and $s(\G_2)$ are conjugate under $k\in\Ffour$.

The maximal connected subgroups of rank 2 in $\G_2$ are 
according to Borel-De Siebenthal \cite{BorelSieb49}
the groups $\SO(4)$ 
(the stabilizer of the quaternion algebra $\HH\SUB\OO$), and $\SU(3)$ 
(the elementwise stabilizer of $\CC\SUB\OO$). According to
Dynkin \cite{DynkinAMST57} and Seitz \cite{Seitz91}
the representation
$^\RR\rho_{6\lambda_1}$ lifts to a maximal connected subgroup of
$\G_2$.
The following table gives all subgroups $H$ of type $\fa_1$ 
in $\G_2$.
\begin{center}
\begin{tabular}{lll}
\multicolumn{3}{c}{\fbox{3-dimensional connected subgroups of $G_2$}} \\
\\
$H$ & $\rho$ & maximal subgroup \\
&& containing $H$ \\
\hline
$\SO(3)$ & $2\cdot{}^\RR\rho_{2\lambda_1}$ & $\SU(3)$ \\
$\SU(2)$ & $^\RR\rho_{\lambda_1}$  & $\SU(3)$ \\
$\SU(2)$ & $^\RR\rho_{\lambda_1}+{}^\RR\rho_{2\lambda_1}$  & $\SO(4)$ \\
$\SO(3)$ & $^\RR\rho_{6\lambda_1}$ & 
\end{tabular}
\end{center}
\end{Num}

%\newpage

%\newpage

\section{The Dynkin index}

Let $\phi:H\too G$ be a homomorphism between compact almost
simple Lie groups. There is a corresponding homomorphism
\[
\pi_3(\phi):\pi_3(H)\too\pi_3(G).
\]
Both groups are infinite cyclic, hence there exists a
number $j=j_\phi\in\NN$ such that the cokernel of
$\pi_3(\phi)$ is cyclic of order $j$. This number is called
the \emph{Dynkin index} of $\phi$. It is clear that the Dynkin
index is multiplicative,
\[
j_{\phi\circ\psi}=j_\phi\cdot j_\psi.
\]
Moreover, if $\rho$ and $\sigma$ are representations, then
\[
j_{\rho+\sigma}=j_\rho+j_\sigma,
\]
cp.~Onishchik \cite{Oni94} Ch.~5, Thm.~2, \S17.2, p.~257.
The following diagrams show the Dynkin indices of various
natural inclusions between compact almost simple Lie groups.
The groups on the right are the stable limits of the classical
groups.
{\small\begin{diagram}
\SO(3) & \rTo_2 & \SO(5) & \rTo_1 & \SO(6) & \rTo_1 & \SO(7) & \rTo_1 &
\SO(8) && \rTo_1 && \SO \\
&&&& \uTo_1 &&&& \uTo_1 &&&& \uTo_1\\
\SU(2) && \rTo_1 && \SU(3) && \rTo_1 && \SU(4) && \rTo_1 && \SU \\
\uTo_1 &&&& &&&& \uTo_1 &&&& \uTo_1 \\
\Sp(1) &&&& \rTo_1 &&&& \Sp(2) && \rTo_1 && \Sp \\
\end{diagram}}%
Here, one can fill in the following subdiagram.
\begin{diagram}
\SO(6) && \rTo_1 && \SO(7) && \rTo_1 && \SO(8) \\
&&& \ruTo_1 &&&& \ruTo_1 \\
\uTo_1 && \G_2 && \rTo_1 && \Spin(7) && \uTo_1 \\
& \ruTo_1 &&&&&& \luTo_1\\
\SU(3) &&&& \rTo_1 &&&& \SU(4) \\
\end{diagram}
Finally, tensoring with $\CC$ and $\HH$ yields a diagram
\begin{diagram}
\SO(3) && \rTo_2 && \SO(5) & \rTo_1 & \SO(6) && \rTo_1 &&
\SO \\
\dTo_4 &&      && \dTo_2 &      & \dTo_2 &&&& \dTo_2 \\
\SU(3) & \rTo_1 & \SU(4) & \rTo_1 & \SU(5) & \rTo_1 & \SU(6) && \rTo_1 
&& \SU \\
\dTo_2 &        & \dTo_2 &        & \dTo_2 &        & \dTo_2 &&&& \dTo_2 \\
\Sp(3) & \rTo_1 & \Sp(4) & \rTo_1 & \Sp(5) & \rTo_1 & \Sp(6) && \rTo_1 
&& \Sp \\
\end{diagram}
Besides this, we use the following result for representations of
$\SU(2)$. The Dynkin index of 
$\rho_{k\lambda_1}:\SU(2) \too \SU(k+1)$ is 
\[
j_{\rho_{k\lambda_1}}=\binom{k+2}{3},
\]
cp.~Onishchik \cite{Oni94} p. 61.

%%%%%%%%%%%%%%%%%%%%%%%%%%%%%%%%%%%%%%%%%%%%%%%%%%%%%%%%%%%%%%%%%%%%%%%%
%                                                                      %
%                                                                      %
%                   Compact homogeneous quadrangles                    %
%                                                                      %
%                          Linus Kramer                                %
%                                                                      %
%                           Memoirs AMS                                %
%                                                                      %
%                         Wuerzburg 2000                               %
%                                                                      %
%                                                                      %
%                            CHQ5.tex                                  %
%                                                                      %
%                                                                      %
%                                                                      %
%%%%%%%%%%%%%%%%%%%%%%%%%%%%%%%%%%%%%%%%%%%%%%%%%%%%%%%%%%%%%%%%%%%%%%%%
\chapter{The case when $G$ is simple}

In this chapter we classify all homogeneous spaces $X=G/H$ with the
following properties.
\begin{enumerate}
\item $G/H$ is compact and 1-connected.
\item $G$ is effective, compact and almost simple.
\item $G/H$ has the same integral cohomology as a product
$\SS^{n_1}\times\SS^{n_2}$, where $3\leq n_1\leq n_2$, and $n_2$ is odd.
\end{enumerate}
There are two cases: either $n_1$ is odd, or $n_1$ is even.

\noindent
\textbf{Case (I):} $n_1$ is odd, $\bfH^\bullet(G/H)\isom\EA_\ZZ(u,v)$, with
$\deg(u)=n_1$ and $\deg(v)=n_2$. In this situation Theorem
\ref{Case(I)General} applies.

\noindent
\textbf{Case (II):} $n_1$ is even, $\bfH^\bullet(G/H)\isom\ZZ[a]/(a^2)\otimes
\EA_\ZZ(u)$, with
$\deg(a)=n_1$ and $\deg(u)=n_2$. In this situation Theorem 
\ref{Case(II)General} applies.

These two cases are considered in the next two sections. Each
section ends with a complete list of all possible pairs $(G,H)$,
see \ref{Case(I)Simple}, \ref{Case(II)Simple}.

\section{Case (I): $\bfH^\bullet(X)=\EA_\ZZ(u,v)$.}

If $n_1$ is odd, then $\bfH^\bullet(X;\QQ)\isom\EA_\QQ(u,v)$.
By Theorem \ref{Case(I)General} $\rk(G)-\rk(H)=2$. Moreover,
$P_H^1\OOT\isom\P_G^1=0$, and
$P^3_H\oot P^3_G$ is a surjection; thus,
$H$ is either trivial or almost simple.
Let $(m_1^G,\ldots,m_n^G)$ denote the degrees of the homogeneous
generators of $G$, and let $(m_1^H,\ldots,m_{n-2}^H)$ denote the
corresponding numbers for $H$. After a suitable permutation of the
indices, $m_i^G=m_i^H$ for $i=1,\ldots,n-2$, cp.~\ref{HowTo(I)}. 
We check the various possibilities.

%\newpage

\begin{Num}\textsc{$G$ of type $\fa_n$}\psn
The homogeneous
generators of the rational cohomology of $\SU(n+1)$ have
degrees 
\[
(3,5,\ldots,2n+1).
\]
If $n\geq8$, then the only
type of simple group which fits in is $\fa_{n-2}$; in low dimensions,
we have also to consider the pairs
$(\fa_7,\fd_5)$,
$(\fa_5,\fb_3)$,
$(\fa_5,\fc_3)$, and
$(\fa_4,\fc_2)$.

\begin{description}
\item[\fbox{$(\fa_n,\fa_{n-2})$}]
For $n\geq4$, it follows from \ref{AnModules} that
$(G,H)=(\SU(n+1),\SU(n-1))$.

Suppose that $n=3$. Then $H$ is of type $\fa_1$.
By \ref{A1Modules}, the representations of $\SU(2)$ on $\CC^4$ are
$\SU(2)$ (i.e.~$\rho_{\lambda_1}$),
$\SO(3)$ (i.e.~$\rho_{2\lambda_1}$),
$2\rho_{\lambda_1}$, and $\rho_{3\lambda_1}$. 
Note that $\SU(4)\cong\Spin(6)$;
the subgroup 
$2\cdot\rho_{\lambda_1}(\SU(2))\SUB\SU(4)$ contains the center, 
and 
$\SU(4)/2\cdot\rho_{\lambda_1}(\SU(2))=\SO(6)/\SO(3)=V_3(\RR^6)$.

\item[\fbox{$(\fa_7,\fd_5)$ and $(\fa_5,\fb_3)$}]
The group $\Spin(10)$ has no non-trivial representation
on $\CC^8$; similarly, there is no non-trivial representation of
$\Spin(7)$ on $\CC^6$, cp.~\ref{BnModules}. This shows that the pairs
$(\fa_7,\fd_5)$ and $(\fa_5,\fb_3)$ do not exist.

\item[\fbox{$(\fa_5,\fc_3)$ and $(\fa_4,\fc_2)$}]
There is a (unique) representation of
$\Sp(3)$ on $\CC^6=\HH^3$, the natural one, and
two representations of $\Sp(2)$ on $\CC^5$, arising from
the inclusions $\SO(5)\SUB\SU(5)$ and $\Sp(2)\SUB\SU(4)\SUB\SU(5)$,
cp.~\ref{BnModules}.
\end{description}
We obtain the following list of groups and homogeneous spaces.
The structure of $\pi_3$ follows from the formula for the Dynkin index of
the various representations of $\SU(2)$.

\begin{center}
\begin{tabular}{llllll} 
\multicolumn{6}{c}{\fbox{$G$ of type $\SU(n+1)$}} \\ 
\\ 
%\hline
$G$             & $H$        & $\Cen_G(H)^\circ$ & $(n_1,n_2)$ & $G/H$ &
Remarks \\ 
\hline
$\SU(n+1)$      & $\SU(n-1)$ & $\U(2)$   & $(2n-1,2n+1)$  & $V_2(\CC^{n+1})$ 
& \\
\hfill$n\geq 3$  \\
\hline
%$\SU(6)$        & $\SU(4)$   & $\U(2)$   & $(9,11)$ & $V_2(\CC^6)$     &  \\
$\SU(6)$        & $\Sp(3)$   & 1         & $(5,9)$  &                  &  
(see below) \\ 
%$\SU(5)$        & $\SU(3)$   & $\U(2)$   & $(7,9)$  & $V_2(\CC^5)$     &  \\
$\SU(5)$        & $\Sp(2)$   & $\U(1)$   & $(5,9)$  &                  &  
(see below) \\
$\SU(5)$        & $\SO(5)$   & 1         & $(5,9)$  &                  &  
$\pi_3=\ZZ/2$ \\
%$\SU(4)$        & $\SU(2)$   & $\U(2)$   & $(5,7)$  & $V_2(\CC^4)$     &  \\
$\SU(4)$        & $\SO(3)$   & $\U(1)$   & $(5,7)$  &                  & 
$\pi_3=\ZZ/4$ \\
$\SU(4)$        & $\rho_{3\lambda_1}$   &  1        & $(5,7)$  &                  &
$\pi_3=\ZZ/10$ \\
$\SO(6)$        & $\SO(3)$   & $\SO(3)$  & $(5,7)$  & $V_3(\RR^6)$     &  
$\pi_3=\ZZ/2$ \\
$\SU(3)$        & 1          & $\SU(3)$  & $(3,5)$  & $V_2(\CC^3)$     &  \\
\end{tabular}
\end{center}
There is one coincidence in this table, arising from the inclusion
$\SU(5)\SUB\SU(6)$: 
\[
\SU(5)/\Sp(2)=\SU(6)/\Sp(3)
\]
\end{Num}

%\newpage

\begin{Num}\textsc{$G$ of type $\fb_n$}\psn
\label{Spin11}%
The homogeneous
generators of the rational cohomology of $\Spin(2n+1)$ have
degrees 
\[
(3,7,\ldots,4n-1).
\]
If $n\geq10$, then the only
types of simple groups which fit in are $\fb_{n-2}$ and
$\fc_{n-2}$. However, $\RR^{2n+1}$ is not a non-trivial
$\Sp(n-2)$-module for $n\geq 5$, so $(\fb_n,\fc_{n-2})$ is
excluded in this range, and $\fb_2\cong\fc_2$ anyway. 
In low dimensions, we have also to consider the pairs
$(\fb_9,\fe_7)$,
$(\fb_6,\ff_4)$, and
$(\fb_4,\fg_2)$.

\begin{description}
\item[\fbox{$(\fb_n,\fb_{n-2})$}]
If $n\geq6$, then $(G,H)=(\SO(2n+1),\SO(2n-3))$, cp.
\ref{BnModules}.

The group $\Spin(7)$ has two non-trivial representations
on $\RR^{11}$, corresponding to the standard inclusion 
$\SO(7)\SUB\SO(9)$,
and to the inclusion $\Spin(7)\SUB\SO(8)$, cp.~\ref{BnModules}.

Similarly, $\Spin(5)$ has two non-trivial representations
on $\RR^9$, one arising from the inclusion $\SO(5)\SUB\SO(9)$, 
and one arising from $\Sp(2)\SUB\SO(8)\SUB\SO(9)$.

The representations of $\SU(2)$ on $\RR^7$ are 
$\SO(3)$ (i.e.~$^\RR\rho_{2\lambda_1}$), $\SU(2)$
\linebreak
(i.e.~$^\RR\rho_{\lambda_1}$),
$^\RR\rho_{\lambda_1}+{}^\RR\rho_{2\lambda_1}$, $^\RR\rho_{4\lambda_1}$, 
and $^\RR\rho_{6\lambda_1}$. We denote the corresponding
homomorphisms into $\Spin(7)$ by a tilde,
\begin{diagram}
 &&&& \Spin(7) \\
 &&& \ruTo(4,2)^{\widetilde\rho} & \dTo \\
\SU(2) &&\rTo^\rho && \SO(7)\rlap{.}
\end{diagram}

\item[\fbox{$(\fb_9,\fe_7)$ and $(\fb_6,\ff_4)$}]
The compact groups $\E_7$ and $\Ffour$  do not have non-trivial 
19- and 13-dimensional representations, respectively, cp.
\ref{F4Modules} and \ref{E7Modules},
thus the pairs $(\fb_9,\fe_7)$ and $(\fb_6,\ff_4)$ do not exist.

\item[\fbox{$(\fb_4,\fg_2)$}]
There is only one 8-dimensional representation of $\G_2$ on
$\RR^8$, given by the action of $\G_2$ on the Cayley algebra
$\OO$, cp.~\ref{G2Modules}.
\end{description}

%\newpage
\begin{center}
{\small
\begin{tabular}{llllll} 
\multicolumn{6}{c}{\fbox{$G$ of type $\SO(2n+1)$}} \\ 
\\
$G$          & $H$           & $\Cen_G(H)^\circ$ & $(n_1,n_2)$ & $X=G/H$ &
Remarks \\ 
\hline
$\SO(2n+1)$  & $\SO(2n-3)$   & $\SO(4)$  & $(4n-5,\ $ %4n-1)$ 
& $V_4(\RR^{2n+1})$ 
& $\pi_{2n-3}=\ZZ/2$ \\
\hfill $n\geq 3$   &&& $\quad4n-1)$  \\ 
\hline
%$\SO(11)$    & $\SO(7)$      & $\SO(4)$  & $(15,19)$  & $V_4(\RR^{11})$ & 
%$\pi_7=\ZZ/2$ \\
$\Spin(11)$  & $\Spin(7)$    &  $\Sp(1)$  & $(15,19)$  &                 &
$\pi_9\neq 0$  \\
%$\SO(9)$     & $\SO(5)$      & $\SO(4)$  & $(11,15)$  & $V_4(\RR^9)$    & 
%$\pi_5=\ZZ/2$  \\
$\Spin(9)$   & $\Sp(2)$      &  $\Sp(1)$ & $(11,15)$  &                 & 
$\pi_5\neq 0$ \\
$\Spin(9)$   & $\G_2$        & $\SO(2)$  & $(7,15)$   & $V_2(\OO^2)$    & \\
%$\SO(7)$     & $\SO(3)$      & $\SO(4)$  & $(7,11)$   & $V_4(\RR^7)$    & 
%$\pi_3=\ZZ/2$ \\
$\Spin(7)$   & $\SU(2)$  & $\SU(2)\times\SO(3)$ & $(7,11)$   & $V_3(\RR^8)$    & 
$\pi_5=\ZZ/2$ \\
$\Spin(7)$   & 
$\widetilde{^\RR\rho_{\lambda_1}+{}^\RR\rho_{2\lambda_1}}$ 
& 1 & $(7,11)$ &                 & 
$\pi_3=\ZZ/3$ \\
$\Spin(7)$   & 
$\widetilde{2\cdot{}^\RR\rho_{2\lambda_1}}$ 
& $\SO(2)$ & $(7,11)$ &                 & 
$\pi_3=\ZZ/4$ \\
$\Spin(7)$   &
$\widetilde{^\RR\rho_{4\lambda_1}}$ 
& $\SO(2)$ & $(7,11)$   &                 & 
$\pi_3=\ZZ/10$ \\
$\Spin(7)$   & 
$\widetilde{^\RR\rho_{6\lambda_1}}$ 
& 1 & $(7,11)$   &                 & 
$\pi_3=\ZZ/28$\\
$\Sp(2)$     & 1              & $\Sp(2)$ & $(3,7)$    & $V_2(\HH^2)$    & \\
\end{tabular}}
\end{center}
The subgroup $\Spin(7)\SUB\SO(8)$ acts transitively on $V_3(\RR^8)$,
hence
\[
\Spin(7)/\SU(2)=\SO(8)/\SO(5).
\]
From the collapsing of the spectral sequence of
\begin{diagram}
\SS^{15}=&\Spin(9)/\Spin(7) & \rTo & \Spin(11)/\Spin(7) \\
&&& \dTo \\
&&& \Spin(11)/\Spin(9) & =V_2(\RR^{11})
\end{diagram}
one sees that 
\begin{align*}
&\bfH^\bullet(\Spin(11)/\Spin(7);\ZZ/2)=
\EA_{\ZZ/2}(x_9,x_{15})\otimes(\ZZ/2)[x_{10}]/(x_{10}^2) \\
\neq &
\bfH^\bullet(V_4(\RR^{11});\ZZ/2)=
\EA_{\ZZ/2}(x_7,x_9)\otimes(\ZZ/2)[x_8,x_{10}]/(x_8^2,x_{10}^2).
\end{align*}
Similarly, the spectral sequence of 
\begin{diagram}
V_2(\RR^7)=&\Spin(7)/\Sp(2) & \rTo & \Spin(9)/\Sp(2) \\
&&& \dTo \\
&&& \Spin(9)/\Spin(7) & =  \SS^{15}
\end{diagram}
collapses, and
\begin{align*}
&\bfH^\bullet(\Spin(9)/\Sp(2);\ZZ/2)=
\EA_{\ZZ/2}(x_5,x_{15})\otimes(\ZZ/2)[x_6]/(x_6^2) \\
\neq &
\bfH^\bullet(V_4(\RR^9);\ZZ/2)=
\EA_{\ZZ/2}(x_5,x_7)\otimes(\ZZ/2)[x_6,x_8]/(x_6^2,x_8^2).
\end{align*}
\end{Num}

%\newpage

\begin{Num}\textsc{$G$ of type $\fc_n$}\psn
The homogeneous
generators of the rational cohomology of $\Sp(n)$ have
degrees 
\[
(3,7,\ldots,4n-1).
\]
If $n\geq10$, then the only
types of simple groups which fit in are $\fb_{n-2}$ and
$\fc_{n-2}$. However, $\HH^n$ is not a non-trivial
$\Spin(2n-3)$-module for $n\geq 5$, cp.~\ref{BnModules}, 
so $(\fc_n,\fb_{n-2})$ is
excluded in this range, and $\fb_2\cong\fc_2$ anyway. 
In low dimensions, we have also to consider the pairs
$(\fc_9,\fe_7)$,
$(\fc_6,\ff_4)$, and
$(\fc_4,\fg_2)$.

\begin{description}
\item[\fbox{$(\fc_n,\fc_{n-2})$}]
If $n\geq4$, then $(G,H)=\Sp(n),\Sp(n-2))$, cp.~\ref{CnModules}.

The group $\Sp(1)$ has several non-trivial representation
on $\HH^3$, corresponding to the inclusions
$\Sp(1)\SUB\Sp(3)$ (i.e.~$^\HH\rho_{\lambda_1}$), 
$\SU(2)\SUB\Sp(2)\SUB\Sp(3)$
(i.e.~$2\cdot{}^\HH\rho_{\lambda_1}$), 
$\SO(3)\SUB\Sp(3)$ (i.e.~$^\HH\rho_{2\lambda_1}$), 
and
$3\cdot{}^\HH\rho_{\lambda_1}$, 
$^\HH\rho_{3\lambda_1}$, 
$^\HH\rho_{3\lambda_1}+{}^\HH\rho_{\lambda_1}$, 
$^\HH\rho_{5\lambda_1}$, cp.~\ref{A1Modules}.

\item[\fbox{$(\fc_9,\fe_7)$, $(\fc_6,\ff_4)$, and $(\fc_4,\fg_2)$}]
The compact groups $\E_7$, $\Ffour$, and $\G_2$  do not have non-trivial 
representations on $\HH^9$, $\HH^6$, and $\HH^4$, respectively,
thus the pairs $(\fc_9,\fe_7)$, $(\fc_6,\ff_4)$, and $(\fc_4,\fg_2)$ 
do not exist, cp.~\ref{E7Modules}, \ref{F4Modules},  \ref{G2Modules}.
\end{description}

\begin{center}
\begin{tabular}{llllll} 
\multicolumn{6}{c}{\fbox{$G$ of type $\Sp(n)$}} \\ 
\\
$G$        & $H$           & $\Cen_G(H)^\circ$ & $(n_1,n_2)$ & $G/H$ & 
Remarks\\
\hline
$\Sp(n)$   & $\Sp(n-2)$    & $\Sp(2)$ & $(4n-5,4n-1)$ & $V_2(\HH^n)$ \\
\hfill $n\geq 3$ \\
\hline
%$\Sp(3)$   & $\Sp(1)$      & $\Sp(2)$ & $(7,11)$ & $V_2(\HH^3)$  \\
$\Sp(3)$   & $\SO(3)$      & 1        & $(7,11)$ & & $\pi_3=\ZZ/8$  \\
$\Sp(3)$   & $2\cdot{}^\HH\rho_{\lambda_1}$  
& $\Sp(1)$ & $(7,11)$      &  & $\pi_3=\ZZ/2$  \\
$\Sp(3)$   & $3\cdot{}^\HH\rho_{\lambda_1}$  
& 1        & $(7,11)$      &  & $\pi_3=\ZZ/3$  \\
$\Sp(3)$   & $^\HH\rho_{3\lambda_1}$  
& $\Sp(1)$ & $(7,11)$      &  & $\pi_3=\ZZ/10$  \\
$\Sp(3)$   & $^\HH\rho_{3\lambda_1}+{}^\HH\rho_{\lambda_1}$  
& 1        & $(7,11)$      &  & $\pi_3=\ZZ/11$  \\
$\Sp(3)$   & $^\HH\rho_{5\lambda_1}$   
& 1        & $(7,11)$      &  & $\pi_3=\ZZ/35$  \\
$\Sp(2)$   & 1             & $\Sp(2)$ & $(3,7)$  & $V_2(\HH^2)$   \\
\end{tabular}
\end{center}
\end{Num}

%\newpage

\begin{Num}\textsc{$G$ of type $\fd_n$}\psn
\label{Spin10}%
The homogeneous
generators of the rational cohomology of $\Spin(2n)$ have
degrees 
\[
(3,7,\ldots,4n-5,2n-1).
\]
If $n\geq5$, then the only
types of simple groups which fit in are $\fb_{n-2}$ and
$\fc_{n-2}$. However, $\RR^{2n}$ is not a non-trivial
$\Sp(n-2)$-module for $n\geq 4$, cp.~\ref{CnModules},
so $(\fd_n,\fc_{n-2})$ is
excluded in this range, and $\fb_2\cong\fc_2$. 
We have also to consider the pair $(\fd_4,\fg_2)$.

\begin{description}
\item[\fbox{$(\fd_n,\fb_{n-2})$}]
If $n\geq6$, then $(G,H)=(\Spin(2n),\Spin(2n-3))$, cp.~\ref{BnModules}.

The two representations of $\Spin(7)$ on $\RR^{10}$ yield
two embeddings \linebreak
$\Spin(7)\SUB\Spin(10)$. Similarly,
the inclusions $\SO(5)\SUB\SO(8)$ and $\Sp(2)\SUB\SO(8)$
yield two inclusions, cp.~\ref{BnModules}. 
The case of $\fd_3\cong\fa_3$ was already considered above.

\item[\fbox{$(\fd_4,\fg_2)$}]
There is a unique representation of $\G_2$ on $\RR^8$, cp.
\ref{G2Modules}.
\end{description}

\begin{center}
\begin{tabular}{llllll} 
\multicolumn{6}{c}{\fbox{$G$ of type $\SO(2n)$}} \\ 
\\
$G$          & $H$         & $\Cen_G(H)^\circ$ & $(n_1,n_2)$ & $G/H$ & 
Remarks \\ 
\hline
$\SO(2n)$    & $\SO(2n-3)$ & $\SO(3)$ & $(2n-1,4n-5)$ & $V_3(\RR^{2n})$ &
$\pi_{2n-3}=\ZZ/2$ \\
\hfill$n\geq 3$ \\
\hline
%$\SO(10)$    & $\SO(7)$    & $\SO(3)$ & $(9,15)$      & $V_3(\RR^{10})$ &
%$\pi_7=\ZZ/2$ \\
$\Spin(10)$  & $\Spin(7)$  & $\SO(2)$ & $(9,15)$      &  &  \\
%$\SO(8)$     & $\SO(5)$    & $\SO(3)$ & $(7,11)$      & $V_3(\RR^8)$    & 
%$\pi_3=\ZZ/2$ \\
$\Spin(8)$   & $\G_2$      & 1        & $(7,7)$       & $\SS^7\times\SS^7$ & 
\\
$\SU(4)$     & $\SU(2)$    & $\U(2)$  & $(5,7)$       & $V_2(\CC^4)$       & 
\\
$\SU(4)$     & $\SO(3)$    & $\U(1)$  & $(5,7)$       &                    &
$\pi_3=\ZZ/4$ \\
$\SU(4)$     & $\rho_{3\lambda_1}$    & 1             & $(5,7)$       &                    &
$\pi_3=\ZZ/10$ \\
%$\SO(6)$     & $\SO(3)$    & $\SO(3)$ &               & $V_3(\RR^6)$       &
%$\pi_3=\ZZ/2$ \\
\end{tabular}
\end{center}
The spectral sequence of the bundle
\begin{diagram}
\SS^{15}=&\Spin(9)/\Spin(7) & \rTo & \Spin(10)/\Spin(7) \\
&&& \dTo \\
&&& \Spin(10)/\Spin(9) & =\SS^9
\end{diagram}
collapses, and thus
$\bfH^\bullet(\Spin(10)/\Spin(7))=\EA_\ZZ(x_9,x_{15})$.
In particular, 
\begin{align*}
\bfH^\bullet(\Spin(10)/\Spin(7);\ZZ/2) & \neq
\bfH^\bullet(V_3(\RR^{10});\ZZ/2) \\
&=\EA_{\ZZ/2}(x_7,x_9)\otimes(\ZZ/2)[x_8]/(x_8^2) 
\end{align*}
\end{Num}

%\newpage

\begin{Num}\textsc{$G$ exceptional}\psn
Finally, we have the exceptional groups. Here, the only pairs
are $(\fe_6,\ff_4)$, $(\ff_4,\fg_2)$, and $(\fg_2,0)$.

\begin{description}
\item[\fbox{$(\fe_6,\ff_4)$}]
By \ref{SubgroupsInExceptionalGroups}, 
there is a unique inclusion $\Ffour\SUB\E_6$; the resulting space 
is Riemannian symmetric.

\item[\fbox{$(\ff_4,\fg_2)$}]
Again by \ref{SubgroupsInExceptionalGroups}, 
there is a unique inclusion $\G_2\SUB\Ffour$.

\item[\fbox{$(\fg_2,0)$}]
\end{description}

\begin{center}
\begin{tabular}{lllll} 
\multicolumn{5}{c}{\fbox{$G$ of exceptional type}} \\ 
\\
$G$          & $H$         & $\Cen_G(H)^\circ$ & $(n_1,n_2)$ &  
Remarks \\ 
\hline
$\E_6$       & $\Ffour$    &        & $(9,17)$  & \\
$\Ffour$     & $\G_2$      &        & $(15,23)$ & $\pi_7=\ZZ/3$\\
$\G_2$       & 1           &        & $(3,11)$  & 2-torsion
\end{tabular}
\end{center}
The $\ZZ/2$-Leray-Serre spectral sequence of the bundle
\begin{diagram}
\SU(2) & \rTo & \G_2 \\
&& \dTo \\
&& \G_2/\SU(2) & = V_2(\RR^7)
\end{diagram}
collapses, and thus
\[
\bfH^\bullet(\G_2;\ZZ/2)\isom\EA_{\ZZ/2}(x_3,x_5)\otimes
(\ZZ/2)[x_6]/(x_6^2).
\]
The inclusion $\G_2\SUB\Ffour$ is analyzed in Borel \cite{Bor54}. 
It is proved in particular, that there is a diagram
\begin{diagram}
\bfH^\bullet(\G_2;\ZZ/3) & \lTo & \bfH^\bullet(\Ffour;\ZZ/3) \\
\dEq && \dEq \\
\EA_{\ZZ/3}(x_3,x_{11}) & \lTo & \EA_{\ZZ/3}(x_3,x_7,x_{11},x_{15})
\otimes(\ZZ/3)[x_8]/(x_8^3)\rlap{,}
\end{diagram}
cp.~\emph{loc.~cit.} 19--23.
Therefore the spectral sequence of the bundle
\begin{diagram}
\G_2 & \rTo & \Ffour \\
&& \dTo \\
&& \Ffour/\G_2
\end{diagram}
collapses, and 
\[
\bfH^\bullet(\Ffour/\G_2;\ZZ/3)\isom
(\ZZ/3)[x_8]/(x_8^3)\otimes\EA_{\ZZ/3}(x_7,x_{15}).
\]
We use this fact later.
\end{Num}

%\newpage
We have classified compact 1-connected homogeneous spaces
whose  rational cohomology is the same as a product of odd-dimensional
spheres. However, the spaces that we are interested in have the
same \emph{integral} cohomology as a product of spheres. 
The following theorem is the main result of this section.

\begin{Thm}
\label{Case(I)Simple}
Let $X=G/H$ be a 1-connected homogeneous space of an effective almost simple
compact Lie group $G$. If $X$ has the same integral cohomology as
a product of odd-dimensional spheres,
\[
\bfH^\bullet(X)\isom\EA_\ZZ(u,v),
\]
$\deg(u)=n_1\geq 3$ odd, $\deg(v)=n_2\geq n_1$ odd, 
then $(G,H)$ is one of the following pairs.

\begin{center}
\begin{tabular}{llllll} 
$G$             & $H$        & $\Cen_G(H)^\circ$ & $(n_1,n_2)$ & $G/H$  \\ 
\hline
$\SU(n+1)$      & $\SU(n-1)$ & $\U(2)$   & $(2n-1,2n+1)$  & $V_2(\CC^{n+1})$ 
\\
\hfill$n\geq 3$  \\
\hline
$\Sp(n)$   & $\Sp(n-2)$    & $\Sp(2)$ & $(4n-5,4n-1)$ & $V_2(\HH^n)$ \\
\hfill $n\geq 3$ \\
\hline
$\E_6$       & $\Ffour$    & $1$      & $(9,17)$      &  \\
$\Spin(10)$  & $\Spin(7)$  & $\SO(2)$ & $(9,15)$      &  \\
$\Spin(9)$   & $\G_2$      & $\SO(2)$ & $(7,15)$      & $V_2(\OO^2)$     \\
$\Spin(8)$   & $\G_2$      & $1$      & $(7,7)$       & $\SS^7\times\SS^7$ 
\\
$\SU(6)$     & $\Sp(3)$    & $1$      & $(5,9)$       &                 \\ 
$\SU(5)$     & $\Sp(2)$    & $\U(1)$  & $(5,9)$       &                 \\
$\SU(3)$     & $1$         & $\SU(3)$ & $(3,5)$       & $V_2(\CC^3)$    \\
$\Sp(2)$     & $1$         & $\Sp(2)$ & $(3,7)$       & $V_2(\HH^2)$    \\
\end{tabular}
\end{center}
All spaces in this list have distinct homotopy types, except for one
coincidence:
the subgroup $\SU(5)\SUB\SU(6)$ acts transitively on $\SU(6)/\Sp(3)$, and 
thus
\[
\SU(5)/\Sp(2)=\SU(6)/\Sp(3).
\]
\begin{proof}
This follows from the tables above. The space $X$ has to be 
$(n_1-1)$-connected and torsion-free.
\end{proof}
\end{Thm}

%\newpage

\section{Case (II): $\bfH^\bullet(X)=\ZZ[a]/(a^2)\otimes\EA_\ZZ(w)$.}

Now we assume that $G/H$ has the same cohomology as a product
$\SS^{n_1}\times\SS^{n_2}$, where $n_1\geq 4$ is even and $n_2>n_1$ is
odd. Thus
\[
\bfH^\bullet(X)\isom\ZZ[a]/(a^2)\otimes\EA_\ZZ(w),
\]
where $\deg(a)=n_1$ and $\deg(w)=n_2$.
If $n_1\geq 6$, then there is an isomorphism
\[
P^3_H\oot P^3_G
\]
and therefore $H$ is almost simple. 
\begin{Lem}
If $n_1=4$, then $G$ is of classical type, with $3\leq \rk(G)\leq 5$.
The group $H$ is semisimple with two almost simple factors of type
$\fa_1,\fa_2,\fb_2$, or $\fg_2$, respectively 
(we show below that a factor $\fg_2$ is in fact not possible).

\begin{proof}
The cokernel of $P^3_H\oot P^3_G$ is 1-dimensional, and
$P^1_H\oot\P^1_G=0$ is an isomorphism. Thus $H$ is semisimple with
two almost simple factors $H_1,H_2$. 
Next, we note that $\rk(\pi_7(X))\in\{1,2\}$. From the
exact sequence
\[
0\too\pi_7(H)\otimes\QQ\too\pi_7(G)\otimes\QQ\too\pi_7(X)\otimes\QQ\too0
\]
we see that 
\[
\rk(\pi_7(G))\geq 1.
\]
Moreover, $\rk(G)=\rk(H)+1\geq 3$. A glance at the table
in \ref{TypeTable} shows that $G$ is of classical type.

We want to show that $\rk(G)\leq 5$.
Assume otherwise. Then $\rk(\pi_7(G))=1$ (because $G$ is of classical type), 
hence $\rk(\pi_7(H))=0$. By the table in
\ref{TypeTable}, the groups $H_1$ and $H_2$ have to be
of type $\fa_1,\fa_2,\fg_2$, or $\fe_6,\fe_7,\fe_8,\ff_4$.

We exclude the four large exceptional groups as follows.
Let $\fh_1$ be one of these four exceptional Lie algebras, and let
$\fg$ be any classical compact
simple Lie algebra. If $\fh_1\SUB\fg$, then 
$\rk(\fg)>2\rk(\fh_1)+1$, 
cp.~\ref{E6Modules}, \ref{E7Modules}, \ref{E8Modules}, \ref{F4Modules}. 
In our situation, $\rk(G)-\rk(H)=1$,
hence $2\rk(H_1)\geq\rk(G)-1$ or $2\rk(H_2)\geq\rk(G)-1$.
This is not possible for these four types of groups.
Therefore $\rk(H)\leq 4$ and thus $3\leq \rk(G)\leq 5$.
\end{proof}
\end{Lem}
A direct computation shows that numerically, the only possibilities are
$(\fa_3,\fa_1\oplus\fa_1)$,
$(\fb_3,\fa_1\oplus\fa_1)$,
$(\fc_3,\fa_1\oplus\fa_1)$,
$(\fa_4,\fa_1\oplus\fa_2)$,
$(\fb_4,\fa_1\oplus\fg_2)$,
$(\fc_4,\fa_1\oplus\fg_2)$,
$(\fd_4,\fa_1\oplus\fb_2)$,
$(\fd_4,\fa_1\oplus\fg_2)$,
$(\fa_5,\fa_2\oplus\fg_2)$, and
$(\fd_5,\fg_2\oplus\fg_2)$.

\begin{Lem}
The pairs
$(\fb_4,\fa_1\oplus\fg_2)$,
$(\fc_4,\fa_1\oplus\fg_2)$,
$(\fd_4,\fa_1\oplus\fg_2)$,
$(\fa_5,\fa_2\oplus\fg_2)$, and
$(\fd_5,\fg_2\oplus\fg_2)$ do not exist.

\begin{proof}
There is only one non-trivial representation of $\G_2$ on
$\RR^k$, $k<14$, the standard one on $\RR^7$, cp.~\ref{G2Modules}. 
Its centralizer
in $\SO(9)$ is $\SO(2)$, and this excludes $(\fb_4,\fa_1\oplus\fg_2)$
and $(\fd_4,\fa_1\oplus\fg_2)$.
Moreover, $\fg_2$ is not contained in $\fc_4$ and $\fa_5$.
Finally, the centralizer of $\G_2$ in $\SO(10)$ is $\SO(3)$, and that
excludes $(\fd_5,\fg_2\oplus\fg_2)$.
\end{proof}
\end{Lem}
The following pairs remain in the case $n_1=4$.
\begin{gather*}
(\fa_3,\fa_1\oplus\fa_1) \\
(\fa_4,\fa_1\oplus\fa_2) \\
(\fb_3,\fa_1\oplus\fa_1) \\
(\fc_3,\fa_1\oplus\fa_1) \\
(\fd_4,\fa_1\oplus\fb_2) 
\end{gather*}

%\newpage

\begin{Num}\textsc{$G$ of type $\fa_n$}\psn
Here, $H$ cannot be almost simple.
The only possibilities are $(\fa_3,\fa_1\oplus\fa_1)$ and
$(\fa_4,\fa_1\oplus\fa_2)$.

\begin{description}
\item[\fbox{$(\fa_3,\fa_1\oplus\fa_1)$}]
The only subgroup of type $\fa_1$ in $\SU(4)$ with large centralizer
is the standard inclusion $\SU(2)\SUB\SU(4)$, cp \ref{A1Modules}. 
Therefore, there is only  the standard inclusion $\SU(2)\times\SU(2)\SUB
\SU(4)$. Note that $\SU(4)/\SU(2)\times\SU(2)=\SO(6)/\SO(4)$.
\item[\fbox{$(\fa_4,\fa_1\oplus\fa_2)$}]
By \ref{AnModules}, we have the standard embedding $\SU(3)\SUB\SU(5)$, 
therefore we obtain
the standard inclusion $\SU(2)\times\SU(3)\SUB\SU(5)$ as the only 
possibility. 
\end{description}

\begin{center}
\begin{tabular}{lllll} 
\multicolumn{5}{c}{\fbox{$G$ of type $\SU(n)$}} \\ 
\\
$G$               & $H$ & $\Cen_G(H)^\circ$ & $G/H$ & $(n_1,n_2)$ \\ 
\hline
$\SU(5)$   & $\SU(3)\times\SU(2)$ & $\U(1)$  & $\widetilde G_3(\CC^5)$ & 
$(4,9)$   \\
$\SO(6)$   & $\SO(4)$             & $\SO(2)$ & $V_2(\RR^6)$ & $(4,5)$  
\\
\end{tabular}
\end{center}
\end{Num}

%\newpage

\begin{Num}\textsc{$G$ of type $\fb_n$}\psn
The only possibilities are $(\fb_n,\fd_{n-1})$, for $n\geq 4$,
$(\fb_3,\fa_2)$, and $(\fb_3,\fa_1\oplus\fa_1)$.

\begin{description}
\item[\fbox{$(\fb_n,\fd_{n-1})$}]
If $n\geq5$, then $(G,H)=(\SO(2n+1),\SO(2n-2))$, cp.
\ref{DnModules}.
For $n=4$ we have the inclusions $\SU(4)\SUB\SO(9)$ and
$\SO(6)\SUB\SO(9)$.

\item[\fbox{$(\fb_3,\fa_2)$}]
By \ref{AnModules}, $(G,H)=(\Spin(7),\SU(3))$.

\item[\fbox{$(\fb_3,\fa_1\oplus\fa_1)$}]
The subgroups of type $\fa_1$ in $\SO(7)$ with large centralizers
are $\SO(3)$, with connected centralizer $\SO(4)$, and
$\SU(2)$, with connected centralizer $\SO(3)\times\SU(2)$, cp.
\ref{A1Modules}.
Hence we obtain subgroups $\SO(3)\times\SO(3)$ and $\SO(3)\times\SU(2)$
in the first case. From $\SU(2)$ we obtain in addition 
$\SU(2)\cdot\SU(2)=\SO(4)$ and another copy of $\SO(4)$, obtained
by mapping $\SU(2)$ diagonally into $\SO(3)\times\SU(2)$. This last
group is the copy of $\SO(4)$ which is contained in $\G_2=\SO(7)
\cap\Spin(7)\SUB\SO(8)$, and $\Spin(7)/\SO(4)=\widetilde G_3(\RR^8)$.
\end{description}

\begin{center}
{\small\begin{tabular}{llllll} 
\multicolumn{6}{c}{\fbox{$G$ of type $\SO(2n+1)$}} \\ 
\\
$G$                     & $H$         & $\Cen_G(H)^\circ$ & $G/H$ &
$(n_1,n_2)$ & Remarks \\ 
\hline
$\SO(2n+1)$     & $\SO(2n-2)$ & $\SO(3)$ & $V_3(\RR^{2n+1})$  & 
$(2n-2,4n-1)$ &  \\
\hfill$n\geq 3$ \\ 
\hline
%$\SO(9)$        & $\SO(6)$    & $\SO(3)$ & $V_3(\RR^9)$  & $(6,15)$ & 
%2-torsion \\
$\Spin(9)$      & $\SU(4)$    & $\U(1)$  &               & $(6,15)$ &  
 \\
$\Spin(7)$      & $\SU(3)$    & $\U(1)$  & $V_2(\RR^8)$  & $(6,7)$  & 
 \\
%$\SO(7)$        & $\SO(4)$    & $\SO(3)$ & $V_3(\RR^7)$  & $(4,11)$ & 
%2-torsion \\
$\SO(7)$        & $\SO(3)\times\SO(3)$ & 1 &             & $(4,11)$ &  
$\pi_2=\ZZ/2$ \\
$\SO(7)$      & $\SO(3)\times\SU(2)$ & 1 & 
& $(4,11)$ & $\pi_2=\ZZ/2$ \\
$\Spin(7)$      & $\SO(4)$    & 1        & $\widetilde G_3(\RR^8)$ 
& $(4,11)$ & $\pi_2=\ZZ/2$ \\
\end{tabular}}
\end{center}
The Leray-Serre spectral sequence of the bundle
\begin{diagram}
 \llap{$\SS^6={}$} \Spin(7)/\SU(4) & \rTo & \Spin(9)/\SU(4) \\
&& \dTo \\
&& \Spin(9)/\Spin(7) \rlap{${}=\SS^{15}$}
\end{diagram}
collapses, and thus
\[
\bfH^\bullet(\Spin(9)/\SU(4))\isom\ZZ[x_6]/(x_6^2)\otimes\EA_\ZZ(x_{15}).
\]
The group $\Spin(7)\SUB\SO(8)$ acts transitively on $V_2(\RR^8)$; thus
\[
\Spin(7)/\SU(3)=\SO(8)/\SO(6).
\]
Similarly, $\Spin(7)$ acts transitively on $\widetilde G_3(\RR^8)$ and
\[
\Spin(7)/\SO(4)=\SO(8)/\SO(5)\times\SO(3).
\]
\end{Num}

%\newpage

\begin{Num}\textsc{$G$ of type $\fc_n$}\psn
\label{CaseSp(n)even}%
The only possibilities are $(\fc_n,\fd_{n-1})$, for $n\geq 4$,
$(\fc_3,\fa_2)$, and $(\fc_3,\fa_1\oplus\fa_1)$.

\begin{description}
\item[\fbox{$(\fc_n,\fd_{n-1})$}]
There is no non-trivial representation of $\Spin(2n-2)$ on
$\HH^n$ for $n\geq5$ by \ref{DnModules}.
Therefore these pairs do not exist.
There is only one non-trivial representation of $\SU(4)$ on $\HH^4$, 
arising from the standard inclusion $\SU(4)\SUB\Sp(4)$,
cp.~\ref{AnModules}.

\item[\fbox{$(\fc_3,\fa_2)$}]
There is only one non-trivial representation of $\SU(3)$ on $\HH^3$, arising
from the standard inclusion $\SU(3)\SUB\Sp(3)$, cp.~\ref{AnModules}.

\item[\fbox{$(\fc_3,\fa_1\oplus\fa_1)$}]
The only representations of $\Sp(1)$ on $\HH^3$ with large
centralizers in $\Sp(3)$ are $\Sp(1)$, with connected
centralizer $\Sp(2)$, and $^\HH\rho_{3\lambda_1}$, with connected 
centralizer $\Sp(1)$, cp.~\ref{A1Modules}. Thus, the possible groups are 
$\Sp(1)\times\Sp(1)$ and $\Sp(1)\times{}^\HH\rho_{3\lambda_1}(\Sp(1))$.
\end{description}

\begin{center}
\begin{tabular}{lllll} 
\multicolumn{5}{c}{\fbox{$G$ of type $\Sp(n)$}} \\ 
\\
$G$      & $H$           & $\Cen_G(H)^\circ$ & $(n_1,n_2)$ & Remarks \\ 
\hline
$\Sp(4)$ & $\SU(4)$      &  1                & $(6,15)$    & $\pi_3=\ZZ/2$ \\
$\Sp(3)$ & $\SU(3)$      &  1                & $(6,7)$     & $\pi_3=\ZZ/2$ \\
$\Sp(3)$ & $\Sp(1)\times\Sp(1)$  & $\Sp(1)$  & $(4,11)$    &  \\
$\Sp(3)$ & $\Sp(1)\times{}^\HH\rho_{3\lambda_1}(\Sp(1))$  
& 1 & $(4,11)$    &  \\
\end{tabular}
\end{center}
The Leray-Serre spectral sequence of the bundle
\begin{diagram}
\llap{$\SS^4={}$} \Sp(2)/\Sp(1)\times\Sp(1) &
\rTo & \Sp(3)/\Sp(1)\times\Sp(1) \\
&& \dTo \\
&& \Sp(3)/\Sp(2) \rlap{${}=\SS^{11}$}
\end{diagram}
collapses, and thus
\[
\bfH^\bullet(\Sp(3)/\Sp(1)\times\Sp(1))\isom\ZZ[x_4]/(x_4^2)\otimes
\EA_\ZZ(x_{11}).
\]
\end{Num}

%\newpage

\begin{Num}\textsc{$G$ of type $\fd_n$}\psn
The only possibilities are $(\fd_n,\fd_{n-1})$, for $n\geq 4$,
$(\fd_5,\fa_4)$, and $(\fd_4,\fa_1\oplus\fb_2)$.

\begin{description}
\item[\fbox{$(\fd_n,\fd_{n-1})$}]
If $n\geq6$, then $(G,H)=(\SO(2n),\SO(2n-2))$, cp.~\ref{DnModules}.
For $\SU(4)$ we obtain $\SU(4)\SUB\SO(8)$ and $\SO(6)\SUB\SO(8)$,
cp.~\ref{AnModules}. Both 
inclusions become equal in $\Spin(8)$ under an automorphism of $\Spin(8)$;
hence we have also only the pair $(\SO(8),\SO(6))$.

\item[\fbox{$(\fd_5,\fa_4)$}]
Then $(G,H)=(\Spin(10),\SU(5))$, cp.~\ref{AnModules}.

\item[\fbox{$(\fd_4,\fa_1\oplus\fb_2)$}]
The non-trivial representations of $\Sp(2)$ on $\RR^8$ are $\Sp(2)$ and 
$\SO(5)$, cp.~\ref{BnModules}. 
Both inclusions become equal in $\Spin(8)$ under an automorphism, hence
we have the pair $(\SO(8),\SO(5)\times\SO(3))$.
\end{description}

\begin{center}
\begin{tabular}{llllll} 
\multicolumn{6}{c}{\fbox{$G$ of type $\SO(2n)$}} \\ 
\\
$G$              & $H$          & $\Cen_G(H)^\circ$ & $(n_1,n_2)$ & $G/H$ &
Remarks \\ 
\hline
$\SO(2n)$        & $\SO(2n-2)$  & $\SO(2)$   & $(2n-2,2n-1)$ & $V_2(\RR^{2n})$ &
 \\
\hfill$n\geq 4$  \\
\hline
$\Spin(10)$      & $\SU(5)$     & $\U(1)$    & $(6,15)$ & &  
(see below) \\
%$\SO(8)$         & $\SO(6)$     & $\SO(2)$   & $(6,7)$  & $V_2(\RR^8)$ &
% \\
$\SO(8)$         & $\SO(5)\times\SO(3)$ & 1  & $(4,11)$ &
$\widetilde G_3(\RR^8)$ & $\pi_2=\ZZ/2$\\
\end{tabular}
\end{center}
The subgroup $\SO(9)\SUB\SO(10)$ 
acts transitively on $\SO(10)/\SU(5)$; consequently,
\[
\Spin(9)/\SU(4)=\Spin(10)/\SU(5).
\]
Similarly,
\[
\Spin(7)/\SU(3)=\SO(8)/\SO(6)=V_2(\RR^8)
\]
and
\[
\Spin(7)/\SO(4)=\SO(8)/\SO(5)\times\SO(3)=\widetilde G_3(\RR^8).
\]
\end{Num}

%\newpage

\begin{Num}\textsc{$G$ exceptional and $H$ almost simple}\psn
The only possibilities are $(\ff_4,\fb_3)$ and $(\ff_4,\fc_3)$.
The maximal connected subgroups of maximal rank of $\Ffour$ are $\Spin(9)$,
$\SU(3)\times\SU(3)/\{\pm1\}$, $\Sp(3)\times\Sp(1)/\{\pm1\}$,
see Borel-De Siebenthal \cite{BorelSieb49}.
According to Dynkin \cite{DynkinAMST57} and
Seitz \cite{Seitz91}, the 
maximal subgroups of strictly smaller rank are of type
$\fa_1$, $\fg_2$ or $\fa_1\oplus\fg_2$; none of these has a subgroup of
type $\fb_3$ or $\fc_3$.
\begin{description}
\item[\fbox{$(\ff_4,\fb_3)$}]
Then $H\SUB\Spin(9)$. There are two conjugacy classes of groups of type
$\Spin(7)$ in $\Spin(9)$; however, both groups are conjugate in $\Ffour$,
as one can see by inspecting their respective fixed point sets in
the Cayley plane $\OO\mathrm{P}^2$.

\item[\fbox{$(\ff_4,\fc_3)$}]
Then $H=\Sp(3)$. There is one conjugacy class of subgroups of this
type.
\end{description}

\begin{center}
\begin{tabular}{lllll} 
\multicolumn{5}{c}{\fbox{$G$ of type $\Ffour$}} \\ 
\\
$G$              & $H$          & $\Cen_G(H)^\circ$ & $(n_1,n_2)$ & 
Remarks \\ 
\hline
$\Ffour$        & $\Spin(7)$    & $\SO(2)$   & $(8,23)$  & (see below) \\
$\Ffour$        & $\Sp(3)$      & $\Sp(1)$   & $(8,23)$  & $\pi_5=\ZZ/2$ \\
\end{tabular}
\end{center}
The inclusion $\G_2\SUB\Spin(7)\SUB\Ffour$ in $\ZZ/3$-cohomology is 
according to Borel \cite{Bor54} given by the diagram
\begin{diagram}
\bfH^\bullet(\G_2;\ZZ/3) & \lTo & 
\bfH^\bullet(\Spin(7);\ZZ/3) & \lTo & 
\bfH^\bullet(\Ffour;\ZZ/3) & {}\\
\dEq && \dEq && \dEq \\
\EA_{\ZZ/3}(x_3,x_{11}) & \lTo &
\EA_{\ZZ/3}(x_3,x_7,x_{11}) & \lTo 
& \EA_{\ZZ/3}(x_3,x_7,x_{11},x_{15})\otimes
(\ZZ/3)[x_8]/(x_8^3)\rlap{.}
\end{diagram}
Therefore the $\ZZ/3$-Leray-Serre spectral sequence of the bundle
\begin{diagram}
\Spin(7) & \rTo & \Ffour \\
&& \dTo \\
&& \Ffour/\Spin(7)
\end{diagram}
collapses, and
\[
\bfH^\bullet(\Ffour/\Spin(7);\ZZ/3)\isom
\EA_{\ZZ/3}(x_{15})\otimes(\ZZ/3)[x_8]/(x_8^3).
\]
\end{Num}

%\newpage
The following theorem summarizes the discussion of case (II).
\begin{Thm}
\label{Case(II)Simple}
Let $X=G/H$ be a 1-connected homogeneous space of an effective almost simple
compact Lie group $G$. If $X$ has the same integral cohomology as
a product of an even- and an odd-dimensional sphere,
\[
\bfH^\bullet(X)\isom\ZZ[a]\otimes\EA_\ZZ(w),
\]
$\deg(a)=n_1\geq 4$, $\deg(w)=n_2>n_1$,
then $(G,H)$ is one of the following pairs.

\begin{center}
\begin{tabular}{lllll} 
$G$              & $H$          & $\Cen_G(H)^\circ$ & $(n_1,n_2)$ & $G/H$ \\ 
\hline
$\SO(2n)$        & $\SO(2n-2)$  & $\SO(2)$   & $(2n-2,2n-1)$ & 
$V_2(\RR^{2n})$
 \\
\hfill$n\geq 3$  \\
\hline
$\Spin(10)$      & $\SU(5)$    & $\U(1)$  & $(6,15)$ &                \\
$\Spin(9)$       & $\SU(4)$    & $\U(1)$  & $(6,15)$ &                \\
$\Spin(7)$       & $\SU(3)$    & $\U(1)$  & $(6,7)$  & $V_2(\RR^8)$    \\
$\Sp(3)$         & $\Sp(1)\times\Sp(1)$ & $\Sp(1)$ & $(4,11)$ & \\
$\Sp(3)$         & $\Sp(1)\times{}^\HH\rho_{3\lambda_1}(\Sp(1))$ & 
$\Sp(1)$ & $(4,11)$ & \\
$\SU(5)$   & $\SU(3)\times\SU(2)$ & $\U(1)$ & $(4,9)$ & 
$\widetilde G_3(\CC^5)$   \\
\end{tabular}
\end{center}
There are two coincidences in this table;
\[
\Spin(9)/\SU(4)=\Spin(10)/\SU(5)
\]
and
\[
\Spin(7)/\SU(3)=\SO(8)/\SO(6)=V_2(\RR^8).
\]

\begin{proof}
This follows from the discussion of the various groups. All other
spaces are either not $(n_1-1)$-connected, or have torsion.
\end{proof}
\end{Thm}

%%%%%%%%%%%%%%%%%%%%%%%%%%%%%%%%%%%%%%%%%%%%%%%%%%%%%%%%%%%%%%%%%%%%%%%%
%                                                                      %
%                                                                      %
%                   Compact homogeneous quadrangles                    %
%                                                                      %
%                          Linus Kramer                                %
%                                                                      %
%                           Memoirs AMS                                %
%                                                                      %
%                         Wuerzburg 2000                               %
%                                                                      %
%                                                                      %
%                            CHQ6.tex                                  %
%                                                                      %
%                                                                      %
%                                                                      %
%%%%%%%%%%%%%%%%%%%%%%%%%%%%%%%%%%%%%%%%%%%%%%%%%%%%%%%%%%%%%%%%%%%%%%%%
\chapter{The case when $G$ is semisimple}

We continue the classification of homogeneous spaces $X=G/H$ which have
the same cohomology as a product of spheres. In this chapter we assume
that $G$ is semisimple with two almost simple factors
(by \ref{LengthAtMost2}, this is the remaining case). More precisely,
our assumptions are as follows.
\begin{enumerate}
\item $G/H$ is compact and 1-connected.
\item $G$ is a product of two compact almost simple Lie groups
$K_1$, $K_2$.
\item the action of $G$ on $X$ irreducible and almost effective.
\item $G/H$ has the same integral cohomology as a product
$\SS^{n_1}\times\SS^{n_2}$, where $3\leq n_1\leq n_2$, and $n_2$ is odd.
\end{enumerate}
%Again, there are two major cases: either $n_1$ is odd, or $n_1$ is even.
Note that $H$ has to be connected.
Put 
\[
H_i=(H\cap K_i)^\circ.
\]
Then $H_1,H_2$ are normal subgroups in $H$, and
thus there exists a compact connected normal subgroup $H_0\SUB H$ such 
that $H$ is an almost direct product
\[
H=(H_1\times H_2)\cdot H_0.
\]
If $H_0=1$, then 
\[
G/H=K_1/H_1\times K_2/H_2
\]
is split. We call this the \textbf{split case}.

If $H_0\neq 1$, then we consider again the two cases

\noindent
\textbf{Case (I):} $n_1$ is odd, $\bfH^\bullet(G/H)\isom\EA_\ZZ(u,v)$, with
$\deg(u)=n_1$ and $\deg(v)=n_2$. In this situation Theorem 
\ref{Case(I)General} applies.

\noindent
\textbf{Case (II):} $n_1$ is even, 
$\bfH^\bullet(G/H)\isom\ZZ[a]/(a^2)\otimes \EA_\ZZ(u)$, with
$\deg(a)=n_1$ and $\deg(u)=n_2$. In this situation Theorem 
\ref{Case(II)General} applies.

These three cases are considered in the next sections. Each
section ends with a complete list of all possible pairs $(G,H)$.

%\newpage

\section{The split case}

Here, 
\[
X=K_1/H_1\times K_2/H_2
\]
is a product. From the K\"unneth theorem we see that
$K_1/H_1$ and $K_2/H_2$ have to be 1-connected cohomology
spheres, say, of dimensions $n_1,n_2$. 
Thus we have to determine all pairs $(K,H)$ with the following properties.
\begin{enumerate}
\item $K$ is a compact almost simple Lie group.
\item $K/H$ is a 1-connected cohomology $m$-sphere.
\item the action of $K$ on $K/H$ is effective.
\end{enumerate}
If $m$ is even, then $K/H$ has Euler characteristic 2, and a
well-known result of Borel-De Siebenthal \cite{BorelSieb49}
and Borel \cite{BorelSphere} says that $K/H$ is
actually a sphere, and that $(K,H)$ is one of the following pairs.

\begin{center}
\begin{tabular}{lllll} 
\multicolumn{5}{c}{\fbox{Homogeneous spaces of Euler characteristic 2}} \\ 
\\ 
%\hline
$K$             & $H$        & $\Cen_K(H)^\circ$ & $m$ & $K/H$ \\ 
\hline
$\SO(2n+1)$     & $\SO(2n)$  & $1$       & $2n$  & $\SS^{2n}$  \\
\hfill$n\geq 1$  \\
\hline
$\G_2$          & $\SU(3)$   & $1$       & $6$   & $\SS^6$     \\
\end{tabular}
\end{center}
We need a similar result for 1-connected homogeneous spaces $K/H$
with 
\[
\bfH^\bullet(K/H;\ZZ)=\EA_\ZZ(w),
\]
where $\deg(w)=m\geq 3$ is odd.
This has been obtained by various authors, cp.~e.g. Bredon \cite{Bre61}, 
Onishchik \cite{Oni94} Ch. 5 \S18 Table 10.
For the sake of completeness, we give a proof. The method is very similar 
to case (I) in the previous chapter. 
The group $G$ has to be almost simple by a similar reasoning as in
\ref{LengthAtMost2}, and $H$ has to be trivial or almost simple, 
because of the 
isomorphism $P_H^1\OOT\isom P_G^1=0$ and the surjection
$P^3_H\oot P_G^3$, cp.~\ref{Case(I)General}.

%\newpage

\begin{Num}\textsc{$K$ of type $\fa_n$}\psn
The only possibilities are $(\fa_n,\fa_{n-1})$ and
$(\fa_3,\fb_2)$.

\begin{description}
\item[\fbox{$(\fa_n,\fa_{n-1})$}]
If $n\geq 3$, then $(K,H)=(\SU(n+1),\SU(n))$ by \ref{AnModules}.
For $n=2$, there is also the pair $(\SU(3),\SO(3))$.

\item[\fbox{$(\fa_3,\fb_2)$}]
By \ref{BnModules}, the only possibility is $(\SO(6),\SO(5))$.
\end{description}

\begin{center}
\begin{tabular}{llllll} 
\multicolumn{6}{c}{\fbox{$K$ of type $\SU(n+1)$}} \\ 
\\ 
$K$             & $H$        & $\Cen_K(H)^\circ$ & $m$ & $K/H$ &
Remarks \\ 
\hline
$\SU(n+1)$      & $\SU(n)$   & $\U(1)$   & $2n+1$  & $\SS^{2n+1}$ 
& \\
\hfill$n\geq 2$  \\
\hline
%$\SU(3)$        & $\SU(2)$   & $\U(1)$   & $7$     & $\SS^7$     &  \\
$\SO(6)$        & $\SO(5)$   & 1         & $5$       & $\SS^5$     & \\
$\SU(3)$        & $\SO(3)$   & 1         & $5$     &             &  
$\pi_3=\ZZ/4$ \\ 
$\SU(2)$        & $1$        & $\SU(2)$  & $3$     & $\SS^3$     &  \\
\end{tabular}
\end{center}
\end{Num}

%\newpage

\begin{Num}\textsc{$K$ of type $\fb_n$}\psn
The only possibilities are $(\fb_n,\fb_{n-1})$, $(\fb_n,\fc_{n-1})$ 
and $(\fb_3,\fg_2)$.

\begin{description}
\item[\fbox{$(\fb_n,\fc_{n-1})$}]
If $n\geq 4$, then there is no non-trivial representation of
$\Sp(n-1)$ on $\RR^{2n+1}$ by \ref{CnModules}, 
hence these pairs do not exist.

\item[\fbox{$(\fb_n,\fb_{n-1})$}]
If $n\geq 5$, then $(K,H)=(\SO(2n+1),\SO(2n-1))$ by \ref{BnModules}.
For $n=2,3,4$, there are also the pairs $(\Spin(9),\Spin(7))$,
$(\Spin(7),\Sp(2))$,
$(\Sp(2),\Sp(1))$, and %\linebreak
$(\Sp(2), {}^\HH\rho_{3\lambda_1}(\Sp(1)))$. 
However, $\Spin(7)$ contains the center of $\Sp(2)$, and thus
%\linebreak
$\Spin(7)/\Sp(2)=\SO(7)/\SO(5)$.

\item[\fbox{$(\fb_3,\fg_2)$}]
By \ref{G2Modules}, $(K,H)=(\Spin(7),\G_2)$.
\end{description}

\begin{center}
\begin{tabular}{llllll} 
\multicolumn{6}{c}{\fbox{$K$ of type $\SO(2n+1)$}} \\ 
\\
$K$          & $H$           & $\Cen_K(H)^\circ$ & $m$ & $K/H$ &
Remarks \\ 
\hline
$\SO(2n+1)$  & $\SO(2n-1)$   & $\SO(2)$  & $4n-1$ & $V_2(\RR^{2n+1})$ 
& $\pi_{2n-1}=\ZZ/2$ \\
\hfill $n\geq 2$  \\ 
\hline
%$\SO(9)$     & $\SO(7)$      & $\SO(2)$  & $15$   & $V_2(\RR^9)$ & 
%2-torsion \\
$\Spin(9)$   & $\Spin(7)$    & $1$       & $15$   & $\SS^{15}$   & \\
%$\SO(7)$     & $\SO(5)$      & $\SO(2)$  & $11$   & $V_2(\RR^7)$ & 
%2-torsion \\
$\Spin(7)$   & $\G_2$        & $1$       & $7$    & $\SS^7$      &  \\
$\Sp(2)$   & $\Sp(1)$        & $\Sp(1)$  & $7$    & $\SS^7$      & \\
%$\SO(5)$   & $\SO(3)$        & $\SO(2)$  & $7$    & $V_2(\RR^5)$ & \\
$\Sp(2)$   & $^\HH\rho_{3\lambda_1}$   & $1$   & $7$  &      & 
$\pi_3=\ZZ/10$ \\
\end{tabular}
\end{center}
\end{Num}

%\newpage

\begin{Num}\textsc{$K$ of type $\fc_n$}\psn
The only possibilities are $(\fc_n,\fc_{n-1})$, $(\fc_n,\fb_{n-1})$ and
$(\fc_3,\fg_2)$.

\begin{description}
\item[\fbox{$(\fc_n,\fb_{n-1})$}]
If $n\geq 4$, then there is no non-trivial representation of
$\Spin(2n-1)$ on $\HH^n$ by \ref{BnModules}, 
hence these pairs do not exist.

\item[\fbox{$(\fc_n,\fc_{n-1})$}]
If $n\geq 3$, then $(K,H)=(\Sp(n),\Sp(n-1))$ by \ref{CnModules}.
For $n=2$, there are also the pairs 
$(\SO(5),\SO(3))$ and $(\Sp(2),{}^\HH\rho_{3\lambda_1}(\Sp(1))$.

\item[\fbox{$(\fc_3,\fg_2)$}]
There is no non-trivial representation of $\G_2$ on $\HH^3$ by
\ref{G2Modules}.
\end{description}

\begin{center}
\begin{tabular}{llllll} 
\multicolumn{6}{c}{\fbox{$K$ of type $\Sp(n)$}} \\ 
\\
$K$        & $H$           & $\Cen_K(H)^\circ$ & $m$ & $K/H$ & 
Remarks\\
\hline
$\Sp(n)$   & $\Sp(n-1)$    & $\Sp(1)$ & $4n-1$ & $\SS^{4n-1}$ \\
\hfill $n\geq 2$ \\
\hline
%$\Sp(2)$   & $\Sp(1)$        & $\Sp(1)$  & $7$    & $\SS^7$      & \\
$\SO(5)$   & $\SO(3)$        & $\SO(2)$  & $7$    & $V_2(\RR^5)$ & 
$\pi_3=\ZZ/2$ \\
$\Sp(2)$   & $^\HH\rho_{3\lambda_1}$ & $1$  & $7$ &   & 
$\pi_3=\ZZ/10$ \\
$\Sp(1)$   & $1$             & $\Sp(1)$  & $3$    & $\SS^3$      & \\
\end{tabular}
\end{center}
\end{Num}

%\newpage

\begin{Num}\textsc{$K$ of type $\fd_n$}\psn
The only possibilities are $(\fd_n,\fb_{n-1})$ and
$(\fd_n,\fc_{n-1})$.

\begin{description}
\item[\fbox{$(\fd_n,\fb_{n-1})$}]
If $n\geq 5$, then $(K,H)=(\SO(2n),\SO(2n-1))$ by \ref{BnModules}.
For $n=4$, there are is also the pair $(\Spin(8),\Spin(7))$. However,
the action of $\Spin(8)$ on $\Spin(8)/\Spin(7)$ is not effective, and
thus $\Spin(8)/\Spin(7)=\SO(8)/\SO(7)$.

\item[\fbox{$(\fd_n,\fc_{n-1})$}]
If $n\geq 4$, then $\RR^{2n}$ is not a non-trivial $\Sp(n-1)$-module
by \ref{CnModules}, hence these pairs do not exist.
\end{description}

\begin{center}
\begin{tabular}{lllll} 
\multicolumn{5}{c}{\fbox{$K$ of type $\SO(2n)$}} \\ 
\\
$K$          & $H$         & $\Cen_K(H)^\circ$ & $m$ & $K/H$ \\ 
\hline
$\SO(2n)$    & $\SO(2n-1)$ & $1$      & $2n-1$ & $\SS^{2n-1}$ \\
\hfill$n\geq 3$ \\
\end{tabular}
\end{center}
\end{Num}

%\newpage

\begin{Num}\textsc{$K$ of exceptional type}\psn
\label{SU(2)inG2}%
The only possibility is $(\fg_2,\fa_1)$. The subgroups of type
$\fa_1$ in $\G_2$ were determined in \ref{SubgroupsInExceptionalGroups}.

\begin{center}
\begin{tabular}{llllll} 
\multicolumn{6}{c}{\fbox{$K$ of exceptional type}} \\ 
\\
$K$          & $H$         & $\Cen_K(H)^\circ$ & $m$ & $K/H$ & Remarks \\ 
\hline
$\G_2$        & $\SU(2)$    & $\SU(2)$   & $11$ & $V_2(\RR^7)$ &
$\pi_5=\ZZ/2$ \\
$\G_2$        & $2\cdot{}^\RR\rho_{2\lambda_1}$  & 1 & $11$ & & 
$\pi_3=\ZZ/4$ \\
$\G_2$        & ${}^\RR\rho_{\lambda_1}+{}^\RR\rho_{2\lambda_1}$
  & $\SU(2)$   & $11$ & & $\pi_3=\ZZ/3$ \\
$\G_2$        & ${}^\RR\rho_{6\lambda_1}$       & 1 & $11$ & & 
$\pi_3=\ZZ/28$ \\
\end{tabular}
\end{center}
Later we will need the $\ZZ/2$-cohomology of $\G_2/H$, where
$H\isom\SU(2)$. There are two such embeddings. In both cases, the
map $\bfH_3(\SU(2))\too\bfH_3(\G_2))$ is multiplication by an odd
number (1 or 3, respectively). Thus we have an isomorphism
$\bfH^3(\SU(2);\ZZ/2)\oot\bfH^3(\G_2;\ZZ/2)$. It follows that the
$\ZZ/2$-Lerray-Serre spectral sequence of the bundle
\begin{diagram}
\llap{$\SU(2)\isom{}$} H & \rTo & \G_2 \\
&& \dTo \\
&& \G_2/H
\end{diagram}
collapses. For $H=\SU(2)$, the base is the Stiefel manifold
$V_2(\RR^7)$. Thus 
\[
\bfH^\bullet(\G_2;\ZZ/2)\isom\EA_{\ZZ/2}(x_3,x_5)\otimes
(\ZZ/2)[x_6]/(x_6^2).
\]
It follows that in both cases
\[
\bfH^\bullet(\G_2/H;\ZZ/2)\isom\EA_{\ZZ/2}(x_5)\otimes
(\ZZ/2)[x_6]/(x_6^2).
\]
\end{Num}

%\newpage

\begin{Thm}
\label{HomogeneousSpheres}
Let $K/H$ be a 1-connected homogeneous space of a compact connected
simple Lie group $K$. Assume that the action of $K$ is effective and
irreducible. Assume that $K/H$ has the same integral cohomology 
as an $m$-sphere, $m\geq 2$.
Then $K/H$ is a (standard) sphere, and $(K,H)$ is one of the
following pairs.
\begin{center}
\begin{tabular}{llll} 
\multicolumn{4}{c}{\fbox{Homogeneous spheres}} \\ 
\\ 
$K$             & $H$        & $\Cen_K(H)^\circ$ & $K/H$ \\ 
\hline
$\SO(2n+1)$     & $\SO(2n)$  & $1$       & $\SS^{2n}$  \\
\hfill$n\geq 1$  \\
\hline
$\G_2$          & $\SU(3)$   & $1$       & $\SS^6$     \\
\hline \hline
$\SU(n+1)$      & $\SU(n)$   & $\U(1)$   & $\SS^{2n+1}$  \\
\hfill$n\geq 2$  \\
\hline
$\Sp(n)$   & $\Sp(n-1)$    & $\Sp(1)$    & $\SS^{4n-1}$ \\
\hfill $n\geq 2$ \\
\hline
$\SO(2n)$    & $\SO(2n-1)$ & $1$         & $\SS^{2n-1}$ \\
\hfill$n\geq 3$ \\
\hline
$\Spin(9)$   & $\Spin(7)$    & $1$       & $\SS^{15}$    \\
$\Spin(7)$   & $\G_2$        & $1$       & $\SS^7$       \\
%$\Sp(2)$     & $\Sp(1)$      & $\Sp(1)$  & $\SS^7$       \\
%$\SU(3)$        & $\SU(2)$   & $\U(1)$   & $\SS^5$       \\
$\SU(2)$        & $1$        & $\SU(2)$  & $\SS^3$       \\
\end{tabular}
\end{center}
\end{Thm}

\begin{Cor}
\label{SplitCase}
Let $G/H$ be a 1-connected homogeneous space of a compact
connected Lie group $G$. Assume that $G=K_1\times K_2$ is a
product of two compact almost simple Lie groups.
Assume moreover that the action is irreducible and split, i.e.~that
$H=(K_1\cap H)\times(K_2\cap H)$.
If $G/H$ has the same integral cohomology as a product of spheres
$\SS^{n_1}\times\SS^{n_2}$, where $n_1\geq 3$ and $n_2\geq n_1$ is odd,
then the factors $K_i/H_i$ are spheres, and the possibilities
of the pairs $(K_1,H_1)$, $(K_2,H_2)$ are given by the table above.
\qed
\end{Cor}
We note also the following.
\begin{Thm}
\label{HomogeneousRealStiefel}
Let $G/H$ be a 1-connected homogeneous space of a compact connected
Lie group $G$. Assume that the action is effective and 
irreducible,
and that $G/H$ has the same rational cohomology as $\SS^{2m-3}$, where
$m\geq 3$. Assume that $G/H$ is $(m-2)$-connected, and that
$\pi_{m-1}(G/H)\isom\ZZ/2$. Then $G/H\isom V_2(\RR^m)$ is a
Stiefel manifold, and there are only the following possibilities.

\begin{center}
\begin{tabular}{llll} 
%\multicolumn{6}{c}{\fbox{$G$ of type $\SO(2n+1)$}} \\ 
%\\
$G$          & $H$           & $\Cen_G(H)^\circ$ &  $X=G/H$ \\
\hline
$\SO(2n+1)$  & $\SO(2n-1)$   & $\SO(2)$  &  $V_2(\RR^{2n+1})$ \\
\hfill $n\geq 2$  \\ 
\hline
$\G_2$  & $\SU(2)$   & $\U(1)$  & $V_2(\RR^7)$ \\
\end{tabular}
\end{center}
\end{Thm}
The following will also be useful.
\begin{Cor}
\label{HomogeneousRationalSphereWithLargeCentralizer}
Let $K/H$ be a 1-connected homogeneous space of a compact connected
Lie group $K$. Assume that the action is effective and 
irreducible,
and that $K/H$ has the same rational cohomology as an odd-dimensional
sphere $\SS^m$ where $m\geq 3$. Assume in addition that $\Cen_K(H)$
contains a subgroup of type $\fa_1$. Then $(K,H)$ is one of the
following pairs.

\begin{center}
\begin{tabular}{ll} 
$K$          & $H$           \\
\hline
$\Sp(n)$     & $\Sp(n-1)$    \\
$n\geq 2$ \\
\hline
$\G_2$       & $\SU(2)$ \\
$\G_2$       & $^\RR\rho_{\lambda_1}+{}^\RR\rho_{2\lambda_1}$ \\
$\Sp(1)$     & $1$           
\end{tabular}
\end{center}
\end{Cor}

%\newpage

\section{The non-split case (I): 
$\bfH^\bullet(X)=\EA_\ZZ(u,v)$.}

The setting is as described at the beginning of this chapter.
Thus $G=K_1\times K_2$ and $H=(H_1\times H_2)\cdot H_0$.
We assume that $n_1$ is odd, and that $H_0\neq 1$.
We first prove some general results.
On the level of Lie algebras we have a direct sum
of ideals $\fh=\fh_1\oplus\fh_2\oplus\fh_0$. Let $\pr_i$ denote
the projection $\fg=\fk_1\oplus\fk_2\too \fk_i$.
The restriction of $\pr_i$ to $\fh_0\oplus\fh_i$ is a monomorphism,
for $i=1,2$. Since $\fh_0$ and $\fh_i$
commute, $[\fh_0,\fh_i]=0$, the same is true for their images under
$\pr_i$. It follows that on the group level the restriction of
$\pr_i$ to $H_0$ has finite kernel, and that $\pr_i(H_0)$ is centralized
by $H_i$,
\[
\pr_i(H_0)\SUB\Cen_{K_i}(H_i).
\]
Put $k_i=\rk(K_i)$ and $h_i=\rk(H_i)$. 
Since $\pr_i(\fh_0)\oplus\fh_i\SUB\fk_i$, 
we have 
\[
k_i-h_i-h_0\geq 0
\]
for $i=1,2$.
Note that $\pr_i(H_0)\neq K_i$, since otherwise
the action would be reducible. In particular, 
\[
k_1,k_2\geq 2
\]
(otherwise $\pr_i$ would be surjective).
Note also that $h_0\geq1$, because we assume that $H_0\neq 1$. Thus 
\begin{align*}
2& = \rk(G)-\rk(H) \\
&= k_1-h_1+k_2-h_2-h_0.
\end{align*}
If we subtract $h_0$, the result is still non-negative, hence 
$1\leq h_0\leq 2$. From the exact sequences
\[
0\oot P^1_H\oot P^1_G\oot0\quad\text{ and }\quad
0\oot P^3_H\oot P^3_G\oot\bfH^3(G/H;\QQ)\oot 0
\]
we see that $H$ is semisimple with at most two almost simple factors.
Therefore we may assume that $H_2=1$ (and thus $h_2=0$).
The equation $2=k_1-h_1-h_0+k_2$ shows that 
\[
k_2=2
\]
and
\[
k_1=h_1+h_0.
\]
Note also that $P_{H_1}\oot P_{K_1}$ is onto, because
$P_H\oot P_G$ is onto.

\begin{Lem}
If $h_0=2$, then $H=H_0$ is of type $\fa_1+\fa_1$ and
$K_1,K_2\in\{\G_2,\Sp(2)\}$. If $K_1=K_2$, then $H_0$ is not
conjugate to a subgroup contained in the diagonal. In particular,
$G\neq\Sp(2)\times\Sp(2)$.

\begin{proof}
Since $\pr_2(H_0)\SUB K_2$ the group $K_2$ has a semisimple subgroup of 
rank 2. The only possibilities for $(K_2,\pr_2(H_0))$ are thus
\[
(\G_2,\SO(4)),(\G_2,\SU(3)),(\Sp(2),\Sp(1)\times\Sp(1)).
\]
Note that $H_1$ has to centralize $\pr_1(H_0)$, and that
$(K_1,H_1)$ is one of the pairs determined in 
\ref{HomogeneousRationalSphereWithLargeCentralizer}.
This shows that $H_1=1$. Since $\pr_1(H_0)\neq K_1$,
we obtain the possibilities $\Sp(2)$ and $\G_2$ for $K_1$.

If $H_0$ is contained in the diagonal, then $\pi_3(G/H_0)\isom\ZZ$, and
thus $\bfH^\bullet(G/H;\QQ)$ contains an element of degree 3, and
thus the structure of the cohomology ring is not the right one
(we need $\bfH^\bullet(G/H;\QQ)\isom\EA_\QQ(P_G/P_G^3)$).

If $H_0=\SU(3)$, then $K_1=K_2=\G_2$. Thus $H_0$ is conjugate to
the copy of $\SU(3)$ which is contained in the diagonal.
\end{proof}
\end{Lem}

\begin{Num}\textsc{The case $\rk(H_0)=2$}\psn
The following possibilities remain: (1) $G=\Sp(2)\times\G_2$, and
$H_0$ is contained in $\Sp(1)\times\Sp(1)\times\SO(4)$. All such
groups which are not contained in one factor are conjugate. The map
$\pi_3(H)\too\pi_3(G)$ is given by the matrix
$\begin{pmatrix} 1 & 3 \\ 1 & 1 \end{pmatrix}$, and thus
$\pi_3(G/H)=\ZZ/2$.

(2) $G=\G_2\times\G_2$ and $H_0\SUB\SO(4)\times\SO(4)$ is a subgroup
which is not contained in one factor. There are two such conjugacy classes,
the diagonal subgroup and the anti-diagonal subgroup. In this case
the map $\pi_3(H)\too\pi_3(G)$ is given by the matrix
$\begin{pmatrix} 1 & 3 \\ 3 & 1 \end{pmatrix}$, and thus
$\pi_3(G/H)=\ZZ/8$.

\begin{center}
\begin{tabular}{lll} 
\multicolumn{3}{c}{\fbox{$H_0$ of type $\fa_1+\fa_1$}} \\ 
\\
$G$            & $(n_1,n_2)$ & Remarks\\
\hline
$\G_2\times\Sp(2)$  & $(7,11)$  & $\pi_3=\ZZ/2$ \\
$\G_2\times\G_2$    & $(11,11)$ & $\pi_3=\ZZ/8$ 
\end{tabular}
\end{center}
\end{Num}

\begin{Lem}
If $h_0=1$, then $(K_1,H_1)$ is one of the pairs
$(\fc_n,\fc_{n-1})$, $n\geq 2$, or $(\fg_2,\fa_1)$.
The group $H_0$ is of type $\fa_1$, and $K_2$ is of type
$\fa_2$, $\fc_2$, or $\fg_2$.

\begin{proof}
The group $K_2$ has rank 2, hence
$K_2$ is one of the groups $\SU(3)$, $\Sp(2)$, $\G_2$, and
$H_0$ is of type $\fa_1$.
Moreover, $k_1-h_1=1$, and $P_{H_1}\oot P_{K_1}$ is a surjection.
Thus $(K_1,H_1)$ is one of the pairs determined in 
\ref{HomogeneousRationalSphereWithLargeCentralizer}.
The only possibilities are the pairs
$(\fc_n,\fc_{n-1})$, $n\geq 2$, and $(\fg_2,\fa_1)$.
\end{proof}
\end{Lem}
Before we consider the various possibilities, we make some more
general observations. It will turn out that in all cases
$\pr_1(H_0)\isom\SU(2)$. Therefore we can write the elements
of $H_0$ as pairs $(h,\phi(h))$, where $\phi:\SU(2)\too \pr_2(H_0)\SUB
K_2$ is some fixed non-trivial homomorphism. Consider the fibre bundle
which arises from the inclusions 
\[
H\SUB K_1\times \phi(H_0)\SUB G
\]
by taking quotients. Note that 
\[
(K_1\times\phi(H_0))/H\isom K_1/(H_1\cdot\ker(\phi))
\]
(because the two spaces have the same dimension) and
\[
G/(K_1\times\phi(H_0))\isom K_2/\phi(H_0).
\]
Hence we obtain a fibre bundle
\begin{diagram}
K_1/(H_1\cdot\ker(\phi)) & \rTo & G/H \\
 & & \dTo \\
 & & K_2/\phi(H_0) \rlap{,}
\end{diagram}
with structure group $K_1\times\phi(H_0)$.

\begin{Num}\textsc{$h_0=1$ and $K_1$ of type $\fc_n$}\psn
Then $K_1=\Sp(n)$, $H_1=\Sp(n-1)$, $\pr_1(H_0)=\Sp(1)$, with
the standard inclusion $\Sp(n-1)\times\Sp(1)\SUB\Sp(n)$, and
$n\geq 2$.
Thus $H_0=\Sp(1)$, and the inclusion $H_0\SUB K_1\times K_2$ is
given by $h\mapstoo(h,\phi(h))$, where $\phi:\Sp(1)\too K_2$ is
a non-trivial representation.

\begin{description}
\item[\fbox{$K_2=\SU(3)$}]
(1)
If $\phi=\rho_{\lambda_1}$ is the standard inclusion $\SU(2)\SUB\SU(3)$,
then we obtain a fibre bundle
\begin{diagram}
\SS^{4n-1} & \rTo & G/H \\
& & \dTo \\
& & \SS^5 \rlap{.}
\end{diagram}
The structure group of this sphere bundle is
$\Sp(n)\times\Sp(1)\SUB\SO(4n)$, hence it is given by a classifying map
\[
\SS^5\too\mathrm{B}(\Sp(n)\times\Sp(1))\too\mathrm{BSO}(4n).
\] 
Since $\pi_5(\mathrm{BSO}(4n))\isom
\pi_4(\SO(4n))\isom\pi_4(\SO)=0$ for $n\geq 2$, the map is homotopic to 
a constant map and the bundle is trivial,
\[
G/H\isom \SS^{4n-1}\times\SS^5.
\]
(2) 
For $\phi={}^\RR\rho_{2\lambda_1}$ we obtain an 
$\RR\mathrm{P}^{4n-1}$-bundle
over $\SU(3)/\SO(3)$. Note that the higher homotopy groups of
the fibre are the same as those of $\SS^{4n-1}$.
This implies that $\pi_3(G/H)\isom
\pi_3(\SU(3)/\SO(3))=\ZZ/4$.

\item[\fbox{$K_2=\Sp(2)$}]
For $\phi:\SU(2)\too\Sp(2)$ we have the three possibilities 
(1) $\phi={}^\HH\rho_{\lambda_1}$,
(2) $\phi=2\cdot{}^\HH\rho_{\lambda_1}$ and
(3) $\phi={}^\HH\rho_{3\lambda_1}$. 
In all three cases the fibre of our bundle is $\SS^{4n-1}$, hence 
$\pi_3(G/H)\isom\pi_3(\Sp(2)/\phi(\Sp(1))=\ZZ/j_\phi$. 

In case (1), the base of the bundle is $\SS^7$. Similarly as above,
\linebreak
$\pi_7(\mathrm{BSO}(4n))\isom\pi_6(\SO(4n))\isom\pi_6(\SO)=0$ for 
$n\geq 2$, hence 
\[
G/H=\SS^{4n-1}\times\SS^7.
\]

\item[\fbox{$K_2=\G_2$}]
We have the three possibilities 
(1) $\phi={}^\RR\rho_{\lambda_1}$,
(2) $\phi={}^\RR\rho_{\lambda_1}+{}^\RR\rho_{2\lambda_1}$,
(3) $\phi=2\cdot{}^\RR\rho_{2\lambda_1}$,
and
(4) $\phi={}^\RR\rho_{6\lambda_1}$.
The fibre of the bundle is $\SS^{4n-1}$ in case (1) and (2), and
$\RR\mathrm{P}^{4n-1}$ in case (3) and (4), hence 
$\pi_3(G/H)\isom\ZZ/j_\phi$.
\end{description}

\begin{center}
\begin{tabular}{llllll} 
\multicolumn{6}{c}{\fbox{$K_1$ of type $\Sp(n)$}} \\ 
\\
$K_2$ & $\phi$           && $(n_1,n_2)$ & $G/H$ & 
Remarks\\
\hline
$\SU(3)$ & $\rho_{\lambda_1}$ & & $(5,4n-1)$ & $\SS^5\times\SS^{4n-1}$ & \\
$\SU(3)$ & $\rho_{2\lambda_1}$ & & $(5,4n-1)$ & & $\pi_3=\ZZ/4$ \\
\hline
$\Sp(2)$ & $^\HH\rho_{\lambda_1}$ & & $(7,4n-1)$ & $\SS^7\times\SS^{4n-1}$ & \\
$\SO(5)$ & $^\RR\rho_{2\lambda_1}$ & & 
$(7,4n-1)$ & & $\pi_3=\ZZ/2$ \\
$\Sp(2)$ & $^\HH\rho_{3\lambda_1}$ & & 
$(7,4n-1)$ & & $\pi_3=\ZZ/10$ \\
\hline
$\G_2$ & $^\RR\rho_{\lambda_1}$ & & $(11,4n-1)$ & & $\pi_5=\ZZ/2$ \\
$\G_2$ & ${}^\RR\rho_{\lambda_1}+{}^\RR\rho_{2\lambda_1}$ & 
& $(11,4n-1)$ & & $\pi_3=\ZZ/3$ \\
$\G_2$ & $2\cdot{}^\RR\rho_{2\lambda_1}$ & 
& $(11,4n-1)$ & & $\pi_3=\ZZ/4$ \\
$\G_2$ & $^\RR\rho_{6\lambda_1}$ & 
& $(11,4n-1)$ & & $\pi_3=\ZZ/28$ \\
\end{tabular}
\end{center}
\end{Num}

\begin{Num}\textsc{$h_0=1$ and $K_1$ of type $\fg_2$}\psn
Then $H_1\SUB\G_2$ is a subgroup of type $\fa_1$ with
large centralizer. The only possibilities are the two
copies of $\SU(2)$ in $\SO(4)\SUB\G_2$. In particular, $H_1\isom\SU(2)$.
For $K_2$ we have the possibilities $\SU(3)$, $\Sp(2)$, and $\G_2$. 
We may write the map $H=H_1\cdot H_0\too G=K_1\times\G_2$ as
$ab\mapstoo(\psi(a)\phi_1(b),\phi(b))$, where
$\psi(a)\phi_1(b)\in\SO(4)\SUB\G_2$ and $\phi:\SU(2)\too G$.
Thus $\{\psi,\phi_1\}=\{{}^\RR\rho_{2\lambda_1},
{}^\RR\rho_{2\lambda_1}+{}^\RR\rho_{\lambda_1}\}$.
It follows that the map $\pi_3(H)\too\pi_3(G)$ is
represented by the matrix
\[
\begin{pmatrix}
j_\psi & j_{\phi_1} \\
0 & j_\phi
\end{pmatrix}.
\]
Note that $j_{({}^\RR\rho_{\lambda_1})}=1$ and
$j_{({}^\RR\rho_{2\lambda_1}+{}^\RR\rho_{\lambda_1})}=3$.
We denote equivalence of matrices over $\ZZ$ by $\sim_\ZZ$.

\begin{description}
\item[\fbox{$K_2=\SU(3)$}]
Then we have the two possibilities $\phi=\rho_{\lambda_1},
\rho_{2\lambda_1}$.
We obtain the matrices
\[
\begin{pmatrix}
1 & 3 \\
0& 1
\end{pmatrix}
\sim_\ZZ
\begin{pmatrix}
1 & 0 \\
0 & 1
\end{pmatrix}
\text{ and }
\begin{pmatrix}
3 & 1 \\
0& 1
\end{pmatrix}
\sim_\ZZ
\begin{pmatrix}
1 & 0 \\
0 & 3
\end{pmatrix}
\qquad(\phi=\rho_{\lambda_1})
\]
\[
\begin{pmatrix}
1 & 3 \\
0& 4
\end{pmatrix}
\sim_\ZZ
\begin{pmatrix}
1 & 0 \\
0 & 4
\end{pmatrix}
\text{ and }
\begin{pmatrix}
3 & 1 \\
0& 4
\end{pmatrix}
\sim_\ZZ
\begin{pmatrix}
1 & 0 \\
0 & 12
\end{pmatrix}
\qquad(\phi=\rho_{2\lambda_1})
\]
\item[\fbox{$K_2=\Sp(2)$}]
The three possibilities for $\phi$ are 
$^\HH\rho_{\lambda_1},2\cdot{}^\HH\rho_{\lambda_1},{}^\HH\rho_{3\lambda_1}$.
Thus we have the matrices
\[
\begin{pmatrix}
1 & 3 \\
0& 1
\end{pmatrix}
\sim_\ZZ
\begin{pmatrix}
1 & 0 \\
0 & 1
\end{pmatrix}
\text{ and }
\begin{pmatrix}
3 & 1 \\
0& 1
\end{pmatrix}
\sim_\ZZ
\begin{pmatrix}
1 & 0 \\
0 & 3
\end{pmatrix}
\qquad(\phi={}^\HH\rho_{\lambda_1})
\]
\[
\begin{pmatrix}
1 & 3 \\
0 & 2
\end{pmatrix}
\sim_\ZZ
\begin{pmatrix}
1 & 0 \\
0 & 2
\end{pmatrix}
\text{ and }
\begin{pmatrix}
3 & 1 \\
0& 2
\end{pmatrix}
\sim_\ZZ
\begin{pmatrix}
1 & 0 \\
0 & 6
\end{pmatrix}
\qquad(\phi=2\cdot{}^\HH\rho_{\lambda_1})
\]
\[
\begin{pmatrix}
1 & 3 \\
0& 10
\end{pmatrix}
\sim_\ZZ
\begin{pmatrix}
1 & 0 \\
0 & 10
\end{pmatrix}
\text{ and }
\begin{pmatrix}
3 & 1 \\
0& 10
\end{pmatrix}
\sim_\ZZ
\begin{pmatrix}
1 & 0 \\
0 & 30
\end{pmatrix}
\qquad(\phi={}^\HH\rho_{3\lambda_1})
\]

\item[\fbox{$K_2=\G_2$}]
In this case we have the possibilities
$\phi={}^\RR\rho_{\lambda_1},
{}^\RR\rho_{\lambda_1}+{}^\RR\rho_{2\lambda_1},
2\cdot{}^\RR\rho_{2\lambda_1},
{}^\RR\rho_{6\lambda_1}$. Accordingly, the matrices are
as follows.
\[
\begin{pmatrix}
1 & 3 \\
0& 1
\end{pmatrix}
\sim_\ZZ
\begin{pmatrix}
1 & 0 \\
0 & 1
\end{pmatrix}
\text{ and }
\begin{pmatrix}
3 & 1 \\
0& 1
\end{pmatrix}
\sim_\ZZ
\begin{pmatrix}
1 & 0 \\
0 & 3
\end{pmatrix}
\qquad(\phi={}^\RR\rho_{\lambda_1})
\]
\[
\begin{pmatrix}
1 & 3 \\
0 & 3
\end{pmatrix}
\sim_\ZZ
\begin{pmatrix}
1 & 0 \\
0 & 3
\end{pmatrix}
\text{ and }
\begin{pmatrix}
3 & 1 \\
0& 3
\end{pmatrix}
\sim_\ZZ
\begin{pmatrix}
1 & 0 \\
0 & 9
\end{pmatrix}
\qquad(\phi={}^\RR\rho_{\lambda_1}+{}^\RR\rho_{2\lambda_1})
\]
\[
\begin{pmatrix}
1 & 3 \\
0& 4
\end{pmatrix}
\sim_\ZZ
\begin{pmatrix}
1 & 0 \\
0 & 4
\end{pmatrix}
\text{ and }
\begin{pmatrix}
3 & 1 \\
0& 4
\end{pmatrix}
\sim_\ZZ
\begin{pmatrix}
1 & 0 \\
0 & 12
\end{pmatrix}
\qquad(\phi=2\cdot{}^\RR\rho_{\lambda_1})
\]
\[
\begin{pmatrix}
1 & 3 \\
0& 28
\end{pmatrix}
\sim_\ZZ
\begin{pmatrix}
1 & 0 \\
0 & 28
\end{pmatrix}
\text{ and }
\begin{pmatrix}
3 & 1 \\
0& 28
\end{pmatrix}
\sim_\ZZ
\begin{pmatrix}
1 & 0 \\
0 & 84
\end{pmatrix}
\qquad(\phi={}^\RR\rho_{6\lambda_1})
\]
\end{description}

\begin{center}
\begin{tabular}{llll} 
\multicolumn{4}{c}{\fbox{$K_1$ of type $\G_2$}} \\ 
\\
$K_2$ & $(\psi,\phi_1,\phi)$ & $(n_1,n_2)$ & Remarks\\
\hline
$\SU(3)$ & 
$({}^\RR\rho_{\lambda_1},{}^\RR\rho_{\lambda_1}+{}^\RR\rho_{2\lambda_1},
\rho_{\lambda_1})$ & $(5,11)$ & $\pi_5=\ZZ/2$ \\
$\SU(3)$ & 
$({}^\RR\rho_{\lambda_1}+{}^\RR\rho_{2\lambda_1},{}^\RR\rho_{\lambda_1},
\rho_{\lambda_1})$ & $(5,11)$ & $\pi_3=\ZZ/4$ \\
$\SU(3)$ & 
$({}^\RR\rho_{\lambda_1},{}^\RR\rho_{\lambda_1}+{}^\RR\rho_{2\lambda_1},
2\rho_{\lambda_1})$ & $(5,11)$ & $\pi_3=\ZZ/3$ \\
$\SU(3)$ & 
$({}^\RR\rho_{\lambda_1}+{}^\RR\rho_{2\lambda_1},{}^\RR\rho_{\lambda_1},
2\rho_{\lambda_1})$ & $(5,11)$ & $\pi_3=\ZZ/12$ \\
\hline
$\Sp(2)$ & 
$({}^\RR\rho_{\lambda_1},{}^\RR\rho_{\lambda_1}+{}^\RR\rho_{2\lambda_1},
{}^\HH\rho_{\lambda_1})$ & $(7,11)$ & $\pi_5=\ZZ/2$ \\
$\Sp(2)$ & 
$({}^\RR\rho_{\lambda_1}+{}^\RR\rho_{2\lambda_1},{}^\RR\rho_{\lambda_1},
{}^\HH\rho_{\lambda_1})$ & $(7,11)$ & $\pi_3=\ZZ/3$ \\
$\Sp(2)$ & 
$({}^\RR\rho_{\lambda_1},{}^\RR\rho_{\lambda_1}+{}^\RR\rho_{2\lambda_1},
2\cdot{}^\HH\rho_{\lambda_1})$ & $(7,11)$ & $\pi_5=\ZZ/2$ \\
$\Sp(3)$ & 
$({}^\RR\rho_{\lambda_1}+{}^\RR\rho_{2\lambda_1},{}^\RR\rho_{\lambda_1},
2\cdot{}^\HH\rho_{\lambda_1})$ & $(7,11)$ & $\pi_3=\ZZ/6$ \\
$\Sp(2)$ & 
$({}^\RR\rho_{\lambda_1},{}^\RR\rho_{\lambda_1}+{}^\RR\rho_{2\lambda_1},
{}^\HH\rho_{3\lambda_1})$ & $(7,11)$ & $\pi_5=\ZZ/10$ \\
$\Sp(3)$ & 
$({}^\RR\rho_{\lambda_1}+{}^\RR\rho_{2\lambda_1},{}^\RR\rho_{\lambda_1},
{}^\HH\rho_{3\lambda_1})$ & $(7,11)$ & $\pi_3=\ZZ/30$ \\
\hline
$\G_2$ & 
$({}^\RR\rho_{\lambda_1},{}^\RR\rho_{\lambda_1}+{}^\RR\rho_{2\lambda_1},
{}^\RR\rho_{\lambda_1})$ & $(11,11)$ & $\pi_5=\ZZ/2\oplus\ZZ/2$ \\
$\G_2$ & 
$({}^\RR\rho_{\lambda_1}+{}^\RR\rho_{2\lambda_1},{}^\RR\rho_{\lambda_1},
{}^\RR\rho_{\lambda_1})$ & $(11,11)$ & $\pi_3=\ZZ/3$ \\
$\G_2$ & 
$({}^\RR\rho_{\lambda_1},{}^\RR\rho_{\lambda_1}+{}^\RR\rho_{2\lambda_1},
{}^\RR\rho_{\lambda_1}+{}^\RR\rho_{2\lambda_1})$ 
& $(11,11)$ & $\pi_3=\ZZ/3$ \\
$\G_2$ & 
$({}^\RR\rho_{\lambda_1}+{}^\RR\rho_{2\lambda_1},{}^\RR\rho_{\lambda_1},
{}^\RR\rho_{\lambda_1}+{}^\RR\rho_{2\lambda_1})$ 
& $(11,11)$ & $\pi_3=\ZZ/9$ \\
$\G_2$ & 
$({}^\RR\rho_{\lambda_1},{}^\RR\rho_{\lambda_1}+{}^\RR\rho_{2\lambda_1},
2\cdot{}^\RR\rho_{2\lambda_1})$ 
& $(11,11)$ & $\pi_3=\ZZ/4$ \\
$\G_2$ & 
$({}^\RR\rho_{\lambda_1}+{}^\RR\rho_{2\lambda_1},{}^\RR\rho_{\lambda_1},
2\cdot{}^\RR\rho_{2\lambda_1})$ 
& $(11,11)$ & $\pi_3=\ZZ/12$ \\
$\G_2$ & 
$({}^\RR\rho_{\lambda_1},{}^\RR\rho_{\lambda_1}+{}^\RR\rho_{2\lambda_1},
{}^\RR\rho_{6\lambda_1})$ 
& $(11,11)$ & $\pi_3=\ZZ/28$ \\
$\G_2$ & 
$({}^\RR\rho_{\lambda_1}+{}^\RR\rho_{2\lambda_1},{}^\RR\rho_{\lambda_1},
{}^\RR\rho_{6\lambda_1})$ 
& $(11,11)$ & $\pi_3=\ZZ/84$ \\
\end{tabular}
\end{center}
\end{Num}

%\newpage

\begin{Thm}
\label{Case(I)SemiSimple}
Let $G/H$ be a 1-connected homogeneous space of a compact connected
Lie group $G$. Assume that 
\[
\bfH^\bullet(G/H)=\EA_\ZZ(u,v)
\]
is an exterior algebra on two generators of odd degrees
$\deg(u)=n_1,\deg(v)=n_2$, with $n_1,n_2\geq 3$.
Assume moreover that the action is irreducible, and not split, and
that $G$ is not almost simple.
Then $(G,H)$ is one of the following pairs.

\begin{center}
\begin{tabular}{llll} 
\multicolumn{4}{c}{\fbox{$G$ semisimple}}  \\ 
\\
$G=K_1\times K_2$ & $H=H_1\cdot H_0$ & $\Cen_G(H)^\circ$ 
& $G/H$ \\
\hline
$\Sp(n)\times\SU(3)$ & $\Sp(n-1)\cdot\Sp(1)$   
& $\U(1)$ & $\SS^{4n-1}\times\SS^5$ \\
\hfill $n\geq 2$ \\
\hline
$\Sp(n)\times\Sp(2)$ & $\Sp(n-1)\cdot\Sp(1)$  
& $\Sp(1)$ & $\SS^{4n-1}\times\SS^7$ \\
\hfill $n\geq 2$ \\
\end{tabular}
\end{center}
In both series, the group $H_0$ is the image of $\Sp(1)$ in
$\Cen_G(H_1)^\circ=\Sp(1)\times K_2$ under the diagonal
embedding $h\mapstoo(h,h)$.
\qed
\end{Thm}

%\newpage

\section{The non-split case (II): 
$\bfH^\bullet(X)=\ZZ[a]/(a^2)\otimes\EA_\ZZ(w)$.}

We use the same conventions  as in the last section.
This time, we have
\begin{align*}
1& = \rk(G)-\rk(H) \\
&= k_1-h_1+k_2-h_2-h_0 \\
0&\leq k_1-h_1-h_0 \\
0&\leq k_2-h_2-h_0,
\end{align*}
whence $h_0=1$. Thus, $H_0$ is either a 1-torus, or of type 
$\fa_1$, and $k_i-h_i-h_0=0$ for $i=1,2$.

\begin{Lem}If $n_1=4$, $H_0$ is of type $\fa_1$, and
$(K_1,H_1)$ is one of the pairs in 
\ref{HomogeneousRationalSphereWithLargeCentralizer}.
Moreover, $(K_2,H_2)=(\Sp(2),\Sp(1))$.

\begin{proof}
From the exact sequence
\[
\underbrace{\pi_2(G/H)}_0\too
\pi_1(H)\too\pi_1(G)\too
\underbrace{\pi_1(G/H)}_0
\]
we see that $H$ is semisimple. Moreover,
\[
\underbrace{\pi_4(G)}_{\text{finite}}\too
\underbrace{\pi_4(G/H)}_\ZZ\too
\pi_3(H)\too\pi_3(G)\too\underbrace{\pi_3(G/H)}_0
\]
implies that $\rk(\pi_3(H))=\rk(\pi_3(G))+1=3$.
Thus, $H$ is semisimple with three almost simple factors. 
One of them has to be $H_0$, hence $H_0$ is
of type $\fa_1$.

Thus, $P^k_{H_i} \oot P^k_{K_i}$ is
a surjection in all degrees $k\geq 5$ (and in fact also in degree 3, since
the map induced on $\pi_3$ is not trivial); 
thus, $(K_1,H_1)$ and $(K_2,H_2)$ are among the pairs determined in 
\ref{HomogeneousRationalSphereWithLargeCentralizer}.
Moreover, $\rk(\pi_7(G/H))\in\{1,2\}$,
cp.~\ref{Case(II)General}, hence
$\rk(\pi_7(G))-\rk(\pi_7(H))=\rk(\pi_7(G/H))\geq 1$.
We may assume that $\rk(\pi_7(K_2))-\rk(\pi_7(H_2))\geq 1$.
Thus $(K_2,H_2)=(\Sp(2),\Sp(1))$.
\end{proof}
\end{Lem}

\begin{Lem}The case $n_1\geq 6$ is not possible.

\begin{proof}
In this case 
\[
\rk(\pi_i(H))=\rk(\pi_i(G))=2\quad\text{ for }i=1,2,3,
\]
hence $H$ is semisimple with two almost simple factors. We may assume that 
$H_2=1$. Note that $h_0=1$, hence $k_2=1$. Therefore $K_2$ is of type 
$\fa_1$. But then $\pr_2(H_0)=K_2$, a contradiction to the irreducibility 
of the action.
\end{proof}
\end{Lem}
We classify the remaining possibilities.

\begin{Num}\textsc{$K_1=\Sp(n)$, $n\geq 2$}\psn
Recall that there is a fibre bundle
\begin{diagram}
\llap{$\SS^{4n-1}={}$} \Sp(n)/\Sp(n-1) & \rTo & G/H \\
&& \dTo \\
&& \Sp(2)/\Sp(1)\times\Sp(1) \rlap{${}=\SS^4$.}
\end{diagram}
In fact, this bundle is the Whitney sum of $n$ copies of
the quaternionic Hopf bundle over $\SS^4=\HH\mathrm{P}^1$.
\end{Num}

\begin{Num}\textsc{$K_1=\G_2$}\psn
There is a similar fibre bundle
\begin{diagram}
\G_2/H_1 & \rTo & G/H \\
&& \dTo \\
&& \Sp(2)/\Sp(1)\times\Sp(1) \rlap{${}=\SS^4$.}
\end{diagram}
By \ref{SU(2)inG2}, 
the fibre has the same $\ZZ/2$-cohomology as the Stiefel manifold
$V_2(\RR^7)$. Therefore the
$\ZZ/2$-Leray Serre spectral sequence of the bundle collapses, and thus
\[
\bfH^\bullet(G/H;\ZZ/2)\isom\EA_{\ZZ/2}(x_5,x_7)
\otimes(\ZZ/2)[x_4,x_6]/(x_4^2,x_6^2).
\]
\end{Num}

%\newpage
\begin{Thm}
\label{Case(II)SemiSimple}
Let $G/H$ be a compact 1-connected homogeneous space of a compact
connected Lie group $G$. Assume that $G/H$ has the same
integral cohomology as $\SS^{n_1}\times\SS^{n_2}$, for 
$n_1\geq 4$ even and $n_2>n_1$ odd.
If the action of $G$ is irreducible and not split, and if 
$G$ is not simple, then $n_1=4$ and $(G,H)$ is one of the following pairs.

\begin{center}
\begin{tabular}{lllll} 
\multicolumn{5}{c}{\fbox{$G$ semisimple}}  \\ 
\\
$G=G_1\times G_2$ & $H=(H_1\times H_2)\cdot H_0$ 
& $\Cen_G(H)^\circ$ & $(n_1,n_2)$ & $G/H$ \\
\hline
$\Sp(n)\times\Sp(2)$ & $(\Sp(n-1)\times\Sp(1))\cdot\Sp(1)$   
& 1 & $(4,4n-1)$ & $S(n\eta_\HH)$ \\
\hfill $n\geq 2$ \\
\end{tabular}
\end{center}
The group $H_0$ is the image of $\Sp(1)$ in
$\Cen_G(H_1\times H_2)^\circ=\Sp(1)\times\Sp(1)$ under the diagonal
embedding $h\mapstoo(h,h)$. The space $G/H$ is the sphere bundle
$S(n\eta_\HH)$ of the Whitney sum $n\eta_\HH$
of $n$ copies of the quaternionic Hopf bundle $\eta_\HH$
over $\SS^4$.
\qed
\end{Thm}

%%%%%%%%%%%%%%%%%%%%%%%%%%%%%%%%%%%%%%%%%%%%%%%%%%%%%%%%%%%%%%%%%%%%%%%%
%                                                                      %
%                                                                      %
%                   Compact homogeneous quadrangles                    %
%                                                                      %
%                          Linus Kramer                                %
%                                                                      %
%                           Memoirs AMS                                %
%                                                                      %
%                         Wuerzburg 2000                               %
%                                                                      %
%                                                                      %
%                            CHQ7.tex                                  %
%                                                                      %
%                                                                      %
%                                                                      %
%%%%%%%%%%%%%%%%%%%%%%%%%%%%%%%%%%%%%%%%%%%%%%%%%%%%%%%%%%%%%%%%%%%%%%%%
\chapter{Homogeneous compact quadrangles}

The principal aim of topological geometry is to classify
reasonable geometries in terms of their automorphism groups.
'Reasonable' geometries are linear spaces (e.g. projective,
affine or hyperbolic geometries, or, more generally,
stable planes),
circle geometries (Laguerre or M\"obius geometries) and finally
Tits buildings, to us the most important class of geometries.

The fundamental theorem of projective geometry asserts that
a projective space of rank at least 3 (i.e.~a projective space
which is not a projective plane) is coordinatized by a
field or skew field. The key step in the proof is to show that
such a projective space satisfies the Desargues condition.
There is a similar result due to Tits \cite{TitsLNM}
for spherical buildings. An irreducible spherical building of rank
at least 3 satisfies the so-called Moufang condition;
furthermore, all irreducible
spherical Moufang buildings of rank at least 2
were determined by Tits \cite{TitsLNM} and Tits-Weiss \cite{Tits-Weiss}.

Moufang buildings are rather important in many branches
of mathematics. Just to mention a few examples,
differential geometers use them in the proof of the rigidity theorems
of Mostow \cite{Mostow}, Gromov \cite{BGS}, Kleiner-Leeb
\cite{KleinerLeeb} and Leeb \cite{Leeb},
and in the classification of isoparametric submanifolds, see
Thorbergsson \cite{Tho91};
they play a r\^ole in connection with $S$-arithmetic groups and
group cohomology, see Rohlfs-Springer \cite{RohlfsSpringer}
and Abramenko \cite{Abramenko},
and in model theory in connection with simple superstable groups,
see Borovik-Nesin \cite{BorovikNesin}
and Kramer-Tent-Van Maldeghem \cite{KTV}.

In Tits' classification,
the assumption that the building has rank at least 3 cannot be dropped;
there exist uncountably many 'wild' buildings of rank 2, even
with large automorphism groups,
see Tits \cite{Tits77} and Tent \cite{Tent99}.
Thus it is desirable to have a criterion which ensures that an
infinite
building of rank 2 is Moufang. Indeed, such a criterion exists
in the topological category.
\psn
\textbf{Theorem A}
{\em  A compact connected irreducible spherical building of rank
at least $2$ with a flag transitive automorphism group is Moufang.}
\medskip

This is proved in a series of papers by
Grundh\"ofer-Knarr-Kramer \cite{GKK95} \cite{GKK98} \cite{GKK00}.
Note that it suffices to consider spherical buildings of rank 2
to prove the theorem.
No similar classification is presently known for finite or
zero-dimensional (totally disconnected) spherical buildings of rank 2.
(The finite Moufang buildings are closely related to finite
simple groups; totally disconnected buildings appear rather naturally
in connection with valuations and algebraic groups over local fields.)
In the course of the proof, the closed connected flag transitive groups
are also explicitly determined. It turns out that such a group $G$
is either the little projective group of the building (the group
generated by the root groups) or a certain compact subgroup of the little
projective group. The compact connected flag transitive groups were
determined by Eschenburg-Heintze \cite{EschHein}; all closed connected
flag transitive groups are determined
in Grundh\"ofer-Knarr-Kramer \cite{GKK95} \cite{GKK98} \cite{GKK00}.

In this chapter we extend Theorem A to
compact connected buildings with a less transitive automorphism group,
that is, an automorphism group which is only transitive on vertices
of a certain type. By Tits' results about buildings of higher rank,
we have only to consider buildings of
rank~2. Spherical buildings of rank two are commonly called
generalized polygons. In building terminology, they are certain
1-dimensional numbered simplicial complexes, i.e.~bipartite graphs:
there are precisely two types of vertices in such
a building which we call points and lines.
The axioms of a building are symmetric in the sense that
we might as well call the lines points and the points lines.
Thus, there is no loss in generality if one studies either
point or line transitive actions.
The main result of this chapter can be stated as follows.
\psn
\textbf{Theorem B}
{\em A compact connected building of type $C_2$ with a
point transitive automorphism group is Moufang, provided that
the point space is $9$-connected (\ref{MainThmCh7-1}, \ref{MainThmCh7-2}}).
\medskip

It is known that there exist line homogeneous compact connected
generalized quadrangles with a 6-connected line space which are
not Moufang, see Kramer \cite{Kr98};
in particular, point or line transitivity does not imply
the Moufang property in general. This indicates already that
the problem is more difficult than the flag transitive case
considered in Grundh\"ofer-Knarr-Kramer \cite{GKK95} \cite{GKK98}.

A result of Knarr \cite{Kn90} 
and the author \cite{Kr94} says that a finite dimensional
compact connected polygon is a building of type $A_2$, $C_2$ or $G_2$,
i.e.~a projective plane, a generalized quadrangle (which is the
same as a polar space of rank 2) or a generalized hexagon.
A point transitive automorphism group implies that the
topological dimension is finite, so this result can be applied.
The three types of buildings  then have to be considered separately.
The $A_2$-case was done by L\"owen and Salzmann about 20 years ago
\cite{Salzmann} \cite{Lowen}.
\psn
\textbf{Theorem C}
{\em  A compact connected projective plane or generalized hexagon
with a
point transitive automorphism group is Moufang, see
Salzmann \cite{Salzmann}, L\"owen \cite{Lowen} or
Salzmann {\sl et al.} \cite{CPP95} 63.8, and Kramer \cite{Kr94} Ch.~5.
A compact connected generalized quadrangle with a point
transitive automorphism group is Moufang, provided that
the point and line space have the same dimension, see Kramer \cite{Kr94}
Ch.~5.}

\medskip
Combining Theorems A, B and C, we obtain the following result.
\psn
\textbf{Theorem D}
{\em Let $\Delta$ be an irreducible spherical compact connected 
building of rank at least $2$. Assume that the automorphism group
of $\Delta$ acts transitively on one type of vertices.
If the building is of type $C_2$, assume in addition that
either
(a) the vertex set in question is $9$-connected, or
(b) the two vertex sets of the building have the same dimension, or
(c) that the action is chamber transitive.
Then $\Delta$ is the Moufang building associated to a non-compact
real simple Lie group.}

\medskip
Before I comment on the proof of Theorem B,
I would like to mention some recent developments.
Biller studied in his Ph.D. thesis \cite{Biller} (not necessarily
transitive) compact group actions on compact generalized quadrangles.
In particular, he characterized the hermitian quadrangles over
the quaternions in terms of the size of their automorphism groups.
His work contains many new results
and will certainly be an important source for further research in this
direction. We use some of his results in Section  4 of this chapter;
my original proof of Theorem \ref{(4,4n-5)} (as given in
\cite{KramerHabil} Section 7.D) is incorrect, and the new proof given here
relies on Biller's ideas.
Recently,
Bletz and Wolfrom started to investigate point transitive actions
on compact connected $(m_1,m_2)$-quadrangles, with $m_1=1,2$.
Finally, Immervoll \cite{StIm} proved that all isoparametric hypersurfaces
with four distinct principal curvatures are generalized quadrangles
(cp.~Chapter 8). This has been a difficult open problem for quite some
time. His proof is a clever combination of transversality arguments
and the algebraic machinery developed by Dorfmeister-Neher in the early 80s.

\medskip
\emph{Some remarks on the proof of Theorem B.}
In order to prove Theorem B, one has to consider generalized
quadrangles whose point and line spaces have different dimensions.
The cohomology ring of the point space of
such a quadrangle is known (Kramer \cite{Kr94}, Strau{\ss} \cite{Strau96}),
and it turns out that it looks
in most cases like the cohomology of a product of spheres.
Hence we can apply our classification of homogeneous
spaces.

Not every homogeneous space in our list \ref{MainTheorem}
corresponds to a compact quadrangle. The point stabilizer
cannot have large normal subgroups, cp.~\ref{NoLargeNormalGroup} below.
This excludes for example the compact symmetric space $\E_6/\Ffour$.

Besides the Stiefel manifold $V_2(\FF^n)$, $\FF=\RR,\CC,\HH$
we have to consider certain homogeneous sphere bundles over spheres, and
also products of homogeneous spheres.

For the Stiefel manifolds we show the following: if the dimension
of the corresponding vector space $\FF^n$ is at least 5
(at least 4 for $\FF=\HH$ ), then
there is a unique quadrangle compatible with the transitive
group action which 'lives' on the Stiefel manifold.

A similar result is proved for one of the three series of homogeneous
sphere bundles. As mentioned before, this proof relies on Biller's
thesis \cite{Biller}.

For the products of homogeneous spheres we show that one of the
two factors has to be $\SS^3$ with the regular $\Sp(1)$-action,
and we obtain some restrictions on the other factor. The examples
of Ferus-Karcher-M\"unzner \cite{FeKaMu81}
show that there is at least one
family of  non-Moufang quadrangles which 'lives' on such homogeneous
products of spheres, see Example \ref{Quad(3,4n-4)} in Chapter 8.

The relevant results about topological generalized quadrangles can
be found in
Grund\-h\"o\-fer-Knarr \cite{GK90},
Grundh\"ofer-Van Maldeghem \cite{GvM90},
Grundh\"ofer-L\"owen \cite{GL95},
Grundh\"ofer-Knarr-Kramer \cite{GKK95} \cite{GKK98},
Knarr \cite{Kn90}, 
Kramer-Van Maldeghem \cite{Ov1} and Kramer \cite{Ov2} \cite{Ov3},
Schroth \cite{Schroth},
and in particular in Kramer \cite{Kr94}; 
Van Maldeghem's monograph \cite{HvM98} is the authoritative source for
geometric properties of generalized polygons; Chapter 9 in his book
summarizes many
results about topological polygons.

%\newpage

\section{Generalized quadrangles and group actions}

A \emph{point-line geometry} is a triple
\[
\frak{G}=\PLF
\]
consisting of a set $\P$ of \emph{points}, a set $\L$ of \emph{lines}
and a set $\F\SUB\P\times\L$ of \emph{flags}. If a pair
$(p,\ell)\in\F$ is a flag, then $p$ and $\ell$ are called
\emph{incident};
one also says that the line $\ell$ passes through the point $p$,
or that $p$ lies on $\ell$. 
(One can turn $\frak G$ into a \emph{bipartite graph} $(V,E)$
with vertex set
$V=\P\cup\L$ and edge set $E=\{\{p,\ell\}|\ (p,\ell)\in\F)\}$;
this is the building point of view.)
A \emph{$k$-chain} joining
$x_0,x_k\in\P\cup\L$ is a sequence
$(x_0,\ldots,x_k)\in(\P\cup\L)^{k+1}$ with the property that $x_i$ is
incident with $x_{i-1}$ for $1\leq i\leq k$. The \emph{distance} 
of $x_0$ and $x_k$ is 
\[
d(x_0,x_k)=k
\]
if there is a $k$-chain joining $x_0$ and $x_k$, but no $j$-chain
joining $x_0$ and $x_k$ for $j<k$. Put
\[
D_k(x)=\{y\in\P\cup\L|\ d(x,y)=k\}.
\]
A point-line geometry is called \emph{thick} if $|D_1(x)|\geq 3$ for
all $x\in\P\cup\L$. We put 
\[
x^\perp=\{x\}\cup D_2(x);
\]
if $x$ is a point, then this is the set of all points which have
a line in common with $x$.

Automorphisms are defined in the obvious way; an automorphism is a
permutation of the point and the line set which preserves the incidence
relation. In graph theoretic terms, an automorphism is a graph automorphism
which preserves the coloring of the bipartite graph.
\begin{Num}\textsc{Generalized quadrangles}\psn
A thick point-line geometry is called a \emph{generalized quadrangle}
if $\F\neq\P\times\L$, and if
\[
|p^\perp\cap D_1(\ell)|=1
\]
holds for all pairs $(p,\ell)\in(\P\times\L)\setminus\F$.
Then there is a unique $3$-chain $(p,h,q,\ell)$ joining $p$ and $\ell$. 
Put 
\[
q=\proj_\ell p\quad\text{ and }\quad h=\proj_p\ell.
\]
The picture below shows the corresponding 3-chain.
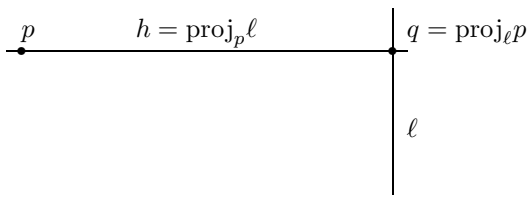
\begin{center}
\unitlength=0.19mm
\begin{picture}(300,150)
\put(10,110){\line(1,0){280}}
\put(20,110){\circle*{5}}
\put(280,110){\circle*{5}}
\put(280,10){\line(0,1){130}}
\put(20,120){$p$}
\put(290,120){$q=\proj_\ell p$}
\put(100,120){$h=\proj_p \ell$}
\put(290,50){$\ell$}
\end{picture}
\end{center}
\end{Num}
Note that a generalized quadrangle contains no digons
(two point are joined by at most one line and two lines intersect
in at most one point) and no triangles. In particular, a line $\ell$ is
uniquely determined by the point row $D_1(\ell)$.
For an example of a generalized quadrangle, see \ref{ClassQuad} below.
\emph{A generalized quadrangle is the same as a spherical building of
type $C_2$.}
\begin{Num}\textsc{Group actions and reconstruction}\psn
\label{Reconstruction}%
An \emph{action} of a group $G$ on a generalized quadrangle
$\frak G$ is a homomorphism 
\[
G\too\Aut(\frak G).
\]
Suppose that 
\begin{description}
\item[\textbf{Rec$_1$}]
$G$ acts transitively on the point set $\P$. 
\end{description}
The question is whether $\frak G$ can be \emph{reconstructed}
from the action of $G$. 
Let $p\in\P$ and assume that 
\begin{description}
\item[\textbf{Rec$_2$}]
the stabilizer $G_\ell$ acts transitively on $D_1(\ell)$ for every
$\ell\in D_1(p)$. 
\end{description}
Then $D_1(\ell)=G_\ell\cdot p\SUB\P$.
A quadrangle satisfying \textbf{Rec$_1$} and
\textbf{Rec$_2$} is an example of what Stroppel \cite{Stro92} \cite{Stro93}
calls a 
\emph{sketched geometry}. Suppose we know the collection
\[
\mathcal{G}=\{G_\ell|\ \ell\in D_1(p)\}
\]
of all stabilizers of lines passing through $p$.
Then we know all point rows
through $p$; they are the sets $G_\ell/G_\ell\cap G_p\SUB G/G_p$.
The collection of their $G$-translates can be canonically
identified with the line set $\L$. Thus, the quadrangle is uniquely 
determined by the triple
\[
(G,G_p,\mathcal{G}).
\]
In fact, put 
\[
\P'=G/G_p
\quad\text{ and }\quad\L'=\bigcup\left\{gHG_p|\ g\in G,\ H\in\mathcal{G}\right\}
\]
and define the incidence as inclusion '$\SUB$'.
It is easy to show that the resulting geometry $(\P',\L',\SUB)$ is
$G$-isomorphic to the original quadrangle $\frak G$.
See Stroppel \cite{Stro92} \cite{Stro93} for more results in this direction.
\end{Num}
We will also use the following elementary result.

\begin{Lem}
\label{TrivialActionLem}
Let $\frak G$ be a point-line geometry, with the property that every line
$\ell$ is uniquely determined by its point row $D_1(\ell)$.
Assume that $d(x,y)<\infty$ holds
for all $x,y\in\P\cup\L$, and that $G\SUB\mathrm{Aut}(\frak G)$
acts transitively on $\P$. Let $p\in\P$ and let $N\SUB G_p$ be a 
subgroup with the following two properties:
\begin{enumerate}
\item $N$ acts trivially on $D_1(p)$
and on $D_1(\ell)$, for all $\ell\in D_1(p)$
(so $N$ fixes $p^\perp$ pointwise).
\item If $gNg^{-1}\SUB G_p$, for some $g\in G$, then
$gNg^{-1}=N$.
\end{enumerate}
Then $N=1$.

\begin{proof}
Let $q$ be a point which is collinear with $p$, i.e.~$q\in D_2(p)$.
Since $G$ acts transitively on $\P$, the stabilizer $G_q$ is conjugate
to $G_p$, say $G_q=gG_p g^{-1}$. Then $gNg^{-1}$ fixes $D_2(q)$
elementwise; in particular, $gNg^{-1}\SUB G_p$.
Thus $N=gNg^{-1}$. It follows inductively that 
$N$ is normal in $G$, and thus $N=1$, because the action of $G$ on $\P$
is effective (here we use that assumption that lines are determined
by their point rows).
\end{proof}
\end{Lem}

%\newpage

\section{Compact quadrangles}

Suppose that $\frak G=\PLF$ is a generalized quadrangle, and that
the sets $\P$ and $\L$ carry Hausdorff topologies. If the map
$(p,\ell)\mapstoo (q,h)=(\proj_\ell p,\proj_p\ell)$ is
continuous on $(\P\times\L)\setminus\F$, then
$\frak G$ is called a \emph{topological quadrangle}; if $\P$ and $\L$ are
in addition compact, then $\frak G$ is called a \emph{compact quadrangle}.

\begin{Prop}
Let $\frak G$ be a generalized quadrangle, and suppose that
$\P$ and $\L$ are compact Hausdorff spaces. Then $\frak G$ is
a compact quadrangle if and only if $\F\SUB\P\times\L$ is
closed (or, equivalently, compact). 

\begin{proof}
See Grundh\"ofer-Van Maldeghem \cite{GvM90} 2.1.
\end{proof}
\end{Prop}
Concerning the set theoretic topology of a compact quadrangle,
the following result is important.
\begin{Prop}
Let $\frak G$ be a compact quadrangle. The point space $\P$ and
the line space $\L$ are second countable, separable and metrizable.
If the point space $\P$ is connected, then $\P$, $\L$, every point
row $D_1(\ell)$ and every line pencil $D_1(p)$ is connected
and locally contractible.

\begin{proof}
See Grundh\"ofer-Knarr \cite{GK90} 3.1 and 4.1, or
Grundh\"ofer-Knarr-Kramer \cite{GKK95} 1.5 and 1.6.
\end{proof}
\end{Prop}

Recall the definition of the \emph{covering dimension}
$\dim(X)$ of a normal space $X$. If every finite open covering
of $X$ has a refinement such that every point of $X$ is contained
in at most $n+1$ sets of the refinement, then $\dim(X)\leq n$,
and $\dim(X)=n$ if $X$ has at most dimension $n$, but not dimension
$n-1$. Equivalently, $\dim(X)\leq n$ holds if and only if the
extension problem
\begin{diagram}[width=3em,height=2em]
A & \rTo^\alpha & \SS^n \\
\dInto & \ruDotsto \\
X \\
\end{diagram}
has a solution for every map $\alpha$ defined on a closed subspace
$A\subseteq X$.
There are other topological notions of dimension,
most notably the small or large inductive dimension, and the
cohomological dimension. 
See Salzmann \emph{et al.} \cite{CPP95}~92 for a discussion of
dimension functions and further references.
For 'good' spaces, these notions of
dimension agree; in particular, an $n$-manifold has dimension $n$.
More generally, an \emph{integral ENR $n$-manifold} has dimension $n$.
\begin{Def}
Let $X$ be an ENR (euclidean neighborhood retract, cp.~Dold \cite{Dold}
IV.8 for properties of such spaces). If there exists a number $n$ such that
\[
\bfH_k(X,X\setminus\{x\})\cong
\begin{cases}\ZZ&\text{ for }k= n\\
0&\text{ for }k\neq n
\end{cases}
\]
holds for all $x\in X$, then $X$ is called an 
\emph{integral ENR manifold}.
Integral ENR manifolds are \emph{generalized manifolds},
i.e.~$cm_R$s and $hm_R$s (for any integral domain $R$) as studied in Bredon
\cite{BredonSheaf}, and it can be shown that $\dim(X)=n$.
This is due to the fact that an ENR is locally contractible
(so singular homology coincides with Borel-Moore homology)
and a result of Bredon \cite{BredonSheaf}. 
(An ENR is locally contractible and in particular $clc_R^\infty$;
since $\dim(X)$ is finite, $\dim_R(X)$ is also finite. 
The local homology groups in Borel-Moore
homology are isomorphic to the local homology groups in singular homology.
The result now follows from Theorem V.16.8 in 
Bredon \cite{BredonSheaf}.)
Integral ENR manifolds
are a very convenient class of spaces,
because they satisfy all standard topological
assumptions in the theory of compact transformation groups.
\end{Def}
Note also that
an ANR of finite covering dimension is the same as an ENR.

A compact quadrangle is called \emph{finite dimensional} if the covering
dimension of $\P$ is finite and positive, or, equivalently, if
the covering dimension of $\L$ is finite and positive.
The following theorem summarizes the most important topological
properties of compact connected finite dimensional quadrangles.

\begin{Thm}
Suppose that $\frak G$ is a finite dimensional compact
quadrangle. Let $(p,\ell)\in\F$ and put 
$(m_1,m_2)=(\dim(D_1(\ell)),\dim(D_1(p)))$.

Then the following hold.
\begin{enumerate}
\item The spaces $\P$, $\L$ and $\F$ are ENRs (euclidean neighborhood
retracts) and in particular ANRs (absolute neighborhood retracts).
\item The spaces $\P$, $\L$ and $\F$, as well as the point rows
and the line pencils are integral (locally and globally homogeneous)
ENR-manifolds.
\item The inclusions $\{p\}\SUB
D_1(\ell)\SUB p^\perp\SUB\P$ are cofibrations
for every pair $(p,\ell)$ (and dually).
\item There are homotopy equivalences
\begin{gather*}
D_1(\ell)\homot\SS^{m_1} \\
p^\perp/D_1(\ell)\homot\SS^{m_1+m_2} \\
\P/p^\perp\homot\SS^{2m_1+m_2}
\end{gather*}
\item $D_1(\ell)$, $p^\perp$ and $\P$ are $(m_1-1)$-connected.
\item Either $m_1=m_2\in\{1,2,4\}$, or $1\in\{m_1,m_2\}$, or
$m_1+m_2$ is odd.
\end{enumerate}
\begin{proof}
This is proved in Knarr \cite{Kn90},
Grundh\"ofer-Knarr \cite{GK90} Sec.~4 and in the present generality
in Kramer \cite{Kr94} Thm.~3.3.6.
\end{proof}
\end{Thm}
The numbers $(m_1,m_2)$ are
called the \emph{topological parameters} of $\frak G$. A compact connected
finite dimensional quadrangle with parameters $(m_1,m_2)$
is called a \emph{compact connected $(m_1,m_2)$-quadrangle} for short.
\begin{Num}\textsc{Example}\label{ClassQuad}
Let $\FF=\RR,\CC,\HH$ and consider the hermitian form in $n+1$ variables
over $\FF$
\[
h(x,y)=-\bar x_0 y_0-\bar x_1y_1+\sum_{k=2}^n\bar x_ky_k.
\]
Let $\P$ denote the collection of all 1-dimensional totally isotropic
subspaces, and $\L$ the collection of all 2-dimensional totally isotropic
subspaces of $\FF^{n+1}$
(a subspace $V$ of $\FF^{n+1}$ is called \emph{totally isotropic}
if the form $h$ vanishes identically
on $V$). If $n\geq 4$, or if $\FF\neq\RR$ and
$n=3$, then 
\[
\mathsf{H}_n(\FF,\RR)=(\P,\L,\subseteq)
\]
is a compact 
connected quadrangle, the \emph{standard hermitian quadrangle} over
$\FF^{n+1}$; one puts
$\mathsf{Q}_n(\RR)=\mathsf{H}_n(\RR,\RR)$
(our notation is essentially the same as in Van Mal\-de\-ghem \cite{HvM98}).
Let $d=\dim_\RR(\FF)=1,2,4$. The quadrangle
$\mathsf{H}_n(\FF,\RR)$ has parameters $(d,d(n-2)-1)$.
Note that there are natural inclusions
\begin{diagram}[width=4em,height=3em]
 &&
\mathsf{Q}_4(\RR) & \rInto &
\mathsf{Q}_5(\RR) & \rInto &
\mathsf{Q}_6(\RR) & \rInto & \cdots \\
&& \dInto && \dInto && \dInto \\
\mathsf{H}_3(\CC,\RR) & \rInto &
\mathsf{H}_4(\CC,\RR) & \rInto &
\mathsf{H}_5(\CC,\RR) & \rInto &
\mathsf{H}_6(\CC,\RR) & \rInto & \cdots \\
\dInto && \dInto && \dInto && \dInto \\
\mathsf{H}_3(\HH,\RR) & \rInto &
\mathsf{H}_4(\HH,\RR) & \rInto &
\mathsf{H}_5(\HH,\RR) & \rInto &
\mathsf{H}_6(\HH,\RR) & \rInto & \cdots 
\end{diagram}
In the limit one obtains non-compact topological
quadrangles with compact point
rows and contractible infinite dimensional
line pencils modeled on $\SS^\infty$.
\end{Num}
Let $\frak G=\PLF$ be a generalized quadrangle, let
$\P'\subseteq\P$ and $\L'\subseteq\L$ be non-empty subsets. If
$\frak G'=(\P',\L',\F\cap(\P'\times\L'))$ is again a generalized
quadrangle, then we call $\frak G'$ a \emph{subquadrangle} of
$\frak G$.
If $D_1(\ell)\subseteq\P'$ holds for all $\ell\in\L'$, then
the subquadrangle is called \emph{full}.
Here is an example of a full and a non-full subquadrangle.
\begin{diagram}[width=6em,height=3em]
\mathsf{H}_3(\CC,\RR) & \rInto^{\text{not full}} & \mathsf{H}_3(\HH,\RR) \\
\dInto_{\text{full}} \\
\mathsf{H}_7(\CC,\RR)
\end{diagram}
The following result is proved in Kramer-Van Maldeghem \cite{Ov1} 
Thm.~4.1.

\begin{Prop}
\label{OvoidProp}
Let $\frak G'\subset\frak G$ be a full, compact, and proper
(i.e.~$\frak G'\neq\frak G$) subquadrangle of the compact connected
$(m_1,m_2)$-quadrangle $\frak G$. Let $(m_1,m_2')$ denote the
parameters of $\frak P'$. Then $m_2'\neq 0$ (so $\frak G'$ is
connected), and
\[
m_1+m_2'\leq m_2.
\]
\begin{proof}
The idea of the proof is to show that $\frak G'$ contains a compact
ovoid $\mathcal{O}$ which injects into a line pencil of $\frak G$.
The dimension of such an ovoid is $m_1+m_2'$, whence $m_1+m_2'\leq m_2$.
For details see Kramer-Van Maldeghem Thm.~4.1 \cite{Ov1}.
The fact that the point rows of $\frak G'$ are connected implies
that the line pencils are also connected, whence $m_2'\neq 0$.
\end{proof}
\end{Prop}
The cohomology of finite dimensional quadrangles is known.

\begin{Num}\textsc{Cohomology of finite dimensional compact quadrangles}\psn
\label{TopologyOfQuad}%
Let $\frak G$ be a compact connected $(m_1,m_2)$-quadrangle. Then
\begin{align*}
\dim(\F)&=2(m_1+m_2) \\
\dim(\P)&=2m_1+m_2 \\
\dim(\L)&=2m_2+m_1.
\end{align*}
If $m_1+m_2$ is odd, then
\begin{align*}
\bfH^\bullet(\P)&\isom\bfH^\bullet(\SS^{m_1}\times\SS^{m_1+m_2}) \\
\bfH^\bullet(\L)&\isom\bfH^\bullet(\SS^{m_2}\times\SS^{m_1+m_2}) \\
\bfH^\bullet(\F)&\isom\bfH^\bullet(\SS^{m_1}\times
\SS^{m_2}\times\SS^{m_1+m_2}).
\end{align*}
If $m_1=1$ and if $m_2>2$ is odd, then
\[
\bfH^\bullet(\L;\QQ)\isom\bfH^\bullet(\SS^{2m_2+1};\QQ).
\]
The group $\pi_{m_2}(\L)$ is generated
by the inclusion $\SS^{m_2}\homot D_1(p)\SUB\L$.
If $m_1+m_2$ is odd, then $\pi_{m_1}(\L)\isom\ZZ$, and
if $m_2$ is odd and $m_1=1$, then $\pi_{m_2}(\L)\isom\ZZ/2$.

These results are proved in Kramer \cite{Kr94} Ch.~3 and 6.4; 
they follow from M\"unzner \cite{Munz81},
cp.~also Strau{\ss} \cite{Strau96}.
\end{Num}
Let $\frak G$ be a compact generalized quadrangle (not necessarily
connected, not necessarily finite dimensional). We endow the group
$\mathrm{AutTop}(\frak G)$ of all continuous automorphisms of 
$\frak G$ with the compact-open topology.
\begin{Thm}[Burns-Spatzier]
The group $\mathrm{AutTop}(\frak G)$ is a locally compact
metrizable group.

\begin{proof}
This follows from Burns-Spatzier \cite{BurSpa87} 2.1, in
combination with Grund\-h\"o\-fer-Knarr-Kramer \cite{GKK95} 1.8 
(there, it is proved that the metrizability assumption of
Burns-Spatzier is superfluous), see also Bletz \cite{Bletz}.
\end{proof}
\end{Thm}
For example, the automorphism group of
$\mathsf{H}_n(\FF,\RR)$ is a finite
extension of the non-compact simple Lie group
$\mathrm{PU}_{2,n-1}(\FF)$.

An \emph{action} of a (Lie) group $G$ on a
compact quadrangle is a continuous homomorphism
\[
G\too\mathrm{AutTop}(\frak G).
\]
For example, the compact Lie group
$\U_2(\FF)\times\U_{n-1}(\FF)$ acts transitively on the points, lines
and flags of $\mathsf{H}_n(\FF,\RR)$.
To get further, we need some facts about (compact)
transformation groups.

\section{Some results about compact transformation groups}

We collect a couple of results which will be used in our
classification of point homogeneous quadrangles.
Let $G$ be a group acting on a set $X$. We denote the image of $G$ in
the symmetric group of $X$ by $G|_X$, and we put $G_X=\{g\in G|\ 
g=\id_X\}$. In other words, $G_X$ is the \emph{kernel of the
action}, $G|_X$ is the \emph{induced group}, and the sequence
\[
1\rTo G_X \rTo G \rTo G|_X \rTo 1
\]
is exact. The $G$-action on $X$ is effective if and only if $G_X=1$.
We put
\[
\mathrm{Fix}(G,X)=\{x\in X|\ G\cdot x=\{x\}\}.
\]
We will use Szenthe's solution of Hilbert's 5th problem:
\begin{Thm}[Szenthe]
\label{Szenthe}
Let $G$ be a locally compact second countable group. Suppose that
$G$ acts effectively and transitively on a locally compact,
locally contractible space $X$. Then $G$ is a Lie group and
$X\cong G/G_x$ is a manifold.

\begin{proof}
This is what is \emph{proved} in Szenthe \cite{Sze74} (there, a more
general result is claimed), cp.~the remarks in Salzmann \emph{et al.} 96.14.
and in Grundh\"ofer-Knarr-Kramer \cite{GKK95} Thm.~2.2.
\end{proof}
\end{Thm}

\begin{Prop}
\label{VerySym}
Let $G$ be a compact connected Lie group acting transitively and
effectively on a manifold $X$. Then 
\[
\dim(G)\leq\binom{\dim(X)+1}2.
\]
If equality holds, then $G=\SO(n+1)$ and $X=\SS^{n}$ or 
$G=\PSO(n+1)$ and $X=\RR\mathrm{P}^n$.

\begin{proof}
The idea of the proof is to introduce a $G$-invariant Riemannian
metric on $X$. The stabilizer $G_x$ of $x\in X$ acts effectively
on $T_xX$, since every point $y\in X$ can be joined by a geodesic with
$x$. Thus $G_x\subseteq\O(T_xX)$, whence $\dim(G_x)\leq\dim(\O(T_x))=
\binom{\dim(X)}2$, see Montgomery-Zippin \cite{MoZi} Thm.~6.2.5
and the following corollary, and Kobayashi-Nomizu \cite{KoNo}
Vol.1 Note 10, Thm.~1, p.~308.
\end{proof}
\end{Prop}
The \emph{type} of a $G$-orbit $G\cdot x$ in a set $X$ is the conjugacy
class of the stabilizer $G_x$. If $G$ is a compact Lie group, then
an orbit $G\cdot x$ is \emph{principal} if the following two conditions
hold.
\begin{description}
\item[\textbf{PO$_1$}]
For every $y\in X$, the stabilizer $G_y$ is conjugate to an overgroup
of $G_x$ (i.e.~there exists an element $g\in G$ such that $G_{g(x)}$ fixes
$y$).
\item[\textbf{PO$_2$}]
The set of all orbits of the same type as $G\cdot x$ is open and
dense in $X$.
\end{description}
It is not difficult to prove that principal orbits exist, provided that
$G$ acts smoothly on a (smooth) manifold, see Bredon \cite{BredTTG} 
IV Thm.~3.1,
tom Dieck \cite{tDieck} Thm.~5.14. However, we need the general result
which is essentially due to Montgomery-Yang \cite{MonYang}
, and which is stated
and proved in the present form in Biller \cite{Biller} Thm.~2.2.3.
\begin{Thm}
Let $G$ be a compact Lie group acting effectively on a connected
integral ENR-manifold $X$. Then there exist principal orbits.
Moreover, $G$ acts effectively on each principal orbit.

\begin{proof}
For $x\in X$ let $d(x)=\dim(G\cdot x)$ and $c(x)=|\pi_0(G_x)|$.
Let $X_r=\{x\in X|\ d(x)=r\}$ and $X_{r,v}=\{x\in X_r|\ c(x)=v\}$.
Choose $k$ as large as possible, such that $X_k\neq\emptyset$, and
let $u=\min\{c(x)|\ x\in X_k\}$. Let $Y=X_{k,u}$. The existence of
slices implies that $Y$ is open, see
Borel \emph{et al.} \cite{BorelTrans} VIII Cor.~3.10. By Montgomery-Yang
\cite{MonYang} Lemma 2, $Y$ is dense in $X$, cp.~Borel \emph{et al.}
\cite{BorelTrans} IX Lem.~3.2. Moreover, $Y/G$ is connected by
Borel \emph{et al.} \cite{BorelTrans} IX Lem.~3.4.
This together with the existence of slices readily implies that
$Y\too Y/G$ is a locally trivial fibre bundle,
and all stabilizers of points in $Y$ are conjugate.

If $g\in G$ fixes all points in some principal orbit, then clearly
$g$ fixes $Y$ elementwise, and thus it fixes $X$, because $Y$ is dense.
Therefore, the
action of $G$ on each principal orbit is effective.
\end{proof}
\end{Thm}

\begin{Cor}
\label{DimensionProp}
Let $G$ be a compact
connected Lie group acting effectively on a connected
$n$-dimensional integral ENR manifold $X$. Then
\[
\dim(G)\leq\binom{n+1}2.
\]
If equality holds, then $G$ acts transitively, and $X=\SS^n$ or
$X=\RR\mathrm{P}^n$, cp.~Proposition \ref{VerySym} above.
If the $G$-action is not transitive, then
\[
\dim(G)\leq\binom{n}2.
\]

\begin{proof}
Let $Z\subseteq X$ be a principal $G$-orbit. By Proposition \ref{VerySym},
$G$ acts effectively on the compact connected manifold $Z$, so
$\dim(G)\leq\binom{\dim(Z)+1}2$. If $Z\neq X$, then $\dim(Z)<\dim(X)$,
since otherwise $Z$ would be open in $X$.
\end{proof}
\end{Cor}
The inequalities can be improved using Mann's result \cite{Man67};
further results are obtained in Biller's thesis \cite{Biller}.

If every point is contained in a principal orbit, then there is
the following useful result due to Borel.
\begin{Prop}[Borel]
\label{BorelStab}
Let $G$ be a compact Lie 
group acting on a space $X$. Suppose that all $G$-orbits
have the same type, i.e.~that all stabilizers $G_x$ are conjugate to
some group $H\SUB G$. The set $\mathcal{H}$
of all $G$-conjugates of $H$ can be identified
with the homogeneous space
$G/\mathrm{Nor}_G(H)$; in this way, it becomes a compact manifold.
Then the map 
\[
X\too \mathcal{H},\quad x\mapstoo G_x
\]
is a continuous surjection (in fact a fibre bundle with $\mathrm{Fix}(H,X)$
as typical fibre).

\begin{proof}
This is proved in Bredon \cite{BredTTG} II.5.9.
\end{proof}
\end{Prop}
Let $p$ be a prime, let $\FF_p$ denote the field with
$p$ elements, and let $A$ be an elementary abelian $p$-group of rank $r$,
i.e.~$A\cong(\ZZ/p)^r$.

\begin{Thm}[Floyd and Borel]
\label{ElementAbelian}
Let $X$ a compact connected $n$-dimensional integral ENR manifold with the
same integral homology as $\SS^n$. Suppose that $A$ acts effectively
on $X$, and let $F=\mathrm{Fix}(A,X)$ denote the fixed point set of $A$.
Then $F$ has the same (sheaf theoretic) $\FF_p$-cohomology as a
sphere $\SS^k$, for $-1\leq k\leq n$ (where $\SS^{-1}=\emptyset$).
Note that $k<n$ if $r\geq 1$.

If $p$ is odd, then $n-k$ is even.

If $r\geq 1$, then there exists a subgroup $A'\subseteq A$ of index $p$
such that $\mathrm{Fix}(A',X)$ is a $\FF_p$-cohomology $k'$-sphere,
for $k'>k$
\begin{proof}
The first and second claim is proved in Borel \emph{et al.}
\cite{BorelTrans} IV 4.3, 4.4, 4.5. To prove the last claim we use
the following result due to Borel. Let $\mathcal{A}$ denote the set of all
subgroups of $A$ of index $p$ (so $|\mathcal{A}|=\frac{p^r-1}{p-1}$).
For each $H\in\mathcal{A}$, let $n_H$ denote the dimension of
the cohomology sphere fixed by $H$. Then
\[
\textstyle n-k=\sum\nolimits_{H\in\mathcal{A}}\left(n_H-k\right),
\]
cp.~Borel \emph{et al.} \cite{BorelTrans} XIII Thm.~2.3.
Put $k'=\max\{n_H|\ H\in\mathcal{A}\}$. Thus $(n-k)(p-1)\leq(p^r-1)(k'-k)$,
and therefore $k'>k$ (note that $k<n$ if $r\geq 1$).
\end{proof}
\end{Thm}

\section{Group actions on compact quadrangles}

\begin{Thm}[Transitive actions on compact quadrangles]\ \psn
\label{Trans1}%
Let $\frak G$ be a compact connected quadrangle. Suppose that
the topological automorphism group $\mathrm{AutTop}(\frak G)$
acts transitively on the point space $\P$. Then the following hold.
\begin{enumerate}
\item The group $\mathrm{AutTop}(\frak G)$ is a Lie group.
\item The quadrangle $\frak G$ has finite dimension; in particular, the 
topological parameters $(m_1,m_2)$ of $\frak G$ are defined, and
\ref{TopologyOfQuad} applies to $\frak G$.
\item If $m_1\geq 2$, then there exists a compact connected Lie subgroup
$G\SUB\linebreak\mathrm{AutTop}(\frak G)$ which acts transitively on $\P$.
\end{enumerate}

\begin{proof}
The first claim follows from Szenthe's solution of Hilbert's 5th
problem, see Theorem \ref{Szenthe} above. Being a homogeneous space,
$\P$ is a manifold and thus finite dimensional.
The last claim follows from Montgomery \cite{Mon50}
Cor. 3, because $\P$ is 1-connected.
\end{proof}
\end{Thm}
The following result is proved in Kramer \cite{Kr94}.
\begin{Thm}
Let $\frak G$ be a compact connected $(m_1,m_2)$-quadrangle and
assume that the topological automorphism group acts transitively on
the points or lines, as in \ref{Trans1}.
If $m_1=m_2$, then $\frak G$ is the real or complex symplectic
quadrangle.
\begin{proof}
See Kramer \cite{Kr94} Thm.~5.2.4 and Thm.~5.2.3.
\end{proof}
\end{Thm}
If $\frak G$ is a compact connected
$(m_1,m_2)$-quadrangle with $m_1\neq m_2$, then
the cohomology of the point space is as in \ref{TopologyOfQuad}.
\begin{Thm}
Let $\frak G$ be a compact connected $(m_1,m_2)$-quadrangle and
assume that the topological automorphism group acts transitively on
the points, as in \ref{Trans1}.
If $m_1\neq m_2$, and if $m_1\geq 3$, then there exists a 
compact connected
Lie group $G\SUB\mathrm{AutTop}(\frak G)$ such that $(G,G_p)$ is
one of the pairs in Theorem \ref{MainTheorem}.

\begin{proof}
Since $m_1\geq 3$, the point space $\P$ is 1-connected. By Theorem
\ref{Trans1}, there exists a compact Lie group $G$ contained in the
automorphism group which acts transitively on $\P$. Moreover,
\ref{TopologyOfQuad} shows that the homogeneous $G$-space $\P$
has the type of cohomology ring which we have classified in
Theorem \ref{MainTheorem}.
\end{proof}
\end{Thm}
The following lemma shows that several of the transitive actions occurring
in Theorem \ref{MainTheorem} are not related to compact quadrangles.

\begin{Lem}
\label{NoLargeNormalGroup}
Let $G$ be a compact Lie group acting effectively on a
compact connected $(m_1,m_2)$-quadrangle $\frak G=\PLF$.
Suppose that $G$ acts transitively on $\P$. Let $N\SUB G_p$
be a normal almost simple subgroup, and assume that $G_p$ has
no other almost simple subgroup which is of the same type as
$N$. Then 
\[
\dim(N)\leq\max\left\{\binom{m_1+1}2,\binom{m_2+1}2\right\}.
\]

\begin{proof}
It follows from the assumptions that $N$ satisfies the
condition (ii) of \ref{TrivialActionLem}. 
Therefore $N$ acts non-trivially on $D_1(p)$
or on $D_1(\ell)$. Since $N$ is almost simple, this non-trivial
action is necessarily almost effective (the kernel is finite).
The claim follows now from \ref{DimensionProp}.
\end{proof}
\end{Lem}
The following results are due to Biller \cite{Biller} Sec.~5.1.
We assume the following.
$\frak G$ is a compact connected
$(m_1,m_2)$-quadrangle, $A$ is an elementary abelian
$p$-group, for some odd prime $p$, and of rank $r\geq 1$, acting effectively
on $\frak G$ and fixing a point row $D_1(\ell)$ pointwise.
\begin{Lem}
If $m_2$ is odd, then the fixed points and lines of $A$ form a 
full, compact, and proper subquadrangle $\frak G'$.

\begin{proof}
Let $h$ be a line which is fixed pointwise by $A$, and
let $q\in D_1(h)$. Then $A$ fixes $h\in D_1(q)$, and in particular
$\mathrm{Fix}(A,D_1(q))$ is non-empty. By Theorem \ref{ElementAbelian},
it is a cohomology sphere
of positive dimension. Iterating this argument, we see that the
substructure of $\frak Q$ which is fixed by $A$ is a full subquadrangle
$\frak G'$ (here we use the structure theorem for weak quadrangles
as in Van Maldeghem \cite{HvM98} Thm.~1.6.2: the fixed structure is
a weak subquadrangle with at least two thick lines and infinitely many
thick points, so it is thick);
since $\frak G$ is compact, the set of all lines and points fixed
by $A$ is compact, so $\frak G'$ is compact. Moreover,
$\mathrm{Fix}(A,\P)\neq\P$, since we assumed that $A\neq 1$ acts
effectively.
\end{proof}
\end{Lem}
Let $(m_1,m_2^{(r)})$ denote the parameters of $\frak G'=\frak G^{(r)}
\subset\frak G=\frak G^{(0)}$.
By Theorem \ref{ElementAbelian}, we find a subgroup $A_1\subseteq A_0=A$ of
index $p$ such that $\mathrm{Fix}(A_1,D_1(p))$ is strictly bigger than
$\mathrm{Fix}(A_0,D_1(p))$, for a point $p\in\mathrm{Fix}(A,\P)$.
If we iterate this process, we obtain
a sequence of full, compact, and proper subquadrangles
$\frak G^{(r)}\subset\frak G^{(r-1)}\subset\cdots\subset\frak G^{(0)}=\frak G$,
with parameters $(m_1,m_2^{(k)})$ for $\frak G^{(k)}$, and
$m_2^{(r)}<m_2^{(r-1)}<\cdots <m_2^{(0)}$.
By Proposition \ref{OvoidProp} we have the inequalities
\begin{align*}
m_1 &\leq m_2^{(0)}-m_2^{(1)}\\
m_1 &\leq m_2^{(1)}-m_2^{(2)}\\
m_1 &\leq m_2^{(2)}-m_2^{(3)}\\
& \ \vdots \\
m_1 &\leq m_2^{(r-1)}-m_2^{(r)}
\end{align*}
Adding these up, we have $rm_1+m_2^{(r)}\leq m_2^{(0)}$.
In particular, $rm_1<m_2$ holds for the parameters of $\frak G$.
\begin{Thm}[Biller]
Let $\frak G$ be a compact connected $(m_1,m_2)$-quadrangle, with 
$m_2$ odd, and let $p$ be an odd prime.
Let $\ell$ be a line, and let
$(\ZZ/p)^r\cong A\subseteq\mathrm{AutTop}(\frak G)$,
be an elementary abelian $p$-group, for some odd prime $p$,
and for some $r\geq 1$.
Suppose that $A$ fixes the point row $D_1(\ell)$ pointwise.
Then $m_1r<m_2$.

\begin{proof}
This follows from the previous discussion.
\end{proof}
\end{Thm}

\begin{Cor}[Biller]
\label{BillerRank}
Let $K$ be a compact connected Lie group acting effectively on a compact
connected $(m_1,m_2)$-quadrangle and fixing a point row $D_1(\ell)$ pointwise.
Assume that $m_2$ is odd. Then $\rk(K)<\frac{m_2}{m_1}$.

In particular, if $m_1=4$ and $m_2=4k-5$, then $\rk(K)\leq k-2$.

\begin{proof}
We choose a maximal torus $T\subseteq K$, of rank, say, $s$.
Thus we have $(\ZZ/3)^s\subseteq K$.
\end{proof}
\end{Cor}

%\newpage

\section{The Stiefel manifolds}
\label{QuadStiefel}

Put $G(n)=\SO(n),\SU(n),\Sp(n)$, let $\FF=\RR,\CC,\HH$ denote the
corresponding skew field, and put $d=\dim_\RR\FF=1,2,4$. 
In this section we classify
those compact quadrangles which admit
a line transitive $G(n)$-action, such that $\L=G(n)/G(n-2)$
is a Stiefel manifold, with $n\geq 5$ for $\FF=\RR,\CC$ and
$n\geq 4$ for $\FF=\HH$. Such a quadrangle has topological
parameters 
\[
(d,d(n-1)-1),
\]
where $d=1,2,4$ for $G=\SO(n),\SU(n),\Sp(n)$, respectively.
\begin{Num}\textsc{Example}
\label{StiefelEx}
Consider the quadrangle $\mathsf{H}_{n+1}(\FF,\RR)=\PLF$ as in \ref{ClassQuad}.
Let $V\in\L$ be a 2-dimensional totally isotropic subspace of $\FF^{n+2}$.
It is not difficult to see that $V$ admits a (unique) basis
$\{u,v\}$ such that 
\[
\sum_{k=2}^{n+1}\bar u_kv_k=0\quad\text{ and }\quad
u_0=v_1=1,\quad u_1=v_0=0.
\]
Thus $((u_2,\ldots,u_{n+1}),(v_2,\ldots,v_{n+1}))$ is an element
of the Stiefel manifold $V_2(\FF^n)$. This correspondence is in fact a
$G(2)\times G(n)$-equivariant homeomorphism $\L\cong V_2(\FF^n)$,
and we may
identify the line space $\L$ of $\mathsf{H}_{n+1}(\FF,\RR)$ with 
$V_2(\FF^n)$.
\end{Num}
Next we note the following. If $k\geq 3$, then 
$\dim(G(k))>\binom{d+1}2$, hence the almost simple group $G(k)$
cannot act non-trivially on a $d$-dimensional point row.
This implies that $G(n)_\ell=G(n)_{p,\ell}$ for all points
$p\in D_1(\ell)$, provided that $n\geq 5$. Thus, 
$G(n)_p\supseteq G_\ell=G(n-2)$.
If $d=4$ and $n=4$, then $G(2)=\Sp(2)\isom\Spin(5)$ could act
non-trivially. However, in this case the action would
necessarily be transitive (we will see that this situation does not occur).
\begin{Lem}
\label{beforeBetweenLem}
Let $\PLF$ be a compact generalized quadrangle, and assume
that $G(n)$ acts transitive on the line space $\L$, such that
there is a $G(n)$-equivariant homeomorphism
$\L\isom G(n)/G(n-2)=V_2(\FF^n)$.
If $n\geq 5$, then the line stabilizer $G(n)_\ell=G(n-2)$
acts trivially on the point row $D_1(\ell)$, whence
\[
G(n)_p\supseteq G(n-2)=G(n)_\ell
\]
for every $p\in D_1(\ell)$.
The same conclusion holds for $n=4$ and $\FF=\HH$.

\begin{proof}
Only the last claim has to be proved. Thus we assume that $n=4$
and $\FF=\HH$.
If $\Sp(2)$ acts non-trivially on $D_1(\ell)$, then it acts
transitively, because $\Sp(2)$ contains no subgroup of codimension
less than 4 (this follows from \ref{DimensionProp}).
Thus $\Sp(4)$ acts transitively on the flags, and
in particular transitively on the points. But we know from our
classification \ref{CaseSp(n)even}
that $\Sp(4)$ cannot act transitively on a
1-connected space with the rational cohomology of $\SS^4\times\SS^{15}$.
\end{proof}
\end{Lem}
We need to know the groups between $G(n-2)$ and $G(n)$.

\begin{Lem}
\label{BetweenLem}
Let $H\SUB G(n)$ be a connected almost simple subgroup such that
$G(n-2)\subsetneq H\subsetneq G(n)$. Then $H$ is conjugate to
$G(n-1)$, provided that $n\geq 5$ for $\FF=\RR,\CC$, and $n\geq 4$
for $\FF=\HH$.

\begin{proof}
Consider the adjoint representation of $G(n-2)$ on the Lie algebra
$\fg(n)$ of $G(n)$. We decompose $\fg(n)$ into irreducible
$G(n-2)$-modules. Then 
\[
\fg(n)=\fg(n-2)\oplus V\oplus V\oplus T,
\]
where $V$ is the natural $G(n-2)$-module, and
$T\isom\RR^{3d-2}$ is a
trivial $G(n-2)$-module. Since $G(n-2)\SUB H$, the Lie
algebra $\fh$ of $H$ is invariant under $G(n-2)$. 

If $\fh\SUB\fg(n-2)\oplus T$, then either $\fh=\fg(n-2)$
or $\fh$ is not almost simple; in either case, we get a contradiction
to our assumptions. Thus there exists an element 
$(0,v_1,v_2,t)\in\fh$
with $\{0\}\neq\{v_1,v_2\}$. We consider the $G(n-2)$-orbit of this element.

If $v_1,v_2$ are $\FF$-linearly independent, then the $G(n-2)$-orbit of
$(0,v_1,v_2,t)$ is a (distorted) Stiefel manifold
$V_2(\FF^{n-2})\isom G(n-2)/G(n-4)$ (here we use the assumption on $n$);
thus, its
$\RR$-span contains $V\oplus V$. It follows that $V\oplus V\SUB\fh$, and, 
consequently, that $\fh=\fg(n)$ and $H=G(n)$, contradicting our
assumptions on $H$.

If $\{v_1,v_2\}\neq\{0\}$ are $\FF$-linearly dependent, then the orbit
of $(0,v_1,v_2,t)$ is a $(d(n-2)-1)$-sphere, and
thus the linear span of the orbit contains an $(n-2)$-dimensional
subspace $W$ of $V\oplus V$ which is $G(n-2)$-isomorphic to $V$.
The Lie algebra $\ff$ generated by $W$ and $\fg(n-2)$ is isomorphic 
and conjugate to $\fg(n-1)$. Thus, $\ff\SUB\fh$, and we have to prove
that equality holds. 

Assume that $(0,v_1',v_2',t')\in\fh\setminus\ff$. If
$\{v_1',v_2'\}\neq \{0\}$, then one sees readily that $V\oplus V\SUB\fh$
and hence $\fh=\fg(n)$, again a contradiction to out assumptions.

We are left with the case that $\fh=\fg(n-2)\oplus W\oplus T'$,
where $W\SUB V\oplus V$ and $T'\SUB T$. We may assume that
$\ff=\fg(n-1)$ and that $W=V\oplus 0$. Thus $\fh\subseteq\fg(n-1)\oplus
T''$, where $T''$ is of dimension $d-1$. Since we assumed that $\fh$
is simple, this implies that $\fh\cong\fg(n-1)$.
\end{proof}
\end{Lem}
Now we go back to the quadrangle. We assume that $n\geq 5$ for 
$\FF=\RR,\CC$ and $n\geq 4$ for $\FF=\HH$. Then
$G(n)_{p,\ell}=G(n)_\ell=G(n-2)$
for all $p\in D_1(\ell)$, and we have inequalities
\begin{align*}
\dim(G(n)_p)&\geq \dim(G(n))-\dim(\P) \\
\dim(G(n)_\ell)&=\dim(G(n))-\dim(\L) \\
\dim(D_1(p))&\geq \dim(G(n)_p)-\dim(G(n)_\ell)
\end{align*}
whence
\[
d(n-2)-1\leq \dim(G(n)_p)-\dim(G(n)_\ell)\leq d(n-1)-1.
\]
Put $(G(n)_p)^\circ=K_1\cdot K_2$, where $K_1$ is almost simple,
$G(n-2)=G(n)_\ell\SUB K_1$, and $K_2\SUB\Cen_{G(n)}(G(n-2))$. 
Note that $\dim(\Cen_{G(n)}(G(n-2))=3d-2$.
If $\dim(G(n)_p)-\dim(G(n)_\ell)\geq \dim(\Cen_{G(n)}(G(n-2))$, then
$K_1\neq G(n-2)$ by the inequalities above.
It follows then from \ref{BetweenLem} that $K_1$ is conjugate to $G(n-1)$.
Thus $\dim(K_1/G(n-2))=d(n-1)-1$. This is the full dimension of a line
pencil. Therefore $G_p$ acts transitively on $D_1(p)$, and thus
$G_p=K_1$ is connected and conjugate to $G(n-1)$.

Now $d(n-2)-1\geq 3d-2$ for $d=1,2,4$, provided that $n\geq 5$.
In these cases, we have proved that $G(n)_p$ is conjugate to
$G(n-1)$, contains $G(n-2)$ and acts transitively on the line pencil
$D_1(p)$. 

Finally, assume that $n=4$ and $\FF=\HH$. Then $21\leq\dim(G_p)\geq 
\dim(\Sp(4))-\dim(\P)=36-19=17$. In fact, $\Sp(4)$ is not transitive
on $\P$ by \ref{CaseSp(n)even}, so we have strict inequality
$\dim(G_p)>17$. Assume that $K_1=G(n-2)$. Then $\dim(K_2)\geq 8$,
and $K_1\SUB\Cen_{\Sp(4)}(Sp(2))^\circ\isom\Sp(2)$.
But $\Sp(2)$ has no subgroup of codimension 2 (by \ref{DimensionProp}).
Thus $K_1$ is strictly bigger than $G(n-2)$, and thus
$K_1\isom G(n-1)=G(n)_p$.

\begin{Thm}
Let $G(n)$ denote one of the groups $\SO(n),\SU(n),\Sp(n)$.
Suppose that $G(n)$ acts line transitively on a compact quadrangle
$\frak G$, such that $\L\isom G(n)/G(n-2)$.
If $\FF=\RR,\CC$ and $n\geq 5$, or if $\FF=\HH$ and $n\geq 4$,
then $\frak G$ is uniquely determined and
continuously $G$-isomorphic to the classical quadrangle 
$\mathsf{H}_{n+1}\FF$.

\begin{proof}
The discussion above shows that $G(n)_p$ acts transitively on
$D_1(p)$ for all $p\in \P$. Put $\mathcal{G}=\{G(n)_p|\ p\in D_1(\ell)\}$. 
The quadrangle is uniquely determined by the triple
\[
(G(n),G(n-2),\mathcal{G}),
\]
cp.~\ref{Reconstruction}. 
It remains to determine $\mathcal{G}$ in group theoretic
terms. Let $\mathcal{G}'$ denote the set of all conjugates
of $G(n-1)$ which contain $G(n-2)$. Clearly, $\mathcal{G}\SUB
\mathcal{G}'$; note that $\mathcal{G}'$ is the collection of subgroups
which arises from the corresponding classical quadrangle. We wish to show
that $\mathcal{G}=\mathcal{G}'$. 

The map $p\mapstoo G(n)_p$ is continuous by \ref{BorelStab}.
If $p,q\in D_1(\ell)$ are distinct points, then 
$G(n)_p\neq G(n)_q$, since otherwise $G(n)_p$ would fix $\ell\in D_1(p)$.
Thus, the map $q\mapstoo G(n)_q$ is injective on the point row $D_1(\ell)$,
and its image is $\mathcal{G}$. Therefore $\mathcal{G}$ is a
homology $d$-sphere inside $\mathcal{G}'$. From the corresponding
classical quadrangle we see that $\mathcal{G}'$ is in fact homeomorphic
to $\SS^d$. It follows that $\mathcal{G}=\mathcal{G}'$.

We have established an abstract isomorphism between $\frak{G}$  and 
a certain classical quadrangle,
\[
\frak G\rTo^\phi\mathsf{H}_{n+1}(\FF,\RR).
\]
The map $\phi$ is $G$-equivariant and continuous on the line space
$\L$, so it is continuous on every line pencil. Every point row
$D_1(\ell)$ can be embedded in $\L$ by the map $q\mapstoo \proj_qh$,
where $h\in D_4(\ell)$ is some line. Thus, $\phi$ is continuous on every
point row and every line pencil. By B\"odi-Kramer \cite{BoeKr95}
Prop.~3.5, the map $\phi$ is continuous (and so is its inverse).
\end{proof}
\end{Thm}
The last step in the proof will be needed again later, so we state
it separately.
\begin{Prop}
\label{Continuity}
Let $\frak G\rTo^\phi\frak G'$ be an abstract isomorphism of generalized
quadrangles. If $\frak G$ and $\frak G'$ are topological quadrangles, and
if $\phi$ is continuous on the line or point set, then $\phi$ is
continuous everywhere.

\begin{proof}
The proof is the same as the one given above.
\end{proof}
\end{Prop}

%\newpage

\section{The $(4,4n-5)$-series}
\label{Quad(4,4n-5)}
Here $G/G_p=\Sp(n)\times\Sp(2)/\Sp(n-1)\cdot\Sp(1)\cdot\Sp(1)$.
The canonical examples for such an action are the quadrangles
$\mathsf{H}_{n+1}(\HH,\RR)$.
\begin{Num}\textsc{Example}
We use the same notation as in \ref{ClassQuad}.
Every $1$-dimensional totally isotropic subspace in $\HH^{n+2}$
is spanned by a vector $(u_0,u_1,\ldots,u_n,u_{n+1})$, with
$|u_0|^2+|u_1|^2=1$. This vector is not unique, but the pair
of vectors
\[
(u_2\bar{u}_0,\ldots,u_{n+1}\bar{u}_0,u_2\bar{u}_1,\ldots,u_{n+1}\bar{u}_1)
\]
in $\HH^{2\times n}$ is. This embedding of the point space $\P$ in
$\HH^{2\times n}$ is $\Sp(2)\times\Sp(n)$-equivariant. The image
is a focal
manifold of the homogeneous isoparametric foliation with $g=4$ and
$(m_1,m_2)=(4,4n-5)$. The other focal manifold is the image of the
line space $\L$ considered in \ref{StiefelEx} before.
\emph{Mutatis mutandis,} similar remarks apply to $\FF=\RR,\CC$
and also to $\FF=\OO$, see Kramer \cite{Ov2} \cite{Ov3} \cite{Kr98}.
\end{Num}
This example gives us a linear model for this homogeneous space, i.e.~an 
equivariant embedding into a vector space. We let $V=\HH^n$ denote the
natural $\Sp(n)$-module and consider $V\oplus V$. Then
\[
\Cen_{\SO(8n)}(\Sp(n))^\circ\isom\Sp(2);
\]
consequently, we have an (irreducible) action of 
$\Sp(n)\times\Sp(2)$ on $V\oplus V$.
We can view $V\oplus V$ as the set of all $2\times n$-matrices
with entries in $\HH$; the group $\Sp(n)$ acts by left multiplication,
and $\Sp(2)$ acts by conjugate transpose right multiplication.
The orbit of the vector $p=(e_1,0)$ is precisely the set
\[
\P=\{(xc,xs)\in V\oplus V|\ x\in V,|x|=1, c,s\in\HH,|c|^2+|s|^2=1\},
\]
and the stabilizer of $p$ is of the form
\[
\left\{
\left.
\begin{pmatrix} a & \\  & A \end{pmatrix}
\times
\begin{pmatrix} a & \\  & b \end{pmatrix}
\right|\ 
A\in\Sp(n-1), a,b\in\SS^3\SUB\HH
\right\}.
\]
Thus $\P\isom\Sp(n)\times\Sp(2)/\Sp(n-1)\cdot\Sp(1)\cdot\Sp(1)$.
Consider the map 
\[
\phi:\P\too\SS^4\SUB\RR\oplus\HH,\quad(xc,xs)\mapstoo(|c|^2-|s|^2,2\bar cs).
\]
The group $G=\Sp(n)\times\Sp(2)$ acts --- via the second factor ---
as $\SO(5)$ on $\RR\oplus\HH$, and the map $\phi$ is $G$-equivariant.

Our aim is to show that $G$ acts transitively on the flag space $\F$
and hence on the line space $\L$. 
Let $G_{[\ell]}$ denote the kernel of the $G_\ell$-action on the
point row $D_1(\ell)$,
\[
1\rTo G_{[\ell]}\rTo G_\ell\rTo G_\ell|_{D_1(\ell)}\rTo1.
\]
Let $H=\Sp(n)\subset G$ denote the normal subgroup acting from the left.
Every $H$-orbit in $\P$ is homeomorphic to a sphere $\SS^{4n-1}$;
these orbits are precisely the fibres of $\phi$ (these fibres turn out
to be ovoids, cp.~Kramer-Van Maldeghem \cite{Ov1} and Kramer
\cite{Ov2} \cite{Ov3}). 
Let $p$ be a point. Then $\mathrm{Fix}(H_p,H\cdot p)\cong\SS^3$, and
$\mathrm{Fix}(H_p,\P)\cong\SS^7$. Thus
$p^\perp\not\subseteq\mathrm{Fix}(H_p,\P)$ if $n\geq 3$. So assume that
$n\geq 3$ and let $q\in D_2(p)$ be a point which is not fixed by $H_p$. Then
$H_{p,q}$ is conjugate to $\Sp(n-2)$. The group $H_{p,q}$ is
almost simple and of dimension at least 21; thus, it fixes the point row
$D_1(\ell)$ pointwise by Proposition \ref{DimensionProp}. 
Now $\rk(G_{[\ell]})\geq \rk(H_{p,q})=n-2$. On the other hand,
$\rk(G_{[\ell]})\leq n-2$ by Corollary \ref{BillerRank}.
The group $H_{p,q}$ is almost simple
and of the same rank as $(G_{[\ell]})^\circ$, so $(G_{[\ell]})^\circ$
itself is almost simple. If $n\geq 4$, then Lemma \ref{BetweenLem}
shows that $(G_{[\ell]})^\circ\cong\Sp(n-2)$.

The dimension of $D_1(\ell)$ is 4, so $G_\ell|_{D_1(\ell)}$ is at most
10-dimensional by Proposition \ref{DimensionProp}, and if the group
is 10-dimensional, then it is
locally isomorphic to $\SO(5)$ and transitive. 
We have
\begin{align*}
\dim(\Sp(n-2))+10 &=\dim(G)-\dim(\L)\\
&\leq\dim(G_\ell)\\
&=\dim(G_\ell|_{D_1(\ell)})+\dim(G_{[\ell]})\\
&\leq\dim(G_{[\ell]})+10\\
&=\dim(\Sp(n-2))+10.
\end{align*}
Therefore $G_{\ell}$ acts transitively on $D_1(\ell)$,
and $G$ acts transitively on $\L$. By our previous classification,
$H$ acts transitively on $\L$, and $\frak G$ is the standard hermitian
quadrangle over the quaternions.
\begin{Thm}
\label{(4,4n-5)}
Suppose that $\Sp(n)\cdot\Sp(2)$ acts point transitively on a compact 
quadrangle $\frak G$, such that $\P\isom \Sp(n)\cdot\Sp(2)/\Sp(n-1)\cdot
\Sp(1)\cdot\Sp(1)$.
If $n\geq 4$, then $\frak G$ is uniquely determined and
continuously $G$-isomorphic to the classical quadrangle 
$\mathsf{H}_{n+1}(\HH,\RR)$. Moreover, $Sp(n)\cdot\Sp(2)$ acts transitively
on $\L$ and on $\F$.
\qed
\end{Thm}

%\newpage

\section{Products of spheres}

We consider the following situation: the compact connected
Lie group $G$ acts on a
compact quadrangle in such a way that
\[
\P=G/H=K_1/H_1\times K_2/H_2=\SS^{m_1}\times\SS^{m_1+m_2}
\]
is a product of homogeneous spheres. 
We use the fact that the inclusion of a point row $D_1(\ell)\too \P$ 
induces a map onto $\pi_{m_1}(\P)$. The composite
\begin{diagram}
\SS^{m_1} & \rTo^\homot &
D_1(\ell) & \rInto & \SS^{m_1}\times\SS^{m_1+m_2} \\
&&& \rdTo & \dTo_{\pr_1} \\
&&&& \SS^{m_1}
\end{diagram}
has non-zero degree, hence it is surjective. This shows that
$D_1(\ell)\cap (\{x\}\times\SS^{m_1+m_2})\neq\emptyset$
for every $x\in\SS^{m_1}$.
The stabilizer $H_1$ fixes a set containing
$\SS^0\times\SS^{m_1+m_2}$. But every point row $D_1(\ell)$
meets this set twice. Therefore $H_1$ fixes every point row,
and $(K_1,H_1)=(\Sp(1),1)$. Thus $m_1=3$.

\begin{Prop}
\label{ProdS3}
Let $G$ be a compact connected Lie group.
If the point space $\P$ of a compact quadrangle is a product
of homogeneous spheres 
\[
\P=G/H=K_1/H_1\times K_2/H_2=\SS^{m_1}\times\SS^{m_1+m_2}
\]
as above, then $m_1=3$, $m_2$ is even, and $(K_1,H_1)=(\Sp(1),1)$.
\qed
\end{Prop}
Now we use Lemma \ref{NoLargeNormalGroup} to exclude some more cases. 
If $H_2$ is almost simple, then
\[
\dim(H_2)\leq\max\left\{\binom42,\binom{m_2+1}2\right\}.
\]
This excludes the pairs $(K_2,H_2)=(\SO(m_2+4),\SO(m_2+3))$, for $m_2\geq 2$, 
and $(K_2,H_2)=(\Spin(7),\G_2)$. We are left with the homogeneous spaces
\begin{align*}
\Sp(n)/\Sp(n-1)&\times\Sp(1) \text{ for $n\geq2$} \\
\SU(n+1)/\SU(n)&\times\Sp(1) \text{ for $n\geq2$} \\
\Spin(9)/\Spin(7)&\times\Sp(1).
\end{align*}
There exist point homogeneous compact connected quadrangles with the first
type of group action, see Example \ref{Quad(3,4n-4)} in Chapter 8.

%\newpage

\section{Summary}

We summarize the main results of this chapter.
\begin{Thm}
\label{MainThmCh7-1}
Let $\frak G=\PLF$ be a compact connected 
quadrangle. Suppose that the topological
automorphism group $\mathrm{AutTop}(\frak G)$ acts transitively
on the point space $\P$. Assume that the topological parameters
$(m_1,m_2)$ of $\frak G$ satisfy the condition $m_1\geq 2$.
Then there exists a compact connected Lie group $G\SUB\mathrm{AutTop}
(\frak G)$ which acts transitively and irreducibly on $\P$. 

\textbf{(A)}
If $m_1=m_2$, then $m_1=m_2=2$, the group $G=\SO(5)$ acts transitively
on the flag space $\F$ and the point space $\P$, and $\frak G$ is 
$G$-equivariantly isomorphic to the complex symplectic quadrangle 
$W(\CC)$ or its dual, the complex orthogonal quadrangle $Q_4(\CC)$.

\textbf{(B)}
If $m_1\geq 3$, then $m_1\neq m_2$. We have classified the following 
subcases.

\textbf{(B1)}
If $(m_1,m_2)=(4,4n-5)$, for $n\geq 4$, then $G=\Sp(n)\cdot\Sp(2)$ and
$\frak G$ is $G$-equivariantly isomorphic to the quaternion hermitian 
quadrangle $\mathsf{H}_{n+1}(\HH,\RR)$.

\textbf{(B2)}
If $(m_1,m_2)=(4n-5,4)$, for $n\geq 4$, then $G=\Sp(n)$ and
$\frak G$ is $G$-equivariantly isomorphic to the dual of the 
quaternion hermitian quadrangle \linebreak $\mathsf{H}_{n+1}(\HH,\RR)$.

\textbf{(B3)}
If $(m_1,m_2)=(2n-3,2)$, for $n\geq 5$, then $G=\SU(n)$ and
$\frak G$ is $G$-equivariantly isomorphic to the dual of the complex 
hermitian quadrangle \linebreak $\mathsf{H}_{n+1}(\CC,\RR)$.

\textbf{(B4)}
If $(m_1,m_2)=(n-1,1)$, for $n\geq 9$, then $\frak G$ is 
$G$-equivariantly isomorphic to the dual of 
the real orthogonal quadrangle $Q_{n+1}(\RR)$.

\textbf{(B5)}
In the remaining cases $(G,\P)$ is one of the pairs determined in
Theorem \ref{MainTheorem} (B) and (C).

\begin{proof}
Case (A) was proved in Kramer \cite{Kr94}. The cases (B1)-B(4)
were considered in the previous sections.
\end{proof}
\end{Thm}
Note that the theorem above does not summarize \emph{all} situations which
were considered in this chapter; we classified some more quadrangles.
\begin{Cor}
\label{MainThmCh7-2}
Let $\frak G=\PLF$ be a compact connected 
quadrangle. Suppose that the topological
automorphism group $\mathrm{AutTop}(\frak G)$ acts transitively
on the point space $\P$. Assume that the topological parameters
$(m_1,m_2)$ of $\frak G$ satisfy the condition
\[
m_1\geq 10.
\]
Then $\frak G$ is $G$-equivariantly isomorphic to the dual of one of
the classical compact connected Moufang
quadrangles $\mathsf{Q}_n(\RR)$, $\mathsf{H}_n(\CC,\RR)$, or
$\mathsf{H}_n(\HH,\RR)$.

\begin{proof}
By \ref{ProdS3}, the point space $\P$ cannot be a product of homogeneous 
spheres. Thus $\P$ is a Stiefel manifold by \ref{HighlyCon}.
\end{proof}
\end{Cor}
Case (B5) in the theorem above consists mainly of infinite series
\begin{gather*}
\Sp(1)\times(\Sp(n)/\Sp(n-1)) \\
\Sp(1)\times(\SU(n+1)/\SU(n)) \\
\Sp(n)\times\SU(3)/\Sp(n-1)\cdot\Sp(1) \\
\Sp(n)\times\Sp(2)/\Sp(n-1)\cdot\Sp(1)
\end{gather*}
and a finite collection of sporadic spaces.
I conjecture that only the first series corresponds to quadrangles.
The following related result is proved by Biller \cite{Biller};
it improves some of our results obtained in this chapter.
\begin{Thm}[Biller]
Let $\frak G$ be a compact connected $(4,4n-5)$-quadrangle, for $n\geq 2$,
and assume that
a compact connected Lie group $G$ of dimension at least
$\binom{2n+1}2+10$ acts effectively on $\frak G$. Then $\frak G$
is continuously isomorphic to $\mathsf{H}_{n+1}(\HH,\RR)$, and
$G\cong\Sp(2)\cdot\Sp(n)$ acts transitively on the flags.

If a compact connected Lie group $G$ of dimension at least 
$22$ acts effectively on a compact connected $(4,5)$-quadrangle
$\frak G$, then
$G$ acts transitively on the flags and $\frak G$ is continuously
isomorphic to a certain classical quadrangle
(the so-called anti-hermitian quadrangle over the quaternions).

\begin{proof}
See Biller \cite{Biller} Sec.~5.3.
\end{proof}
\end{Thm}
Some more cases (e.g. $\E_6/\Ffour$ or $\SU(6)/\Sp(3)$) can be excluded 
by \ref{NoLargeNormalGroup}. 
More cases can be excluded using the 
following deep result of Stolz and Markert.

\begin{Thm}[Markert]
Let $\frak G$ be a compact connected $(m_1,m_2)$-quadrangle,
with $2\leq m_1<m_2$. Put $k=\min\{m_2-m_1,m_1-1\}$. Then
\[
2^{\phi(k)}\text{ divides }m_1+m_2+1,
\]
where $\phi(k)=|\{i|\ 1\leq i\leq k\text{ and } i\equiv 0,1,2,4\pmod 8\}|$.

\begin{proof}
This follows from results of Stolz \cite{Stolz98}. However, the
proof requires several non-trivial modifications of Stolz' paper, see
Markert \cite{Mark}.
\end{proof}
\end{Thm}

A complete classification of all compact connected quadrangles
which admit a point or line transitive automorphism group seems
to be difficult, but not impossible.
\begin{Num}\textsc{Conjecture}
\emph{
A compact connected quadrangle which admits a point or line transitive
automorphism group is either a Moufang quadrangle or a quadrangle 
associated to an isoparametric hypersurface of 
Clifford type (with parameters $(3,4k)$ or $(8,7)$).}
\end{Num}

%%%%%%%%%%%%%%%%%%%%%%%%%%%%%%%%%%%%%%%%%%%%%%%%%%%%%%%%%%%%%%%%%%%%%%%%
%                                                                      %
%                                                                      %
%                   Compact homogeneous quadrangles                    %
%                                                                      %
%                          Linus Kramer                                %
%                                                                      %
%                           Memoirs AMS                                %
%                                                                      %
%                         Wuerzburg 2000                               %
%                                                                      %
%                                                                      %
%                            CHQ8.tex                                  %
%                                                                      %
%                                                                      %
%                                                                      %
%%%%%%%%%%%%%%%%%%%%%%%%%%%%%%%%%%%%%%%%%%%%%%%%%%%%%%%%%%%%%%%%%%%%%%%%
\chapter{Homogeneous focal manifolds}

We apply our classification of homogeneous spaces to obtain a classification
of homogeneous focal manifolds of isoparametric hypersurfaces. 
A closed submanifold $M^n\SUB\RR^{n+k}$
is called isoparametric if its normal bundle is (globally)
flat, and if the eigenvalues of the Weingarten map along any
parallel normal field are constant. Typical examples are
principal orbits of isotropy representations of Riemannian
symmetric spaces. However, Ferus-Karcher-M\"unzner showed that
there are infinitely many isoparametric hypersurfaces which
do \emph{not} arise from isotropy representations of symmetric
spaces. On the other hand, Thorbergsson proved that an irreducible
full isoparametric submanifold with codimension strictly bigger
than 2 (i.e.~$M$ is not congruent to a product, the codimension is
$k\geq 3$, and $M$ is not contained in an affine hyperplane)
is a principal orbit of the isotropy
representation of a Riemannian symmetric space. This should
be compared to Tits' result about buildings of higher rank
mentioned in the previous chapter.

This correspondence is not accidental.
The rank of an isoparametric submanifold $M^n\SUB\RR^{n+k}$
is defined as
$\dim(\bra{M}_{\text{aff}})-\dim(M)$, where $\bra{M}_{\text{aff}}$
denotes the affine span of $M$, i.e.~the intersection of all affine
hyperplanes containing $M$. It is called irreducible if it is not
congruent to a product $M_1\times M_2\SUB\RR^{m_1}\times\RR^{m_2}$
of isoparametric submanifolds $M_1\SUB\RR^{m_1}$, $M_2\SUB\RR^{m_2}$
under an isometry of $\RR^{m_1+m_2}$.
Thorbergsson \cite{Tho91}
proved that an irreducible isoparametric submanifold $M^n\SUB\RR^{n+k}$
of rank $k\geq3$ carries in a natural way the structure of an
irreducible spherical building of rank $k$. Using Tits' classification 
\cite{TitsLNM} (combined with a result by Burns-Spatzier \cite{BurSpa87}),
he deduces that the
building is Moufang and arises as the spherical building at infinity
of a Riemannian symmetric space $G/K$. The isotropy representation of
$K$ on $G/K$ yields an isoparametric family
in $\mathfrak p$, where $\mathfrak g=\mathfrak k\oplus\mathfrak p$
is a Cartan decomposition of $G$. A closer inspection shows that there
is an isometry $\RR^{n+k}\isom\mathfrak p$ (up to some real scaling factor)
which  carries $M$ onto an isoparametric submanifold in $\mathfrak p$
of the family belonging to $G/K$.
Thus he obtains the following classification \cite{Tho91}.
\psn
\textbf{Theorem (Thorbergsson)}
{\em A closed irreducible isoparametric submanifold of rank at least $3$
arises from the isotropy representation of an irreducible Riemannian
symmetric space of non-compact type.}
\medskip

An isoparametric submanifold of rank 2 is (up to some real scalar) congruent
to an isoparametric hypersurface in the unit sphere of the ambient
vector space.
Ferus-Karcher-M\"unzner \cite{FeKaMu81}
constructed an infinite series of families of isoparametric hypersurfaces
in spheres which are not congruent (not even homeomorphic) to principal
orbits of isotropy representations of Riemannian symmetric spaces.
So Thorbergsson's Theorem does not carry over to the rank 2 case.
However, a weak form of his result (which is due to Hsiang-Lawson 
\cite{HsLaw71} and
much older than Thorbergsson's result) is still valid in the rank 2 situation.
\psn
\textbf{Theorem (Hsiang-Lawson)}
{\em A closed irreducible isoparametric submanifold of rank $2$
arises from the isotropy representation of an irreducible Riemannian
symmetric space of non-compact type and rank $2$, provided that it admits a
transitive isometry group.}
\medskip

Table II in Hsiang-Lawson \cite{HsLaw71}
is incomplete (the hypersurface with 6 distinct
principal curvatures and multiplicities $(1,1)$ is
erroneously omitted, cp.~Uchida \cite{Uch80}); 
a corrected table can be found in Takagi-Takahashi \cite{TkgTsh72}. 

Note also that in the case of higher rank Thorbergsson \emph{proves}
the existence
of a transitive isometry group (the isotropy group $K$ of the
symmetric space $G/K$).
Recently, Heintze-Liu \cite{HeLi} gave a new proof for the existence
of a transitive isometry group in the case of rank at least 3;
this yields also a new proof of Thorbergsson's Theorem. Yet another
proof is due to Olmos \cite{Olm}.
Most of the new examples discovered by Ferus-Karcher-M\"unzner admit no
transitive isometry groups.

Generalizing the result of Hsiang-Lawson \cite{HsLaw71}, 
we consider the following situation: $M$ is an irreducible
isoparametric hypersurface whose isometry group acts transitively on one
of the focal manifolds of $M$. Some (but not all) of the Ferus-Karcher-M\"unzner
examples have this property. By M\"unzner \cite{Munz81},
the spectrum of the
Weingarten map (i.e.~the set of all principal curvatures)
of an isoparametric hypersurface has $g=1,2,3,4$ or 6 elements.
The cases $g=1,2$ are easy to classify, so the interesting cases are
$g=3,4,6$. For $g=3,6$ one can show that a group
which acts transitively on one focal manifold acts also transitively on the
isoparametric hypersurface itself. Also, the case $g=3$ was completely
classified by Cartan, and the case $g=6$ was partially classified
by Dorfmeister-Neher (see Section \ref{IsoHyp} below for references).
Thus, the interesting case is $g=4$.
In this chapter we prove the following result.
\psn
\textbf{Theorem A}
{\em Let $M$ be a closed isoparametric hypersurface,
with $4$ distinct principal curvatures. Assume that the isometry group
of $M$ acts transitively on one of the two focal manifolds of $M$, and
that this focal manifold is $2$-connected. Then either $M$ itself is
homogeneous, or $M$ is of Clifford type, with multiplicities
$(8,7)$ or $(3,4k)$.}
\medskip

Wolfrom \cite{Wol} recently extended this result to the case
where one multiplicity is $m_1=2$ (for all $g$), showing that
no new examples occur. If $g=4,6$, and if one of the multiplicities is
$m_1=1$, then the results of Takagi \cite{Takagi} and
Dorfmeister-Neher \cite{DorNeh85} apply. Combining Theorem A with
these results, we have the following theorem.
\psn
\textbf{Theorem B}
{\em Let $M$ be a closed irreducible
isoparametric submanifold of rank at least $2$.
Assume that the isometry group
of $M$ acts transitively on one of the focal manifolds of $M$.
Then either $M$ itself is
homogeneous and arises from the isotropy representation of an irreducible
Riemannian symmetric space (of non-compact type),
or $M$ is of Clifford type, with multiplicities
$(8,7)$ or $(3,4k)$.}
\medskip

Finally, I would like to mention Immervoll's new and beautiful result
\cite{StIm}.
\psn
\textbf{Theorem (Immervoll)}
{\em Let $M$ be a closed irreducible
isoparametric submanifold of rank at least $2$; if the rank is $2$,
assume we are not in the situation $(g,m_1,m_2)=(6,2,2)$.
Then the simplicial complex associated to $M$ is a compact
connected Tits building.}
\medskip

Probably, the case $(g,m_1,m_2)=(6,2,2)$ is not really an exception;
it is conjectured that only one isoparametric hypersurface
with these parameters exists (and this known hypersurface
is a building, the so-called split Cayley hexagon over the
complex numbers).

\emph{The strategy of the proof of Theorem A}
is as follows. In the situation of 
Theorem A, the focal manifold has the same integral cohomology as
a product of spheres. Thus we can apply our classification of
homogeneous spaces. However, not every homogeneous space in
\ref{MainTheorem}
is a focal manifold. The fact that the homogeneous space
admits an equivariant embedding into a vector space of the right
dimension rules out many of the 'wrong' manifolds. If the homogeneous
space is one of the known examples of focal manifolds, then
it is in most cases easy to see that this equivariant embedding
is unique. Finally, some low-dimensional cases and some exceptional
spaces need special attention.

Our method does not use very much Riemannian geometry; what we use
is the global behavior of isoparametric submanifolds, and also the
point-line geometry associated to such a submanifold. In a sense,
we treat these submanifolds as if they were nicely embedded
Tits buildings; many proofs are similar to the proofs in the
last chapter. I think that both the effectiveness of this method
and the close relation between the known isoparametric submanifolds
and compact buildings justify this approach.

\medskip
There are many papers and several books on isoparametric hypersurfaces and
submanifolds, we just mention (in alphabetic order)
Abresch \cite{Abrsch83},
Cecil-Ryan \cite{CeRy85},
Console-Olmos \cite{ConOlm},
Dorfmeister-Neher \cite{DorNeh83a}, \cite{DorNeh83b},
\cite{DorNeh84}, \cite{DorNeh85},
Ferus-Karcher-M\"unz\-ner \cite{FeKaMu81},
Grove-Halperin \cite{GroHal87},
Heintze-Liu \cite{HeLi},
Hsiang-Palais-Terng \cite{HsPaTe88},
Knarr-Kramer \cite{KnaKra95},
M\"unzner \cite{Munz80} \cite{Munz81},
Olmos \cite{Olm},
Palais-Terng \cite{PaTe88},
Stolz \cite{Stolz98},
Str\"ubing \cite{Strub86},
Terng \cite{Te85},
Thorbergsson \cite{Tho91}, \cite{Tho92}, \cite{ThoSur},
Wang \cite{Wan88}.
Thorbergsson's survey \cite{ThoSur} is a good introduction to
the subject and contains many references; a somewhat older
bibliography is K\"uhnel-Cecil \cite{KuehCec}.
%\newpage

\section{Isoparametric hypersurfaces}
\label{IsoHyp}
Let $X\SUB\SS^{r+1}$ be a submanifold. We denote the normal
bundle of $X\SUB\SS^{r+1}$ by ${\perp}X$.

\begin{Num}\textsc{Isoparametric hypersurfaces}\psn
A compact connected hypersurface $\F\SUB\SS^{r+1}$ is called
\emph{isoparametric} if its principal curvatures are constant.
Fix a unit normal field $N$ and consider the normal exponential map
\[
\RR\times\F\too\SS^{r+1},\qquad
(t,x)\mapstoo\exp_x(tN).
\]
If $t$ is small enough, then the endpoint map $\eta_t:x\mapstoo \exp_x(tN)$
is a diffeomorphism, and the image of this map is a parallel hypersurface
which is also isoparametric. More precisely, there are real numbers 
$\theta_1<0<\theta_2$ such that $\eta_t$ is a diffeomorphism for 
$\theta_1<t<\theta_2$;
if $t=\theta_1,\theta_2$, then $\eta_t$ is a submersion onto a manifold
of strictly smaller dimension. Put 
\[
\P=\eta_{\theta_2}(\F)\text{ and }\L=\eta_{\theta_1}(\F)
\]
and 
\[
\pr_\P=\eta_{\theta_2}\text{ and }\pr_\L=\eta_{\theta_1}.
\]
The parallel hypersurfaces of $\F$ are again isoparametric; thus,
we have a foliation of the sphere by isoparametric hypersurfaces
and the two focal manifolds, a so-called \emph{isoparametric foliation}.

Let $g$ denote the number of distinct principal curvatures of
$\F$. A remarkable result of M\"unzner \cite{Munz81} says that 
$g\in\{1,2,3,4,6\}$. The isoparametric
hypersurfaces with $g=3$ were classified by Cartan \cite{Cart39}, 
cp.~also Karcher \cite{Karch88}, Knarr-Kramer \cite{KnaKra95}
and Console-Olmos \cite{ConOlm}. The 6-dimensional isoparametric
hypersurfaces with $g=6$ were classified by
Dorfmeister-Neher \cite{DorNeh85}.

It is not difficult to prove that
for $g=1$ one has (up to a scalar factor) an $r$-sphere 
$\mathrm{S}^r_\eps\SUB\SS^{r+1}$ of radius $\eps$, and
for $g=2$ one has a Clifford torus $\mathrm{S}^{m_1}_\eps\times
\mathrm{S}^{m_2}_{\delta}\SUB\SS^{m_1+m_2+1}$, for
$0<\eps<1$ and $\eps^2+\delta^2=1$, cp.~Palais-Terng \cite{PaTe88}.

The results of M\"unzner also show the following. Let $m_i$ denote
the multiplicity of the $i$th principal curvature. After rearranging
the indices one has $m_i=m_{i+2}$ for $i=1,\ldots,g$ (indices mod $g$). 
The numbers $m_1$, $m_2$
are called the \emph{multiplicities} of $\F$. If $g=1,3,6$, then
$m_1=m_2$. If $g=4$, then either $m_1=m_2\in\{1,2\}$, or
$1\in\{m_1,m_2\}$, or $m_1+m_2$ is odd, cp.~M\"unzner \cite{Munz81},
Abresch \cite{Abrsch83}.

Stolz \cite{Stolz98} recently proved the following deep result. 
Suppose that $g=4$ and that $2\leq m_1<m_2$. Then 
\[
(m_1,m_2)=(4,5)\quad\text{ or }\quad
2^{\phi(m_1-1)}\text{ divides }m_1+m_2+1,
\]
where $\phi(m)=|\{i|\ 1\leq i\leq m\text{ and } i\equiv 0,1,2,4\pmod 8\}|$.
\end{Num}

\begin{Num}\textsc{The global geometry}\psn
\label{IsoGlobal}%
The fibres of $\eta_{\theta_i}$ are spheres of
dimension $m_i$, for $i=1,2$. 
Conversely, $\F$ can be identified
with the normal sphere bundle $\mathrm{S}_{2\sin(\theta_1/2)}({\perp}\P)$ 
consisting of all normal vectors of length $2\sin(\theta_1/2)$ (and similarly,
$\F\isom\mathrm{S}_{2\sin(\theta_2/2)}({\perp}\L))$. In fact,
\[
\F=\{x\in\SS^{r+1}|\ \mathrm{dist}(x,\P)=2\sin(\theta_1/2)\}
\]
and
\[
\L=\{x\in\SS^{r+1}|\ \mathrm{dist}(x,\P)=2\sin((\theta_2-\theta_1)/2)\}.
\]
In particular, the isoparametric foliation is completely determined
by one focal manifold.
\end{Num}
We need the following result.
\begin{Lem}
\label{WeingartenLemma}
Let $p\in\P$ and let $N\in{\perp}_p\P$ be a non-zero normal vector.
Then the kernel of the Weingarten map $A_N$ is $m_2$-dimensional.
A similar result holds for $\L$; here, the kernel is $m_1$-dimensional.

\begin{proof}
This is proved e.g. in Cecil-Ryan \cite{CeRy85} Ch.~2, Cor.~2.2.
\end{proof}
\end{Lem}

\begin{Num}\textsc{The corresponding point-line geometry}\psn
Let $\F$ be an isoparametric hypersurface with focal manifolds
$\P$ and $\L$ and $g\geq 2$. Call the elements of $\P$ \emph{points} and the
elements of $\L$ lines. A point $p$ and a line $\ell$ are \emph{incident}
if and only if
\[
\mathrm{dist}(p,\ell)=2\sin((\theta_2-\theta_1)/2).
\]
In that case
the geodesic arc between $p$ and $\ell$ meets $\F$ in a unique element 
$x\in\F$; thus, we can identify $\F$ with the flag set of this geometry.
We denote the resulting point-line geometry by $\frak{G}(\P,\L)$.

This geometry was introduced in Thorbergsson \cite{Tho91};
there, and in Knarr-Kramer \cite{KnaKra95} it is used to classify
certain isoparametric submanifolds. The same geometry is also
used in Eschenburg-Schr\"oder \cite{EschSchr91} App. 
If the isoparametric foliation arises from a non-compact
Riemannian symmetric space, then this is the same as
the building at infinity in the compactification of the
symmetric space.

One can show that $\frak{G}(\P,\L)$ is connected; in fact,
$d(x,y)\leq g$ for all elements $x,y\in\P\cup\L$, 
cp.~Eschenburg-Schr\"oder \cite{EschSchr91}. Also, one can
show that $\frak{G}(\P,\L)$ contains no digons, provided that $g\geq 3$
(Thorbergsson, unpublished; see Knarr-Kramer \cite{KnaKra95} for a proof).
Let $\ell\in\L$. The point row corresponding to $\ell$ has a simple
description: the sphere of radius $2\sin((\theta_2-\theta_1)/2)$ around
$\ell$ touches $\P$ along a sphere $S$, and this set $S$ is precisely
the set of all points which are incident with $\ell$.
In differential-geometric terms, $S\SUB\P$ is a \emph{curvature sphere}.

For all known examples of isoparametric hypersurfaces (with $g\geq 3$),
the geometry $\frak{G}(\P,\L)$ is a spherical building of rank 2;
more precisely, it is a
compact generalized $g$-gon, cp.~Thorbergsson \cite{Tho92}, and,
more generally, Immervoll \cite{StIm}.
\end{Num}

\begin{Num}\textsc{Isometry groups}\psn
Let $\mathrm{Isom}(\F)$ denote the group of all isometries of $\F$. 
Isoparametric
hypersurfaces are \emph{rigid}: every isometry of $\F$ is induced by
an isometry of $\SS^{r+1}$. Thus $\mathrm{Isom}(\F)\SUB\O(r+2)$.
Hsiang-Lawson classified all isoparametric submanifolds where
$\mathrm{Isom}(\F)$ acts transitively on $\F$. We call such a hypersurface
\emph{homogeneous}.

The connected component of $\mathrm{Isom}(\F)$ 
acts also on the focal manifolds $\P$, $\L$; in particular,
\[
\mathrm{Isom}(\F)^\circ\SUB\Aut(\frak{G}(\P,\L)).
\]
Conversely, if $G\SUB\O(r+2)$ is a group which leaves $\P$ or $\L$
invariant, then $G\SUB\mathrm{Isom}(\F)$, cp.~\ref{IsoGlobal}.
As the examples by Ferus-Karcher-M\"unzner show, it is possible that
one of the focal manifolds is homogeneous, while $\F$ itself is
not homogeneous.
\end{Num}
We derive some more general results about transitive actions
on focal manifolds.
\begin{Lem}
\label{NoTrivialFactor}
Assume that $g\geq 3$.
Let $G\SUB\SO(r+2)$ be a compact connected subgroup. 
If $G$ acts transitively on the focal manifold $\P$ (or $\L$), then
there is no invariant subspace with trivial $G$-action.

\begin{proof}
Assume otherwise. Let $U\oplus V=\RR^{2(m_1+m_2)+2}$
be a decomposition into $G$-modules such that $U$ is a trivial $G$-module
of positive dimension.
Let $(u,v)\in\P$. Then $\P=G\cdot(u,v)=u\times 
G\cdot v$ is contained in 
the proper affine subspace $u+V$. But for $g\geq 3$ the focal submanifolds 
are full, i.e.~they are not contained in any proper affine subspace of
$\RR^{r+2}$.
\end{proof}
\end{Lem}
Next we note the following. Let $x\in\F$. Then $G\cdot x$
surjects onto $G\cdot\pr_\P(x)$. Thus, if $\P$ is $G$-homogeneous and
if $y\in\SS^{r+1}\setminus\L$, then the orbit 
$G\cdot y \SUB\SS^{r+1}$ has at least dimension $\dim(\P)$.

\begin{Lem}
Suppose that $\P=G/G_p\SUB\SS^{r+1}$ is a homogeneous focal manifold. If
$x\in\SS^{r+1}$ is a point with $\dim (G\cdot x)<\dim(\P)$, then $G\cdot x$
is contained in the other focal manifold $\L$.
\qed
\end{Lem}
Let $\P=G/G_p$ be a homogeneous focal manifold. Then $G_p$ acts on
the normal space ${\perp}_p\P\SUB T_p\SS^{r+1}$. 
This is the \emph{normal isotropy
representation} of the isotropy group $G_p$. The normal sphere bundle
$S({\perp}\P)$ can be identified with the isoparametric hypersurface $\F$, 
and the action of $G$ on this sphere bundle coincides with the action
on $\F$. Let $\L$ denote the other focal manifold, and let $S$ denote the
curvature sphere corresponding to $p$. The action of $G_p$ on $S$
coincides with the action of $G_p$ on the normal sphere
in ${\perp}_p\P$. Note that \ref{TrivialActionLem}
applies also to isoparametric hypersurfaces.

\begin{Lem}
\label{NoLargeNormalGroup2}
Let $N$ be a normal almost simple subgroup of $G_p$, and assume that
$G_p$ has no other closed subgroup isomorphic to $N$. Then
\[
\dim(N)\leq \max\left\{\binom{m_1+1}2,\binom{m_2+1}2\right\}.
\]
\begin{proof}
The proof is the same as in \ref{NoLargeNormalGroup}.
\end{proof}
\end{Lem}
The problem we will consider is the following. Let $G\SUB\SO(r+2)$
be a compact connected subgroup which leaves $\P$ or $\L$ 
invariant, and which acts transitively on $\P$ or on $\L$.
We wish to determine all possibilities for $(G,\P,\L)$.

The cases $g=1,2,3$ are rather simple. In fact, for $g=3,6$, 
or if $g=4$ and $m_1=m_2$, one
can show that transitivity of $\mathrm{Isom}(\F)$ on one focal
manifold implies transitivity on $\F$ itself; the methods developed in
Kramer \cite{Kr94} for compact polygons carry over to isoparametric 
hypersurfaces (with some modifications in the smallest case $m_1=m_2=1$),
see Wolfrom \cite{Wol}.

\begin{Num}\textsc{The topology for $g=4$}\psn
We consider the remaining case where $g=4$ and $m_1\neq m_2$.
Then
\begin{align*}
\dim(\F)&=2(m_1+m_2) \\
\dim(\P)&=2m_1+m_2 \\
\dim(\L)&=2m_2+m_1
\end{align*}
If $m_1\neq m_2$, and $2\leq m_1,m_2$, then $m_1+m_2$ is odd.
If $m_1+m_2$ is odd, then
\begin{align*}
\bfH^\bullet(\P)&\isom\bfH^\bullet(\SS^{m_1}\times\SS^{m_1+m_2}) \\
\bfH^\bullet(\L)&\isom\bfH^\bullet(\SS^{m_2}\times\SS^{m_1+m_2}) \\
\bfH^\bullet(\F)&\isom\bfH^\bullet(\SS^{m_1}\times
\SS^{m_2}\times\SS^{m_1+m_2}).
\end{align*}
If $m_1=1$ and if $m_2$ is odd, then
\[
\bfH^\bullet(\L;\QQ)\isom\bfH^\bullet(\SS^{2m_2+1};\QQ).
\]
Moreover, the following is true. The focal manifold $\L$ is
$(m_2-1)$-connected, and $\P$ is $(m_1-1)$-connected. 
Let $S\SUB\L$ be a curvature
sphere of dimension $m_2$, corresponding to $p\in\P$. 
Then $\pi_{m_2}(\L)$ is generated
by the inclusion $\SS^{m_2}\isom S\SUB\L$; a similar result holds for $\P$ 
(note that $\pi_{m_2}(\L)\isom\ZZ$ if $m_1<m_2$ and $m_1+m_2$ is odd and
that $\pi_{m_2}(\L)\isom\ZZ/2$ if $m_1=1$ and $m_2$ is odd).

This follows from M\"unzner \cite{Munz81}. M\"unzner determined
structure constants for the cohomology rings, although he did
not write down the resulting rings. 
See Strau\ss\ \cite{Strau96} for a thorough discussion.
\end{Num}
By the result above, we can apply our classification of homogeneous
spaces to focal manifolds.

%\newpage

\section{The Stiefel manifolds}
\label{IsoStiefel}

Let $\FF\in\{\RR,\CC,\HH\}$. The Stiefel manifolds $V_2(\FF^n)$
can be embedded into spheres as focal manifolds of isoparametric
hypersurfaces.
Let $G(n)$ be one of the groups $\SO(n)$, $\SU(n)$ or $\Sp(n)$,
let $V$ denote the natural $G(n)$-module, and let $\FF$ denote the
corresponding skew field, $\FF=\RR,\CC,\HH$. 
Consider the action of $G(n)$ on the unit sphere 
$\SS^{2dn-1}\SUB V\oplus V$, where $d=\dim_\RR\FF=1,2,4$.
There are precisely two orbit types. Let $(x,y)\in\SS^{2dn-1}
\SUB V\oplus V$. If $x,y$ are $\FF$-linearly
independent, then $G(n)_{(x,y)}\isom G(n-2)$, and if $x,y$ are
$\FF$-linearly dependent, then $G_{(x,y)}\isom G(n-1)$.

The union 
\[
X=\{(x,y)\in\SS^{2dn-1}|\ G(n)\cdot x\isom\SS^{dn-1}\}
\]
of the singular orbits is one focal manifold
of the homogeneous isoparametric hypersurface corresponding to this
situation, cp.~Section \ref{QuadStiefel} and 
Section \ref{Quad(4,4n-5)} in the
previous chapter. The other focal manifold is the orbit
$Y=G(n)\cdot(x,y)$, where $|x|^2=|y|^2=1/2$ and $(x|y)=0$. However,
this orbit is (topologically) not unique; all principal orbits
are homeomorphic.

Suppose now that $\L=G(n)/G(n-2)\SUB\RR^{2dn}$ is a homogeneous focal 
manifold.
Assume in addition that the action of $G(n)$ is as above, i.e.~that 
$\RR^{2dn}\isom V\oplus V$. The union $X$ of the singular orbits has to be 
contained in the other focal submanifold $\P$. Now $X$ is a manifold of 
the same dimension as $\P$, hence $\P=X$. But the isoparametric foliation 
is uniquely determined by one focal manifold, so $\L=Y$.

\begin{Lem}
The homogeneous spaces $G(n)/G(n-2)$ can be realized as homogeneous focal
manifolds. If $\RR^{2dn}=V\oplus V$, then there is a unique isoparametric
foliation in $\SS^{2dn-1}$ corresponding to this action.
\qed
\end{Lem}
To obtain a complete classification, it remains to show that 
$\RR^{2dn}\isom V\oplus V$. If $n$ is large enough,
then there is only one irreducible representation of $G(n)$ of dimension
at most $2dn$, the natural one on $V=\FF^n$. From 
\ref{BnModules}, \ref{DnModules}, \ref{AnModules}, \ref{CnModules}, 
we see that
the precise numbers are as follows: for $\FF=\RR,\CC,\HH$ we need
$n\geq 10,6,5$, respectively.
For the low-dimensional cases we use the following fact.
According to Lemma \ref{NoTrivialFactor}, the $G(n)$-module $\RR^{2dn}$
cannot have any trivial factors.

\begin{description}
\item[\fbox{$G(5)=\SU(5)$, $2dn=20$}]
By \ref{AnModules}, the semisimple 20-dimensional
$\SU(5)$-$\RR$-\linebreak mod\-ules without trivial factors are  $V\oplus V$
and $^\RR X_{\lambda_2}$.
The orbits in $^\RR X_{\lambda_2}$ yield
isoparametric hypersurfaces, but with other orbit types
(the multiplicities are $(4,5)$).
Thus in our setting we have $V\oplus V$ as the only possibility.

\item[\fbox{$G(4)=\SU(4)$, $2dn=16$}]
By \ref{AnModules}, there are no semisimple 16-dimensional
\linebreak
$\SU(4)$-$\RR$-modules without trivial factors, except for $V\oplus V$.

\item[\fbox{$G(3)=\SU(3)$, $2dn=12$}]
By \ref{AnModules}, the semisimple 12-dimensional
$\SU(3)$-$\RR$-\linebreak
mod\-ules without trivial factors are  $V\oplus V$
and $^\RR X_{2\lambda_1}=\mathrm{S}^2\CC^3$.
I am indebted to R. Bryant for pointing out that every complex
symmetric $3\times 3$-matrix can be diagonalized under the action
of $\SU(3)$. Therefore, the isotropy group of every element in
$\mathrm{S}^2\CC^3$ contains $\ZZ/2\oplus\ZZ/2$, and thus
there are no orbits of type $\SU(3)/1$ in $\mathrm{S}^2\CC^3$.
Thus we have $V\oplus V$ as the only possibility.

\item[\fbox{$G(9)=\SO(9)$, $2dn=18$}]
By \ref{BnModules}, there are no semisimple 18-dimensional
\linebreak
$\SO(9)$-$\RR$-modules without trivial factors, except for $V\oplus V$.

\item[\fbox{$G(7)=\SO(7)$, $2dn=14$}]
By \ref{BnModules}, there are no semisimple 14-dimensional
\linebreak
$\SO(7)$-$\RR$-modules without trivial factors, except for $V\oplus V$.

\item[\fbox{$G(5)=\SO(5)$, $2dn=10$}]
By \ref{BnModules}, the semisimple 10-dimensional
$\SO(5)$-$\RR$-\linebreak
mod\-ules without trivial factors are  $V\oplus V$
and $\mathrm{Ad}$.
In $\mathrm{Ad}$, the principal orbits are $\SO(5)/\TT^2$. These
orbits belong to a homogeneous  isoparametric hypersurface with 
multiplicities $(2,2)$, and this excludes this module in our situation. 
Thus, we we have $V\oplus V$. 

\item[\fbox{$G(4)=\Sp(4)$, $2dn=32$}]
By \ref{CnModules}, there are no semisimple 32-dimensional
\linebreak
$\Sp(4)$-$\RR$-modules without trivial factors, except for $V\oplus V$.

\item[\fbox{$G(3)=\Sp(3)$, $2dn=24$}]
By \ref{CnModules}, there are no semisimple 24-dimensional
\linebreak
$\Sp(3)$-$\RR$-modules without trivial factors, except for $V\oplus V$.

\item[\fbox{$G(2)=\Sp(2)$, $2dn=16$}]
By \ref{BnModules}, there are no semisimple 16-dimensional
\linebreak
$\Sp(2)$-$\RR$-modules without trivial factors, except for $V\oplus V$.

\item[\fbox{$G(8)=\SO(8)$, $2dn=16$}]
By \ref{BnModules}, there are no semisimple 16-dimensional
\linebreak
$\SO(8)$-$\RR$-modules without trivial factors, except for $V\oplus V$.

\item[\fbox{$G(6)=\SO(6)$, $2dn=12$}]
By \ref{BnModules}, there are no semisimple 12-dimensional
\linebreak
$\SO(6)$-$\RR$-modules without trivial factors, except for $V\oplus V$.
\end{description}
Thus we have a complete classification for classical groups acting
on Stiefel manifolds. However, there are also the exceptional actions
$\G_2/\SU(2)=V_2(\RR^7)$,\linebreak $\Spin(7)/\SU(3)=V_2(\RR^8)$, and 
$\Spin(9)/\G_2=V_2(\OO^2)$.
We consider these cases, and some more, in the next section.

%\newpage

\section{Some sporadic cases}

Now we consider some homogeneous spaces of almost simple Lie
groups which do not belong to any series. Nevertheless, the
ideas are very similar as in the last section.
\begin{description}
\item[\fbox{$\L=\Spin(9)/\G_2=V_2(\OO^2)$}]
By \ref{BnModules}, the only 32-dimensional
semisimple \linebreak $\Spin(9)$-$\RR$-module 
without trivial factors is $V\oplus V$, where $V=\OO\oplus\OO$
is the affine Cayley plane. 
The orbit types are $\Spin(9)/\Spin(7)$ (15-dimensional orbits),
$\Spin(9)/\Spin(6)$ (21-dimensional orbits), and 
$\Spin(9)/\G_2$ (22-di\-men\-sio\-nal orbits), see
Salzmann \emph{et al.} \cite{CPP95}, Chapter 1. 
The union $X$ of the singular
orbits is the non-homogeneous focal manifold of the Clifford hypersurface
with multiplicities $(8,7)$ (the definite case,
cp.~Ferus-Karcher-M\"unzner 
\cite{FeKaMu81} 6.6). Since $X\SUB\P$ and $\dim(X)=\dim(\P)$, we have
uniqueness.

\item[\fbox{$\P=\Spin(10)/\SU(5)$ or $\L=\Spin(10)/\Spin(7)$}]
By \ref{DnModules}, the only 32-dimension\-al
semisimple
$\Spin(10)$-$\RR$-module  without trivial factors is
$^\RR X_{\lambda_4}=\CC^{16}$. 
Inspecting the homogeneous hypersurface with multiplicities $(6,9)$, 
we see
that there are exactly two singular orbits, of dimensions
$24$ and $21$, respectively; the principal orbits have
codimension 1. Thus, the singular orbits are the focal manifolds.

\item[\fbox{$\P=\SU(5)/\SU(2)\times\SU(3)$ or $\L=\SU(5)/\Sp(2)$}]
By \ref{AnModules}, the only 20-di\-men\-sional
semisimple $\SU(5)$-$\RR$-modules  without trivial factors are
$2\cdot{}^\RR X_{\lambda_1}=\CC^5\oplus\CC^5$ and
${}^\RR X_{\lambda_2}=\EA^2\CC^5$. 
In both cases, the orbit structure is known. In our present
situation, we have the 20-dimensional simple
module. There are exactly 2 singular orbits of type
$\P=\SU(5)/\SU(2)\times\SU(3)$ and
$\L=\SU(5)/\Sp(2)$.

\item[\fbox{$\P=\Spin(9)/\SU(4)$}]
By \ref{BnModules}, the only 32-dimensional
semisimple $\Spin(9)$-$\RR$-module 
without trivial factors is $V\oplus V$, where $V=\OO\oplus\OO$
is the affine Cayley plane. As observed above, the orbits
have dimension $21=\dim(\P)$, 22, or 15. The 15-dimensional orbits
have to be contained in the other focal manifold $\L$, and
(by counting dimensions) $\L$ is the union of these 15-spheres.
Thus we know $\L$, and this determines the isoparametric
foliation uniquely.

\item[\fbox{$\L=\Spin(7)/\SU(3)=V_2(\RR^8)$}]
By \ref{BnModules}, the only 16-dimensional
semisimple \linebreak $\Spin(9)$-$\RR$-module 
without trivial factors is $\RR^8\oplus\RR^8$.
The orbits are either 7-spheres (for pairs $(cx,sx)$ with
$|x|^2=1=c^2+s^2$) or Stiefel manifolds $V_2(\RR^8)$. As before,
the union $X$ of the 7-spheres has to be contained in the other focal
manifold $\P$, and thus $X=\P$.

\item[\fbox{$\L=\G_2/\SU(2)=V_2(\RR^7)$}]
By \ref{G2Modules}, the only 14-dimensional
semisimple $\G_2$-$\RR$-modules without
trivial factors are $\RR^7\oplus\RR^7$ and $\mathrm{Ad}$.
The orbits in $\mathrm{Ad}$ belong to an isoparametric hypersurface 
$G_2/\TT^2$ with 6 distinct principal curvatures and multiplicities 
$(2,2)$. This excludes this module.
Therefore, we have the standard action of $\G_2$ on pairs of pure
octonions. The orbits are either 6-spheres, or Stiefel manifolds
$V_2(\RR^7)$. As above, the union $X$ of the 6-spheres in $\SS^{13}$ 
is precisely the other focal manifold $\P$.
\end{description}
There are some more pairs $(G,H)$ where $G$ is almost simple.
We show that they do not belong to isoparametric foliations.

\begin{description}
\item[\fbox{$\Sp(3)/\Sp(1)\times\Sp(1)$ 
and $\Sp(3)/\Sp(1)\times{}^\HH\rho_{3\lambda_1}(\Sp(1))$ are not possible.}]
\ {}
The \linebreak
multiplicities would be $(4,7)$. However, there is no 24-dimensional 
semisimple $\Sp(3)$-$\RR$-module without trivial factors by \ref{CnModules}.

\item[\fbox{$\P=\E_6/\Ffour$ is not possible.}]
This follows from Lemma \ref{NoLargeNormalGroup2}, because 
$\dim(\Ffour)=52>\max\left\{
\binom92,\binom82\right\}=36$. It follows also from Stolz \cite{Stolz98}: 
the multiplicities $(8,9)$ are not possible.

\item[\fbox{$\P=\SU(6)/\Sp(3)$ is not possible.}]
This follows from Lemma \ref{NoLargeNormalGroup2}: 
$\dim(\Sp(3))=21>\max\left\{
\binom42,\binom52\right\}=10$.

Note however that $\SU(6)/\Sp(3)$ and $\SU(5)/\Sp(2)$ are 
$\SU(5)$-equivariant\-ly
homeomorphic, and that $\SU(5)/\Sp(2)$ \emph{is} a focal manifold.
\end{description}

%\newpage

\section{The semisimple case}

We begin with the split case.
Let $\P=K_1/H_1\times K_2/H_2$ be a product of homogeneous
spheres. We show that in 'most' cases, $\P$ cannot be a focal
submanifold of an isoparametric foliation.
The idea is exactly the same as in the last chapter.
Let $S\SUB\P=\SS^{m_1}\times\SS^{m_1+m_2}$ be an
$m_1$-dimensional curvature sphere. Then the composite
\begin{diagram}
S & \rTo & \P \\
&& \dTo_{\pr_1} \\
&&\SS^m_1
\end{diagram}
is a homotopy equivalence, hence
$S\cap (\{x\}\times\SS^{m_1+m_2})\neq\emptyset$ for all 
$x\in\SS^{m_1}$. The fixed point set $\mathrm{Fix}(H_1,\P)$ contains
(at least) $\SS^0\times\SS^{m_1+m_2}\SUB\P$. It follows that
$H_1$ fixes all curvature spheres of the type of $S$ globally.
But this implies that $H_1$ fixes $\P$ pointwise, hence 
$(K_1,H_1)=(\Sp(1),1)$.
\begin{Lem}
If $\P=K_1/H_1\times K_2/H_2$ is a product of homogeneous spheres,
then one of the factors is $\SS^3$ with the regular $\Sp(1)$-action.
\qed
\end{Lem}
To exclude the homogeneous spheres
$K_2/H_2=\SO(m_2+4)/\SO(m_2+3)$, for $m_2\geq 2$, we
apply Lemma \ref{NoLargeNormalGroup2} to $H_2$.
In fact $\dim(H_2)\leq\max\left\{\binom42,\binom{m_2+1}2\right\}$.
Note that $\SS^3\times\SS^4$ cannot be a focal manifold, since
$\pi_3(\SS^3\times\SS^4)\neq\ZZ/2$, so $m_2=1$ is also excluded.

Suppose that $\P=(\SU(r+1)/\SU(r))\times\SS^3$. Then $m_2=2r-2$.
Here we cannot use \ref{NoLargeNormalGroup2}. 
However, $\SU(r)$ cannot act non-trivially on 
$\RR^{2r-1}$, provided that $r\neq 2,4$, cp.~\ref{AnModules}.
Thus, $\SU(r)$ cannot act non-trivially on ${\perp}_p\P$
for $n\neq 2,4$, and this excludes these groups.
Similarly, the cases $\P=(\G_2/\SU(3))\times\SS^3$ or $\P=
(\Spin(7)/\G_2)\times\SS^3$ are not possible.
We are left with the cases
\begin{gather*}
(\Sp(n)/\Sp(n-1))\times\SS^3, \text{ \rlap{$n\geq 2$}} \\
(\Spin(9)/\Spin(7))\times\SS^3 \\
(\SU(5)/\SU(4))\times\SS^3 \\
(\SU(3)/\SU(2))\times\SS^3
\end{gather*}
In the non-split case 
we have also to consider the pairs $(\Sp(n)\times\Sp(2),\Sp(n-1)\cdot
\Sp(1))$ and $(\Sp(n)\times\SU(3),\Sp(n-1)\cdot\Sp(1))$.
By the same reasoning as above, one sees that these groups cannot occur 
for $n\geq 4$, since then $\Sp(n-1)$ cannot act non-trivially on
real vector spaces of dimension less than $4(n-1)$.
We have to consider the pairs
\begin{gather*}
(\Sp(3)\times\Sp(2),\Sp(2)\cdot\Sp(1)) \\
(\Sp(3)\times\SU(3),\Sp(2)\cdot\Sp(1)) \\
(\Sp(2)\times\SU(3),\Sp(1)\cdot\Sp(1)) \rlap{,}
\end{gather*}
and some more cases.
\begin{description}
\item[\fbox{$\P=\Spin(9)/\Spin(7)\times\SS^3$ is not possible.}]
Otherwise we have a 32-dimensional $\Spin(9)$-$\RR$-module. The 
low-dimensional irreducible representations have dimension $9$ and $15$,
so $\RR^{32}=V\oplus W$ splits off the $\Spin(9)$-module $V=\OO\oplus\OO$.

Suppose that $W$ splits off the natural $\SO(9)$-$\RR$-module
$X=\RR^9$. The remaining 7-dimensional $\Spin(9)$-module has to be trivial,
cp.~\ref{BnModules}. Moreover, $\Cen_{\SO(32)}(\Spin(9))\isom\SO(7)$.
Therefore we see that $\P=\SS^{15}\times\SS^3$ factors in such a
way that $\SS^{15}\SUB V\oplus W$ and $\SS^3\SUB\RR^7$.
Let $p=(v,w,x)\in\P\SUB\V\oplus W\oplus\RR^7$ and let 
$N\in{\perp}_p\P\cap\RR^7$ be a non-zero normal vector
(i.e. choose a non-zero vector in $\RR^7$ perpendicular to the
$\SS^3$-orbit of $p$).
Then $N_{(v',w',x)}=N$ defines a normal vector field for 
all $(v',w')\in\SS^{15}$. The Weingarten map is $A_NX=0$ for 
all $X\in T_{(v,w)}\SS^{15}\SUB T_(v,w,x)\P$. 
Thus, the kernel of $A_N$ is at least 15-dimensional, contradicting 
\ref{WeingartenLemma}.

If $W$ is a trivial $\Spin(9)$-module, then by a similar argument $A_N$
has a kernel of dimension at least 15, provided that $N\SUB W$ is a
normal vector.

Finally, suppose that $W\isom V$. Then $\Cen_{\SO(32)}(\Spin(9))^\circ
\isom\SO(2)$; thus, there is no room left for $\Sp(1)$.

\item[\fbox{$\P=\SU(5)/\SU(4)\times\SS^3$ is not possible.}]
The corresponding isoparametric hypersurface would have
multiplicities $(3,6)$, contradicting Stolz' result.

\item[\fbox{$\P=\SU(3)/\SU(2)\times\SS^3$ is not possible.}]
Suppose otherwise. We have a $12$-di\-men\-sional $\SU(3)$-$\RR$-module.
From \ref{AnModules} we see that $\RR^{12}\isom V\oplus W$, where 
$V=\CC^3$ is the natural module.

If $\RR^{12}\isom V\oplus V$, then $\Cen_{\SO(12)}(\SU(3))\isom\SU(2)$;
thus, we know the action of $\SU(3)\times\SU(2)$ on $\RR^{12}$.
However, the principal orbits of this action are isoparametric hypersurfaces
with multiplicities $(2,1)$, and this excludes the module $V\oplus V$.

Otherwise, $W$ splits off a non-zero trivial $\SU(3)$-module.
By the same argument as above, we see that then the Weingarten map
of $\P$ has a 5-dimensional kernel. However, the multiplicities
in this case are $(3,2)$, and thus we have again a contradiction to
\ref{WeingartenLemma}.

\item[\fbox{$\P=\Sp(3)\times\Sp(2)/\Sp(2)\cdot\Sp(1)$ is not possible.}]
Otherwise we have a 24-di\-men\-sional
$\Sp(3)$-$\RR$-module. Then
$\RR^{24}$ splits off the natural $\Sp(3)$-module $V$ by \ref{CnModules}.

If $\RR^{24}\isom V\oplus
V$ (as a $\Sp(3)$-module), then $\Cen_{\SO(24)}(\Sp(3))^\circ
\isom\Sp(2)$, thus we know the orbit structure. The
principal orbits are isoparametric hypersurfaces with
multiplicities $(4,7)$, but both focal manifolds have other
orbit types than the one which we consider here.

Otherwise, $\RR^{24}$ splits off a 12-dimensional
trivial module, cp.~\ref{CnModules}, and therefore
$\Cen_{\SO(24)}(\Sp(3))^\circ
\isom\Sp(1)\times\SO(12)$. However, $\RR^{12}$ cannot be a semisimple
$\Sp(2)$-module without trivial factors, cp.~\ref{BnModules}.
Thus, $\RR^{24}$ cannot be a semisimple $\Sp(3)\times\Sp(2)$-module 
without trivial factors.

\item[\fbox{$\P=\Sp(3)\times\SU(3)/\Sp(2)\cdot\Sp(1)$ is not possible.}]
\ 
Here, the multiplicities \linebreak
would be $(5,6)$, contradicting Stolz 
\cite{Stolz98}.

\item[\fbox{$\P=\Sp(2)\times\SU(3)/\Sp(1)\cdot\Sp(1)$ is not possible.}]
\ Here, the multiplicities \linebreak would be $(5,2)$.
We decompose $\RR^{16}$ as a semisimple $\Sp(2)$-$\RR$-module.
From \ref{BnModules} we see that $\RR^{16}$ splits off the natural
$\Sp(2)$-module $V$. If $\RR^{16}\isom V\oplus V$, then
$\Cen_{\SO(16)}(\Sp(2))\isom\Sp(2)$; however, there is no
embedding $\SU(3)\SUB\Sp(2)$, cp.~\ref{AnModules}.
Put $\RR^{16}=V\oplus W$. Then $\SU(3)\SUB\Cen_{\SO(16)}(\Sp(2))^\circ
\SUB\Sp(1)\times \SO(8)$. Thus, $V$ is a trivial $\SU(3)$-module.
Now $W$ is an $8$-dimensional non-trivial $\SU(3)$-$\RR$-module.
By \ref{AnModules}, $W$ splits off the natural $\SU(3)$-module.
We have a remaining 2-dimensional $\Sp(2)\times\SU(3)$-$\RR$-module, which
has to be trivial. Thus, $\RR^{16}$ cannot be a semisimple
$\Sp(2)\times\SU(3)$-module without trivial factors.

\end{description}
In the split case, only the infinite series 
$\SS^{4n-1}\times\SS^3=\Sp(n)/\Sp(n-1)\times\SS^3$
with multiplicities $(3,4n-4)$ remains; in the non-split case, the series
$\Sp(n)\times\Sp(2)/\Sp(n-1)\cdot\Sp(1)\cdot\Sp(1)$ with multiplicities
$(4,4n-5)$ remains.

%\newpage

\section{The $(4,4n-5)$- and the $(3,4n-4)$-series}

Now we consider the remaining two infinite series.
We start with the series
\[
\P=\Sp(n)\times\Sp(2)/\Sp(n-1)\cdot\Sp(1)\cdot\Sp(1)),
\]
cp.~Section \ref{Quad(4,4n-5)}.
Here, the multiplicities are $(4,4n-5)$ and
$G=\Sp(n)\times\Sp(2)$. 
Consider the action of $\Sp(n-1)$ on the normal space
${\perp}_p\P\isom\RR^{4n-4}$. 

If this action is trivial, then
$\Sp(n-1)$ fixes a curvature sphere in the other focal manifold $\L$.
The normal space of $\L$ is 5-dimensional. If $n-1\geq 3$, then
$\Sp(n-1)$ acts also trivially  on the normal space ${\perp}_\ell\L$, 
since there is no 5-dimensional non-trivial $\Sp(n-1)$-$\RR$-module, cp.
\ref{CnModules}. It follows then that
$\Sp(n-1)$ fixes every point in $\P$, which is absurd.

Thus, $\Sp(n-1)$ acts non-trivially on ${\perp}_p\P\isom\RR^{4n-4}$.
If $n-1\geq 3$, then this is the natural module by \ref{CnModules}, and 
thus the action of $G$ on the normal sphere bundle 
$\mathrm{S}({\perp}\P)\isom\F$ is transitive. It follows that the 
action of $G$ on $\L$ is also transitive. In Section \ref{IsoStiefel}
we classified all homogeneous focal manifolds $\L$ belonging to
isoparametric foliations with multiplicities $(4,4n-5)$, for $n\geq 3$;
then $\L$ is a Stiefel manifold $V_2(\HH^n)$. Note that the action of $G$ on
$\L$ is not irreducible; the normal subgroup $\Sp(n)$ acts already 
transitively on $\L$. We have established the uniqueness of the
isoparametric foliation for $n\geq 4$. We consider the cases $n=2,3$
separately.

\begin{description}
\item[\fbox{$n=3$}]
Then we have a 24-dimensional $\Sp(3)$-module.
Thus $\RR^{24}=V\oplus W$, and $V=\HH^3$ is the natural $\Sp(3)$-module,
cp.~\ref{CnModules}. 

If $W$ is a trivial $\Sp(3)$-$\RR$-module, then
$\Cen_{\SO(24)}(\Sp(3))\isom\Sp(1)\times\SO(12)$.
Thus $\Sp(2)\SUB\SO(12)$, and $W\isom\RR^{12}$ has to be a semisimple
$\Sp(2)$-module without trivial factors. By \ref{BnModules},
this is not possible.

Thus $\RR^{24}=V\oplus V$. Then $\Cen_{\SO(24)}(\Sp(3))\isom\Sp(2)$.
Thus, we have uniqueness of the $G$-action. The principal $G$-orbits are
isoparametric.

\item[\fbox{$n=2$}]
Then we have a 16-dimensional $\Sp(2)$-$\RR$-module.
Thus $\RR^{16}=V\oplus W$, and $V=\HH^2$ is the natural $\Sp(2)$-module,
cp.~\ref{BnModules}. 

If $W$ is a trivial $\Sp(2)$-$\RR$-module, then
$\Cen_{\SO(16)}(\Sp(2))\isom\Sp(1)\times\SO(8)$.
Then the other factor $\Sp(2)$ is contained in $\SO(8)$; as in the case
$n=3$ above, we conclude that it acts as $\Sp(2)$ on $\HH^2\isom\RR^8$.
Then the $G$-orbits are either products of spheres $\SS^7\times\SS^7$
or spheres $\SS^7$; thus, $\P$ cannot be an orbit in this module.

Suppose that $W=\RR^5\oplus\RR^3$ decomposes into the natural
$\SO(5)$-$\RR$-module and a trivial 3-dimensional module.
Then $\Cen_{\SO(16)}(\Sp(2))\isom\Sp(1)\times\SO(3)$, and there is
no more room for the other factor $\Sp(2)$.

Finally, suppose that $\RR^{16}=V\oplus V$. Then
$\Cen_{\SO(16)}(\Sp(2))\isom\Sp(2)$, and we have uniqueness of the
action of $G$. The principal orbits are isoparametric.
\end{description}

\begin{Prop}
Each of the homogeneous spaces in the $(4,4n-5)$-series, $n\geq 2$,
is in a unique way a focal manifold of an isoparametric foliation.
\qed
\end{Prop}
Now we consider the series
\[
\P=\Sp(n)/\Sp(n-1)\times\Sp(1).
\]
The multiplicities are $(3,4n-4)$.
First we describe the known example.
\begin{Num}\textsc{Example}
\label{Quad(3,4n-4)}\

\psn
Let
\begin{align*}
\P&=\{(x,xa)\in\HH^n \oplus\HH^n|\ x\in\HH^n,\ a\in\HH,\ 
|x|^2=1/2,\ |a|^2=1\} \\
\L&=\{(u,v)\in\HH^n \oplus\HH^n|\ (u|v)=0,\ |u|^2+|v|^2=1\}
\end{align*}
(here $(-|-)$ denotes the standard positive definite hermitian
form on $\HH^n$).
These are focal manifolds of an isoparametric hypersurface with
$g=4$ distinct principal curvatures. The multiplicities are
$(3,4n-4)$. The incidence in the corresponding point-line
geometry is as follows. A point $p=(x,xa)$ is incident with the
line $\ell=(u,v)$ if and only if 
\[
x=\frac{1}{\sqrt 2}(u+v\bar a).
\]
The group $\Sp(n)\times\Sp(1)\times\Sp(1)$ acts as a group of
automorphisms and isometries on the isoparametric foliation
(the group $\Sp(n)$ as a matrix group from the left, the
group $\Sp(1)\times\Sp(1)$ as multiplication by pairs of
scalars from the right), and the action is
transitive on the point space $\P$.
The subgroup $\Sp(n)\times\Sp(1)$ (where $\Sp(1)$ is any of the
two factors) acts still transitively on the points. 
\end{Num}
Now we get back to the unknown homogeneous focal manifold.
If the action of $\Sp(n-1)$ on ${\perp}_p\P\isom\RR^{4n-3}$ is trivial, 
then $\Sp(n-1)$ has to act on the space ${\perp}_\ell\L\isom\RR^4$,
which is impossible for $n\geq 3$. 
Thus we have a non-trivial action of $\Sp(n-1)$ on 
${\perp}_p\P\isom\RR^{4n-3}$, provided that $n\geq 3$.
If $n\geq 4$, then the only possibility is that 
$\RR^{4n-4}\isom V\oplus \RR$ is the natural $\Sp(n-1)$-module $V$
plus a 1-dimensional trivial module. Thus, we see that
the normal bundle of $\P$ is
\[
{\perp}\P=(\Sp(n)\times\Sp(1))\times_{\Sp(n-1)}(V\oplus\RR).
\]
The isoparametric hypersurface is the unit sphere bundle of
the normal bundle ${\perp}\P$.

Consider the $\Sp(n)$-$\RR$-module $\RR^{8n}$. From \ref{CnModules}
we see that $\RR^{8n}\isom V\oplus W$ splits off the natural
module $V=\HH^n$ (this holds for all $n\geq 2$).

Suppose that $W$ is the trivial module. Let $(v,w)\in\P\SUB V\oplus W$.
Then $w\neq 0$ and $\Sp(n)\cdot(v,w)$ is a sphere $\SS^{4n-1}$.
The normal isotropy representation of $\Sp(n-1)$ on
${\perp}_p(\Sp(n)\cdot(v,w))\isom\RR\oplus W$ is trivial; it follows that
the normal isotropy representation of $\Sp(n-1)$ on
${\perp}_p\P$ is trivial as well. But this is impossible for $n\geq 3$.

If $n\geq 3$, then $W$ has to be trivial or isomorphic to the natural
module $V$, cp.~\ref{CnModules}.
Thus, $\RR^{8n}\isom V\oplus V$ as a $\Sp(n)$-module, provided that
$n\geq 3$.

\begin{description}
\item[\fbox{$n=2$}]
Then it is also possible that $W=X\oplus\RR^3$ splits off the
5-dimensional $\SO(5)$-$\RR$-module $X$. 
Suppose that we have this module. Then \linebreak
$\Cen_{\SO(16)}(\Sp(2))\isom
\Sp(1)\times\SO(3)$.
Let $p=(w,x,z)\in\P\SUB\V\oplus W\oplus\RR^3$. If $w\neq 0\neq x$, then
the $\Sp(2)$-stabilizer of $p$ is trivial. But for $p\in\P$, the 
stabilizer is $\Sp(1)$. Thus $x=0$, and $\P\SUB V\oplus\RR^3$, which is
impossible. Thus $\RR^{16}=V\oplus V$ as a $\Sp(2)$-module.
\end{description}
We have established for all $n\geq 2$ that
$\RR^{8n}=V\oplus V$ as a $\Sp(n)$-module, and
\[
\SS^3\SUB\Cen_{\SO(8n)}(\Sp(n))\isom\Sp(2).
\] 
Let $p=(x,y)\in\P\SUB V\oplus V$. 
The $\Sp(n)$-orbit of $p$ is a sphere, therefore $(x,y)$ have to be 
$\HH$-linearly dependent. Put 
\[
\widehat\P=\{(xc,xs)|\ x\in\SS^{4n-1},c,s\in\HH,|c|^2+|s|^2=1\}.
\]
Then $\P\SUB\widehat\P$. 
As in Section \ref{Quad(4,4n-5)}, consider the map 
\[
\phi:\widehat\P\too\SS^4, \quad (xc,xs)\mapstoo(|xc|^2-|xs|^2,2\bar cs).
\]
The fibres of $\phi$ are precisely
the $\Sp(n)$-orbits in $\widehat\P$. The group 
$\Cen_{\SO(8n)}(\Sp(n))\isom\Sp(2)$ permutes these fibres and acts as
$\SO(5)$ on the image $\SS^4$ of $\phi$; in other words, $\phi$ is
$\Sp(2)$-equivariant. Now $\SS^3\SUB\Sp(2)$
has an orbit on $\P$ which meets every $\Sp(n)$-orbit precisely once, 
and which is
homeomorphic to $\SS^3$. Thus we have to consider subgroups isomorphic 
to $\SS^3$ in $\SO(5)$. By \ref{AnModules}, there is just
one conjugacy class of such groups, $\SU(2)\SUB\SO(4)\SUB\SO(5)$.
We lift this group into $\Sp(2)$ to obtain $\SS^3\SUB\Sp(2)$.

So far, we have determined the subgroup $\SS^3\SUB\Sp(2)$; it is
conjugate to $\Sp(1)\SUB\Sp(2)$, and the action of $G=\Sp(n)\times\Sp(1)$ 
on $\RR^{8n}$. The problem is that there are many $G$-orbits in 
$\SS^{8n-1}$ which are homeomorphic to $\P$, and we have to find the right 
one. Note that for each orbit homeomorphic to $\P$, the normal bundle has
the form described above. Therefore
we know the orbits which are contained in the hypersurfaces, or in $\P$; 
they are either homeomorphic to $\P$ (and thus $(4n+2)$-dimensional), or
$(8n-3)$-dimensional (for pairs $(x,y)$, where $x,y$ are linearly
independent, but $(x|y)\neq 0$). However, there are also
$(8n-6)$-dimensional orbits (for pairs $(u,v)$, where $u,v$
are linearly independent, with $(u|v)=0$)
and $(4n-1)$-dimensional orbits (for $x=0$ or $y=0$). These
orbits have to be contained in the other focal manifold $\L$.
In fact, $\L$ consists entirely of these orbits, because the
union of these orbits is a manifold of the same dimension as $\L$
--- here we look at our model as in \ref{Quad(3,4n-4)}. Thus, $\L$ is
unique, and this establishes the uniqueness of the isoparametric
foliation.

\begin{Prop}
Each of the spaces in the $(3,4n-4)$-series, $n\geq 2$,
is in a unique way a focal manifold of an isoparametric foliation.
\qed
\end{Prop}

%\newpage

\section{Summary}

In this chapter, we have established the following result.

\begin{Thm}
\label{MainIsoTheorem}
Let $\F\SUB\SS^{r+1}$ be an isoparametric hypersurface with four distinct
principal curvatures and multiplicities $(m_1,m_2)$.
Let $\P$, $\L$ denote the focal manifolds of $\F$. 
Suppose that $G\SUB\SO(r+2)$ is a compact connected subgroup which acts 
transitively on the focal manifold $\P$ and assume that
$m_1\geq 3$. 

There are the following possibilities for the isoparametric
foliation (and no others).
\begin{description}
\item[\textbf{$\F$ is homogeneous.}]
Then $\F$ belongs to the homogeneous series with multiplicities
$(n-2,1)$, for $n\geq 5$, or
$(2n-3,2)$, for $n\geq 3$, or
$(4n-5,4)$, or dually $(4,4n-5)$, for $n\geq 2$, or has
multiplicities $(4,5)$, $(5,4)$, $(6,9)$, or $(9,6)$ .
\item[\textbf{$\F$ is not homogeneous.}]
Then $\F$ is of Clifford type. 
Either it belongs to the non-ho\-mo\-ge\-ne\-ous Clifford
series with multiplicities $(3,4n-4)$, for $n\geq 2$, or $(m_1,m_2)=
(7,8)$, and $(\P,\L,\F)=(M_{+},M_{-},M)$ is dual to the definite
isoparametric foliation of Clifford type with multiplicities $(8,7)$
(cp.~Ferus-Karcher-M\"unzner \cite{FeKaMu81} for the terminology).
\end{description}
\end{Thm}

\begin{Cor}
Let $\F$ be an isoparametric hypersurface with four
distinct principal curvatures. If both focal manifolds
$\P$ and $\L$ are homogeneous, then $\F$ is homogeneous.
\qed
\end{Cor}
The following result was recently proved by Wolfrom.
\begin{Thm}[Wolfrom]
\label{WolfromTheorem}
Let $\F$ be an isoparametric hypersurface with $g=3,4,6$, and assume that
$m_1=2$. If the focal manifold $\P$ is homogeneous, then $\F$ is
homogeneous \cite{Wol}.
\end{Thm}

\begin{Cor}
Let $\F$ be an isoparametric hypersurface, and assume that one focal
manifold is homogeneous. Then either $\F$ is homogeneous, or
$\F$ is of Clifford type, with $g=4$ and $(m_1,m_2)=(8,7)$ or
$(m_1,m_2)=(3,4k)$.

\begin{proof}
We may assume that $\P$ is homogeneous, and that $g=4,6$.
If $m_1\geq 3$, then the result
follows from our Theorem \ref{MainIsoTheorem} above, and for $m_1=2$
it follows from Wolfrom's
Theorem \ref{WolfromTheorem}. If $m_1=1$,
then the results of Takagi \cite{Takagi} and
Dorfmeister-Neher \cite{DorNeh85} apply.
\end{proof}
\end{Cor}
Finally, we state Immervoll's theorem.
\begin{Thm}[Immervoll]
The point-line geometry associated to an isoparametric hypersurface
with $g=4$ is a compact connected and smooth generalized quadrangle
\cite{StIm}.
\end{Thm}
The following difficult problem is open.
\begin{Num}\textsc{Conjecture}
{\em Every isoparametric hypersurface is either homogeneous, or
of Clifford type.}
\end{Num}
Looking back at our proof of Theorem \ref{MainIsoTheorem}, we note
that we have obtained 
the following classification of compact connected
transitive groups.
Let $\F$ be an isoparametric hypersurface as above, and let 
$G\SUB\SO(r+2)$ be a compact connected group which acts
transitively on $\P$. There are precisely the following 
possibilities for $G$. 

\begin{Num}\textsc{Multiplicities $(n-2,1)$, for $n\geq 5$}\psn
Then $\P=V_2(\RR^n)$ is a Stiefel
manifold, and $G=\SO(n)$ or $G=\SO(n)\cdot \SO(2)$, 
where $\SO(n)\cap\SO(2)=\SO(n)\cap\{\pm1\}$.

For $n=8$, there are the additional possibilities
$G=\Spin(7)$ and $G=\Spin(7)\cdot\SO(2)$, where
$\Spin(7)\cap\SO(2)=\{\pm1\}$, and for
$n=7$ there is the additional possibility
$G=\G_2$ or $G=\G_2\cdot\SO(2)$, where $\G_2\cap\SO(2)=1$.

The action of $G$ is transitive on $\L$ and $\F$ if and only if $G$ 
contains the second factor $\SO(2)$.
\end{Num}

\begin{Num}\textsc{Multiplicities $(2n-3,2)$, for $n\geq 3$}\psn
Then $\P=V_2(\CC^n)$ is a Stiefel
manifold, and $G$ is one of the groups
$\SU(n)$, $\U(n)$, $\SU(n)\cdot\SU(2)$, $\U(n)\cdot\U(2)$, where
$\SU(n)\cap\SU(2)=\SU(n)\cap\{\pm1\}$ and
$\U(n)\cap\SU(2)=\{\pm1\}$.

The group $G$ acts transitively on $\L$ or $\F$ only if $G$ contains
the $\SU(2)$ factor.
\end{Num}

\begin{Num}\textsc{Multiplicities $(4,4n-5)$ or $(4n-5,4)$, for $n\geq 2$}\psn
Then $\L=V_2(\HH^n)$ is a Stiefel manifold, and $G$ is of the form
$G=\Sp(n)\cdot K$, where $K$ is a connected subgroup of $\Sp(2)$,
and $\Sp(2)\cap\Sp(n)=\{\pm1\}$ (it is not difficult to determine
all connected subgroups of $\Sp(2)$). The action of $G$ is
transitive on the focal manifold $\L$, and
transitive on $\P$ or $\F$ if and only if $K=\Sp(2)$.
\end{Num}

\begin{Num}\textsc{Multiplicities $(4,5)$ or $(5,4)$}\psn
The group is either $G=\SU(5)$ or $G=\U(5)$. In each case,
$G$ acts transitively on $\P$, $\L$ and $\F$.
\end{Num}

\begin{Num}\textsc{Multiplicities $(9,6)$ or $(6,9)$}\psn
If $G$ acts transitively on $\L$, then $G=\Spin(10)$ or
$G=\Spin(10)\cdot\U(1)$, and $G$ acts transitively on $\L$ and $\F$.
The groups $G=\Spin(9)$ and $G=\Spin(9)\cdot\U(1)$ act
transitively on $\P$, but not on $\L$ or $\F$.
\end{Num}

\begin{Num}\textsc{Multiplicities $(3,4n-4)$, for $n\geq 2$}\psn
Then $\P=\SS^3\times\SS^{4n-5}$, and $G$ is of the form
$G=\Sp(n)\cdot K$, where $K$ is a connected subgroup of 
$\Sp(1)\times\Sp(1)$,
containing a subgroup isomorphic to $\Sp(1)$.
The group $G$ acts transitively on $\P$, but not on $\L$ or $\F$.
\end{Num}

\begin{Num}\textsc{Multiplicities $(7,8)$}\psn
If $G$ acts transitively on $\P$, then $G=\Spin(9)$ or
$G=\Spin(9)\cdot\SO(2)$, and $G$ does not act transitively on $\L$ or 
$\F$. Here, $\P=V_2(\OO^2)=M_{+}$.
\end{Num}

\backmatter
%%%%%%%%%%%%%%%%%%%%%%%%%%%%%%%%%%%%%%%%%%%%%%%%%%%%%%%%%%%%%%%%%%%%%%%%
%                                                                      %
%                                                                      %
%                   Compact homogeneous quadrangles                    %
%                                                                      %
%                          Linus Kramer                                %
%                                                                      %
%                           Memoirs AMS                                %
%                                                                      %
%                         Wuerzburg 2000                               %
%                                                                      %
%                                                                      %
%                           CHQlit.tex                                 %
%                                                                      %
%                                                                      %
%                                                                      %
%%%%%%%%%%%%%%%%%%%%%%%%%%%%%%%%%%%%%%%%%%%%%%%%%%%%%%%%%%%%%%%%%%%%%%%%

\include{index}

\begin{thebibliography}{95}
%\addcontentsline{toc}{chapter}{Bibliography}

\bibitem{Abramenko}
P.~Abramenko,
\textit{Twin buildings and applications to $S$-arithmetic groups.}
LNM 1641, Springer-Verlag, Berlin (1996) x+123 pp.
\textsf{MR~99k:20060}

\bibitem{Abrsch83}
U.~Abresch,
\textit{Isoparametric hypersurfaces with four or six distinct principal
  curvatures.} 
Math. Ann. 264 (1983) 283--302.
\textsf{MR~85g:53052b}

\bibitem{BGS}
W.~Ballmann, M.~Gromov, and V.~Schr\"oder,
\textit{Manifolds of nonpositive curvature.} 
Birkh\"auser, Boston, MA. (1985) vi+263 pp.
\textsf{MR~87h:53050}

\bibitem{Biller}
H.~Biller,
\textit{Actions of compact groups on spheres and on
generalized quadrangles.}
Ph.D. Thesis, Stuttgart: {M}ath. {F}ak., {U}niv. Stuttgart,
(1999) xxvi+195 pp.

\bibitem{Bletz}
O.~Bletz,
\textit{Ein Beweis f\"ur die Lokalkompaktheit der
 Automorphismengruppen kompakter Polygone.}
\textit{Diplomarbeit}, W\"urzburg: {M}ath. {F}ak., {U}niv. W\"urzburg,
(1999) 40 pp.
 

\bibitem{BoediJoswig93}
R.~B{\"o}di and M.~Joswig,
\textit{Tables for an effective enumeration of real representations of
  quasi-simple {L}ie groups.} 
Sem. Sophus Lie 3 (1993) 239--253.
\textsf{MR~95f:22003}

\bibitem{BoeKr95}
R.~B{\"o}di and L.~Kramer,
\textit{On homomorphisms between generalized polygons.} 
Geom. Dedicata 58 (1995) 1--14.
\textsf{MR~96k:51017}

\bibitem{BorelSphere}
A.~Borel,
\textit{Some remarks about {L}ie groups transitive on spheres and tori}. 
Bull. Amer. Math. Soc. 55 (1949) 580--587.
\textsf{MR~10,680c}

\bibitem{BorelThesis}
A.~Borel,
\textit{Sur la cohomologie des espaces fibr\'es principaux et des espaces
  homog\`enes de groupes de {L}ie compacts.} 
Ann. of Math. 57 (1953) 115--207.
\textsf{MR~14,490e}

\bibitem{Bor54}
A.~Borel,
\textit{Sur l'homologie et la cohomologie des groupes de {L}ie compacts
  connexes.} 
Amer. J. Math. 76 (1954) 273--342.
\textsf{MR~16,219b}

\bibitem{BorelTrans}
A.~Borel,
\textit{Seminar on transformation groups.}
With contributions by G.~Bredon, E.~E.~Floyd, D.~Montgomery, R.~Palais,
Princeton Univ. Press, Princeton, NJ. (1960) vii+245 pp.
\textsf{MR~22\#7129}

\bibitem{BorelLNM}
A.~Borel,
\textit{Topics in the homology theory of fibre bundles.} 
LNM 36, Springer-Verlag, Berlin-New York (1967) 95 pp.
\textsf{MR~36\#4559}

\bibitem{BorelSieb49}
A.~Borel and J.~De~Siebenthal,
\textit{Les sous-groupes ferm\'es de rang maximum des groupes
 de {L}ie clos.} 
Comment. Math. Helv. 23 (1949) 200--221.
\textsf{MR~11,326d}

\bibitem{BorovikNesin}
A.~Borovik and A.~Nesin,
\textit{Groups of finite Morley rank.} 
Oxford University Press, New York (1994) xviii+409 pp.
\textsf{MR~96c:20004}

\bibitem{Bre61}
G.~E.~Bredon,
\textit{On homogeneous cohomology spheres.}
Ann. of Math. 73 (1961) 556--565.
\textsf{MR~23\#A243}

\bibitem{BredTTG}
G.~E.~Bredon,
\textit{Introduction to compact transformation groups.} 
Academic Press, New York (1972) xiii+459 pp.
\textsf{MR~54\#1265}

\bibitem{BredonSheaf}
G.~E.~Bredon,
\textit{Sheaf theory. Second ed.}
Springer GTM 170,
Springer Verlag, NewYork (1997) xii+502 pp.
\textsf{MR~98g:55005, MR~36\#4552}

\bibitem{BurSpa87}
K.~Burns and R.~Spatzier,
\textit{On topological {T}its buildings and their classification.} 
Inst. Hautes \'Etudes Sci. Publ. Math. 65 (1987) 5--34.
\textsf{MR~88g:53049}

\bibitem{Cart39}
E.~Cartan,
\textit{Sur des familles remarquables d'hypersurfaces isoparam{\'e}triques
  dans les espaces sph{\'e}riques.} 
Math. Z. 45 (1939) 335--367.
\textsf{MR~1,28f}

\bibitem{CarSer52} 
H.~Cartan and J.-P. Serre,
\textit{Espaces fibr\'es et groupes d'homotopie. {I}{I}. {A}pplications.} 
C. R. Acad. Sci. Paris 234 (1952) 393--395.
\textsf{MR~13,675b}

\bibitem{CeRy85}
T.~E. Cecil and P.~J. Ryan,
\textit{Tight and taut immersions of manifolds}. 
Pitman Research Notes in Mathematics 107,
Pitman, Boston, MA. (1985). vi+335 pp.
\textsf{MR~87b:53089}

\bibitem{ConOlm}
S.~Console and C.~Olmos,
\textit{Clifford systems, algebraically constant second fundamental
form and isoparametric hypersurfaces.} 
Manuscripta Math. 97 (1998) 335--342.  
\textsf{MR~99j:53072}

\bibitem{tDieck}
T.~tom Dieck,
\textit{Transformation groups.}
de Gruyter Studies in Mathematics, 8, Berlin (1987) x+312 pp.
\textsf{MR~89c:57048}

\bibitem{Dold}
A.~Dold,
\textit{Lectures on algebraic topology.}
Springer Verlag, New York-Berlin (1972) xi+377 pp.
\textsf{MR~54\#3685}

\bibitem{DorNeh83a}
J.~Dorfmeister and E.~Neher,
\textit{Isoparametric triple systems of {F}{K}{M}-type. {I}.} 
Abh. Math. Sem. Univ. Hamburg 53 (1983) 191--216.
\textsf{MR~85i:17002a}

\bibitem{DorNeh83b}
J.~Dorfmeister and E.~Neher,
\textit{Isoparametric triple systems of {F}{K}{M}-type. {I}{I}.} 
Manuscripta Math. 43 (1983) 13--44.
\textsf{MR~85i:17002b}

\bibitem{DorNeh84}
J.~Dorfmeister and E.~Neher,
\textit{Isoparametric triple systems of {F}{K}{M}-type. {I}{I}{I}.} 
Algebras Groups Geom. 1 (1984) 305--343.
\textsf{MR~86f:17005}

\bibitem{DorNeh85}
J.~Dorfmeister and E.~Neher,
\textit{Isoparametric hypersurfaces, case $g=6,\;m=1$.}
Comm. Algebra 13 (1985) 2299--2368.
\textsf{MR~87d:53096}

\bibitem{DynkinAMST57}
E.~B. Dynkin,
\textit{Semisimple subalgebras of semisimple {L}ie algebras.}
Amer. Math. Soc. Transl. (2) 6 (1957) 111--244.
\textsf{MR~13,904c}

\bibitem{Eng89}
R.~Engelking,
\textit{General topology.}
Heldermann Verlag, Berlin, second edition (1989) viii+529 pp.
\textsf{MR~91c:54001}

\bibitem{EschHein}
J.~H. Eschenburg and E.~Heintze,
\textit{On the classification of polar representations.}
Math. Z.  232 (1999) 391--398.
\textsf{MR~2001g:53099}

\bibitem{EschSchr91}
J.-H. Eschenburg and V.~Schroeder,
\textit{Tits distance of {H}adamard manifolds and isoparametric
  hypersurfaces.}
Geom. Dedicata 40 (1991) 97--101.
\textsf{MR~92i:53038}

\bibitem{FHT}
Y.~F\'elix, S.~Halperin and J.-C.~Thomas,
\textit{Rational homotopy theory.}
Springer GTM 205,
Springer Verlag, New York (2001) xxv+535 pp.

\bibitem{FeKaMu81}
D.~Ferus, H.~Karcher, and H.F. M{\"u}nzner,
\textit{Cliffordalgebren und neue isoparametrische {H}yperfl{\"a}chen.}
Math. Z. 177 (1981) 479--502.
\textsf{MR~83k:53075}

\bibitem{FFG86}
A.~T. Fomenko, D.~B. Fuchs, and V.~L. Gutenmacher,
\textit{Homotopic topology.}
Akad\'emiai Kiad\'o, Budapest (1986) 310 pp.
\textsf{MR~88f:55001}

\bibitem{GroHal87}
K.~Grove and S.~Halperin,
\textit{Dupin hypersurfaces, group actions and the double mapping cylinder.}
J. Differential Geom. 26 (1987) 429--459.
\textsf{MR~89h:53113}

\bibitem{GK90}
Th.~Grundh{\"o}fer and N.~Knarr,
\textit{Topology in generalized quadrangles.}
Topology Appl. 34 (1990) 139--152.
\textsf{MR~91c:51025}

\bibitem{GKK95}
Th.~Grundh{\"o}fer, N.~Knarr, and L.~Kramer,
\textit{Flag-homogeneous compact connected polygons.}
Geom. Dedicata 55 (1995) 95--114.
\textsf{MR~96a:51009}

\bibitem{GKK98}
Th.~Grundh{\"o}fer, N.~Knarr, and L.~Kramer,
\textit{Flag-homogeneous compact connected polygons, II.}
Geom. Dedicata 83 (2000) 1--29.

\bibitem{GKK00}
Th.~Grundh{\"o}fer, N.~Knarr, and L.~Kramer,
\textit{The classification of compact buildings.}
Preprint W\"urzburg (1999).

\bibitem{GL95}
Th.~Grundh{\"o}fer and R.~L{\"o}wen,
\textit{Linear topological geometries.}
In: \textit{Handbook of incidence geometry}, p. 1255--1324. 
North-Holland, Amsterdam (1995).
\textsf{MR~97c:51008}

\bibitem{GvM90}
Th.~Grundh{\"o}fer and H.~Van~Maldeghem,
\textit{Topological polygons and affine buildings of rank three.}
Atti Sem. Mat. Fis. Univ. Modena 38 (1990) 459--479.
\textsf{MR~92c:51028}

\bibitem{HeLi}
E.~Heintze and X.~Liu,
\textit{Homogeneity of infinite-dimensional isoparametric submanifolds.}
Ann. of Math. 149 (1999) 149--181.
\textsf{MR~2000c:58007}

\bibitem{HoMo}
K.~H.~Hofmann and S.~A.~Morris,
\textit{The structure of compact groups.}
De Gruyter, Berlin (1998) xviii+835 pp.
\textsf{MR~99k:22001}

\bibitem{HsLaw71}
W.-y.~Hsiang and H.~B.~Lawson, Jr.,
\textit{Minimal submanifolds of low cohomogeneity.}
J. Differential Geometry 5 (1971) 1--38.
\textsf{MR~45\#7645}

\bibitem{HsPaTe88}
W.-y.~Hsiang, R.~S.~Palais, and C.-L.~Terng,
\textit{The topology of isoparametric submanifolds.}
J. Differential Geom. 27 (1988) 423--460.
\textsf{MR~89m:53104}

\bibitem{HsSu68}
W.-y.~Hsiang and J.~C.~Su,
\textit{On the classification of transitive effective actions on {S}tiefel
  manifolds.}
Trans. Amer. Math. Soc. 130 (1968) 322-336.
\textsf{MR~36\#4581}

\bibitem{StIm}
S.~Immervoll,
\textit{Smooth projective planes, smooth generalized quadrangles, and
  isoparametric hypersurfaces.}
Ph.D. {T}hesis, {T\"u}bingen: {M}ath. {F}ak., {U}niv. {T\"u}bingen,
in preparation (2001).

\bibitem{Karch88}
H.~Karcher,
\textit{A geometric classification of positively curved symmetric spaces
 and the isoparametric construction of the {C}ayley plane.}
Ast\'erisque (163-164) 6 (1988) 111--135.
\textsf{MR~90g:53063}

\bibitem{KleinerLeeb}
B.~Kleiner and B.~Leeb,
\textit{Rigidity of quasi-isometries for symmetric spaces
  and Euclidean buildings.}
Inst. Hautes \'Etudes Sci. Publ. Math. 86 (1997) 115--197.
\textsf{MR~98m:53068}

\bibitem{Kn90}
N. Knarr,
\textit{The nonexistence of certain topological polygons.}
Forum Math. 2 (1990) 603--612.
\textsf{MR~92a:51018}

\bibitem{KnaKra95}
N.~Knarr and L.~Kramer,
\textit{Projective planes and isoparametric hypersurfaces.}
Geom. Dedicata 58 (1995) 193--202.
\textsf{MR~96i:53059}

\bibitem{KoNo}
S.~Kobayashi and K.~Nomizu,
\textit{Foundations of differential geometry, Vol.~I.}
Interscience Publishers, New York-London (1963) xi+329 pp.
\textsf{MR~27\#2945}

\bibitem{Kr94}
L.~Kramer,
\textit{Compact polygons.}
{Ph.D.} {T}hesis, {T\"u}bingen: {M}ath. {F}ak., {U}niv. {T\"u}bingen,
(1994) vi+72 pp. 
\textsf{Zbl 844.51006} \\
available as math.DG/0104064 at the Mathematics ArXiv,

{\tt http://front.math.ucdavis.edu/math.DG/0104064}


\bibitem{Kr98}
L.~Kramer,
\textit{Octonion Hermitian quadrangles.}
Bull. Belg. Math. Soc. Simon Stevin 5 (1998) 353--362.
\textsf{MR~99g:51003}

\bibitem{KramerHabil}
L.~Kramer,
\textit{Compact homogeneous quadrangles and focal manifolds.}
\textit{Habilitationsschrift,} W\"urzburg: {M}ath. {F}ak., {U}niv.
W\"urzburg,
(1998) 138 pp.

\bibitem{Ov2}
L.~Kramer,
\textit{Compact ovoids in quadrangles II. The classical quadrangles.}
Geom. Dedicata 79 (2000) 179--188.
\textsf{MR~2001b:51007}

\bibitem{Ov3}
L.~Kramer,
\textit{Compact ovoids in quadrangles III. Clifford algebras
 and isoparametric hypersurfaces.}
Geom. Dedicata 79 (2000) 321--339.

\bibitem{Ov1}
L.~Kramer and H.~Van Maldeghem,
\textit{Compact ovoids in quadrangles I. Geometric constructions.}
Geom. Dedicata 78 (1999) 279--300.
\textsf{MR~2000m:51008}

\bibitem{KTV}
L.~Kramer, K.~Tent, and H.~Van Maldeghem,
\textit{Simple groups of finite Morley rank and Tits buildings.}
Israel J. Math. 109 (1999) 189--224.
\textsf{MR~2000f:51022}
	        
\bibitem{KS2}	 
M.~Kreck and S.~Stolz,
\textit{A diffeomorphism classification of $7$-dimensional homogeneous
 {E}instein manifolds with
${\rm SU}(3)\times{\rm SU}(2)\times{\rm U}(1)$-symmetry.}
Ann. of Math. 127 (1988) 373--388.
\textsf{MR~89c:57042}

\bibitem{KS1}
M.~Kreck and S.~Stolz,
\textit{Some nondiffeomorphic homeomorphic homogeneous $7$-manifolds with
  positive sectional curvature.}
J. Differential Geom. 33 (1991) 465--486.
\textsf{MR~92d:53043}

\bibitem{KS1a}
M.~Kreck and S.~Stolz,
\textit{A correction on:
'Some nondiffeomorphic homeomorphic homogeneous $7$-manifolds with
positive sectional curvature'.}
J. Differential Geom. 49 (1998) 203--204.
\textsf{MR~99h:53069}

\bibitem{KuehCec}
W.~K\"uhnel and T.~Cecil,
\textit{Bibliography on tight, taut and isoparametric submanifolds.}
In: \textit{Tight and taut submanifolds, Berkeley 1994}, p.~307--339,
MSRI Publ. 32,
Cambridge Univ. Press, Cambridge (1997).

\bibitem{Leeb}
B.~Leeb,
\textit{A characterization of irreducible symmetric spaces
 and Euclidean buildings of higher rank by their
 asymptotic geometry.}
\textit{Habilitationsschrift,} Bonn (1997) 42 pp.

Bonner Math. Schriften 326 (2000).

\bibitem{Lowen}
R.~L\"owen,
\textit{Homogeneous compact projective planes.}
J. Reine Angew. Mathematik 321 (1981) 217--220.
\textsf{MR 82c:51023}

\bibitem{Man67}
L.~N. Mann,
\textit{Dimensions of compact transformation groups.}
Michigan Math. J. 14 (1967) 433--444.
\textsf{MR~36\#3916}

\bibitem{Mark}
E.~Markert,
\textit{Isoparametric hypersurfaces and generalized quadrangles.}
{D}iplomarbeit, {W\"u}rzburg: {M}ath. {I}nst., {U}niv. {W\"u}rzburg
(1999) 45 pp.

\bibitem{McCl85}
J.~McCleary,
\textit{User's guide to spectral sequences.}
Publish or Perish Inc., Wilmington, DEL. (1985) xiv+423 pp.
\textsf{MR~87f:55014}

\bibitem{Mi95}
M.~Mimura,
\textit{Homotopy theory of {L}ie groups.}
In: \textit{Handbook of algebraic topology}, p. 951--991.
  North-Holland, Amsterdam (1995).
\textsf{MR~97c:57038}

\bibitem{MiTo91}
M.~Mimura and H.~Toda,
\textit{Topology of {L}ie groups. {I}, {I}{I}.}
Transl. of Math. Monographs 91,
American Mathematical Society, Providence, RI. (1991) iv+451 pp.
\textsf{MR~92h:55001}

\bibitem{Mon50}
D.~Montgomery,
\textit{Simply connected homogeneous spaces.}
Proc. Amer. Math. Soc. 1 (1950) 467--469.
\textsf{MR~12,242c}

\bibitem{MonYang}
D.~Montgomery and C.~T.~Yang,
\textit{Orbits of highest dimension.}
Trans. Amer. Math. Soc. 87 (1958) 284--293.
\textsf{MR~20\#6705}

\bibitem{MoZi}
D.~Montogmery and L.~Zippin,
\textit{Topological transformation groups.}
Interscience Publishers, New York-London (1955) xi+282 pp.
\textsf{MR~17,282b}

\bibitem{Mostow}
G.~D.~Mostow,
\textit{Strong rigidity of locally symmetric spaces}.
Princeton University Press,
Princeton, NJ. (1973) v+195 pp.
\textsf{MR~52\#5874}

\bibitem{Munz80}
H.~F. M{\"u}nzner,
\textit{Isoparametrische {H}yperfl\"achen in {S}ph\"aren, I.}
Math. Ann. 251 (1980) 57--71.
\textsf{MR~82a:53058}

\bibitem{Munz81}
H.~F. M{\"u}nzner,
\textit{Isoparametrische {H}yperfl\"achen in {S}ph\"aren, II.}
Math. Ann. 256 (1981) 215--232.
\textsf{MR~82m:53053}

\bibitem{FN}
F.~Nietzsche,
\textit{Also sprach Zarathustra}.
Deutscher Taschenbuch Verlag/de Gruyter (1988) 420~pp.

\bibitem{Olm}
C.~Olmos,
\textit{Isoparametric submanifolds and their homogeneous structures.}
J. Differential Geom. 38 (1993) 225--234.
\textsf{MR~95a:53102}

\bibitem{Oni94}
A.~L.~Onishchik,
\textit{Topology of transitive transformation groups}. 
Johann Ambrosius Barth Verlag GmbH, Leipzig (1994) xvi+300 pp.
\textsf{MR~95e:57058}

\bibitem{OniVin90}
A.~L. Onishchik and {\`E}.~B. Vinberg,
\textit{Lie groups and algebraic groups}.
Springer-Verlag, Berlin (1990). xx+328 pp.
\textsf{MR~91g:22001}

\bibitem{PaTe88}
R.~S.~Palais and C.-L.~Terng,
\textit{Critical point theory and submanifold geometry}. 
LNM 1353 Springer-Verlag, Berlin (1988). x+272 pp.
\textsf{MR~90c:53143}

\bibitem{RohlfsSpringer}
J.~Rohlfs and T.~A.~Springer,
\textit{Applications of buildings.}
In: \textit{Handbook of incidence geometry}, p.~1085--1114.
North-Holland, Amsterdam (1995).
\textsf{MR~97b:20041}

\bibitem{Salzmann}
H.~Salzmann,
\textit{Homogene kompakte projektive Ebenen.}
Pacific J. Math. 60 (1975) 217--234.
\textsf{MR 53\#3878}

\bibitem{CPP95}
H.~Salzmann, D.~Betten, Th.~Grundh{\"o}fer, 
H.~H{\"a}hl, R.~L{\"o}wen, and  M.~Stroppel,
\textit{Compact projective planes}.
Walter de Gruyter \& Co., Berlin (1995) xiv+688 pp.
\textsf{MR~97b:51009}

\bibitem{Scheerer}
H.~Scheerer,
\textit{Transitive actions on Hopf homogeneous spaces.}
Manuscripta Math. 4 (1971) 99--134.
\textsf{MR 45\#2086}

\bibitem{Schneider}
V.~Schneider,
\textit{Transitive actions on highly connected spaces.}
Proc. Amer. Math. Soc. 38 (1973) 179--185.
\textsf{MR~47\#9658}

\bibitem{Schroth}
A.~E.~Schroth,
\textit{Topological circle planes and topological quadrangles.}
Pitman Research Notes in Mathematics 337,
Longman, Harlow (1995) x+155 pp. 
\textsf{MR~97b:51010}

\bibitem{Seitz91}
G.~M.~Seitz,
\textit{Maximal subgroups of exceptional algebraic groups.} 
Mem. Amer. Math. Soc. 90 (1991) iv+197 pp.
\textsf{MR~91g:20038}

\bibitem{Spa66}
E.~H.~Spanier,
\textit{Algebraic topology}.
Springer-Verlag, New York-Berlin (1981). xvi+528 pp.
\textsf{MR~96a:55001}

\bibitem{Stolz98}
S.~Stolz,
\textit{Multiplicities of {D}upin hypersurfaces.}
Invent. Math. 138 (1999) 253--279.
\textsf{MR~2001d:53065}

\bibitem{Strau96}
Th.~Strau\ss,
\textit{Cohomology rings, sphere bundles, and double mapping cylinders.} 
{D}iplomarbeit, {W\"u}rzburg: {M}ath. {I}nst., {U}niv. {W\"u}rzburg,
(1996) 45 pp.

\bibitem{Stro92}
M.~Stroppel,
\textit{Reconstruction of incidence geometries from groups of automorphisms.}
Arch. Math. 58 (1992) 621--624.
\textsf{MR~93e:51026}

\bibitem{Stro93}
M.~Stroppel,
\textit{A categorical glimpse at the reconstruction of geometries.}
Geom. Dedicata 46 (1993) 47--60.
\textsf{MR~94c:51036}

\bibitem{Strub86}
W.~Str{\"u}bing,
\textit{Isoparametric submanifolds.}
Geom. Dedicata 20 (1986) 367--387.
\textsf{MR~87j:53084}

\bibitem{Sze74}
J.~Szenthe,
\textit{On the topological characterization of transitive {L}ie group
  actions.}
Acta Sci. Math. (Szeged) 36 (1974) 323--344.
\textsf{MR~50\#13368}

\bibitem{Takagi}
R.~Takagi,
\textit{A class of hypersurfaces with constant principal curvatures
 in a sphere.}
J. Differential Geom. 11 (1976) 225--233.
\textsf{MR 54\#13798}

\bibitem{TkgTsh72}
R.~Takagi and T.~Takahashi,
\textit{On the principal curvatures of homogeneous hypersurfaces in a sphere.}
In: \textit{Differential geometry in honor of K. Yano}, p. 469--481, 
Kinokuniya, Tokyo (1972).
\textsf{MR~48\#12413}

\bibitem{Tent99}
K.~Tent,
\textit{Very homogeneous generalized polygons of finite Morley rank.}
J. London Math. Soc. 62 (2000) 1--15.
\textsf{MR~2001f:51011}

\bibitem{Te85}
C.-L.~Terng,
\textit{Isoparametric submanifolds and their {C}oxeter groups.}
J. Differential Geom. 21 (1985) 79--107.
\textsf{MR~87e:53095}

\bibitem{Tho91}
G.~Thorbergsson,
\textit{Isoparametric foliations and their buildings.}
Ann. of Math. 133 (1991) 429--446.
\textsf{MR~92d:53053}

\bibitem{Tho92}
G.~Thorbergsson,
\textit{Clifford algebras and polar planes.}
Duke Math. J. 67 (1992) 627--632.
\textsf{MR~93i:51033}

\bibitem{ThoSur}
G.~Thorbergsson,
\textit{A survey on isoparametric hypersurfaces and their
 generalizations.}
In: \textit{Handbook of differential geometry, I}, P. 963--995
North Holland, Amsterdam (2000).
\textsf{MR~2001a:53097}

\bibitem{TitsTabellen}
J.~Tits,
\textit{Tabellen zu den einfachen {L}ie {G}ruppen und ihren
  {D}arstellungen.}
LNM 40, Springer-Verlag, Berlin (1967). v+53 pp.
\textsf{MR~36\#1575}

\bibitem{Tits71}
J.~Tits,
\textit{Repr\'esentations lin\'eaires irr\'eductibles d'un groupe
 r\'eductif sur un corps quelconque.}
J. Reine Angew. Math. 247 (1971) 196--220.
\textsf{MR~43\#3269}

\bibitem{TitsLNM}
J.~Tits,
\textit{Buildings of spherical type and finite BN-pairs}. 
LNM 386,
Springer-Verlag, Berlin (1974). x+299 pp.
\textsf{MR~57\#9866}
   
\bibitem{Tits77}
J.~Tits,
\textit{Endliche Spiegelungsgruppen, die als Weylgruppen auftreten.}
Invent. Math. 43 (1977) 283--295. 
\textsf{MR~57\#478}

\bibitem{Tits-Weiss}
J.~Tits and R.~Weiss,
\textit{Moufang polygons}.
Book in preparation (2001).

\bibitem{Uch80}
F.~Uchida,
\textit{An orthogonal transformation group of $(8k-1)$-sphere.}
J. Differential Geom. 15 (1980) 569--574.
\textsf{MR~83a:57056}

\bibitem{HvM98}
H.~Van~Maldeghem,
\textit{Generalized polygons}.
Birkh\"auser Verlag, Basel (1998). xvi+502 pp.
\textsf{MR~2000k:5104}

\bibitem{Wang}
H.~C. Wang,
\textit{Homogeneous spaces with non-vanishing Euler characteristics.}
Ann. Math. 50 (1949) 925--953.
\textsf{MR 11,326c}

\bibitem{Wan88}
Q.~M.~Wang,
\textit{On the topology of {C}lifford isoparametric hypersurfaces.}
J. Differential Geom. 27 (1988) 55--66.
\textsf{MR~89e:53093}

\bibitem{Whi78}
G.~W.~Whitehead,
\textit{Elements of homotopy theory}.
Springer-Verlag, New York-Berlin (1978). xxi+744 pp.
\textsf{MR~80b:55001}

\bibitem{Wol}
M.~Wolfrom,
{Ph.D.} {T}hesis, W\"urzburg: {M}ath. {F}ak., {U}niv. W\"urzburg,
in preparation (2001).

\end{thebibliography}
\end{document}